\documentclass[12pt,leqno]{amsart}
\usepackage{amsmath}
\usepackage{amssymb}
\usepackage{amsthm} 
\usepackage{epsf}
\usepackage{graphicx}
\usepackage{esint} 
\usepackage{dsfont}
\usepackage[pagebackref]{hyperref} 
\hypersetup{backref,  pdfpagemode=FullScreen,  colorlinks=true} 
\usepackage{color}
\usepackage{enumerate}

\numberwithin{equation}{section}

\newcommand{\nn}{\nonumber}
\newcommand{\ms}{\medskip}

\newcommand{\R}{\mathbb{R}}
\renewcommand{\H}{\mathcal H}

\newcommand{\bN}{\mathbb{N}}

\renewcommand{\d}{\partial}
\newcommand{\dist}{\,\mathrm{dist}}

\newcommand{\sm}{\setminus}
\DeclareMathOperator{\supp}{supp}
\DeclareMathOperator{\diam}{diam}

\newcommand{\wt}{\widetilde}
\newcommand{\wh}{\widehat}
\newcommand{\ol}{\overline}
\newcommand{\ub}{\underbar}

\newcommand{\cW}{{\mathcal  W}}

\newcommand{\cD}{{\mathcal  D}}

\newcommand{\cR}{{\mathcal  R}}

\newcommand{\1}{{\mathds 1}}

\newcommand{\dr}{\partial}
\DeclareMathOperator*{\osc}{osc}
\newcommand{\ds}{\displaystyle}
\newcommand{\WW}{W_r}

\newcommand{\g}{\mathfrak g}

\newcommand{\D}{\mathbb D}

\theoremstyle{plain}
\newtheorem{theorem}[equation]{Theorem}
\newtheorem{lemma}[equation]{Lemma}
\newtheorem{corollary}[equation]{Corollary}
\newtheorem{proposition}[equation]{Proposition}

\theoremstyle{definition}
\newtheorem{definition}[equation]{Definition}

\theoremstyle{remark}
\newtheorem{remark}[equation]{Remark}

\newcommand{\RR}{{\mathbb{R}}}

\newcommand{\ZZ}{{\mathbb{Z}}}

\newcommand{\eps}{\varepsilon}

\newcommand{\A}{\mathcal{A}}
\newcommand{\C}{\mathcal{C}}

\newcommand{\po}{{\partial\Omega}}

\newcommand{\bp}{\noindent {\it Proof}.\,\,}
\newcommand{\ep}{\hfill$\Box$ \vskip 0.08in}

\DeclareMathOperator{\Tr}{Tr}
\DeclareMathOperator{\Ext}{Ext}
\def\div{\mathop{\operatorname{div}}}

\newcounter{hcountno}
\newcommand{\hcount}[1]{\refstepcounter{hcountno} \label{#1}}

\begin{document}

\title{Elliptic theory in domains with boundaries of mixed dimension}

\author[David]{Guy David}
\address{Guy David. Universit\'e Paris-Saclay, CNRS, Laboratoire de math\'ematiques d'Orsay, 91405 Orsay, France}
\email{guy.david@universite-paris-saclay.fr}

\author[Feneuil]{Joseph Feneuil}
\address{Joseph Feneuil. Mathematical Sciences Institute, Australian National University, Acton, ACT 2601, Australia}
\email{joseph.feneuil@anu.edu.au}

\author[Mayboroda]{Svitlana Mayboroda}
\address{Svitlana Mayboroda. School of Mathematics, University of Minnesota, Minneapolis, MN 55455, USA}
\email{svitlana@umn.edu}

\thanks{David was partially supported by the European Community H2020 grant GHAIA 777822, and the Simons Foundation grant 601941, GD. Mayboroda was supported in part by the Alfred P. Sloan Fellowship, the NSF grants DMS 1344235, DMS 1839077, and the Simons foundation grant 563916, SM}
\date{}

\maketitle

 
\begin{abstract}
Elliptic theory in domains with boundaries of mixed dimension.

Take an open domain $\Omega \subset \R^n$ whose boundary may be composed of pieces of different dimensions.  For instance, $\Omega$ can be a ball on $\R^3$, minus one of its diameters $D$, 
or  
a so-called saw-tooth domain, with a boundary consisting of pieces of 1-dimensional curves intercepted by 2-dimensional spheres. It could also be a domain with a fractal (or partially fractal) boundary. Under appropriate geometric assumptions, essentially 
the existence of doubling measures on $\Omega$ and $\d \Omega$ with appropriate size conditions-- 
we construct a class of second order degenerate elliptic operators $L$ adapted to the geometry, and establish key estimates of elliptic theory associated to those operators. This includes boundary Poincar\'e and Harnack inequalities, maximum principle, and H\"older continuity of solutions at the boundary. We introduce Hilbert spaces naturally associated to the geometry,  construct appropriate trace and extension operators, and use them to define weak solutions to $Lu=0$. Then we prove De Giorgi-Nash-Moser estimates inside $\Omega$ and on the boundary, solve the Dirichlet problem and thus construct an elliptic measure $\omega_L$ associated to $L$. We construct Green functions and use them to prove a comparison principle and the doubling property for $\omega_L$. 
Since our theory emphasizes measures, rather than the geometry per se, the results are new even in the classical setting of a half-plane $\RR^2_+$ when the boundary $\partial \RR^2_+= \RR$ is equipped with a doubling measure $\mu$ singular with respect to the Lebesgue measure on $\RR$. Finally, the present paper provides a generalization of the celebrated Caffarelli-Sylvestre extension operator from its classical setting of $\RR^{n+1}_+$ to general open sets, and hence, an extension of the concept of fractional Laplacian to Ahlfors regular boundaries and beyond. 

\smallskip\noindent
R\'ESUM\'E. 
Th\'eorie elliptique dans des domaines \`a fronti\`eres de dimension mixte.

Soit $\Omega \subset \R^n$ un domaine dont la fronti\`ere peut contenir des morceaux de dimensions diff\'erentes. Par exemple, $\Omega$ peut \^etre une boule de $\R^3$, moins l'un de ses diam\`etres $D$, ou un domaine en dents de scies, avec une fronti\`ere compos\'ee de morceaux de courbes et de morceaux de sph\`eres. Ou encore, un domaine avec une fronti\`ere (partiellement) fractale. Avec des hypoth\`eses g\'eom\'etriques convenables, essentiellement 
l'existence de mesures doublantes sur $\Omega$ et $\d\Omega$ de tailles appropri\'ees, on construit une classe d'op\'erateurs elliptiques d'ordre $2$ d\'eg\'en\'er\'es de mani\`ere adapt\'ee \`a la g\'eom\'etrie, et on prouve les estimations cl\'e associ\'ees \`a ces op\'erateurs $L$. Ceci inclue des in\'egalit\'es de Poincar\'e et de Harnack, le principe du maximum, et la continuit\'e H\"old\'erienne \`a la fronti\`ere des solutions. On introduit les espaces de Hilbert naturellement associ\'es \`a la g\'eom\'etrie, on construit les op\'erateurs de trace et d'extension associ\'es, on les utilise pour d\'efinir les solutions faibles de $Lu=0$, puis on prouve les in\'egalit\'es de De Giorgi-Nash-Moser dans $\Omega$ et \`a la fronti\`ere, on r\'esout le probl\`eme de Dirichlet, qu'on utilise pour construire une mesure elliptique $\omega_L$ associ\'ee \`a $L$. On construit les fonctions de Green et on les utilise pour obtenir le principe de comparaison et la propri\'et\'e doublante pour $\omega_L$. 
Notre th\'eorie \'etant centr\'ee sur les mesure, en pas seulement sur la geometrie de $\Omega$, les r\'esultats sont nouveaux m\^eme  dans le cas classique du demi-plan $\RR^2_+$, mais o\`u la fronti\`ere 
$\partial \RR^2_+= \RR$ est munie d'une mesure doublante $\mu$ singuli\`ere par rapport \`a la mesure de Lebesgue sur $\RR$. Finalement, ce papier donne une g\'en\'eralisation du c\'el\`ebre op\'erateur d'extension de  Caffarelli-Sylvestre, depuis son cadre classique de $\RR^{n+1}_+$ vers des ouverts plus g\'en\'eraux, et donc une extension du concept de Laplacien fractionnaire \`a des fronti\`eres Ahlfors reguli\`eres et au del\`a. 

\smallskip\noindent
SHORT ABSTRACT.
We study the initial regularity results (H\"older continuity, De Giorgi-Nash-Moser inequalities, 
maximum principle, existence and doubling property for the elliptic measure, and estimates for the 
Green function) for a class of second order elliptic operators associated to the geometry of a domain,
whose boundary can have pieces of different dimensions, but where we have two related doubling measures, one on the domain and one on the boundary.

\smallskip\noindent
R\'ESUM\'E COURT. 
On \'etudie les r\'esultats de r\'eglarit\'e initiaux (continuit\'e H\"old\'erienne, in\'egalit\'es de
De Giorgi-Nash-Moser, principe du maximum, existence et propri\'et\'e doublante de la mesure harmonique, estim\'ees pour la fonction de Green) pour une classe d'op\'erateurs elliptiques du second degr\'e associ\'ee  \`a la g\'eom\'etrie d'un domaine dont la fronti\`ere peut avoir des morceaux de dimensions diverses, mais avec deux mesures doublantes li\'ees, l'une sur le domaine et l'autre sur la fronti\`ere.

\end{abstract}

\ms\noindent{\bf Key words:} 1-sided NTA domains, boundaries of mixed dimensions, homogeneous weighted Sobolev spaces, Trace theorem, Extension theorem, Poincar\'e inequalities, degenerate elliptic operators, De Giorgi-Nash-Moser estimates, harmonic measure, Green function, comparison principle.

\ms\noindent{\bf Mots cl\'e:} domaines \`a acc\`es non tangentiel, fronti\`eres de dimensions mixtes, espace de Sobolev homog\`enes \`a poids, th\'eor\`eme de trace, th\'eor\`eme d'extension, in\'egalit\'es de Poincar\'e, op\'erateurs elliptiques d\'eg\'en\'er\'es, estimations de De Giorgi-Nash-Moser, mesure harmonique, fonction de Green, principe de comparaison.

\ms\noindent
AMS classification (MSC2020): 28A15, 28A25, 31B05, 31B25, 35J25, 35J70, 42B37.

\tableofcontents

\section{Motivation and a general overview of the main results}


%
%
%
%
%


\subsection{Motivation}
Massive efforts of the past few decades at the intersection of analysis, PDEs, and geometric measure theory have recently culminated in a comprehensive understanding of the relationship between the absolute continuity of the harmonic measure with respect to the Hausdorff measure and rectifiability of the underlying set \cite{AHM3TV, AHMMT}. Even more recently, in 2020, we could identify a sharp class of elliptic operator for which the elliptic measure behaves similarly to that of the Laplacian in the sense that analogues of the above results could be obtained, at least under mild additional topological assumptions \cite{HMMTZ2}, \cite{KePiDrift}.

Unfortunately, all of those results have been restricted to the case of $n$-dimensional domains with $n-1$ dimensional boundaries, and as such, left completely beyond the scope of the discussion a higher co-dimensional case, such as, for example, a complement of a curve in $\RR^3$. The authors of the present paper have recently launched a program investigating the latter, which we will partially review below, and which quite curiously brought a completely different level of understanding of $n-1$ dimensional results and a plethora of open problems, again, relevant even in the context of ``classical" geometries, e.g., simply connected planar domains or even a half-space. What is the role of measure on the boundary and given a rough measure, possibly singular with respect to the Hausdorff measure, can we define an elliptic operator whose solutions would be well-behaved near the boundary? What is the role of the dimension, especially when fractional dimensions are allowed? Even in the case of the Laplacian the dimension of the harmonic measure is a mysterious and notoriously difficult subject with scarce celebrated results due to Makarov, Bourgain, Wolff, and many problems open to this date, but what if we step out of the context of the Laplacian and similar operators? Closely related to this question is another one: what is the role of degeneracy, that is, where are the limits of the concept of ``ellipticity" which could still carry reasonable PDE properties. This brought us, in particular, to a new version of the Caffarelli-Sylvestre extension operator and hence, a new fractional Laplacian (or, one could say, a certain form of  differentiation) on general Ahlfors regular sets. Let us discuss all this in more details.

As we mentioned above, this project started as a continuation of efforts in \cite{DFM-CRAS}, 
\cite{DFMprelim}, \cite{DFMDahl},  \cite{DFMKenig}, \cite{DEM}, \cite{MZ}, \cite{FMZ} to define an analogue of harmonic measure on domains with lower dimensional boundaries and ultimately to develop a PDE theory comparable in power and scope to that of $n-1$-rectifiable sets. 
Initially, we focused on domains $\Omega \subset \R^n$ 
whose boundary $\Gamma = \d \Omega$ is Ahlfors regular of dimension $d < n-1$ 
(see \eqref{defADR} below). When $d \leq n-2$, such sets would not be recognized by harmonic functions, and we were led to a class of degenerate elliptic operators $L$ adapted to the dimension. Taking the coefficients of $L$ to be, roughly speaking, of the order of $\dist(x, \po)^{-(n-d-1)}$, we managed to define a well behaved elliptic measure $\omega_L$ associated to $L$ and $\Omega$ and prove the estimates for $\omega_L$ and for the Green functions, similar to the classical situation where $d=n-1$ and $L$ is elliptic. Furthermore, we proved in \cite{DFMDahl} that $\omega_L$ is absolutely continuous with respect to the Hausdorff measure $\mu = \H^d_{\vert \Gamma}$, with an $A_\infty$ density, when $\Gamma$ is a Lipschitz graph with a small Lipschitz constant and the coefficients of $L$ are proportional to $D(x, \po)^{-(n-d-1)}$, where $D$ is a carefully chosen, appropriately smooth, distance function. However, in an effort to extend these results to the context of uniformly rectifiable domains we faced some fundamental problems which bring us to the setting of the present paper. 

A key feature of (uniformly) rectifiable sets is the fact that at every scale a significant portion of such a set can be suitably covered by well-controlled Lipschitz images. To take advantage of this, one has to develop an intricate procedure which allows one to ``hide the bad parts" and more precisely, it is absolutely essential to be able to consider suitable subdomains of an initial domain which carry similar estimates on harmonic measure, within the scales under consideration. The latter are referred to as the saw-tooth domains and the reader can imagine ``biting off" from the initial domain a ball, or rather a cone, surrounding a bad subset of the boundary. The problem is that when the initial domain is, say, the complement of a curve in $\RR^3$, any subdomain would have a boundary of a mixed dimension and the specific procedure that we are describing yields pieces of one-dimensional curves intercepted by 2-dimensional spheres, or more precisely, 2-dimensional Lipschitz images. We will give in Section~\ref{SExa} 
a careful description of this example. Similarly, any attempt to localize a problem on a set with lower dimensional boundary (e.g., $\RR^n\setminus \RR^d$) yields a new domain, given by an $n$ dimensional ball minus a $d$-dimensional curve, which now has a union of an $n-1$ dimensional sphere and a $d$-dimensional surface as its boundary. These challenges led us to a necessity to develop a meaningful elliptic theory in the presence of the mixed-dimensional boundaries.

This immediately raises a question: what are the appropriate elliptic operators, as our favorite choice $L=-\div D(x, \po)^{-(n-d-1)} \nabla $, and similar ones, carry a power which depends on the dimension of the boundary $d$. To some extent, this is necessary: as we mentioned above, the Laplacian would not see very low-dimensional sets and this argument can be generalized. But to which extent? Can  $L=-\div D(x, \po)^{-(n-d-1+\beta)}\, \nabla $ be allowed for some $\beta$?  Can  $L=-\div D(x, \po)^{-\beta}\, \nabla $ be allowed for some $\beta$ for ``classical" domains with $n-1$ dimensional boundaries? 

Even for the Laplacian these issues are extremely challenging. Fundamental results of Makarov \cite{Mak1}, \cite{Mak2} establish that on the plane, the Hausdorff dimension $\dim_\H \omega$ is equal to $1$ if the set $\partial\Omega$ is connected. More generally, for any domain $\Omega$ on the Riemann sphere whose complement has positive logarithmic capacity there exists a subset of $E\subset \po$ which supports harmonic measure in $\Omega$ and has Hausdorff dimension at most 1 \cite{JW1}. In particular, the supercritical regime is fully characterized on the plane: if $s \in (1, 2)$, $0 < \H^s(E) < \infty$, then $\omega$ is always singular with respect to $\H^s|_E$. However, for $n>2$ the picture is far from being well-understood. On the one hand, Bourgain \cite{Bo} proved that the dimension of harmonic measure always drops: $\dim_\H \omega < n$. On the other hand, even for connected $E=\partial\Omega$, it turns out that $\dim_\H \omega$ can be strictly bigger than $n-1$, due to a celebrated counterexample of Wolff \cite{W}. Some recent efforts in this direction include, e.g.,  \cite{Azzam}, but overall the problem of the dimension of the harmonic measure remains open, and to the best of our knowledge there exist no results for other elliptic operators, with the only exception of \cite{Sw}. Definitely we have not encountered any results of this type for degenerate elliptic operators. 

On the other hand, in a more benign geometric setting degenerate operators have of course been studied in the literature. The most obvious example resonating with our setting is the celebrated Caffarelli-Sylvestre extension operator. In \cite{CS} the authors proposed that the fractional Laplacian $(-\Delta)^\alpha$, $\alpha\in (0,1)$, on $\RR^d$ can be realized as a Dirichlet-to-Neumann mapping corresponding to the operator $L=-\div \dist (\cdot, \RR^d)^{-\beta}\,\nabla$ on $\RR^{d+1}$ with $\beta=2\alpha-1$, $\beta\in (-1,1).$ This turned out to be an extremely fruitful idea, facilitating many properties of the fractional Laplacian and similar operators, and was extended to other $\alpha$ by A. Chang and R. Yang in \cite{CY}. One of the outcomes of the present paper is an extension of the elliptic theory to the complement of any $d$-dimensional Ahlfors regular set for  $L=-\div A(x) \,\dist(x, \po)^{-(n-d-1+\beta)} \nabla$, $\beta\in (-1, 1)$, including the Caffarelli-Sylvestre extension operator and generalizing it to extremely rough geometric situations, fractal sets, mixed dimensions, etc. 

We point out, parenthetically, that while this paper concerns the fundamental elliptic estimates, we plan to address also absolute continuity of elliptic measure for this type of operators in the forthcoming publications. It is slightly surprising that such a study has not been pursued before even in the $\RR^{d+1}$ setting, but this seems to be the case. The only known results pertain to the degenerate operators with weights independent of the distance to the boundary (see, e.g., \cite{ARR}).  

Returning to the general elliptic theory, a search for the appropriate assumptions on the boundary and the coefficients of the corresponding allowable elliptic operators have quickly revealed that the key players are the measure $\mu$ on $\po$ and the corresponding measure $m$ in $\Omega$ which will define the ``ellipticity" of $L$. This brings out two more issues. First, even in the simple case of the half-plane $\RR^{2}_+$ there is another layer of complexity possibly introduced by the boundary measure. Specifically, one can ask whether there exists an elliptic operator which is well-behaved with respect to an arbitrary doubling measure $\mu$ on the boundary, for instance, the one furnished using the Riesz products on $\RR$. In the present paper we allow the measure $\mu$ on $\Gamma=\dr \Omega$ to be wild: we will present in Section \ref{S3.6} an example where, around any point in $\Gamma$, $\mu$ is not absolutely  continuous with respect to the surface measure, and yet the corresponding elliptic operator has well-behaved solutions. Secondly, one can encounter a matrix of coefficients which is (also) degenerate at interior points of $\Omega$. 

\subsection{Additional historical comments} Obviously this paper is not the first one where degenerate elliptic operators were studied. Perhaps most closely related to our questions are the works of Fabes, Jerison, Kenig, Serapioni, \cite{FKS,FJK,FJK2}, Maz'ya \cite{Mazya11}, Heinonen, Kipel\"ainen, Martio \cite{NDbook}, Ammann, Bacuta, Mazzucato, Nistor, Zikatonov \cite{AN,BMNZ}, and, as far as Sobolev-Poincar\'e inequalities and similar questions are concerned, Haj\l asz and Koskela  \cite{HK2,HK00}. One can see our operators as elliptic operators with drifts, and special cases arising in modeling population genetics were studied in \cite{EM1,EM2,EP1,EP2,EP3}. While different from the scope of this paper, in some respects they guided our intuition, and there is even some overlap with our results. However, typically their stress is rather on the singularities of the weight inside the domain; here we emphasize its behavior near the boundary $\Gamma$, and, respectively, the behavior of solutions near $\po$ depending on the geometry and the underlying measure on $\po$ and on $\Omega$. When the impact of the boundary is considered, the aforementioned works concentrate on the Wiener criterion and surrounding questions, often of a qualitative nature, while we aim at the uniform scale-invariant quantitative results. 
And even more, the boundary results in \cite{FKS,FJK,FJK2} are stated for 2-sided NTA domains, 
which forces the existence of  
a big portion of the complement $\Omega^c$ around any point of the boundary $\Gamma$, 
a condition that we do not want to impose when  
a part of $\Gamma$ has codimension higher than 1 (like for instance when $\Omega$ is a ball deprived of a diameter).
 Also, on a more technical side, the estimates of \cite{FKS,FJK,FJK2} or \cite{NDbook} would be hard to use here, because we need to be able to consider unbounded domains and boundaries. Finally, once again, our coverage, including boundaries of mixed-dimension equipped with possibly complicated doubling measures, and the overall point of view of designing elliptic operators which respect the geometric and measure theoretic setting of the problem, ends up in a different range of results.

It is interesting to point out that an alternative route to generalization of elliptic theory to sets with lower dimensional boundaries consists of studying the $p$-Laplacian operator for a suitable range of values of $p$. This approach has been developed by Lewis, Vogel, Nystr\"om, and others -- see, e.g., \cite{LN} for boundary Harnack estimates; however, to the best of our knowledge, it did not yield the absolute continuity results for the underlying elliptic measure on uniformly rectifiable sets (not to mention that the role of the elliptic measure for a non-linear PDE is quite different) and for that reason we pursue a different route. 

\subsection{A rough outline of the main assumptions and results}
The general assumptions of this paper will be described precisely in the next section,
but let us give a first overview right now. We are given a domain $\Omega \subset \R^n$,
and a doubling measure $\mu$ on its boundary $\Gamma$. 
We are also given a doubling measure
$m$ on $\Omega$, which is assumed to be absolutely continuous with respect to the Lebesgue
measure (that is, $dm(X) = w(X) dX$ for some weight $w$). The reader can think of 
$w(X) = \dist(X,\Gamma)^{-\gamma}$ with  $\gamma \in (n-d-2, n-d)$, or even more general weights.
The key assumptions is a relation between our two doubling measures, that says that on balls $B(x,r)$ centered on $\Gamma$, one measure does not grow much faster than the other: 
\begin{equation} \label{HH5b-intro}
\frac{m(B(x,r) \cap \Omega)}{m(B(x,s) \cap \Omega)} 
\leq C \left(\frac rs\right)^{2-\epsilon} \frac{\mu(B(x,r))}{\mu(B(x,s))}
\qquad \text{ for $x\in \Gamma$, $0<s<r$,}
\end{equation}
for suitable $C, \epsilon>0$. This is the condition (\hyperref[H5]{H5}) below (or rather \eqref{HH5b}), 
and it is responsible for our requirement that $n-d-2 <\gamma< n-d$ above. It is somewhat surprising perhaps that we only need an estimate from above. 

We also have a requirement on $w$, related to its behavior far from $\Gamma$;
in the same spirit as in \cite{NDbook}, we demand a weak Poincar\'e estimate for the space $(\Omega, m)$ (with the usual metric),
which is explained below as (\hyperref[H6]{H6}). 
If $w$ is regular enough away from $\Gamma$, for instance if 
$\sup_{X \in B} w(X) \leq C \inf_{X \in B} w(X)$ for all balls $B$ such that  
$2B \subset \Omega$ and for some constant $C$ that does not depend on $B$, 
then (\hyperref[H6]{H6}) is automatically satisfied. This is the case in our previous papers,
and in the context of sawtooth domains which have been alluded to above. 

With these preliminaries, the operator $L = -\div A \nabla$ can be any elliptic operator as long as the ellipticity condition is satisfied with respect to our measure $m$; that is, we simply require that $w(X)^{-1} A(X)$ satisfy the standard  boundedness and ellipticity conditions on $\Omega$. 

The final set of assumptions pertains to connectivity. When $\Gamma$ is an Ahlfors regular set of dimension $d<n-1$, we do not need to add any topological conditions ensuring the (quantitative) connectedness of $\Omega$, because they are automatically satisfied. Here our setting allows boundaries of all dimensions, and in such a setting some topological restrictions are necessary \cite{BJ}. In line with many antecedents,  we require that $\Omega$ satisfy the ``one sided NTA conditions''. That is, we demand the existence of corkscrew balls
and Harnack chains in $\Omega$; see the conditions (\hyperref[H1]{H1}) and (\hyperref[H2]{H2}) 
in the next section, and the discussion that follows them.

\ms
All these assumptions will be described in detail in the next section,
and we will then give examples in Section \ref{SExa}.
Under these assumptions, we will be able to define an elliptic measure associated to $L$ and establish the fundamental properties for solutions.

First of all, in Section \ref{SecW}, we will define an energy space $W$, the Hilbert space of functions on $\Omega$
with a derivative in $L^2(m) = L^2(w(X) dX)$. In this section, $\Omega$ and $\Gamma$ are unbounded, so the space $W$ is a homogeneous space. This is the most useful scenario for our applications  but we also treat the extension of all our results to the case where $\Gamma$ or $\Omega$ are bounded in Section \ref{SBounded}.

An important tool in our theory, which allows us to dispense with the existence of
large balls in the complement of $\Omega$, or barrier functions, is a Poincar\'e
estimate at the boundary. The next two sections aim for that result.

In Section \ref{Scones}, we introduce some technical material, such as the dyadic pseudocubes
$Q$, $Q \in \D$, on $\Gamma$, an analogue of the Whitney cubes in $\Omega$ (the sets $\cW_Q$
of \eqref{defWQ}), and some non-tangential access regions $\gamma(x)$
and their truncated versions $\gamma_Q(x)$;
see near \eqref{defgamma}. The sets $\gamma^\ast(x)$ and $\gamma_Q^\ast(x)$ are analogue to  $\gamma(x)$ and $\gamma_Q(x)$, and are obtained from the latter by fattening them a bit so that $\gamma^\ast(x)$ and $\gamma_Q^\ast(x)$ are well connected sets. We rely on (\hyperref[H2]{H2}) for this procedure. The access cones $\gamma_Q^\ast(x)$ are used to constructs well connected tent sets $T_{2Q}$, which will advantageously substitute the sets $B\cap \Omega$ where $B\subset \R^n$ is a ball centered on $\Gamma$ (indeed, the sets $B\cap \Omega$ have no reason to be connected).
We use those tent sets to extend Poincar\'e inequalities given in (\hyperref[H6]{H6}) to sets that actually get close to the boundary. In particular, we prove in  Theorem \ref{Poincare1} that

\begin{theorem} \label{Main1}
There exists $k>1$ such that for any $u\in W$ and any $Q\in \D$,
\[\left( \fint_{T_{2Q}} |u- u_{T_{2Q}}|^{2k} \right)^{1/2k} \leq C \diam(Q) \left(  \fint_{T_{2Q}} |\nabla u|^2 \right)^{1/2},\]
where $u_{T_{2Q}}:= \fint_{T_{2Q}} u \, dm$.
\end{theorem}

Our next goal is to obtain a variant of the above theorem, where $u_{T_{2Q}}$ is removed but we assume that $u =0$ on the boundary $2Q$. To this end, we need a notion of trace. We define then a Sobolev space $H$ as the set of $\mu$-measurable functions $g$ on $\Gamma$ such that
\[\|g\|_H := \left( \int_\Gamma \int_\Gamma
\frac{\rho(x,|x-y|)^2 |g(x) - g(y)|^2}{m(B(x,|x-y|) \cap \Omega)} \, d\mu(y) d\mu(x)\right)^\frac12\]
is finite. As the reader can see, $H$ depends on both $\mu$ and $m$, and the dependence on $m$ is a bit surprising at first; but if one recalls that our objective is to construct a bounded trace from $W$ (that depends on $m$) and $H$, it makes sense. Here you can see $\rho(x,|x-y|)$ as a corrective term that takes into account how far $m(B(x,r) \cap \Omega)$ is from $r \mu(B(x,r))$. Of course, if $\mu$ and $m$ are intertwined so that $m(B(x,r) \cap \Omega) \approx r \mu(B(x,r))$ for $x\in \Gamma$ and $r>0$ - which will be the most natural situation - then the strange term $\rho(x,|x-y|)$ disappears and the space $H$ does not depend on $m$ anymore.

With the space $H$ at hand, we construct in Section \ref{STrace}  a bounded trace operator $\Tr$ from $W$ to $H$. We later build in Section \ref{Sext} a nice extension operator $\Ext : H \to W$, such that $\Tr \circ \Ext = I$. Those results are given in Theorem \ref{TraceTh} and Theorem \ref{ThExt}, which are summarized below.

\begin{theorem} \label{Main2}
There exists two bounded linear operators $\Tr :\, W \to H$ and $\Ext: \, H \to W$ such that for $u\in W$ and $\mu$-almost every $x\in \Gamma$,
\[\Tr u(x) = \lim_{\begin{subarray}{c} X \in \gamma(x) \\ \delta(X) \to 0 \end{subarray}} \fint_{B(X,\delta(X)/2)} u \, dm\]
and such that for $g \in H$ and $\mu$-almost every $x\in \Gamma$,
\[\Tr \circ \Ext g(x) = g(x).\]
\end{theorem}

By combining the trace which was just introduced with Theorem \ref{Main1}, we established that

\begin{theorem} \label{Main3}
There exists $k>1$ such that for $Q\in \D$ and for $u\in W$ such that $\Tr u = 0$ $\mu$-a.e. on $2Q$,
\[\left( \fint_{T_{2Q}} |u|^{2k} \right)^{1/2k} \leq C \diam(Q) \left(  \fint_{T_{2Q}} |\nabla u|^2 \right)^{1/2}.\]
\end{theorem}

The theorem above is a particular case of Theorem \ref{PoincareTh2}. 

Our next big objective is to get estimates on solutions to appropriate degenerate elliptic operator. To prepare for this, in Sections \ref{Sext} and \ref{SComplete}, we check some density and stability results for our spaces; these should not be surprising but they are very useful for our later arguments. In Section \ref{SLocal} we add one last bit of functional analysis, which is the definition of localized versions $W_r(E)$ of our space $W$, and the way they co-operate with the trace operator (Lemma~\ref{defTrWE}).

We start the study of our degenerate operators $L = \div A \nabla$
and their solutions in Section~\ref{SSolutions}.
We require $w(X)^{-1} A(X)$ to satisfy the usual
ellipticity conditions, so the bilinear form naturally associated to $L$ is coercive
on $W$, and getting weak solutions in $W$ with a given trace on $H$
is rather easy, with the help of the Lax-Milgram Theorem.

We define weak subsolutions, supersolutions, and solutions in our local $W_r(\Gamma)$ spaces,
and start studying their regularity properties.  
We first prove an interior Caccioppoli inequality
(Lemma~\ref{CaccioI}), then extend it to the boundary (Lemma \ref{CaccioB}), then
prove interior Moser estimates (Lemma \ref{MoserI}), and extend them to the boundary
(Lemma \ref{MoserB}).
The next step is to prove interior H\"older estimates (Lemma \ref{HolderI})
and Harnack inequalities (Lemmas \ref{HarnackI} and~\ref{HarnackI2}). Some  
of the proofs in this section are just sketched, since they use the same arguments as, e.g., in \cite{DFMprelim}. The reader may be interested in some of the results and not the others, and we do not want to state all of them here. The theory was developed with boundary estimates in our mind, so they are the ones that we shall first present here.

\begin{theorem}[Moser estimates on the boundary] \label{Main4}
Let $B$ a ball centered on $\Gamma$ and $u$ be a non-negative subsolution to $Lu=0$ in $2B \cap \Omega$ such that $\Tr u = 0$ $\mu$-a.e. on $2B$. Then
\[\sup_{B \cap \Omega} u \leq C \fint_{2B \cap \Omega} |u|\, dm.\]
\end{theorem}

\begin{theorem}[H\"older estimates on the boundary] \label{Main5}
Let $x\in \Gamma$ and $r>0$. Assume that $u$ be a solution to $Lu=0$ in $B(x,r) \cap \Omega$. Then for $0<s<r$,
\[\osc_{B(x,s) \cap \Omega} u \leq C \left( \frac{r}{s} \right)^\alpha \osc_{B(x,r) \cap \Omega} u + C \osc_{B(x,\sqrt{sr}) \cap \Gamma} \Tr u,\]
where $C$ and $\alpha$ are positive constants independent of $x$, $s$, $r$, and $u$.
\end{theorem}

Of course, the Harnack inequality below, in particular when the weight is {\bf not} bounded from above or below by a positive constant, is interesting on its own right.

\begin{theorem}[Harnack inequality] \label{Main6}
Let $B$ be a ball such that $2B \subset \Omega$ and let $u$ be a non-negative solution to $Lu=0$ in $2B$. Then
\[\sup_B u \leq C \inf_B u.\]
\end{theorem}

We continue our article with a construction of the harmonic measure. The construction is classical, and relies on the maximum principle (Lemma \ref{lMPs}).

\begin{theorem}[Maximum principle] \label{Main7}
Let $u \in W$ be a solution to $Lu=0$ in $\Omega$. Then
\[\sup_\Omega u \leq \sup_\Gamma \Tr u \quad \text{ and } \quad \inf_\Omega u \geq \inf_\Gamma \Tr u.\]
\end{theorem}

The maximum principle combined to the Lax-Milgram theorem allows us to
solve the Dirichlet problem for compactly supported continuous functions on $\Gamma$ (Lemma \ref{ldefhm}), and thus define the desired harmonic measure $\omega_L^X$ with the Riesz representation theorem (Lemma~\ref{ldefhm2}).

\begin{theorem} \label{Main8}
For any $X \in \Omega$, there exists a unique positive Borel measure $\omega^X:=\omega^X_L$ on $\Gamma$ such that for any continuous and compactly supported $g \in H$, we have
\[u_g(X) = \int_\Gamma g(y) d\omega^X(y),\]
where $u_g$ is the solution in $W$ given by the Lax-Milgram theorem to $Lu = 0$ in $\Omega$ and $\Tr u = g$.

Furthermore, $\omega^X$ is a probability measure, that is $\omega^X(\Gamma) = 1$.
\end{theorem}


We end the article by building Green functions and 
using them to prove the non-degeneracy and the doubling property of the harmonic measure, as well as a comparison principle (which, applied to the harmonic measure, is also called change of pole property). The Green functions and the comparison principle have been 
companions of the mathematicians for ages. Maybe the first people to intensively study the Green functions in the case of general (non-degenerate) elliptic operators are Littman, Stampacchia, and Weierberger \cite{LSW}. Gr\"uter and Widman deepened the analysis of Green functions and established a comparison principle \cite{GW}. Fabes, Jerison, and Kenig worked with degenerate 
operators in \cite{FJK}, \cite{FJK2}, and some of their results are very similar to ours. However, those authors worked with bounded and 2-sided Non Tangentially Accessible domains, while we are interested in unbounded and weaker 1-sided NTA domains.
The Green functions were studied for systems in \cite{HoK} and \cite{DK} by assuming only De Giorgi-Nash-Moser estimates (and in particular not the maximum principle). We do not follow this route since the maximum principle is a prerequisite for the construction of the harmonic measure, which is the object that we are particularly interested in.

The Green function $g(x,y)$ is function on $\Omega \times \Omega$ such that
\begin{equation} \label{GreenIntro1}
\left\{ \begin{array}{l} 
L g(.,y) = \delta_y \text{ in } \Omega \\
g(.,y) \equiv 0 \text{ on } \Omega,
\end{array}\right.
\end{equation}
where $\delta_y$ is the delta distribution centered on $y$. We follow the strategy of \cite{GW}, in particular, we define the Green functions $g(.,y)$ as a limit of solutions in $W^{1,2}(\Omega,m)$ given by the Lax-Milgram theorem (and not by  taking the inverse of the operator $L$ on measures as in \cite{FJK}). The properties of the Green functions are given in Theorem \ref{GreenEx}, Lemma \ref{GreenSym}, 
Lemma~\ref{GreenDi}, Lemma~\ref{GreenFS} 
and Lemma \ref{GreenUn}. For instance we have the following pointwise bounds.

\begin{proposition} 
For $x,y\in \Omega$ such that $|x-y| \geq \dist(y,\Gamma)/10$, 
\[0 \leq g(x,y) \leq C\frac{|x-y|^2}{m(B(y,|x-y|) \cap \Omega)},\]
and for $x,y\in \Omega$ such that $|x-y|\leq \dist(y,\Gamma)/2$,
\[C^{-1} \int_{|x-y|}^{\dist(y,\Gamma)} \frac{s^2}{m(B(y,s))} \frac{ds}{s} \leq g(x,y) \leq C \int_{|x-y|}^{\dist(y,\Gamma)} \frac{s^2}{m(B(y,s))} \frac{ds}{s},\]
where the constant $C$ is of course independent of $x$ and $y$.
\end{proposition}

The harmonic measure is non-degenerate, 
in the following sense.  

\begin{theorem} \label{NDIntro}
Let $B$ be a ball centered on $\Gamma$ . Then
\[\omega^X(B \cap \Omega) \geq C^{-1} \qquad \text{ for } X \in \frac12 B \cap \Omega.\]
\end{theorem}

We prove a comparison principle between harmonic measures and Green functions. We need to define corkscrew points: $X_0 \in \Omega$ is a Corkscrew point associated to a ball $B=B(x,r)$ if
\[ |X_0-x| \leq r \ \text{ and } \ \dist(X_0,\Gamma) \geq \epsilon r,\]
for some $\epsilon>0$ that depends only on $\Omega$. 
We will assume in (\hyperref[H1]{H1}) that such points always exists.

\begin{theorem} \label{GHMIntro}
Let $B$ be a ball centered on $\Gamma$ and let $X_0$ be a corkscrew point associated to B. Then
\[C^{-1} \frac{m(B \cap \Omega)}{r^2} g(X,X_0) \leq \omega^X(B \cap \Gamma) \leq C \frac{m(B \cap \Omega)}{r^2} g(X,X_0) 
\ \ \ \text{ for }
X\in \Omega \setminus 2B.\]
\end{theorem}

With this comparison in hand, we show that the harmonic measure is doubling.

\begin{theorem} \label{DoubIntro}
Let $B$ be a ball centered on $\Gamma$. Then
\[\omega^X (2B\cap \Omega) \leq C \omega^X(B\cap \Gamma) \ \ \ \text{ for }
X\in \Omega \setminus 4B.\]
\end{theorem}

The change of pole property comes next.

\begin{theorem}
Let $B$ be a ball centered on $\Gamma$ and $X_0$ be a corkscrew point associated to B. Let $E,F \subset \Gamma \cap B$ be two Borel subsets of $\Gamma$ such that $\omega^{X_0}(E)$ and $\omega^{X_0}(F)$ are positive. 
Then
\[C^{-1}\frac{\omega^{X_0} (E)}{\omega^{X_0} (F)} \leq \frac{\omega^{X} (E)}{\omega^{X} (F)} \leq C \frac{\omega^{X_0} (E)}{\omega^{X_0} (F)} \ \ \ \text{ for }
X\in \Omega \setminus 2B.\]
\end{theorem}

At last, we give properties on the harmonic measure analogous to Theorems \ref{NDIntro}, \ref{GHMIntro}, \ref{DoubIntro} but for $\omega^X(\Gamma \setminus B)$ instead. We use them to prove a comparison principle for positive local solutions.

\begin{theorem}
There exists a large $K\geq 2$ that depends on how well $\Omega$ is connected (if $\Omega$ is well connected, we can take $K=2$)  such that the following holds.

Let $B$ be a ball centered on $\Gamma$, and let $X_0 \in \Omega$ be a corkscrew point associated to $B$. 
Let $u, v$ be two non-negative,
not identically zero, solutions to $Lu = Lv = 0$ in $KB \cap \Omega$ which are zero on the large boundary ball 
$KB \cap \Gamma$. Then
\[
C^{-1} \frac{u(X_0)}{v(X_0)} \leq \frac{u(X)}{v(X)} \leq C \frac{u(X_0)}{v(X_0)}
\ \ \ \text{ for } X \in \Omega \cap B.
\]
\end{theorem}

We have to be a little careful with this theorem, because we will allow 
the case when $\Gamma$ is bounded and $\Omega$ is the unbounded connected component of $\R^n \setminus \Gamma$. Then for large balls $B$, it will happen that
$\omega^X(\Gamma \setminus B) =0$ for $X\in \Omega$, and then we do not have 
bound on $\omega^X(\Gamma \setminus B)$ as in Theorem \ref{NDIntro}

\medskip

In this paper we only worry about the properties of $\omega_L$ in a  general setting. Then one may  
continue the study with more specific situations and carefully chosen operators.
Some of our earlier results, such as the extension of Dahlberg's result
in \cite{DFMDahl}, namely the $A_\infty$ absolute continuity of $\omega_L$ when $\Gamma$ is a 
Lipschitz graph with small constant and $L$ is well chosen, also work when 
$w(X) = \dist(X,\Gamma)^{-\gamma}$, $\gamma \in (n-d-2, n-d)$. 
This too will be studied more systematically in upcoming publications.

\section{Our assumptions}
\label{S2}

We aim to develop an elliptic theory on an open domain $\Omega$, equipped with a measure $m$. 
We are particularly interested in 
boundary estimates, and we want to be able to deal with a large class of measures $\mu$ (supported) on the boundary $\Gamma := \dr \Omega$.

Let us review previously known settings  
that we aim to generalize. 
Recall that a $d$-dimensional Ahlfors regular set $E \subset \R^n$ - denoted $d$-AR for short - is a set for which there exists a measure $\sigma$ supported on $E$ and a constant $C>0$ such that 
\begin{equation} \label{defADR}
C^{-1} r^d \leq \sigma(B(x,r)) \leq Cr^d \qquad \text{ for } x\in E, \, r>0.
\end{equation}
It is well known that if \eqref{defADR} holds for a measure $\sigma$ and a constant $C$, then \eqref{defADR} also holds for $\sigma' = \mathcal H_{|E}^d$ - the $d$-dimensional Hausdorff measure restricted to $E$ - and another constant $C'$.

The ``classical setting'' consists in taking 
an open domain $\Omega \subset \R^n$ such that its boundary $\Gamma$ is a $(n-1)$-AR set. The measure $m$ is taken as the $n$-dimensional Lebesgue measure and we choose $\mu$ as the surface measure on $\Gamma$, or in fact any measure satisfying \eqref{defADR}. 
The properties of elliptic PDEs in this context on relatively nice (e.g., Lipschitz) domains have been 
studied for 50 years;  see \cite{Stampacchia65,GW,CFMS} to cite a few,  and \cite{KenigB} for a extended presentation  of the properties. A more challenging setting of domains with uniformly rectifiable or even general AR boundaries came to the focus of development in the last 20 years. Unfortunately, there is no good single reference reviewing the underlying elliptic theory, but we can generally point the reader to recent works of Hofmann, Martell, Toro, Tolsa, and their collaborators. In addition to boundary regularity, typically, some mild topological assumptions (such as one-sided non-tangential accessibility or a weak local John condition) are needed for satisfactory PDE results.

In \cite{DFMprelim}, the authors developed an elliptic theory when $\Gamma \subset \R^n$ is a $d$-AR set, $d<n-1$, and $\Omega := \R^n \setminus \Gamma$. When $d \leq n-2$, 
the boundary is too thin to be seen by a solution of an elliptic PDE in the classical sense 
(for instance by solutions of the Laplacian), and the authors worked with degenerate elliptic operators $-\div A \nabla$, 
such that $w(X)^{-1} A(X)$ satisfy the standard  boundedness and ellipticity conditions on $\Omega$ with the weight $w(X)=\dist(X, \po)^{-(n-d-1)}$. 
This can be reformulated by saying that the underlying measure on $\Omega$ 
is given by $dm(x) := w(x) dx$, and the measure $\mu$ on the boundary $\Gamma$ 
is given by \eqref{defADR}.

In the present article, we give a very large range of choice for $\Omega \subset \R^n$, $m$, and $\mu$, pushing the limits of geometric and measure-theoretic assumptions as well as degeneracy of coefficients of the operators. 
In the rest of the section, we introduce the hypotheses 
on $\Omega$, $m$, and $\mu$, that we shall use for 
most of
the rest of our paper. 

Let us denote
\begin{equation}
\delta(X) := \dist(X,\Gamma)
\end{equation}
for $X \in \Omega$. Since we  allow $\Omega$ to have  boundaries containing pieces of dimension $n-1$ and even higher, we shall need to  deal with connectedness issues that didn't appear in \cite{DFMprelim}.  To this end, we start with standard quantitative connectedness assumptions on $\Omega$,
the Corkscrew and Harnack Chain conditions.

\medskip
\begin{enumerate}[(H1)] 
\item \hcount{H1}
There exists for any $x\in \Gamma= \dr \Omega$ and any $r >0$ 
--- or $r \in (0,\diam \Omega)$ if $\Omega$ is bounded ---
we can find $X \in B(x,r)$ such that 
$B(X,C_1^{-1}r) \subset \Omega$. 
\end{enumerate}

The assumption (\hyperref[H1]{H1}) is widely known as the Corkscrew point condition, and can be seen as quantitative openness. When we say that $Y$ is a Corkscrew point associated to the couple $(y,s)$, we mean that $Y$ is a 
point $X$ given by (\hyperref[H1]{H1}) with  
$x =y$ and $r=s$.

\begin{enumerate}[(H1)] \addtocounter{enumi}{1}
\item \hcount{H2}
There exists a positive integer $C_2 = N+1$
such that if $X,Y \in \Omega$ satisfy $\delta(X)>r$, $\delta(Y) >r$, and $|X-Y| \leq 7C_1r$, 
then we can find $N+1$ points $Z_0:=X,Z_2,\dots,Z_N = Y$ such that for any $i\in \{ 0,\dots,N-1\}$, 
$|Z_i - Z_{i+1}| < \frac12\delta(Z_i)$. 
\end{enumerate}

\medskip

The assumption (\hyperref[H2]{H2}) is a condition of 
quantitative connectedness, and is a slightly weaker way to state the usual
Harnack chain condition. 
We shall discuss (\hyperref[H2]{H2}) more at the end of the section, and in particular prove that 
together with (\hyperref[H1]{H1}), it implies a stronger version of (\hyperref[H2]{H2}), but let us 
first describe the conditions on the two measures $\mu$ (supported on $\Gamma$)
and $m$ (supported on $\Omega \subset \R^n$). 

\medskip

\begin{enumerate}[(H1)] \addtocounter{enumi}{2}
\item \hcount{H3}
The support of $\mu$ is $\Gamma$ and 
$\mu$ is doubling, i.e., there exists $C_3 > 1$ 
such that \[\mu(B(x,2r)) \leq C_3 \mu(B(x,r)) \qquad \text{ for } x\in \Gamma, \, r>0 .\]

\medskip

\item \hcount{H4}
The measure $m$ is mutually absolutely continuous with respect to the Lebesgue measure; 
that is, 
there exists a weight $w \in L^1_{loc}(\Omega)$ such that
\[m(A) = \int_A w(X) \, dX \qquad \text{ for any Borel set } A \subset \Omega\]
and such that $w(X) >0$ for (Lebesgue) almost every $X\in \Omega$.
In addition, $m$ is doubling, i.e. there exists $C_4 \geq 1$ 
such that 
\begin{equation} \label{HH4a}
m(B(X,2r) \cap \Omega) \leq C_4 m(B(X,r) \cap \Omega) 
\qquad \text{ for } X\in \overline{\Omega}, \, r>0 .
\end{equation}
\end{enumerate}

We  included also $X \in \po$, because this is often easier to use, and anyway the version
of (\hyperref[H4]{H4}) with $X\in \overline{\Omega}$ follows easily from the version with $X\in \Omega$.

In some cases, we can get (\hyperref[H4]{H4}) as a consequence of the fact that $m$ is the restriction to $\Omega$
of a doubling measure on $\R^n$. That is, let us say that (H4') holds when 
$m = m'_{\vert \Omega}$ for some absolutely continuous measure $m'$ supported on $\R^n$, 
and which is doubling, i.e.,
\begin{equation} \label{HH4b}
m'(B(X,2r)) \leq C_{4'} \, m'(B(X,r) \cap \Omega) 
\qquad \text{ for } X\in \R^n, \, r>0 
\end{equation}
and some $C_{4'} \geq 1$.

We claim that (H4') and (\hyperref[H1]{H1}) imply (\hyperref[H4]{H4}).
Indeed, take $X\in \overline{\Omega}$ and
separate the two cases $\delta(X) > r/2$ and $\delta(X) \leq r/2$. In the first case,
\[ 
m(B(X,2r) \cap \Omega) \leq m'(B(X,2r)) 
\leq C_{4'}^2 m'(B(X,r/2)) = C_{4'}^2 m(B(X,r)).
\]
In the second case, take $x$ such that $\delta(X) = |X-x|$, and then let $X'$ be a  
Corkscrew point associated to $x$ and $r/2$. 
  Thus 
$B(X',\frac{r}{2C_1}) 
\subset B(X,r) \subset \Omega$ and $B(X',4r) \supset B(X,2r)$, hence 
if $\kappa$ denotes the smaller integer such that $2^\kappa \geq 8C_1$, 
\[
m(B(X,2r) \cap \Omega) \leq m'(B(X',4r)) \leq C_{4'}^{\kappa} \, m'(B(X',\frac{r}{2C_1})) 
\leq C_{4'}^{\kappa} \, m(B(X,r)).
\]
The claim follows.

\smallskip
It is classical (and easy to prove) that the condition (\hyperref[H4]{H4}) is equivalent to the apparently stronger following condition: there exists $d_m>0$ 
and $C>0$, both depending only on $C_4$, such that
\begin{equation} \label{mdoublingG}
m(B(X,\lambda r) \cap \Omega) \leq C \lambda^{d _m} 
m(B(X,r) \cap \Omega) 
\qquad \text{ for } X\in \overline{\Omega}, \, \lambda \geq 1, \, r>0.
\end{equation}

We now state a condition (H5) on the compared growths of $m$ and $\mu$.

\begin{enumerate}[(H1)] \addtocounter{enumi}{4}
\item \hcount{H5}
The quantity $\rho$ defined for
$x\in \Gamma$ and $r>0$ by
\begin{equation} \label{defrho}
\rho(x,r) := \frac{m(B(x,r) \cap \Omega)}{ r\mu(B(x,r))}
\end{equation}
satisfies 
\begin{equation} \label{HH5a}
\frac{\rho(x,r)}{\rho(x,s)} \leq C_5 \left(\frac rs\right)^{1-\epsilon} 
\qquad \text{ for $x\in \Gamma$, $0<s<r$,}
\end{equation}
for some constants $C_5>0$ and $\epsilon := C_5^{-1}$.
\end{enumerate}

\ms
We like $\rho(x,r)$ because it is a dimensionless quantity, but we may 
also write \eqref{HH5a} as 
\begin{equation} \label{HH5b}
\frac{m(B(x,r) \cap \Omega)}{m(B(x,s) \cap \Omega)} 
\leq C_5 \left(\frac rs\right)^{2-\epsilon} \frac{\mu(B(x,r))}{\mu(B(x,s))}
\qquad \text{ for $x\in \Gamma$, $0<s<r$,}
\end{equation}
with a slightly more surprising exponent $2-\epsilon$ due to a different scaling. 

The condition (\hyperref[H5]{H5}) means that the two measures $\mu$ and $m$ need to be intertwined,
in a more precise way that we would get from merely the doubling conditions. That is, we
require that $m(B(x,r) \cap \Omega)$ does not grow much faster than $\mu(B(x,r))$,
with a precise limitation on the exponent.
It is not shocking that something like this is needed. Indeed we 
require for our theory to have a trace theorem (see Section \ref{STrace}), that says that
the functions in the weighted Sobolev space $W^{1,2}(\Omega,m)$ have a trace on $(\Gamma,\mu)$. The function $\rho$ is used in the definition of the 
space of traces, and quantifies the ``deviation'' of the measure of tent 
sets $m(B(x,r) \cap \Omega)$ from the measure on $B(x,r) \cap \Omega$ induced by $\mu$, 
which is $r \mu(B(x,r))$.
It is perhaps more surprising that we only need an upper bound in \eqref{HH5a} and \eqref{HH5b}.

\medskip

Our last condition (H6) requires that the measure $m$ be regular enough, 
and satisfy a weak Poincar\'e inequality.

\medskip
\begin{enumerate}[(H1)] \addtocounter{enumi}{5}
\item \hcount{H6}
If $D \Subset \Omega$ is open and $u_i \in C^\infty(\overline {D})$ is a sequence of functions such that $\int_D |u_i| \, dm \to 0$ and $\int_D |\nabla u_i - v|^2\, dm \to 0$ as $i\to +\infty$, where $v$ is a vector-valued function in $L^2(D,m)$, then $v\equiv 0$.

In addition, there exists $C_6$ such that for any ball $B$ satisfying $2B \subset \Omega$ and any function $u\in C^\infty(\overline{B})$, one has
\begin{equation} \label{defPoincare}
\fint_B |u - u_B| \, dm \leq C_6 r \left( \fint_B |\nabla u|^2 \, dm \right)^\frac12,
\end{equation}  
where $u_B$ stands for $\fint_B u\, dm$ and $r$ is the radius of $B$. 
\end{enumerate}

\medskip

The first part of the condition is technical; it is here to make sure that when we define an appropriate 
notion of gradient for functions that are not smooth, we still have that the convergence in $L^1$
implies a convergence of the gradients. With this property, we shall be able to show the completeness of the weighted Sobolev space we shall work with, which is also essential to be able to get weak solutions. One can find the same condition in \cite{NDbook}.

Using the theory of Haj{\l}asz and Koskela \cite{HK2,HK00}, we will be able to improve the second part of (\hyperref[H6]{H6}) into a Sololev-Poincar\'e inequality. Furthermore, because we can prove a trace theorem, we will also be able to get Poincar\'e inequalities on the sets $B\cap \Omega$, where $B$ is a ball centered on the boundary $\Gamma$. This Poincar\'e inequality at the boundary will be crucial for our proof of the boundary De Giorgi-Nash-Moser estimates in Section \ref{SSolutions}.

The condition (\hyperref[H6]{H6}) will be sometimes replaced by the much stronger condition (\hyperref[H7]{H6'}), defined as

\medskip

\begin{enumerate}[(H1)] \addtocounter{enumi}{5}
\item[(H6')] \hcount{H7}
There exists $C_{6'}$ such that for any ball $B\subset \R^n$ satisfying $2B \subset \Omega$, one has the following condition on the weight $w$:
\begin{equation} \label{HH6'}
\sup_{B} w \leq C_{6'} \inf_B w.
\end{equation}
\end{enumerate}

\medskip

The proof of the fact that (\hyperref[H7]{H6'}) implies (\hyperref[H6]{H6}) goes as follows. 
The second part of (\hyperref[H6]{H6}) is a consequence of the classical Poincar\'e inequality 
(with the Lebesgue measure) and the fact that $w(Z) \approx m(B)/(\diam B)^n$ for all $Z$ 
in a ball $B$ such that $2B \subset \Omega$, which is an easy consequence of (\hyperref[H7]{H6'}).
We turn to the first part of (\hyperref[H6]{H6}). Take $D$, $u_i$, $v$ as in (\hyperref[H6]{H6}). We can cover $D$ by a finite number of balls $B$ satisfying $2B \subset \Omega$, so by (\hyperref[H7]{H6'}), we can find a constant $C_D$ such that $C_D^{-1} \leq w(X) \leq C_D$ for any $X\in D$. We have thus $\int_D |u_i| dx \to 0$, which means that $u_i$ converges to 0 in $L^1(D)$ and hence $\nabla u_i$ converges to 0 in the distributional sense. Since we have also $\int_D |\nabla u_i- v|^2 dx \to 0$, which implies that $\nabla u_i$ converges to $v$ in the distributional sense, we also
have $v \equiv 0$. 

\ms
This completes our list of assumptions concerning the measures $\mu$ and $m$. 
Once we have them, the results of this paper hold for any divergence form
operator $L = \div A \nabla$, where
\begin{equation} \label{2.11}
w(X)^{-1} A(X) \text{ satisfies the standard 
boundedness and ellipticity conditions on $\Omega$,}
\end{equation}
namely, there exists a constant $C_A>0$ such that
\begin{equation} \label{defEllip-intro}
A(X) \xi \cdot \xi \geq C_A^{-1} w(x) |\xi|^2 \quad \text{ for } X\in \Omega \text{ and } \xi \in \R^n
\end{equation}
and
\begin{equation} \label{defBdd-intro}
A(X) \xi \cdot \zeta \leq C_A w(x) |\xi||\zeta| \quad \text{ for } X\in \Omega \text{ and }  \xi,\zeta \in \R^n.
\end{equation}
Of course in practice we may have $L$ and $A$ initially, and then this more or less forces 
the definition of $w$ and $m$.

\ms
We end this section with a further discussion of the Harnack chain condition (\hyperref[H2]{H2}). In the following arguments 
we shall write $Z_i[X,Y]$, $0\leq i \leq N$, when we want to specify the endpoints of the sequence 
given by (\hyperref[H2]{H2}). 
The number $N=C_2-1$ is independent of $X,Y\in \Omega$ as long as $X,Y \in \Omega$ satisfy $\delta(X)>r$, $\delta(Y) >r$, and $|X-Y| \leq 7C_1r$; indeed, even if the sequence is shorter, we can repeat a point $Z_i$ as many times as we want to match the proper length. At last, the ``chain'' in the Harnack chain condition (\hyperref[H2]{H2}) is given by the balls
\begin{equation} \label{HarnackBalls}
B_i =B_i[X,Y] := B(Z_i[X,Y],\delta(Z_i[X,Y])/2).
\end{equation}
From (\hyperref[H2]{H2}), we can easily see that $B_0[X,Y]$ is $B(X,\delta(X)/2)$, 
$B_N[X,Y]$ is $B(Y,\delta(Y)/2)$. 
Moreover, $Z_{i+1}[X,Y] \in B_i[X,Y]$, from which we deduce that 
\begin{equation} \label{propZi}
\delta(Z_{i}[X,Y]) \geq \frac12 \delta(Z_{i-1}[X,Y]) \geq 2^{-N} \delta(X) \geq 2^{-N} r,
\end{equation}
\begin{equation} \label{propZi3}
\delta(Z_{i}[X,Y]) \leq \frac32 \delta(Z_{i-1}[X,Y]) \leq \left(\frac32\right)^i \delta(X) \leq \left(\frac32\right)^N\delta(X),
\end{equation}
and
\begin{equation} \label{propZi2}
|X-Z_{i}[X,Y]| \leq \frac12 \sum_{j=1}^{i-1} \left(\frac32\right)^j \delta(X) \leq \left(\frac32\right)^N\delta(X) \leq 2^N \delta(X),
\end{equation}
that is all the balls $B_i[X,Y]$ from the Harnack chain linking $X$ to $Y$ have equivalent radii, don't get close to the boundary, and are all included in $B(X,2^{N+1}\delta(X))$.  

With the help of (\hyperref[H1]{H1}), the condition (\hyperref[H2]{H2}) self improves, as shown by the result below.

\begin{proposition} \label{propHarnack}

Let $\Omega$ 
satisfy (\hyperref[H1]{H1}) and (\hyperref[H2]{H2}). 
There exists $C:=C(C_1,C_2)>0$ such that if $X,Y \in \Omega$ satisfy 
$\min\{\delta(X),\delta(Y)\} > r$ and $|X-Y| \leq \Lambda r$
(for some choice of $r > 0$ and $\Lambda \geq 1$), 
then we can find $N_\Lambda:= \lceil C\ln(1+\Lambda)\rceil$ points $Z_0:=X,Z_1,\dots,Z_{N_\Lambda} = Y$ such that for any $i\in \{ 0,\dots,N_\Lambda-1\}$, 
\begin{enumerate}[(i)]
\item $|Z_i - Z_{i+1}| \leq \frac12\delta(Z_i)$, 
\item $\delta(Z_i) \geq 2^{-N}r$,
\item $Z_i \in B(X,C_12^{N+4}\Lambda r)$,
\end{enumerate}
where $N:=C_2-1$ 
comes from (\hyperref[H2]{H2}) 
and $C_1$ comes from (\hyperref[H1]{H1}).
\end{proposition}

\bp Let $X,Y \in \Omega$ satisfy $\min\{\delta(X),\delta(Y)\} > r$ and $|X-Y| \leq \Lambda r$. 
Observe first that if 
$\Lambda r \leq \delta(X)$ or $\Lambda r \leq \delta(Y)$ there is no need for (\hyperref[H2]{H2}).
Indeed, the segment $[X,Y]$ is included in $\Omega$, and one can construct the chain recursively as: $Z_0 = X$, $Z_{i+1}$ is either (if it exists) the only point further from $X$ than $Z_i$ at the intersection of $[X,Y]$ and the sphere centered at $Z_i$ and radius $\delta(Z_i)$, or simply $Y$ if 
this point 
doesn't exist. 

\medskip

In the remaining 
case where $\delta(X),\delta(Y) \leq \Lambda r $ (which forces $\Lambda \geq 1$), the idea is to use the condition (\hyperref[H1]{H1}) to find enough points between $X$ and $Y$ to be able to split the distance $|X-Y|$ into small jumps where we can use (\hyperref[H2]{H2}).

Let $x,y\in \Gamma$ be such that $|X-x| = \delta(X)$ and $|Y-y| = \delta(Y)$, so that $X,Y$ are Corkscrew points associated to respectively $(x,C_1\delta(X))$ and $(x,C_1\delta(Y))$. 
We define the points $X_j,Y_j$ as follows: $X_0=X$, $Y_0=Y$ and, if 
$j\geq 1$, $X_j$ is a Corkscrew point associated to $(x,C_12^j\delta(X))$ and $Y_j$ is a Corkscrew point associated to $(y,C_12^j\delta(Y))$. Then we set $j_x$ and $j_y$ as the smallest values of $j\in \bN$ such that $2^j\delta(X) \geq \Lambda r$ and $2^j\delta(Y) \geq \Lambda r$ respectively. It is easy to check that by construction, 
\begin{equation} \label{Lambdalength}
j_x,j_y \leq 1+ \ln_2(\Lambda) \leq C \ln(1+\Lambda),
\end{equation}
where $C$ is a universal constant. We set $M_\Lambda = j_x+j_y+1 \leq C\ln(1+\Lambda)$, we and construct a first sequence of points $(Z^j)_{0\leq j \leq M_\Lambda}$ as $Z^j = X_j$ when $j\leq j_x$ and $Z^j = Y_{M_\Lambda-j}$ when $j>j_x$. We want to verify that two successive elements of $(Z^j)$ satisfy the assumptions for the use of (\hyperref[H2]{H2}). Indeed, $X_j$ and $X_{j+1}$ are such that $\min\{\delta(X_j),\delta(X_{j+1})) > 2^j \delta(X)$ and 
\[|X_j - X_{j+1}| \leq |X_j - x| + |X_{j+1} - x| \leq C_12^{j+2}\delta(X);\]
the same kind of estimates holds between $Y_j$ and $Y_{j+1}$; 
the points $X_{j_x}$ and $Y_{j_y}$ are such that 
\[
\min\{\delta(X_{j_x}),\delta(Y_{j_y})\} 
\geq \Lambda r\]
and, since $\delta(X),\delta(Y) \leq \Lambda r$,
\[\begin{split}
|X_{j_x} - Y_{j_y}|  
& \leq |X_{j_x} - x| + |x-X| + |X-Y| + |Y-y| + |y - Y_{j_y}| \\
& \leq 2C_1\Lambda r + \delta(X) + \Lambda r + \delta(Y) + 2C_1\Lambda r \\
& \leq 7C_1\Lambda r.
\end{split}\]
The fact that two consecutive points of $(Z^j)_j$ satisfies the assumption of the Harnack chain condition follows. Now, $N$ stands for $C_2-1$, and the sequence $(\mathcal Z_i)_{0 \leq i \leq N.M_\Lambda}$ is built such that $\mathcal Z_i = Z_i[Z^j,Z^{j+1}]$ if $jN \leq i \leq (j+1)N$. 
The conclusion (i) is given by the definition/construction of the points $Z_i[X',Y']$; 
the conclusion (ii) is immediate from \eqref{propZi} since all $Z^j$ are such that 
$\delta(Z_i) \geq \min\{\delta(X),\delta(Y)\} \geq r$. As for (iii), we estimate the distance between $X$ and the $Z^j$ rather brutally and we let the reader check that $|X-Z^j| \leq 14C_1\Lambda r$ for any $j\in \bN$, which, combined with \eqref{propZi2}, gives (iii).
\ep

\ms
In the rest of the paper, the notation $u\lesssim v$ means $u \leq C v$,
where $C>0$ is a constant that depends on parameters which will be 
either obvious from the
context or recalled. The expression $u \approx v$ is used when 
$u \lesssim v$ and $u \gtrsim v$.

\section{Some examples where our assumptions hold}
\label{SExa}

The assumptions of the previous section may still look complicated to the reader, 
so let us mention some situations where they are satisfied, and hence we can define an
elliptic measure $\omega_L$ with the properties described in the introduction.

\subsection{Classical elliptic operators}
We start with the classical elliptic operators $L = \div A \nabla$, 
where $A(X)$ satisfies the standard boundedness and ellipticity conditions 
\eqref{defBdd-intro} and \eqref{defEllip-intro} on $\Omega$. 
In view of \eqref{2.11}, 
$w=1$. We also require the one-sided
NTA conditions (\hyperref[H1]{H1}) and (\hyperref[H2]{H2}), which happen to hold automatically when
$\Gamma = \d\Omega$ is Ahlfors regular of dimension $d < n-1$, but not in general. 
Then our additional assumptions are the existence of a doubling measure $\mu$
on $\Gamma$ (as in (\hyperref[H3]{H3})) that satisfies (\hyperref[H5]{H5}); the other conditions,
including (\hyperref[H7]{H6'}) are trivially satisfied. In particular if $\Gamma$ is Ahlfors regular of dimension
$d \in (n-2,n)$, it is easy to check that $\mu = \H^d_{\vert \Gamma}$ satisfies (\hyperref[H5]{H5}).

\subsection{Ahlfors regular sets}\label{ARex}
Our next example is when $\Gamma = \d\Omega$ is an Ahlfors regular set
of dimension $d$.
When $n-1 \leq d < n$, we also require one-sided NTA conditions (\hyperref[H1]{H1}) and (\hyperref[H2]{H2}).
The simplest option is to take $w(X) = \dist(X,\Gamma)^{-\gamma}$ for some
$\gamma \in (n-d-2,n-d)$. Then $w$ is locally integrable, by \eqref{defADR}, 
because $\gamma < n-d$, and by a simple estimate on the measure of the 
$\varepsilon$-neighborhoods of $\Gamma$. The same estimates yield that 
$m(B(X,r)) \approx r^{n-\gamma}$ when $\delta(X) < 4r$ (with a lower bound that uses (\hyperref[H1]{H1})) and $m(B(X,r)) \approx r^n \delta(X)^{-\gamma}$ when $\delta(X) >2r$. This proves that $m$ is doubling; then \eqref{HH5b} holds as soon as $\gamma > n-d-2$.
The other conditions, including (\hyperref[H7]{H6'}), are easy to check, and so our results apply
to operators $L = \div A \nabla$, where $\dist(X,\Gamma)^{\gamma} A(X)$ satisfies 
the standard boundedness and ellipticity conditions 
\eqref{defBdd2} and \eqref{defEllip2} on $\Omega$.

In this Ahlfors regular setting, we can also deal with more general weights $w$,
that would also have mild local singularities in the middle of $\Omega$, and then
the corresponding classes of degenerate elliptic operators $L = \div A \nabla$, 
where $w(X)^{-1} A(X)$ is bounded elliptic, but then we have to check 
(\hyperref[H6]{H6}) too, in addition to (\hyperref[H5]{H5}).

\subsection{Caffarelli and Sylvestre fractional operators}  \label{S3CS}
A special case of the weight $w(X) = \dist(X,\Gamma)^{-\gamma}$ for an Ahlfors 
regular boundary was considered by L. Caffarelli and L. Sylvestre \cite{CS}, 
although in a very different context.
Take $n=d+1$, $\Omega = \R^{d+1}_+$ (the upper half space), and $\Gamma = \R^{d}$.
Choose $\mu$ to be the Lebesgue measure on $\R^{d}$, and for $m$ take the 
weight $w(X) = \dist(X,\Gamma)^{-\gamma} = t^{-\gamma}$,
where we write $(x,t)$ the coordinates of $X$ in $\R^{d+1}$. As before, we restrict to
$\gamma \in (-1,1)$. 

Caffarelli and Sylvestre considered the fractional operator
$T = (-\Delta)^{s}$ on $\R^d$, with $s = \frac{1+\gamma}{2} \in (0,1)$, and 
proved that for $f$ defined on $\R^d$, in the appropriate space,
$Tf$ can also be written as $Tf(x) = - C \lim_{t \to 0+} t^{-\gamma} \frac{\d u}{\d t}$,
where $u$ is the solution of $Lu=0$, with $L = \div t^{-\gamma} \nabla$, whose trace 
on $\R^d$ is $f$. This point of view turned out to be a very useful way to deal with unpleasant
aspects of the non-local character of $T$. 

We can generalize some of this to the context of Ahlfors regular sets, as above, with 
$L = -\div \dist(X,\Gamma)^{-\gamma} \nabla$ (or a similar operator). 
When $f$ lies in our Sobolev space $H = H(\Gamma)$, the results of the present paper 
allow us to solve the Dirichlet problem for $f$, i.e., find a function 
$u \in W^{1,2}(\Omega, wdX) =  W^{1,2}(\Omega,\dist(X,\Gamma)^{-\gamma}dX)$ 
such that $Lu=0$ and ${\rm Tr}(u) = f$. Then we can also define an operator $T$, 
that generalizes the fractional operator of \cite{CS}, 
by saying that $Tf$ is a distribution on $\Gamma$ (or a continuous 
linear operator on $H$), such that
\begin{equation} \label{Syl}
\langle Tf, \varphi \rangle 
= \int_\Omega \dist(X,\Gamma)^{-\gamma} \, \nabla u(X) \cdot \nabla \varphi F(X) dX
\end{equation}
when $F$ is a function of $W = W^{1,2}(\Omega, dm)$ such that ${\rm Tr}(F) = \varphi$.
For example, we could take for $F$ the extension of $f$ given in Section \ref{Sext},
but taking another extension $F'$ should give the same result, because $F-F'$
lies in the space $W_0$ of functions of $W$ with a vanishing trace 
(see Definition \ref{defW0}), and because $u$ is a weak solution (see
the definition \eqref{defsol} and the last item of Lemma \ref{rdefsol}).
Now $Tf$ can be seen as a weak limit of normalized derivatives $w(X)\frac{\d u}{\d \nu}(X)$
in the normal direction, as above: we can integrate by parts on a smaller 
domain $\Omega_\varepsilon$ and try to take a limit. 
Ultimately, it would be nice to have a more precise and constructive 
definition of $T$, with estimates in a better space than $H^{-1/2}$ (the dual of $H$);
however this will require a better analysis of $L$, and quite probably stronger assumptions on 
$\Gamma$. Yet the similarity with the situation in \cite{CS} is intriguing.
 
 \subsection{Sawtooth domains} \label{S3s}
We now turn to Sawtooth domains. 
Let us first describe the simpler case of Lipschitz graphs, and then later comment
on Ahlfors regular boundaries.
Let us assume that $0 < d < n-1$ (the case when $d=n-1$ is simpler, and well known)
and that $\Gamma$ is the graph of some Lipschitz function $A : \R^{d} \to \R^{n-d}$,
where we identify $\R^d$ and $\R^{n-d}$ to the obvious coordinate subspaces of $\R^n$.
Also let $E \subset \R^{d}$ be a given subset of $\R^d$, for which we want to hide
$\wt E = \big\{ (x,A(x) \, ; \, x \in E \big\}$. We assume that $F = \R^d \sm E$ is 
not empty (otherwise, there is no point in the construction), and we set
$\wt F = \big\{ (x,A(x) \, ; \, x \in F \big\} = \Gamma \sm \wt E$.
Let us define a sawtooth domain
$\Omega_s \subset \Omega$ such that $\wt F \subset \d\Omega_s$ and
$\wt E \subset \R^n \sm \overline{\Omega_s}$.
We choose $M$ larger than the Lipschitz norm of $A$, and set 
\begin{equation} \label{st1}
\Omega_s = \big\{ (x,t) \in \R^n \times \R^{n-d} \, ; \, |t-A(x)| > M \dist(x,F) \big\}.
\end{equation}
Thus we are removing from $\Omega$ some sort of a conical tube
around $\wt E \subset \Gamma$. In co-dimension $1$, we would proceed
similarly, but restrict to the part of $\wt\Omega$ that lies above $\Gamma$, for instance.
We can forget about this case because it is very classical anyway.

It is clear from the definition that $\Omega_s$ is an open set that does not meet
$\Gamma$, and its boundary consists in the closure of $\wt F$, plus
the conical piece 
\begin{equation} \label{st2}
Z = \big\{ (x,t) \in \R^n \times \R^{n-d} \, ; \, |t-A(x)| = M \dist(x,F) > 0 \big\},
\end{equation}
which nicely surrounds $\wt E$.

The verification that $\Omega_s$ contains Corkscrew balls and Harnack chains
(as in (\hyperref[H1]{H1}) and (\hyperref[H2]{H2})) is rather easy, because we can always escape
in a direction opposite from $\Gamma$ to find extra room; we skip the verification.
In this type of situation, we probably want to be able to use the same operators $L$
as we had on $\Omega$, so let us consider the restriction to $\Omega_s$
of our earlier weight $w(X) = \dist(X,\Gamma)^{-\gamma}$, with some $\gamma \in (n-d-2,n-d)$.
As usual, this defines a class of matrices $A$.

We have to construct a new measure $\mu_s$ on $\Gamma_s := \d\Omega_s$,
and we choose
\begin{equation} \label{st3}
\mu_s = \mu_{\vert \wt F} + \dist(X,\Gamma)^{d+1-n} \,\H^{n-1}_{\vert Z},
\end{equation}
where $\mu$ is an Ahlfors regular measure on $\Gamma$ that we like, such as 
$\H^d_{\vert \Gamma}$ or the image of $\H^d_{\vert \R^d}$
by the mapping $x \to (x,A(x))$. There may be locally more subtle choices, but
we do not worry too much here because for our purpose $\mu_s$ only needs to be known within
bounded multiplicative errors. We mostly care that 
if $\pi$ denotes the orthogonal projection on $\R^d$ and $\pi_\ast\mu_s$
is the push-forward image of $\mu$ by $\pi$, then 
\begin{equation} \label{st4}
C^{-1} \H^d_{\vert \R^d} \leq  \pi_\ast\mu_s \leq C \H^d_{\vert \R^d}.
\end{equation}
This last is easy to prove, because when $x\in E$ is such that 
$d = \dist(x, F) > 0$, and $r > 0$ is much smaller than $d$,
the surface measure of $Z \cap \pi^{-1}(B(x,r))$ is comparable to
$r^{d} (Md)^{n-d-1}$, so $\mu_s(\pi^{-1}(Z \cap B(x,r))) \approx r^{d}$
(the dependence on $M$ does not interest us).

As in the previous examples, $w$ is essentially constant on the balls $B$ such
that $2B \subset \Omega_s$, so (\hyperref[H7]{H6'}) and (\hyperref[H6]{H6}) hold; the verification
of the doubling property (\hyperref[H4]{H4}) for $m$ is the same as for the initial open 
set $\Omega$, and we could also use (H4') directly; so we only need
to check the doubling property (\hyperref[H3]{H3}) for $\mu_s$ and the intertwining
growth condition (\hyperref[H5]{H5}). 

To this effect, let us first show that 
\begin{equation} \label{st5}
C^{-1} r^d \leq \mu_s(B(x,r)) \leq C r^d 
\qquad \text{ for $x\in \wt F$ and $r > 0$}.
\end{equation}
The second inequality, which incidentally holds when $x\in Z$ too, follows from
\eqref{st4} because $B(x,r) \subset \pi^{-1}(\R^d \cap B(\pi(x),r))$.
For the first inequality, set $r_1 = (1+M)^{-1}r$. 
We claim that $\Gamma_s \cap \pi^{-1}(\pi(B(x,r_1))) \subset B(x,r)$; once we prove
the claim, \eqref{st5} will follow from \eqref{st4} because
$\H^d(\pi(B(x,r_1))) \geq C^{-1}r_1^d$. 
Now let $z\in \Gamma_s \cap \pi^{-1}(\pi(B(x,r_1)))$ be given. 
If $z\in \Gamma$, then $|z-x| < (1+M)r_1 = r$, so $z\in B(x,r)$.
Otherwise $z\in Z$, and let $y\in \Gamma$ be such that $\pi(y) = \pi(z)$. Notice
that $|\pi(y)-\pi(x)| = |\pi(z)-\pi(x)| < r_1$, so 
$\dist(\pi(y), F) \leq |\pi(y)-\pi(x)| < r_1$ and now 
\eqref{st2} says that $|z-y| \leq M \dist(\pi(y), F) < Mr_2$. Again $z\in B(x,r)$
and the claim follows. This proves \eqref{st5}.
Next we check that 
\begin{equation} \label{st6}
C^{-1} r^d \leq \mu_s(B(z,r)) \leq C r^d 
\qquad \text{ for $z\in Z$ and $r > (2+2M) \dist(\pi(z), F)$}.
\end{equation}
Recall that the second inequality always holds. 
For the first one, set $r_2 = (2+2M)^{-1}r >\dist(\pi(z), F)$, 
choose $p \in F$ such that $|p-\pi(z)| < r_2$, and then let $x\in \Gamma$
be such that $\pi(x)=p$. Observe that $x\in \wt F$.

Also let $y \in \Gamma$ be such that $\pi(y)=\pi(z)$. Then 
$|y-x| \leq (1+M)|\pi(y)-\pi(w)| = (1+M)|\pi(z)-p| < (1+M) r_2$,
and now \eqref{st2} yields $|z-y| \leq M \dist(\pi(y), F) 
= M \dist(\pi(z), F) < M r_2$; altogehter 
$|z-x| \leq |z-y| + |y-x| \leq (1+2M) r_2$, so $B(w,r_2) \subset B(z,r)$, and 
\eqref{st6} follows from \eqref{st5} because $x\in \wt F$.

We are left with the case when $z\in Z$ and $r < (2+2M) \dist(\pi(z), F)$;
we claim that then
\begin{equation} \label{st7}
C^{-1} r^{n-1} \dist(\pi(z), F)^{d+1-n} \leq \mu_s(B(z,r)) 
\leq C  r^{n-1} \dist(\pi(z), F)^{d+1-n}.
\end{equation}
Set $d(z)=\dist(\pi(z), F)$. When $10^{-1} d(z) < r \leq (2+2M) d(z)$,
$r^{n-1} d(z)^{d-n-1}$ is roughly the same as the $r^d$ that we had before, 
so there is some continuity in our estimate. Also, the upper bound stays true as before,
and the lower bound will follow as soon as we prove it for $r = 10^{-1} d(z)$.
Finally, the remaining case when $r \leq 10^{-1} d(z)$ is easy, because in 
$B(z,r)$, $\wt \Gamma$ coincides with $Z$ and looks like the product of $\R^d$
(or $\Gamma$) with an $(n-d-1)$-sphere, with $d(z)^{d+1-n}$ times the surface
measure.

The doubling condition (\hyperref[H3]{H3}) for $\wt \mu$ follows at once from
the estimates above, so let us just check the intertwining growth condition (\hyperref[H5]{H5}).
We start with the case when $x\in \wt F$. Then 
$m(B(x,r)) \approx r^{n-\gamma}$ as in the standard Ahlfors regular case, 
$\mu_s(B(x,r)) \approx r^{d}$ by \eqref{st5}, so \eqref{HH5a} also holds as in the 
Ahlfors regular case. 

Next assume that $x\in Z$. 
When $r < 10^{-1}d(x)$, $m(B(x,r)) \approx d(x)^{-\gamma} r^n$ because
$w(X) = \dist(X,\Gamma)^{-\gamma} \approx d(x)^{-\gamma}$ in $B(x,r)$, and
$\mu_s(B(x,r)) \approx r^{n-1} d(x)^{d+1-n}$ by \eqref{st7}. 
Thus $\rho(x,r) = m(B(x,r)) r^{-1} \mu_s(B(x,r))^{-1} \approx d(x)^{-\gamma - d - 1 + n}$.

When $10^{-1}d(x) \leq r \leq (2+2M) d(x)$, none of these number changes too much,
so $\rho(x,r) \approx d(x)^{-\gamma - d - 1 + n}$ as well. 

Finally, when $r > (2+2M) d(x)$, $m(B(x,r)) \approx r^{n-\gamma}$ as in the standard 
Ahlfors regular case, $\mu_s(B(x,r)) \approx r^d$ by \eqref{st6}, and so 
$\rho(x,r) \approx r^{-\gamma - d - 1 + n}$. Now \eqref{HH5a} and (\hyperref[H3]{H3})
follow because $|n-d-1-\gamma| < 1$.

Thus, in the case of Lipschitz graphs, the sawtooth regions $\Omega_s$ that we constructed,
together with the measures $m$ on $\Omega_s$ and $\mu_s$ on $\d\Omega_s$,
satisfy the requirements of Section \ref{S2}.

\medskip
There is a more general construction of sawtooth regions, 
adapted to the case when $\Omega$ is a one-sided NTA domain (i.e., (\hyperref[H1]{H1}) 
and (\hyperref[H2]{H2}) hold) with an Ahlfors regular boundary $\Gamma$ --
see for instance \cite{HM1} for a first occurrence. It was used in quite a few papers later, at least when $\Gamma$ is of co-dimension $1$.

Before we start with a very rough description of a sawtooth construction, it is convenient 
to take a collection of dyadic pseudocubes $Q$, like the one that will be described in
Proposition~\ref{ChristQ} below. This is a collection of sets $Q \subset \Gamma$,
$Q \in \D$, that have roughly the same covering and inclusion properties as
the usual dyadic cubes in $\R^d$. 

Then, to each pseudocube $Q \in \D$, we can also associate a 
Whitney region $\cW(Q) \subset \Omega$, such that for some large $C \geq 1$, 
$C^{-1} \diam(Q) \leq \dist(X,\Gamma)$ and $\dist(X,Q) \leq C \diam(Q)$ for $X \in \cW(Q)$.
We make sure to take the $\cW(Q)$ to be sufficiently large, so that they cover $\Omega$.
And also, for the construction below to work, one should choose them carefully,
so that they have sufficiently simple boundaries (for instance, by requiring that each
$\cW(Q)$ is composed of a finite union of cubes in a sufficiently sparse collection), and  
that they are 
well connected with each other. 

Then, we are given a one-sided NTA domain (i.e., with (\hyperref[H1]{H1}) and (\hyperref[H2]{H2})),
with an Ahlfors regular boundary of any dimension $d < n$. We are also given a
stopping time region, some times also called regime,
where one starts from a top cube $Q_0$, and one keeps a collection $\mathcal S$
of subcubes of $Q_0$, with some coherence conditions. For instance, if
$R \in \mathcal S$, then all the cubes $S \in \D$ such that $R \subset S \subset Q_0$
lie in $\mathcal S$ too. The set we want to keep access too is the set $F$ of points
of $Q_0$ such that all the cubes $Q$ such that $x \in Q \subset Q_0$ lie in
$\mathcal S$. And the corresponding sawtooth region is the union of all the Whitney
sets $\cW(Q)$, $Q \in \mathcal S$. We claim that if the sets
$\cW(Q)$, $Q \in \D$, are carefully chosen, then the assumptions of the current paper
are satisfied. But we do not check this
here, because we intend to do this in a next paper,
where this will be useful for a comparison of elliptic measures. 
The general idea of the verification is the same as for Lipschitz graphs, but 
the details are a little more complicated. 

\subsection{Balls minus an Ahlfors regular set of low dimension}

Let $\Gamma \subset \R^n$ be an Ahlfors regular set of dimension $d<n-1$, that is a set that satisfies \eqref{defADR}. Consider any ball $B \subset \R^n$. We want to show that the theory developed in this article applies to $B \setminus \Gamma$. Of course, by taking the particular case where $\Gamma = \R$ and $B \subset \R^3$ centered on $\R$, we see that balls deprived of one diameter - as claimed in the abstract - are included in our theory.

First of all, by translation and scale invariance of the problem, we can assume that $B$ is the ball centered at 0 and of radius 1. In this subsection, we plan to give measures on $\Omega := B \setminus \Gamma$ and $\dr \Omega = \dr B \cup (\Gamma \cap B)$, and establish (\hyperref[H1]{H1}), (\hyperref[H2]{H2}), (\hyperref[H3]{H3}), (H4'), (\hyperref[H5]{H5}), and (\hyperref[H7]{H6'}). The weights on $\Omega$ that we choose are the same as the ones in Subsection \ref{ARex}, that is
$w(X) = \dist(X,\Gamma)^{-\gamma}$ for some $\gamma \in (n-d-2,n-d)$,
 which is natural because we expect the present domains to appear when we want to localize problems and properties from the situation given in Subsection \ref{ARex}.
The boundary $\dr \Omega$ is divided into $\Gamma_1:= \Gamma \cap B$ and $\Gamma_2 := \dr B$. And then in the spirit of what we did in Subsection \ref{S3s}, the measure $\mu$ on $\dr\Omega$ is
\begin{equation}
\mu = \mu_1+ \mu_2 := \sigma_{|\Gamma_1} + \dist(X,\Gamma)^{d+1-n} \mathcal H_{|\Gamma_2}^{n-1},
\end{equation}
where $\sigma$ is an Ahlfors regular measure that satisfies in \eqref{defADR} and $\mathcal H_{|\Gamma_2}^{n-1}$ is the surface measure on the sphere $\Gamma_1$.

\medskip

The conditions (H4') and (\hyperref[H7]{H6'}) are the same as in Subsection \ref{ARex}. For (\hyperref[H3]{H3}) and (\hyperref[H5]{H5}), we want to check that the transition between two parts of the boundary with different dimensions goes smoothly. We only need to prove our hypotheses for $r\leq 2$, which is the diameter of our domain, and so we assume $r\leq 2$ for the rest of the subsection. The proofs of (\hyperref[H3]{H3}) and (\hyperref[H5]{H5}) works a bit like when we have sawtooth domains, in particular we have the same type of estimates. First, for $x\in \dr\Omega$ and $\dist(x,\Gamma) \geq 2r >0$, we prove that
\begin{equation} \label{ss6c}
C^{-1} r^{n-1}\dist(x,\Gamma)^{d+1-n} \leq \mu(B(x,r)) \leq C r^{n-1} \dist(x,\Gamma)^{d+1-n}.
\end{equation}
Of course, our assumption on $x$ and $r$ forces $\mu_1(B(x,r)) = \sigma(B(x,r)) = 0$. Moreover, the weight $\dist(X,\Gamma)^{d+1-n}$ used to define $\mu_2$ is essentially constant on $B(x,r)$, and added to the facts that $x$ has to belong to $\Gamma_2$ and $\mathcal H_{|\Gamma_2}^{n-1}$ is an Ahlfors regular measure of dimension $n-1$, we deduce that $\mu_2(B(x,r)) \approx r^{n-1} \dist(x,\Gamma)^{d+1-n}$. The claim \eqref{ss6c} follows.

Next, we want
\begin{equation} \label{ss6d}
C^{-1} r^{d} \leq \mu(B(x,r))  \qquad \text{ for $x\in \dr \Omega$ and $\dist(x,\Gamma) < 2r$}.
\end{equation}
We need to distinguish two cases. Either $\dist(x,\Gamma_2) \geq r/2$, and then $\mu(B(x,r)) \geq \sigma(B(x,r/2)) \geq C^{-1} r^d$. Or $\dist(x,\Gamma_2) < r/2$ and we can find $x' \in \Gamma_2$ such that $|x-x'| < r/2$, from which we deduce \[\mu(B(x,r)) \geq \mu_2(B(x',r/2)) \geq C^{-1} r^{d+1-n}\mathcal H_{|\Gamma_2}^{n-1}(B(x',r/2)) \geq C^{-1} r^{d}\]
because $H_{|\Gamma_2}^{n-1}$ is a $(n-1)$-Ahlfors regular measure.  

The last inequality that we need is
\begin{equation} \label{ss6e}
\mu(B(x,r)) \leq C r^{d} \qquad \text{ for $x\in \dr \Omega$ and $\dist(x,\Gamma) < 2r$}.
\end{equation}
We take a point $x' \in \Gamma$ so that $|x-x'| < 2r$. The above claim will be a consequence of the fact that $\mu(B(x',3r)) \lesssim r^d$. The inequality $\mu_1(B(x',3r)) \lesssim r^d$ is a free consequence of the fact that $\Gamma$ is $d$-Ahlfors regular. As for $\mu_2$, we divide $B(x',3r)$ into the strips
\[S_k := \{y\in B(x',3r), \, 2^{2-k}r \geq \dist(y,\Gamma) \geq 2^{1-k}r\}.\]
We use \eqref{defADR} to cover $\Gamma \cap B(x',10r)$ with less than $C2^{kd}$ balls $\{B_j\}$ of radius $2^{-k}r$ and $S_k$ is contains in the union of the $5B_j$.
We deduce that $\mathcal H^{n-1}_{|\Gamma_2}(S_k) \leq C2^{kd} (2^{-k}r)^{n-1}$ and then
\[\begin{split}
\mu_2(B(x',3r)) & \leq \sum_{k\in \bN} \mu_2(S_k) \leq C(2^{-k}r)^{d+1-n} \mathcal H^{n-1}_{|\Gamma_2}(S_k) \leq Cr^d
\end{split}\]
as desired.
We let the reader check that the estimates \eqref{ss6c}, \eqref{ss6d}, and \eqref{ss6e} easily imply (\hyperref[H3]{H3}) and (\hyperref[H5]{H5}).

\medskip

The last conditions are (\hyperref[H1]{H1}) and (\hyperref[H2]{H2}). The proof of Lemma 2.1 in \cite{DFMprelim} - which treats the case $\R^n \setminus \Gamma$ - can actually be repeated in our case without any changes. We obtain that two points can be linked by 3 consecutive tubes that don't intersect $\Gamma$ and stay in $B$. The assumption (\hyperref[H2]{H2}) follows by taking an appropriate sequence of points in those tubes.

(\hyperref[H1]{H1}) is just a bit more complicated, because it requires to distinguish cases, and still relies on what we did for $\R^n \setminus \Gamma$ in \cite{DFMprelim}. Take $x\in \dr \Omega$ and $0<r\leq 2$. We can find $X'$ such that $B(X',r/2) \subset B \cap B(x,r)$. Indeed, if $x\in \Gamma_2 = \dr [B(0,1)]$, take $X':= (1-r/2) x$, and if $x\in \Gamma_1$, observe that we are further from the sphere $\dr B$ so it is easier to find such $X'$ (for instance $X' = \frac{x}{|x|} \min\{(1-r/2),|x|\}$).  We look then at $\Gamma$ inside $B(X',r/2) \subset B$. If $B(X',r/4) \cap \Gamma = \emptyset$, then $X'$ is our point for (\hyperref[H1]{H1}). Otherwise, we can find $y\in \Gamma \cap B(X',r/4)$, and Lemma 11.6 in \cite{DFMprelim} gives us $X$ such that $B(X,C^{-1}r) \subset B(y,r/4) \setminus \Gamma \subset B \setminus \Gamma = \Omega$.

\medskip

Maybe the reader will be interested to observe that we could replace the ball by other sets, like cubes. We claim that we can replace $B$ by any 1-sided NTA domain $D$ - that is $D$ to satisfies (\hyperref[H1]{H1}) and (\hyperref[H2]{H2}) - such that $\dr D$ is a $(n-1)$-Ahlfors regular set, and we let the reader verify that all the computations above can be adapted.

\subsection{Nearly $t$-independent $A_2$-weights}  \label{SS3.5}

The $t$-independent elliptic operators have a special status among divergence form operators,
in particular, because some control of behavior of the coefficients in the direction transversal to the boundary is necessary for absolute continuity of harmonic measure with respect to the Lebesgue measure -- see \cite{CFK}. 

We start with the simplest case in co-dimension $1$.
Let $\omega : \R^d  \to \R_+$ be any $A_2$-weight on $\R^d$ (see \cite{Jou, Rub} for details)
and use it to define a weight $w$ on $\R^{d+1}_+ = \R^d \times \R$ by 
$w(x,t) = \omega(x)$. Then set $dm(x,t) = w(x,t) dx dt$ as usual. 
This is a doubling measure on $\Omega$ because $\omega(x)dx$ is doubling 
on $\R^d$ for any $A_\infty$ weight $\omega$.

On $\Gamma = \R^d$, we simply put the measure $d\mu = \omega(x)dx$.
With the mere assumption that $\omega$ is doubling, we immediately
get the one-sided NTA conditions (\hyperref[H1]{H1}) and (\hyperref[H2]{H2}),
the doubling conditions (\hyperref[H3]{H3}) and (\hyperref[H4]{H4}), and even the
intertwining condition (\hyperref[H5]{H5}), because $\mu$ is doubling and
$m(B(x,r)) \approx r \mu(B(x,r) \cap \R^d)$, so $\rho \approx 1$ in \eqref{defrho}
and (\hyperref[H5]{H5}). 

We are left with the last condition (\hyperref[H6]{H6}), and 
this is where we really use our assumption that $\omega \in A_2(\R^d)$.
It is easily checked that $w \in A_2(\R^n)$ too, and now we can use Theorem 1.2 
in \cite{FKS} to deduce the density and the Poincar\'e results of (\hyperref[H6]{H6}). 
The reader may be worried about a minor point. 
In the case of $A_2$ weights, the authors of \cite{FKS} first claim
the first part of (\hyperref[H6]{H6}) only when $\int_D |u_i|^2 dm \to 0$. 
This is then easy; one applies Cauchy-Schwarz and uses the fact that $1/w$ is locally 
integrable. But then the slightly stronger version stated in (\hyperref[H6]{H6}), 
where we only assume the $L^1$ convergence, follows:
by the discussion above \eqref{NablaLip}, our functions $u_i$ actually
lie in $W \subset L^1_{loc}$ and their gradient of Definition \ref{defW} is also
their distribution gradient; then the uniqueness of (\hyperref[H6]{H6}) comes as in the case of 
(\hyperref[H7]{H6'}) (see below (\hyperref[H6]{H6})).

\smallskip
The previous verification of (\hyperref[H1]{H1})-(\hyperref[H2]{H2}) can easily be extended to the 
case when $\Gamma \subset \R^{d+1}$ 
is the graph of a Lipschitz function $A : \R^d \to \R$, and $\omega \in  A_2(\R^d)$.
We take $w(x,t) = \omega(x)$ as before, and for $\mu$ the image of $\omega(x)dx$
by the mapping $x \to (x,A(x)$ from $\R^d$ to $\Gamma$. Not much changes, because
our conditions are essentially invariant under bilipschitz mappings; we only need to check
that the Poincar\'e estimate \eqref{defPoincare} away from $\Gamma$ stays true,
with merely $B$ in its right hand side, but this is all right. 

Finally, we can further generalize all this to higher co-dimensions, except that we replace 
$t$-invariance by a more reasonable homogeneity. Let us now take integers $d < n-1$,
$\Gamma = \R^d \subset \R^n$ and $\Omega = \R^n \sm \Gamma$ for simplicity, 
and as before $\omega \in A_2(\R^d)$. We keep $d\mu(x) = \omega(x) dx$ on $\Gamma$,
but now use $w(x,t) = |t|^{-\gamma} \omega(x)$ on $\Omega$, with as usual
$\gamma \in (n-d-2, n-d)$. Again the results of this paper apply in this context, 
and the verification is the same as in the first case. In particular, observe that now 
$m(B(x,r)) \equiv r^{n-d - \gamma} \mu(B(x,r) \cap \R^d)$, so 
$\rho \approx r^{n-d-1 - \gamma}$ in \eqref{defrho}, with an exponent smaller than $1$, and
for (\hyperref[H6]{H6}), that in a ball $B$ such that $2B \subset \Omega$, we multiply $\omega(x)$ by a
function $|t|^{-\gamma}$ which is roughly constant.

Again, in such circumstances, the results of this paper are all valid, but more
precise results on the corresponding elliptic measures would need more precise assumptions on
the operators $L = \div A \nabla$.

\subsection{Stranger measures $\mu$}
\label{S3.6}

Even when $\Gamma = \d\Omega$ is a nice hypersurface, the measure $\mu$
on $\Gamma$ does not need to be as simple as surface measure; the next example
shows that it does not need to be absolutely continuous with respect to surface measure.
In dimension $d=1$, it could be given by a Riesz product, for instance, and hence one-dimensional
yet singular with respect to the Lebesgue measure.

Let $\Omega$ satisfy (\hyperref[H1]{H1}) and (\hyperref[H2]{H2}); we need to ask this because the
next assumption does not really say anything nice on the geometry of $\Omega$. 
Then let $\mu$ be any doubling measure whose support is $\Gamma$ (so (\hyperref[H3]{H3}) holds). 
We shall now define a measure $m$ on $\Omega$ such that all the other assumptions 
(\hyperref[H4]{H4})-(\hyperref[H6]{H6})
are satisfied. In view of \eqref{defrho} and the intertwining condition, it is reasonable to
take
\begin{equation} \label{st9}
w(X) = \delta(X)^{1-n} \mu(\Gamma \cap B(X,2\delta(X))),
\end{equation}
where we recall that $\delta(X) = \dist(X,\Gamma)$ for $X\in \Gamma$.
Let us first show that
\begin{equation} \label{st10}
m(B(x,r)) := \int_{\Omega \cap B(x,r)} w(X) dX \approx r \mu(B(x,r)) 
\qquad \text{ for $x\in \Gamma$ and $r > 0$.}
\end{equation}
To this effect, cover $B(x,r)$ by the sets
$R_k = \big\{ X\in B(x,r)  \, ; \, 2^{-k-1}r \leq \delta(X) \leq 2^{-k} r \big\}$,
$k \geq 0$, and further cover each $R_k$ by the balls $B_{k,l} = B(z_{k,l},2^{-k+2}r)$,
where $\{ z_{k,l} \}$ is a maximal collection of points of $\Gamma \cap B(x,4r)$
that lie at mutual distances larger than $2^{-k-1} r$.
Notice that $B(X,2\delta(X)) \subset B(z_{k,l}, 2^{-k+3}r)$
for $X \in R_k \cap B_{k,l}$, so
\[
m(R_k \cap B_{k,l}) \leq (2^{-k-1}r)^{1-n} \mu(B(z_{k,l},2^{-k+3}r)) \, |B_{k,l}| 
\leq C 2^{-k} r \mu(B(z_{k,l},2^{-k+3}r)).
\]
For each $k$, the $B_{k,l}$ have bounded overlap, so
\[
m(R_k) \leq \sum_l m(R_k \cap B_{k,l}) \leq C 2^{-k} r \, \sum \mu(B(z_{k,l},2^{-k+3}r))
\leq C 2^{-k} r \, \mu(\Gamma \cap B(x,12r)).
\]
We sum over $k \geq 0$ and get the upper bound in \eqref{st10}.
For the lower bound, we use (\hyperref[H1]{H1}) to select a corkscrew ball 
$B \subset \Omega \cap B(x,r)$. Observe that $w(X) \geq C^{-1} r^{1-n} \mu(B(x,r))$
for $X\in B$ (because $\mu$ is doubling), so $m(B(x,r)) \geq m(B) \geq C^{-1} r \mu(B(x,r))$, 
which completes the proof \eqref{st10}.

It follows from \eqref{st10} that $m$ is locally finite. The intertwining property (\hyperref[H5]{H5})
follows at once from \eqref{st10}, which says that $\rho(x,r) \approx 1$, and (\hyperref[H6]{H6})
holds because of (\hyperref[H7]{H6'}), by \eqref{st9} and because $\mu$ is doubling. We are left with the
doubling property (\hyperref[H4]{H4}) for $m$.

So let $X \in \Omega$ be given, and choose $x\in \Gamma$ such that $|X-x| = \delta(X)$.
For $R < \delta(X)/2$, $m(B(X,R)) \approx \delta(X) \mu(B(X,2\delta(X)))
\approx \delta(X) \mu(B(x,\delta(X))$ because $\mu$ is doubling.
For $R$ larger than $2\delta(X)$, $B(x,R/2) \subset B(X,R) \subset B(x,3R)$, so
$m(B(X,R)) \approx m(B(x,R) \approx R \mu(B(x,R))$ by \eqref{st10}. 
The doubling property follows easily.

So in this setting too, our assumptions hold and the rest of the paper will show that
there is a well behaved elliptic measure associated to each operator $L = \div A \nabla$
such that $w(x)^{-1}A$ satisfies the usual boundedness and ellipticity properties.

\section{The definition of the space $W$}
\label{SecW}

We want to define the space $W$ as the space of functions $u$ such that $\nabla u$ is in $L^2(\Omega,m)$, and we wish to prove that this space is complete; 
more precisely that the quotient of $W$ by constants, 
equipped with the quotient of the semi-norm $\|\nabla \cdot\|_{L^2(\Omega,m)}$, 
is complete.

However, 
$W$ is not entirely defined by the fact that $\nabla u$ is in $L^2(\Omega,m)$, because we don't explain where the functions $u$ are taken from. The first natural space where we could 
take $u$ from is $L^1_{loc}(\Omega,m)$, 
but in this case 
recall that we do not assume enough regularity on $w$ to make sure that 
$u \in L^1_{loc}(\Omega,dx)$, and then maybe $u$ does not define a distribution and
we do not know the meaning of $\nabla u$.
The second choice would be to take $u$ is the space of distributions, or in $L^1_{loc}(\Omega,dx)$, but nothing guarantees that the quotient of the constructed space by constants will be complete.

To solve this problem, we use the strategy from \cite{NDbook}, which consists in completing the smooth functions with respect to an appropriate norm. Our spaces shall be homogeneous, 
while the ones in \cite{NDbook} are inhomogeneous. 
Homogeneous spaces are slightly 
trickier, because we need to 
quotient by constant functions  
to get a Hilbert space. 

\medskip

\begin{definition} \label{defW}
A function $u$ belongs to $W$ if $u \in L^1_{loc}(\Omega,m)$ and there exists a vector valued function $v \in L^2(\Omega,m)$ such that for some sequence $\{\varphi_i\}_{i\in \bN} \in C^\infty(\overline{\Omega})$, one has
\begin{equation} \label{defW9}
\int_\Omega |\nabla \varphi_i|^2 \, dm < +\infty \qquad \text{ for any $i\in \bN$,}
\end{equation}
\begin{equation} \label{defW1}
\lim_{i\to \infty} \int_B |\varphi_i-u| \, dm = 0 \qquad \text{ for any ball $B$ satisfying $2B \subset \Omega$}
\end{equation}
and
\begin{equation} \label{defW2}
\lim_{i\to \infty} \int_\Omega |\nabla \varphi_i-v|^2 \, dm = 0.
\end{equation}
\end{definition}

Observe that if $u\in W$, then, assuming the first part of (\hyperref[H6]{H6}), the vector $v$ 
from the definition is unique. In the rest of the article, if $u\in W$, we use the notation 
$\nabla u$ (or $\nabla_W u$ when the notion of derivative we are talking is not obvious) for the 
unique vector valued function $v$ given by Definition \ref{defW}. In particular, we can equip $W$ with the semi-norm
\[\| u \|_W := \| \nabla u \|_{L^2(\Omega,m)} \qquad \text{ for $u\in W$}.\]
We want to highlight that $\nabla \cdot$ is a linear operator, but is not (necessarily) the gradient in the sense of distribution. An example where the two notions of derivative are different is given page 13 of \cite{NDbook}. 

Let us recall now some cases where the two notions of derivative actually coincide. First, if $L^1_{loc}(\Omega,dx) = L^1_{loc}(\Omega,dm)$ - which is the case for instance under the assumption (\hyperref[H7]{H6'}) - then \eqref{defW1} implies the convergence of the $\varphi_i$ to $u$ in $L^1_{loc}(\Omega,dx)$, which in turn implies the convergence $\varphi_i \to u$ in the sense of distribution. So the only possible limit of $\nabla \varphi_i$ is $\nabla u$, where $\nabla u$ is the derivative taken in the sense of distributions.

Let us present another case, given on
page 14 of \cite{NDbook}. Let the measure $m$ be absolutely continuous with respect to the Lebesgue measure, so that there exists a weight $w$ satisfying $dm(x) = w(x) dx$. 
Assume in addition that 
$w$ belongs to the Muckenhoupt class 
$\A_2$. Then 
the measure $m$ satisfies (\hyperref[H4]{H4}) and $\nabla_W u$ is 
the distribution gradient of $u$ in $\Omega$.
See Subsection~\ref{SS3.5} for a bit more information. 

In general,
Lemma 1.11 in \cite{NDbook} shows that 
\begin{equation} \label{NablaLip}
\text{ if
$u$ is compactly supported and Lipschitz, then 
$u\in W$ and $\nabla_W u$ is the usual gradient.}
\end{equation}
The proof of \eqref{NablaLip} uses the fact that $m$ is absolutely continuous 
with respect to the Lebesgue measure.

\medskip

We now show 
that the Poincar\'e inequality given as (\hyperref[H6]{H6}) extends to all functions in $W$.

\begin{lemma} \label{Poincare0}
Let $(\Omega,m)$ satisfy (\hyperref[H6]{H6}). Then for any ball $B$ such that 
$2B \subset \Omega$ and any $u\in W$, one has
\[\fint_B |u - u_B| \, dm \leq C_6 r \left( \fint_B |\nabla u|^2 \, dm \right)^\frac12,\]
where $u_B$ stands for $\fint_B u\, dm$ and $r$ is the radius of $B$. 
\end{lemma}

\bp
By definition of $W$, we can find a sequence of functions $\varphi_i \in C^\infty(\overline{\Omega})$ such that
\[\lim_{i\to \infty} \fint_B |\varphi_i - u| \, dm = 0\]
 and
\[\lim_{i\to \infty} \left(\fint_B |\nabla \varphi_i - \nabla u|^2 \, dm\right)^\frac12 = 0.\]
 Let $\varphi_{i,B}$ stand for $\fint_B \varphi_i \, dm$; then for $i\geq 0$
\[\begin{split}
\fint_B |u - u_B| \, dm & \leq \fint_B |\varphi_i - \varphi_{i,B}| \, dm + \fint_B |\varphi_i - u| \, dm  + |\varphi_{i,B} - u_B| \\
& \leq C_6 r \left(\fint_B |\nabla \varphi_i|^2 \, dm\right)^\frac12 + 2  \fint_B |\varphi_i - u| \, dm  \\
& \leq C_6 r \left(\fint_B |\nabla u|^2 \, dm\right)^\frac12 + C_6 r \left(\fint_B |\nabla \varphi_i - \nabla u|^2 \, dm\right)^\frac12 + 2 \fint_B |\varphi_i - u| \, dm, 
\end{split}\]
where we use (\hyperref[H6]{H6}) in the second inequality. Taking the limit as $i\to \infty$ gives the desired result.
\ep

We have the following nice improvement of Lemma \ref{Poincare0} by  
Keith and Zhong \cite{KZ}, where it is enough to control $|\nabla u|$ in some $L^p$ norm, $p < 2$.

\begin{lemma} \label{Poincare00}
Let $(\Omega,m)$ satisfy (\hyperref[H4]{H4}) and (\hyperref[H6]{H6}). 
There exists $p_0 \in [1,2)$ such that for any $p\in [p_0,2]$, any ball $B$ satisfying 
$2B \subset \Omega$, and any $u\in W$, 
\begin{equation} \label{defpPoincare}
\fint_B |u - u_B| \, dm \leq C r \left( \fint_B |\nabla u|^p \, dm \right)^\frac1p,
\end{equation}
where $u_B$ stands for $\fint_B u\, dm$ and 
$r$ is the radius of $B$. The parameter $p_0$ and the constant $C$ depends only on $C_4$ and $C_6$.
\end{lemma}

\begin{remark}
If \eqref{defpPoincare} is true for some $p_0<2$, then it holds 
for all $p\in [p_0,2]$, 
by Jensen's inequality.
\end{remark}

\bp The Poincar\'e inequality \eqref{defPoincare} holds for all locally Lipschitz functions according to Lemma \ref{Poincare0}, \eqref{NablaLip}, and the fact that \eqref{defPoincare} is a local property. 
We deduce that our metric measured spaces $(B,|.|,m)$ are doubling spaces that admit a $(1,2)$-Poincar\'e inequality in the sense of \cite{KZ}, and the doubling constant and the Poincar\'e constant are uniform on the balls $B$. 
Theorem 1.0.1 in \cite{KZ} applies, so we have the existence of $p\in(1,2)$ and $C>0$ independent of $B$ such that our spaces $(B,|.|,m)$ admit a $(1,p)$-Poincar\'e inequality with constant $C$, which means in the terminology of \cite{KZ} that \eqref{defpPoincare} holds for any locally Lipschitz function.
The proof of \eqref{defpPoincare} for all $W$ then follows 
from the same density argument as 
in the proof of Lemma \ref{Poincare0}.
\ep

We end the section with a simple but useful lemma.

\begin{lemma} \label{Nablau}
Let $(\Omega,m)$ satisfy (\hyperref[H1]{H1}), (\hyperref[H2]{H2}), and (\hyperref[H6]{H6}). 
Take $u\in W$. Then  
$\|u\|_W = 0$
if and only if $u$ is $m$-almost everywhere equal to a constant function.
\end{lemma}

\bp
First, constant functions are in $C^\infty(\Omega)$. So if $u$ matches a constant function $c$ except maybe on a set of $m$-measure $0$, we can take $v=0$ and $\varphi_i = c$ in the Definition \ref{defW}. By the uniqueness of $v = \nabla u$, we deduce that 
$\nabla u = 0$. 

Conversely, let $u \in W$ be such that $\|u\|_W = 0$.
By Lemma \ref{Poincare0}, for any ball $B$ such that  
$2B \subset \Omega$, we have that $\int_B |u-u_B| \, dm = 0$, which implies immediately that $u \equiv u_B$ $m$-a.e. on $B$. Yet, $\Omega$ is connected (and can even be connected by a chain of balls $\{B_i\}_i$ satisfying $2B_i \subset \Omega$, 
thanks to Proposition \ref{propHarnack}), so $u$ is $m$-a.e. equal to a constant function.
\ep

\section{The access
cones and their properties}

\label{Scones}

In all this section, we assume that $\Omega$ satisfies (\hyperref[H1]{H1})--(\hyperref[H2]{H2}) and that the measures $\mu$ and $m$ satisfy 
(\hyperref[H3]{H3}) and (\hyperref[H4]{H4}). 
We also choose to take $\Gamma$ (and thus $\Omega$) to be infinite. This assumption is not part of (\hyperref[H1]{H1})--(\hyperref[H4]{H4}), and is not even necessary for our proofs to work. The proofs of the bounded and unbounded cases only differ slightly, but will require us to separate cases. We will present the infinite case - which we plan to use in future articles and which we believe is less common - and we shall discuss the differences in Section \ref{SBounded}.
We first describe the dyadic decomposition of $(\Gamma,\mu)$ of M. 
Christ (see \cite[Theorem 11]{Christ}). 

\begin{proposition} \label{ChristQ}
There exists a collection of measurable 
subsets - we call them cubes by comparison with the Euclidean case - $\{Q_j^k, \, k\in \mathbb Z, \, j\in \mathcal J_k\}$, and some constants $\eta$, $a_0$, $C$ - all of them depending only on $C_3$ - such that
\begin{enumerate}[(i)]
\item $\Gamma = 
\bigcup_{j \in \mathcal J_k} 
Q^k_j$ for all $k\in \mathbb Z$.
\item If $\ell \geq k$, then either $Q_i^\ell \subset Q^k_j$ or $Q_i^\ell \cap Q^k_j = \emptyset$.
\item For each pair $(k,j)$ and each $\ell < k$, there exists a unique $i$ such that $Q^k_j \subset Q^\ell_i$.
\item $\diam \, Q^k_j \leq 2^{-k}$.
\item 
$Q^k_j$ contains some surface ball $B(z^k_j,a_02^{-k}) \cap \Gamma$.
\item $\mu(\{x\in Q^k_j, \, \dist(x,\Gamma\setminus Q^k_j) \leq \rho 2^{-k}\}) \leq C\rho^\eta \mu(Q^k_j)$ for all $k\in \mathbb Z$, 
$j\in \mathcal J_k$, and some constant $\rho >0$.
\end{enumerate}
We shall denote by $\mathbb D_k$ the collection 
\[\mathbb D_k:= \{Q^k_j,\, j \in \mathcal J_k\}\]
and by $\D$ the collection
\[\mathbb D := \bigcup_{k\in \ZZ} \mathbb D_k.\]
\end{proposition}

\begin{remark}
An element of $\D$ is given by a subset $Q$ of $\Gamma$ and a generation $k$. Indeed, if we only know the set $Q$, contrary to dyadic cubes in $\R^n$, we cannot be sure of the generation. 

Despite the above comment, 
we shall abuse notation and use 
the term $Q$ for both an element of $\D$ and the corresponding subset of $\Gamma$. We write $k(Q)$ when we want to refer to the ``dyadic generation'' of the cube $Q\in \D$, that is the \underline{only} integer $k$ such that $Q \in \D_k$. The length of a dyadic cube is $\ell(Q) = 2^{-k(Q)}$.
\end{remark}

The conclusion $(vi)$ will not be used in this article, but we wanted to state the complete result of Christ nevertheless. Moreover, properties $(iv)$ and $(v)$ of the decomposition implies the existence of $z_Q \in \Gamma$ such that 
\begin{equation} \label{defzQrQ}
B(z_Q,r_Q) \cap \Gamma \subset Q \subset B(z_Q,R_Q),
\quad \text{with $r_Q= a_0 \ell(Q)$ and $R_Q = \ell(Q)$.}
\end{equation}
When $Q\in \D$ and $\lambda \geq 1$, we also use the notation $\lambda Q$ for the set $\{x\in \Gamma:\, \dist(x,Q) \leq (\lambda-1)\ell(Q)\}$. 
As a consequence, if $Q$ and $Q'$ are from the same generation, i.e., 
$k(Q) = k(Q')$, and $Q$ and $Q'$ are adjacent, i.e., 
$\dr Q \cap \dr Q' \neq \emptyset$, then $Q' \subset 2Q$.

\medskip

Also, as in the first pages of \cite{Stein}, we can define a Whitney decomposition of $\Omega \subset \R^n$ made by (true) dyadic cubes. To do this, take a dyadic decomposition of $\R^n$ by cubes $I$, ordered by inclusion, and we define $\mathcal W$ as a the set of dyadic cubes $I \subset \Omega$ for which $4 \, \diam \, I \leq \dist(4I,\Gamma)$ but the parent $I'$ of $I$ - that is the only dyadic cube $I'\supset I$ satisfying $\ell(I') = 2\ell(I)$ - doesn't satisfy  $4 \, \diam \, I' \leq \dist(4I',\Gamma)$. It is easy to check that $\mathcal W$ is a non-overlapping covering of $\Omega$, that for $I\in \mathcal W$
\begin{equation} \label{WQ11}
4 \, \diam \, I \leq \dist(4I,\Gamma) \leq \dist(I,\Gamma) \leq 12 \diam(I)
\end{equation}
and if $I_1,I_2 \in \mathcal W$ are two adjacent cubes 
\begin{equation} \label{WQ12}
\frac{\diam \, I_1}{\diam \, I_2} \in \Big\{\frac12,1,2 \Big\}.
\end{equation}
Let us write $X_I$ for the center of $I \in \mathcal W$, $\ell(I)$ for its side length 
(thus $\ell(I) \approx \diam\, I$), and $k(I)$ for the integer $k$ that satisfies $\ell(I) = 2^{-k}$.

\medskip

Now, let us match the dyadic decomposition $\mathbb D$ of $\Gamma$ with 
the Whitney decomposition $\mathcal W$ of $\Omega$. 
For each $Q\in \mathbb D$, we define $\mathcal W_Q$ as
\begin{equation} \label{defWQ}
\mathcal W_Q:= \{I \in \mathcal W, \, C_a^{-1}\ell(Q) \leq \ell(I) \text{ and } \dist(I,Q) \leq 2\ell(Q)\},
\end{equation}
where $C_a^{-1} = C(C_1,n) > 1$ is chosen 
in the following next lines. 
Set $X_Q$ as a Corkscrew point associated to a point $x_Q \in 2Q$ and a distance $\ell(Q)$, 
that is $X_Q \in B(x_Q,\ell(Q))$ and $B(X_Q,\ell(Q)/C_1) \subset \Omega$. 
The point $X_Q$ belongs to some $I_Q \in \mathcal W$. Observe that 
\[\dist(I_Q,Q) \leq |X_Q - x_Q| + \dist(x_Q,Q) \leq 2\ell(Q)\]
and
\[\ell(I_Q) \geq \frac1{16\sqrt{n}} \dist(X_Q,\Gamma) \geq \frac{1}{16C_1\sqrt{n}} \ell(Q);\]
we can pick 
the constant $C_a$ in \eqref{defWQ} as  for instance  $1000C_1\sqrt n$, so that $I_Q \in \mathcal W_Q$. But the choice of 
$C_a$ doesn't really matter (as long as it is big enough); 
we can choose it as an additional parameter and make the future results depend on $C_a$ too. 
Now define 
the associated Whitney region
\begin{equation} \label{defUQ}
U_Q := \bigcup_{I \in \mathcal W_Q} I,
\end{equation}
which contains by construction of $\cW_Q$ all the Corkscrew points associated 
to a point $x\in 2Q$ and the 
distance $\ell(Q)$. We also define, for each $x\in \Gamma$, the 
``dyadic access'' cone
\begin{equation} \label{defgamma}
\gamma(x) := \bigcup_{Q\in \mathbb D:\, Q\ni x} U_Q.
\end{equation}

\medskip

We also need cones with a ``larger aperture''. We consider 
 the collection $\mathcal W_Q^0$ of dyadic
cubes that meet $B(X,\delta(X)/2)$ for some $X \in U_Q \cup U_{Q'}$, where $Q'$ is the parent of $Q$. 
Thus, when $I \in \mathcal W_Q^0$,
$\delta(X_I) \approx \ell(Q)$ with constants that depends only on $C_a$ (i.e., $n$ and $C_1$), so each couple of centers $X_I,X_{I'}$,  $I,I' \in \mathcal W_Q^0$, can be linked by a Harnack chain (see Proposition~\ref{propHarnack}). We define 
$\mathcal W_Q^*$ as the collection of cubes in $\cW$ that meet at least one of those Harnack chains
from 
\eqref{HarnackBalls}, and finally define
\begin{equation} \label{defU*Q}
U_Q^* := \bigcup_{I \in \mathcal W_Q^*} I 
\end{equation}
and, for $x\in \Gamma$, the cone 
\begin{equation} \label{defgamma*}
\gamma^*(x) := \bigcup_{Q\in \mathbb D:\, Q\ni x} U^*_Q.
\end{equation}
We shall also need the truncated cone
\begin{equation} \label{defgammaQ*}
\gamma^{*}_Q(x) := \bigcup_{Q'\in \mathbb D:\, x \in Q' \atop \ell(Q') \leq \ell(Q)} U^*_{Q'},
\end{equation} 
and the ``tent sets'' 
\begin{equation} \label{defTQ}
T_Q:= \bigcup_{x\in Q} \gamma^{*}_{Q}(x) \quad \text{ and } \quad T_{2Q}:= \bigcup_{x\in 2Q} \gamma^{*}_{Q}(x). 
\end{equation}

\medskip
The following standard properties of the sets above are easy to check.
The cones $\gamma(x), \gamma^*(x)$ are such that $\gamma(x) \subset \gamma^*(x)$ and
\begin{equation} \label{prgamma}
\delta(X) > c|X-x| \quad \text{ for } X \in \gamma^*(x).
\end{equation}
The Whitney regions $U_Q$ and $U^*_Q$ are such that $U_Q \subset U^*_Q$ and
\begin{equation} \label{prU}
\ell(Q) \lesssim \dist(U^*_Q, Q) \leq \dist(U_Q, Q) \lesssim \diam\, U_Q \leq \diam\, U^*_Q 
\lesssim \ell(Q), 
\end{equation}
where the constants depends only on $n$, $C_1$, and $C_2$. 
The tent sets $T_Q$ and $T_{2Q}$ satisfy
\begin{equation} \label{propTQ}
B(z_Q,r'_Q) \cap \Omega \subset T_Q \subset T_{2Q} \subset B(z_Q,R'_Q),
\end{equation}
where $z_Q$ is as 
in \eqref{defzQrQ}, and $r'_Q,R'_Q \approx \ell(Q)$. 
Indeed, the second inclusion is easy; 
for the first one, observe that if $Z \in B(z_Q,r'_Q)$ with $r'_Q$ small enough, 
then any point in $\Gamma$ such that 
$|z-Z| = \delta(Z)$ lies in $B(z_Q,r_Q)$, 
where $r_Q = a_0 \ell(Q)$ as 
in \eqref{defzQrQ}. 
The point $Z$ is a Corkscrew point for $z$, so $Z \in \gamma(z)$, and as long as $r'_Q$ is small enough, it is also in $\gamma_Q^*(z) \subset T_Q$. The measure of the various sets that we just introduced 
are given by the following lemma. 

\begin{lemma} \label{lemmQ}
Let $Q\in \mathbb D$ and $x\in Q$. 
Then 
\begin{enumerate}[(i)]
\item $\mu(Q) \approx \mu(B(x,\ell(Q)))$,
\item $m(U_Q) \approx m(U^*_Q) \approx m(B(x,\ell(Q)) \cap \Omega)$,
\item $\rho(x,\ell(Q)) \approx \frac{m(U^*_Q)}{\mu(Q) \ell(Q)}$.
\end{enumerate}
In (i), the constants depends only on $C_3$, and in (ii) and (iii), the constants depend also
on $n$, $C_1$, $C_2$, and $C_4$.

\smallskip
In particular, we can define $\rho(Q)$ as 
\begin{equation} \label{defrhoQ}
\rho(Q) := \frac{m(U^*_Q)}{\mu(Q) \ell(Q)}, 
\end{equation}
and if (\hyperref[H5]{H5}) if satisfied, we have 
\begin{equation} \label{H5Q}
\frac{\rho(Q^*)}{\rho(Q)} \leq C \left(\frac{\ell(Q^*)}{\ell(Q)}\right)^{1-\epsilon},
\end{equation}
where $C>0$ depends on $n$, $C_1$ to $C_5$.
\end{lemma}

\bp Let us prove (i). By \eqref{defzQrQ} and (\hyperref[H3]{H3}),  
\[\mu(Q) \leq \mu(B(z_Q,R_Q)) \leq \mu(B(x,2R_Q)) \lesssim \mu(B(x,\ell(Q)))\]
and
\[\mu(B(x,\ell(Q)) \leq \mu(B(z_Q,2R_Q)) \lesssim \mu(B(z_Q,r_Q)) \leq \mu(Q).\]
The assertion (i) follows. As for (ii), since $U_Q,U_Q^*$ are Whitney regions associated to $Q$, \eqref{prU} shows that we can find $K>1$ and $X \in U_Q$ such that 
\[B(X,K^{-1}\ell(Q)) \subset U_Q \subset U^*_Q\subset B(x,K\ell(Q)) \cap \Omega \subset B(X,K^2\ell(Q)).\]
The assertion (ii) is now an immediate consequence of (\hyperref[H4]{H4}), the doubling measure property for $m$. The conclusion (iii) is no difficulty from (i) and (ii).
\ep

One can also easily check that the number of dyadic cubes in $\cW^*_Q$ is uniformly bounded. Indeed, 
the cubes in $\cW^*_Q$ are pairwise disjoint, and their diameters 
are all equivalent to the diameter of $U^*_Q$ - which is their union. One can also easily check that $U^*_Q$ is connected (by construction, we linked the points in $U_Q \cup U_{Q'}$ by Harnack chains). So since $W^*_Q$ is only constituted of dyadic cubes, for any couple $I,I' \in \mathcal W^*_Q$, we can find a sequence of cubes in $\mathcal W^*_Q$ linking $I$ to $I'$, where two consecutive cubes are adjacent; the sequence has uniformly bounded length because there is a bounded number 
of cubes in $\mathcal W^*_Q$. We summarize these conclusions in the following lemma.

\begin{lemma} \label{lemUQchain}
There exists $N_0:= N_0(n,C_1,C_2) \in \bN$ such that for $Q\in \mathbb D$ and 
$I,I' \in \mathcal W^*_Q$, we can find a collection $\{I_i\}_{0\leq i \leq N_0}$ of cubes in $\mathcal W^*_Q$ such that 
\begin{enumerate}[(i)]
\item $I_0 = I$, $I_{N_0} = I'$,
\item for any $i\in \{1,\dots,N_0\}$, $I_{i-1}$ and $I_i$ are adjacent or equal.
\end{enumerate}
\end{lemma}

As a corollary, we get the following result with balls instead of cubes.

\begin{lemma} \label{lemUQchain2}
There exists $N_0:= N_0(n,C_1,C_2) \in \bN$ such that for $Q\in \mathbb D$ and for $I,I' \in \mathcal W^*_Q$, we can find a collection $\{B_i\}_{0\leq i \leq N_0}$ of balls such that 
\begin{enumerate}[(i)]
\item $2B_i \subset \Omega$ and $B_i \subset U^*_Q$,
\item $B_0$ is $B(X_I,\ell(I)/2)$ and $B_{N_0}$ is $B(X_{I'},\ell(I')/2)$,
\item for any $i\in \{1,\dots,N_0\}$, we have $r_i \approx \ell(I)$, where $r_i$ is the radius of $B_i$,
\item for any $i\in \{1,\dots,N_0-1\}$, one has $|X_{i+1} - X_i| \leq r_i$, where $X_{i}$ is the center of $B_{i}$.
\end{enumerate}
\end{lemma}

\bp
We construct the sequence of balls $\{B_i\}$ from the sequence of dyadic cubes $\{I_i\}_{0\leq i \leq N_0}$ as follows. We replace each cube $I_i$, $i<N_0$, by $n+2$ balls $\{B^j\}_{0\leq j \leq n+1}$, according to the following procedure:
\begin{itemize}
\item If $I_{i+1}$ is smaller than $I_i$, then since $I_i$ and $I_{i+1}$ are adjacent, 
hence $\ell(I_i) = 2\ell(I_{i+1})$ by \eqref{WQ12}.
So up to translation, rotation, and dilatation, $I_{i}$ is the cube $[0,4]^n$ and $I_{i+1}$ is the cube $[-2,0]\times [0,2]^{n-1}$. In this case, we take $B^0$ as the ball with center at $(2,\dots,2)$ - the center of $I_i$ - and radius $2$, the balls $B^j$, $1\leq j\leq n$, are centered on
\[(\underbrace{2,\dots,2}_{n-j},\underbrace{1,\dots,1}_{j})\]
and are of radius 1, the ball $B^{n+1}$ is centered on $(0,1,\dots,1)$ and again of radius 1.
\item If $I_{i+1}$ has the same size of $I_i$, yet is different from $I_i$, then up to rotation, 
translation and dilatation, $I_i=[0,4]^n$ 
and $I_{i+1}=[-4,0]\times [0,4]^{n-1}$. The  
$B^j$ have the same 
radius $2$, $B^0$ is the ball centered on $X_{I_{i}} = (2,\dots,2)$, and all 
the other balls $B^j$ are equal and centered on $(0,2,\dots,2)$.
\item If $I_{i+1}$ is bigger than $I_i$, then as 
before we necessary have $2\ell(I_i) = \ell(I_{i+1})$. So up to translation, rotation, and dilatation, 
$I_{i}=[-2,0] \times [0,2]^{n-1}$ and $I_{i+1}=[0,4]^n$. 
All the balls have but the last one have radius $1$ and $B^{n+1}$ has radius 2; $B^0$ is centered on $(-1,1\dots,1)$, $B^1$ is centered on $(0,1,\dots,1)$, and for $2\leq j \leq n+1$, 
$B_j$ is centered on 
\[(\underbrace{1,\dots,1}_{n+2-j},\underbrace{2,\dots,2}_{j-1}).
\]
\item If $I_{i+1} = I_i$, then $B^j$ is 
always the same ball $B(X_{I_i},\ell(I_i)/2)$.
\end{itemize}
We replace $I_{N_0}$ by the ball $B(X_{I_{N_0}}, \ell(I_{N_0})/2)$.

The balls that we constructed satisfy (i), because first 
$B^j \subset I_{i} \cup I_{i+1}$ 
and second, the Whitney cubes $I_i$  satisfy \eqref{WQ11}, which ensures that 
$2B^j \subset \Omega$; 
(ii) and (iv) are not hard to check by construction, (iii) comes from the fact that all $I_i$ have similar radius (equivalent to the diameter of $U^*_Q$). The lemma follows.
\ep

\medskip

We shall 
use the last lemma to prove quantitative connectedness on the sets $U^*_Q$, $\gamma^*_Q(x)$, 
and $T_Q$. We start with a definition. 

\begin{definition} \label{defPchain}
We say that a (bounded) 
set $D \subset \Omega$ satisfies the chain condition $C(\kappa,M)$, where $\kappa \in [1/2,1)$, if there exists a distinguished ball $B_0 \subset D$ such that for every $x\in D$, there exists an infinite sequence of balls $B_0,B_1,\dots$ (called chain) with the following properties:
\begin{enumerate}[(i)]
\item for $i\in \bN$, we have $B_i \subset D$ and $2B_i \subset \Omega$;
\item for $i\geq 0$, $x\in MB_i$;
\item for $i\geq 0$, one has 
$$M^{-1} (\diam \,D) \kappa^{i} \leq r_i \leq  M (\diam \,D) \kappa^{i},$$
where $r_i$ is the radius of $B_i$;
\item for $i\geq 0$, if $X_j$ denotes the center of $B_j$, we have $|X_{i+1} - X_i| \leq r_i$
\end{enumerate}
\end{definition}

\begin{remark}
The definition above is shamelessly inspired by the $C(\lambda,M)$ condition in \cite{HK2}. 
Notice
that $\kappa$ in our condition doesn't correspond to $\lambda$ in the chain condition of \cite{HK2}. Indeed, $\kappa$ is fixed equal to 1/2 in \cite{HK2}, while the $\lambda$ in \cite{HK2} 
doesn't really have an equivalent in our condition. However, these technicalities don't really 
change the core the proofs.
\end{remark}

\begin{lemma} \label{lemchain}
For every
$\kappa \in [1-n^{-1/2},1)$, there exists $M:=M(\kappa,n,C_1,C_2)$ such that each Whitney cube $I \in \cW$, and each set $U^*_Q$, $Q\in \mathbb D$, satisfies the chain condition $C(\kappa,M)$.

There exists $\kappa \in [1/2,1)$ and $M \geq 1$ - both depending only on $n$, $C_1$, $C_2$, and $C_4$ - such that for any $Q\in \mathbb D$ and any $x\in 2Q$, the sets $\gamma^*_Q(x)$, 
$T_Q$, and $T_{2Q}$, satisfy the chain condition $C(\kappa,M)$.
\end{lemma}

\bp
We start with an (open) Whitney cube $I\in \cW$. Take $\kappa \in [1-n^{-1/2},1)$. We choose the distinguish ball associated to $I$ as $B_0 := B(X_I,\ell(I)/2)$. Then we take $X\in I$ and we construct the chain of balls $\{B_i\}_{i\geq 0}$ as follows. For $i\geq 1$, the ball $B_i$ has radius $r_i=\kappa^i \ell(I)/2$ and its center $X_i$ is the closest point to $X$ on the segment $[X_I,X]$ which satisfies $|X_i - X_{i-1}| \leq r_{i-1}$ and $\dist(X_i,\dr I) \leq r_i$. If $M = \sqrt n$, the points (iii) and (iv) of Definition \ref{defPchain} are true by construction, as well as the fact that $B_i \subset I$. The condition $2B_i \subset \Omega$ is true because we have $B_i \subset I$ and \eqref{WQ11}. The condition (ii) of Definition \ref{defPchain} holds because we chose $\kappa$ large enough to ensure that we can get (at least infinitely close) to $X$ at some point.

\smallskip

Now let $\kappa \in [1-n^{-1/2},1)$ and $Q\in \D$ be given.
We want to prove that the sets $U^*_Q$ 
satisfy the chain condition $C(\kappa,M)$ for some $M$. We choose $I_0$ as any dyadic cube in 
$\mathcal W^*_Q$ (the choice is not important here), and then we choose the distinguished
ball $B_0$ as $B(X_{I_0},\ell(I_0)/2)$. 
Take then $X\in U^*_Q$. There exists $I\in \cW^*_Q$ such that $X\in I$. The balls $B_i$ are constructed as follows: $\{B_i\}_{0\leq i \leq N_0}$ is the collection of balls linking the center of $I_0$ to the center of $I$ given by Lemma~\ref{lemUQchain2}, and the balls $\{B_i\}_{i>N_0}$ are the chain associated to the cube $I$ and the point $x$ that we constructed above. We can
check that the chain satisfies all the conditions 
of Definition~\ref{defPchain} when $M$ is large enough.

\smallskip

We turn to the proof of the chain condition for the sets $\gamma^*_Q(x)$.
For each $j\in \bN$, we define $Q_j$ as the dyadic cube in $\mathbb D_{j+k(Q)}$ that contains $x$. 
Choose for $X_j$
a Corkscrew point associated to $x$ and $\ell(Q_j)=2^{-j}\ell(Q)$. 
By construction of $\mathcal W_{Q_j}$, we can find a dyadic cube $I^j \in \mathcal W_{Q_j}$ that contains $X_j$. We construct the chain $\{\mathcal B_i\}$ as follows: 
for $j\in \mathbb N$, $\{\mathcal B_i\}_{jN_0 \leq i \leq (j+1)N_0}$ is the collection linking the center of $I^j$ to the center of $I^{j+1}$ given by Lemma \ref{lemUQchain2} (recall that both $I^j$ and $I^{j+1}$ are in $\mathcal W^*_{Q_j}$).

Now let us take $X \in \gamma^*_Q$. By construction, $X$ lies
in $U^*_{Q_{j(X)}}$ for some $j(X) \in \bN$. We construct the chain $\{B_i\}_{i\geq 0}$ as follows: 
if $i\leq jN_0$, then $B_i  = \mathcal B_i$; and then
the chain $\{B_i\}_{i\geq jN_0}$ is the one used to prove that $U_{Q_j}^*$ satisfies the chain condition $C(\kappa,M)$ with $\kappa=2^{-1/N_0}$. 
\medskip

At last, we shall prove that $T_Q$ and $T_{2Q}$ satisfy the chain condition $C(\kappa,M)$ for $\kappa:=2^{-1/N_0}$ and for some $M$ independent of $Q$. We only prove it for 
$T_{2Q}$, since $T_Q$ is very similar. It is actually an easy consequence of the chain condition 
of $\gamma_Q^*(x)$ and of $U_Q^*$. Indeed, we chose the distinguish ball $B_0^x$ 
of $\gamma_Q^*(x)$ as a ball centered on a dyadic cube $I_0^x$ containing a Corkscrew point 
associated to $(x,\ell(Q))$. However, by construction of $\cW_Q$, all the
$B_0^x \subset I_0^x$ are subsets of the same $U_Q \subset U^*_Q$. 
So we take any  $I_0\in \cW_Q^*$, we chose $B_0:= B(X_{I_0},\ell(I_0)/2)$ as the distinguish ball. Take then $X\in T_{2Q}$, and pick $x \in 2Q$ so that $X \in \gamma^*_Q(x)$. We construct the chain between the distinguish cube $B_0$ and $X$ as the concatenation of the chain (of finite length) linking $B_0$ to $B_0^x$ given by Lemma \ref{lemUQchain} and the one linking $B_0^x$ to $X$ given by the fact that $\gamma^*_Q(x)$ satisfies the $C(\kappa,M)$ chain condition.  
The lemma follows.
\ep

We may now
extend the Poincar\'e inequality given in (\hyperref[H6]{H6}) to domains that are not balls.

\begin{theorem} \label{Poincare1}
Assume that $(\Omega,m,\mu)$ satisfies (\hyperref[H1]{H1})--(\hyperref[H4]{H4}) and (\hyperref[H6]{H6}). 
Let $p_0\in (1,2)$ be as in Lemma \ref{Poincare00}, and take $p\in [p_0,2]$.

Let $M>1$ and $\kappa \in (1/2,1)$. Assume that $D \subset \Omega$ satisfies the chain condition $C(\kappa,M)$. Then there exists $k>1$,
that depends only on $C_4$, 
such that, for any $u \in W$, 
\begin{equation} \label{Poincare1a}
\left( \fint_D |u- \bar u|^{pk} \, dm \right)^{1/pk} \leq C \diam(D) \left(  \fint_D |\nabla u|^p \, dm \right)^p,
\end{equation}
where $\bar u$ is the average of $u$ on any set $E \subset D$ satisfying $m(E) \geq c m(D)$, and where $C>0$ depends only on $\kappa$, $M$, $C_4$, $C_6$, and $c$.

In particular, for any cube $Q\in \D$ and any $x\in 2Q$, 
\eqref{Poincare1a} holds for $D = U^*_Q$, 
$\gamma^*_Q(x)$, $T_Q$, or $T_{2Q}$, and the constant $C$ depends now (only) on $n$, $C_1$, $C_2$, $C_4$, $C_6$, and $c$.
\end{theorem}

\begin{remark} \label{rmkL1}
The theorem gives in particular that any function $u\in W$ lies 
in $L^1(D)$, where $D$ is any 
domain that satisfies the $C(\kappa,M)$ condition for some $\kappa$ and $M$. In particular $D$ can stand for $\gamma^*_Q(x)$, $T_Q$, or $T_{2Q}$, despite the fact that none of these domains are relatively compact in $\Omega$.
\end{remark}

\begin{remark} \label{rmkL2g}
We can apply the theorem when $D = 2B$, where $B$ is a ball such that $2B \subset \Omega$,
and $u \in W$ vanishes a.e. on $2B \sm B$; then we can take $E = 2B \sm B$ and \eqref{Poincare1a} becomes
\begin{equation} \label{Poincare1ag}
\left( \fint_B |u|^{pk} \, dm \right)^{1/pk} \leq C \diam(B) \left(  \fint_B |\nabla u|^p \, dm \right)^p,
\end{equation}
because $u=|\nabla u| = 0$ a.e. on $D \sm B = 2B \sm B$ anyway.
\end{remark}

\bp Let us not lie, our proof is the one of \cite{HK2} with very small
modifications. But we write it for completeness (and since it is quite short and fun). Also, in all the proof, if $S \subset D$, then $u_S$ denotes $\fint_S u \, dm$.

\medskip

Let $B_0 \subset D$ be the distinguished ball given by the $C(\kappa,M)$ condition. 
Also
write $r$ for the diameter of $D$. 
From (ii) of Definition \ref{defPchain}, the radius $r_0$ of $B_0$ is equivalent to $r$, so we deduce from 
(\hyperref[H4]{H4}) that $m(I_0) \approx m(D)$. As a consequence,
\[\begin{split}
\left(\fint_{D} |u- \bar u|^{kp} \, dm \right)^{1/kp} & \leq \left( \fint_{D} |u- u_{B_0}|^{kp} \, dm \right)^{1/kp} + |\bar u - u_{B_0}| \\
& \leq \left( \fint_{D} |u- u_{B_0}|^{kp} \, dm \right)^{1/kp} +  \fint_E |u - u_{B_0}| \, dm \\
& \lesssim \left( \frac1{m(D)} \int_{D} |u- u_{B_0}|^{kp} \, dm \right)^{1/kp} 
\end{split}\] 
by the H\"older inequality and the fact that $m(E) \approx m(D)$.

\medskip

So it is enough to prove the theorem when $\bar u = u_{B_0}$. Besides, without loss of generality, we can assume that $u_{B_0} = 0$. Our goal is to establish a weak-type $L^q-L^p$ estimate for $q>p$ that will be improved into a strong $L^{q'}-L^p$ estimate for $q' \in (p,q)$ by a standard argument.

Let $Z\in A_t:=\{|u|>t\}$ be a Lebesgue point for 
$u$, i.e., a point $Z$ such that 
\[ \lim_{r \to 0} \sup_{B_r \text{ ball of radius } r \atop
\text{and }
x\in MB_r} \fint_{B_r} |u(X) - u(Z)| \, dm(X) = 0.\]
It is well known 
that the Lebesgue points have full measure, i.e. $m(A_t) = m(\{Z\in A_t, \, Z \text{ is a Lebesgue point}\})$.

Let $B_0,B_1,\dots$ be the chain assigned to $Z$ and given by Definition \ref{defPchain}, and write $r_i$ for the radius of $B_i$. Pick a ball $B'_i \subset B_i \cap B_{i+1}$ with radius comparable to $r_i$ 
(and $r_{i+1}$); it is indeed possible since $r_i \approx r_{i+1}$ and, thanks to (iv) of 
Definition \ref{defPchain}, the center $X_{i+1}$ of $B_{i+1}$ belongs to $\overline{B_i}$. 
Since 
$Z$ is a Lebesgue point of $u$ and since the chain $\{B_i\}_{i\geq 0}$ satisfies (ii) and (iii) of Definition \ref{defPchain}, 
\[\begin{split}
t & < |u(Z) - u_{B_0}| \leq \sum_{i\in \bN} |u_{B_i} - u_{B_{i+1}}| \leq \sum_{i\in \bN} \left( |u_{B_i} - u_{B'_i}| + |u_{B_{i+1}} - u_{B'_i }| \right) \\
& \leq \sum_{i\in \bN} \fint_{B'_i} [|u-u_{B_i}| + |u-u_{B_{i+1}}|] \, dm \\
& \lesssim \sum_{i\in \bN}   \left(   \fint_{B_i} |u-u_{B_i}| \, dm +  \fint_{B_{i+1}} |u-u_{B_{i+1}}| \, dm \right).
 \end{split}\]
Poincar\'e's inequality (\hyperref[H6]{H6}) implies that
\[\begin{split}
t & \lesssim r \sum_{i\in \bN}  \kappa^i \left[  \left( \fint_{B_i} |\nabla u|^p \, dm \right)^\frac1p +  \left(\fint_{B_{i+1}} |\nabla u|^p \, dm\right)^\frac1p \right]\\
& \lesssim r \sum_{i\in \bN}  \kappa^i  \left(\fint_{B_i} |\nabla u|^p \, dm\right)^\frac1p,
 \end{split}\]
which can be written, when $\epsilon>0$, as
\begin{equation}
r \sum_{i\in \bN} \kappa^i\left(\fint_{B_i} |\nabla u|^p\, dm \right)^\frac1p \gtrsim t \gtrsim t \sum_{i\in \bN} \kappa^{i\epsilon}.
\end{equation}
The estimate above proves that there exists $i_Z$ such that 
\[ r \kappa^{i_Z}\left(\fint_{B_{i_Z}} |\nabla u|^p\, dm \right)^\frac1p \gtrsim t \kappa^{i_Z\epsilon}\]
hence, taking the power $p$ and writing the average explicitly,
\begin{equation} \label{PoincareA}
\kappa^{i_Zp(\epsilon-1)}  m(B_{i_Z})  \lesssim  \left(\frac rt\right)^p \int_{B_{i_Z}} |\nabla u|^p\, dm.
\end{equation}
Condition (ii) of Definition \ref{defPchain} gives that $Z\in MB_{i_Z}$. 
Another way to say this
is that $B_{i_Z} \subset B_Z:=B(Z,r_Z)$ for some $r_Z \approx r_{i_Z} \approx \kappa^{i_Z}r$. 
Moreover, due to (\hyperref[H4]{H4}), 
$m(B_{i_Z}) \approx m(B_Z \cap \Omega)$ and 
\[\kappa^{-i_Zd} \gtrsim \left( \frac{r}{r_Z} \right)^d 
\gtrsim \frac{m(B(Z,r) \cap \Omega)}{m(B_Z \cap \Omega)} 
\gtrsim \frac{m(D)}{m(B_Z \cap \Omega)}, 
\]
where $d$ is the exponent $d_m$ given in \eqref{mdoublingG}, and where we recall that $r:=\diam D$. 
We can freely assume that $\epsilon <1$,
and \eqref{PoincareA} becomes
\begin{equation} \label{lastestimate}
m(B_Z \cap \Omega)^{1+(\epsilon-1)p/d} m(D)^{(1-\epsilon)p/d} 
\lesssim  \left(\frac rt\right)^p \int_{B_Z \cap D} 
|\nabla u|^p\, dm.
\end{equation}
The balls $B_Z$, where $Z\in A_t$ is a Lebesgue point, 
cover almost all of $A_t$.
Hence
the Vitali covering lemma entails that there exists a collection of pairwise disjoints balls 
$B_{Z_j}$, $j\in J$,
 such that $A_t \subset \Omega \cap \big(\bigcup_{j\in J} 5B_{Z_j} \big)$
modulo a negligible set.
 We fix $\epsilon$ such that $1+(\epsilon-1)p/d = 1-p/(d+p) = d/(d+p)$, that is 
 $(\epsilon-1)p/d = -p/(d+p)$. Then 
\[\begin{split}
m(A_t)^{d/(d+p)} & \leq \Big[ \sum_{j\in J} m(5B_{Z_j} \cap \Omega)
\Big]^{d/(d+p)} \leq \sum_{j\in J} m(B_{Z_j} \cap \Omega)^{d/(d+p)} \\
& \lesssim m(D)^{p/(d+p)} \left(\frac rt\right)^p \sum_{j\in J}  \int_{B_{Z_j} \cap D} |\nabla u|^p\, dm \leq m(D)^{p/(d+p)} \left(\frac rt\right)^p \int_{D} |\nabla u|^p\, dm
\end{split}\]
by the covering property, because $d/(d+p) < 1$, then by \eqref{lastestimate}, 
our choice of $\varepsilon$,
because the exponent for $m(D)$ is $(1-\epsilon)p/d = -1 + [1+(1-\epsilon)p/d]
= 1 - d/(d+p) = p/(d+p)$, and finally because the $B_{Z_j}$ are disjoint.
Written differently, we proved that
\[
\frac{m(A_t)}{m(D)} 
\leq C \left(\frac rt\right)^{p(d+p)/d} \Big\{\int_{D} |\nabla u|^p\, dm\Big\}^{p(d+p)/d},
\]
or in other words $u$ lies in the weak Lebesgue space $L^{p(d+p)/d}_{w}(D)$.
We can use this and the Cavalieri formula to estimate $|| u ||_{L^q(D)}$ for any $q<p(d+p)/d$, 
and get that
\[\left( \fint_D |u-u_{B_0}|^q \, dm\right)^\frac1q 
= \left( \fint_D |u|^q \, dm\right)^\frac1q \leq C_q r \left( \fint_D |\nabla u|^p \, dm\right)^\frac1p;
\]
Theorem \ref{Poincare1} follows.
\ep

\begin{remark} \label{R5.32}
A careful inspection on the proof would show that we can prove 
\[\left( \frac1{m(D)}\int_D |u- \bar u|^{q} \, dm \right)^{1/q} \leq C_q \diam(D) \left( \frac1{m(D)} \int_D |\nabla u|^p \, dm \right)^{1/p} ,\] 
for every $q<+\infty$ if $p\geq d$ and every $q< \frac{pd}{d-p}$ if $p<d$. 
So in Theorem \ref{Poincare1}, if $2\geq d$, we can take for $k$ every positive value, and if $2 < d$, $k$ can take every value smaller than $d/(d-2)$.
\end{remark}

\section{The Trace Theorem}

\label{STrace}

As in the previous section, we assume that $\Gamma$ and $\Omega$ are infinite, but  
the results of this section, in particular Theorem \ref{TraceTh}, and Lemma \ref{lmult}, 
still hold when $\Gamma$ and/or $\Omega$ are finite. 
We shall discuss this again 
in Section \ref{SBounded}.

Let us first play 
a bit with the dyadic decomposition $\mathbb D$,
H\"older, and Fubini.

\begin{lemma}\label{HardyInequality}
Assume that $(\Omega,m,\mu)$ satisfies (\hyperref[H1]{H1})--(\hyperref[H3]{H3}) and let 
$q>1$. For $g\in L^q(\Omega,m)$,
\[ \sum_{Q\in \mathbb D} m(U^*_Q)^{1-q} \left( \int_{U^*_{Q}} g\, dm \right)^q \leq C \int_\Omega |g|^q \, dm,\]
where $C$ depends only on constants $C_1$ to $C_5$, $n$ and $q$.
\end{lemma}

\bp
The H\"older inequality implies that
for every $Q\in \mathbb D$,
\[\begin{split}
\left( \int_{U^*_Q} g\, dm \right)^q &  \leq  \left(\int_{U^*_Q} |g(Z)|^q  dm(Z) \right) m(U^*_Q)^{q-1}.
\end{split}\]
We sum over 
the dyadic cubes $Q$ to get that
\[\begin{split}
\sum_{Q\in \mathbb D} m(U^*_Q)^{1-q}  \left( \int_{U^*_{Q}} g\, dm \right)^q
& \lesssim \sum_{Q\in \mathbb D} \int_{U^*_Q} |g(Z)|^q \, dm(Z) \\
& \lesssim \int_{\Omega} |g(Z)|^q h(Z) \, dm(Z)\\
\end{split}\]
by Fubini's lemma, and where
\[h(Z) = \sum_{Q\in \mathbb D} \1_{U^*_Q}(Z).\]
The sets $U^*_Q$ are Whitney regions 
associated to the cubes $Q$, so
$Z \in U^*_Q$ implies that $\delta(Z) \approx \ell(Q) \approx \dist(Z,Q)$, and 
for each $Z$
there can be only a bounded 
number of such dyadic cubes in $\mathbb D$ (the number depends only on $n$, $C_1$, $C_2$, $C_3$). 
Hence
$h(Z) \lesssim 1$ and 
\[\sum_{Q\in \mathbb D} m(U^*_Q)^{1-q} \left( \int_{U^*_{Q}} g\, dm \right)^q \lesssim \int_{\Omega} |g(Z)|^q \, dm(Z).\]
The lemma follows.
\ep

We also need the following Hardy inequality.

\begin{lemma} \label{rmkHardyIneq}
Let $q>1$. Assume that $\{s_i\}_{i\in \mathbb Z}$ is a weight on $\mathbb Z$ that satisfies
\begin{equation} \label{condsi}
\frac{s_i}{s_j} \leq C_s2^{(j-i)\epsilon} \quad 
\text{ for } i > j,
\end{equation}
for some positive constants $C_s$ and $\epsilon$. 
Then, for $ [g_i]_{i\in \mathbb Z} \in \ell^q(\mathbb Z,s_i)$,
\[ \sum_{k\in \mathbb Z} s_k^{1-q} \left( \sum_{i>k} s_i g_i \right)^q 
\leq C \sum_{i\in \mathbb Z} s_i |g_i|^q ,\]
where $C$ depends only on $q$, $\epsilon$ and $C_s$.
\end{lemma}

\begin{remark} \label{rmkHardy2}
If $g_i = 0$ for $i > i_0$,
then we only need to require
\eqref{condsi}  for $i\leq i_0$.
\end{remark}

\bp Let $\alpha = \epsilon/2 >0$. Then by H\"older's inequality
\[\begin{split}
\left( \sum_{i>k} s_i g_i \right)^q 
& = \left( \sum_{i>k} 2^{-i\alpha} s_i2^{i\alpha}  g_i \right)^q \\
& \leq \left( \sum_{i>k}  2^{-i\alpha} |s_i2^{i\alpha} g_i|^q \right) \left(\sum_{i>k} 2^{-i\alpha} \right)^{q-1} \\
& \lesssim 2^{-k\alpha(q-1)} \sum_{i>k} 2^{-i\alpha} |s_i2^{i\alpha} g_i|^q
\end{split}\]
because $\alpha >0$. We sum in $k \in \mathbb Z$ 
and then apply Fubini's lemma
to get
\[ \begin{split}
\sum_{k\in \mathbb Z} s_k^{1-q} \left( \sum_{i>k} s_i g_i \right)^q 
& \lesssim \sum_{k\in \mathbb Z} (2^{k\alpha}s_k)^{1-q} \sum_{i>k} 2^{-i\alpha} |s_i2^{i\alpha} g_i|^q 
\\
& \lesssim \sum_{i\in \mathbb Z} 2^{-i\alpha} |s_i2^{i\alpha} g_i|^q  \sum_{k<i} (2^{k\alpha}s_k)^{1-q}.
\end{split}\]
By \eqref{condsi}, $2^{i\alpha}s_i = 2^{i(\alpha-\varepsilon)} 2^{i\varepsilon} s_i
\lesssim 2^{i(\alpha-\varepsilon)} 2^{k\varepsilon}s_k 
= 2^{(i-k)(\alpha-\varepsilon)} 2^{k\alpha}s_k$ for $k < i$; then
\[\sum_{k<i} (2^{k\alpha}s_k)^{1-q} 
\lesssim \sum_{k<i}(2^{i\alpha}s_i)^{1-q} 2^{(k-i)(\epsilon-\alpha)(q-1)}
\lesssim (2^{i\alpha}s_i)^{1-q},\]
because $q > 1$ and $\alpha < \epsilon$. This yields
\[\sum_{k\in \mathbb Z} s_k^{1-q} \left( \sum_{i>k} s_i g_i \right)^q \lesssim \sum_{i\in \mathbb Z} s_i |g_i|^q \] 
and the lemma follows.
\ep

The aim of the section is to show that the functions in $W$ have a trace, 
and that the traces lie in the space $H$ defined as
\begin{equation} \label{defH}
H:= \Big\{g: \Gamma \to \R \, ; \,  \text{$g$ is  $\mu$-measurable and} \, 
\int_\Gamma \int_\Gamma \frac{\rho(x,|x-y|)^2 |g(x) - g(y)|^2}{m(B(x,|x-y|) \cap \Omega)}  
d\mu(y) d\mu(x)  < +\infty \Big\},
\end{equation}
where $\rho$ is as in \eqref{defrho}.
The space $H$ is equipped with the semi-norm 
\[\|g\|_H := \left( \int_\Gamma \int_\Gamma 
\frac{\rho(x,|x-y|)^2 |g(x) - g(y)|^2}{m(B(x,|x-y|) \cap \Omega)} \, d\mu(y) d\mu(x)\right)^\frac12
\]
(adding a constant to $g$ keeps $g$ in $H$ and does not change $\|g\|_H$).

The existence of traces is given by the following result.
Recall the nontangential cones $\gamma(x)$, $x\in \Gamma$, from \eqref{defgamma}.

\begin{theorem} \label{TraceTh}
Assume that $(\Omega,m,\mu)$ satisfies (\hyperref[H1]{H1})--(\hyperref[H6]{H6}). There exists a bounded linear operator $\Tr: \, W \to H$ (a trace operator) with the following properties. The trace of $u\in W$ is such that, for $\mu$-almost every $x\in \Gamma$,
\begin{equation} \label{defTr}
\Tr u(x) = \lim_{X \in \gamma(x) \atop {\delta(X) \to 0}} \fint_{B(X,\delta(X)/2)} u \, dm
\end{equation}
and even, 
analogously to the Lebesgue density property,
\begin{equation} \label{defLebesgue}
\lim_{X \in \gamma(x) \atop {\delta(X) \to 0}} \fint_{B(X,\delta(X)/2)} |u(Z) - \Tr u(x)|  \, dm(Z) =  0.
\end{equation}
\end{theorem}

\bp For $x\in \Gamma$ and $k \in \mathbb Z$, we write $\Tr_k u(x)$ for any quantity
$$\Tr_k u(x) := \fint_{B(X,\delta(X)/2)} u(Z) \, dm(Z),$$
where 
$X$ is picked in $U_{Q^k(x)}$
and $Q^k(x)$ is the only set in $\mathbb D_k$ containing $x$. Keep in mind that 
$\Tr_k u(x)$ is not uniquely defined, but the estimates on $\Tr_k u$ that will be proven here hold with a constant independent of the choice of $X\in U_{Q^k(x)}$. For the rest of the proof, we also write $B^k_x$ for $B(X,\delta(X)/2)$ when $X\in U_{Q^k(x)}$. For any couple of integers $k<j$, one has
\[
\left|\Tr_j u(x) - \Tr_k u(x)\right| \leq \sum_{k < i \leq j} \left|\Tr_{i-1} u(x) - \Tr_{i} u(x)\right|
\]
and then for all $i \in \mathbb Z$,  since both $B^{i-1}_x$ and $B^{i}_x$ belong to 
$U^*_{Q^i(x)}$ by construction,
\[\begin{split}
\left|\Tr_{i-1} u(x) - \Tr_{i} u(x)\right| 
& \leq \left|\fint_{B^{i-1}_x} u \, dm - \fint_{U^*_{Q^i(x)}} u\, dm \right| 
+ \left|\fint_{B^{i}_x} u\, dm - \fint_{U^*_{Q^i(x)}} u\, dm \right| \\
& \leq \fint_{U^*_{Q^i(x)}} \left( \left| u(Z) -  \fint_{B^{i-1}_x} u\, dm \right| + \left| u(Z) -  \fint_{B^{i}_x} u\, dm \right| \right) dm(Z) \\
& \lesssim 2^{-i} \left(\fint_{U^*_{Q^i(x)}} |\nabla u(Z)|^p dm(Z)\right)^\frac1p,
\end{split}\]
where the last inequality and the parameter $p \in (1,2)$ are given by the Poincar\'e inequality 
(Theorem \ref{Poincare1}). The combination of the two
proves that
\begin{equation} \label{claimTrth1}
\left|\Tr_j u(x) - \Tr_k u(x)\right| \lesssim \sum_{i=k+1}^{j}  2^{-i} \left(\frac{1}{m(U^*_{Q^i(x)})} \int_{U^*_{Q^i(x)}} |\nabla u|^p dm\right)^\frac1p.
\end{equation}

Take $Q^* \in \D$
and write $k^*$ for $k(Q^*)$. Let us prove that $(\Tr_k u)_{k\geq k^*}$ is a Cauchy sequence in $L^2(Q^*,\mu)$. We integrate in $x$ to get
\begin{equation} \label{Trthc} 
\begin{split}
\int_{Q^*} \left|\Tr_j u - \Tr_k u\right|^2 d\mu & \lesssim \int_{Q^*} \left(  \sum_{i=k+1}^{j} 2^{-i} \left(\frac{1}{m(U^*_{Q^i(x)})} \int_{U^*_{Q^i(x)}} |\nabla u|^p dm\right)^\frac1p \right)^2 dx \\
& \hspace{-2cm} \lesssim \int_{Q^*} \left[\sum_{i=k+1}^{j}  \frac{2^{-i}\rho(x,2^{-i})}{m(U^*_{Q^i(x)})^{2/p}} \left( \int_{U^*_{Q^i(x)}} |\nabla u(Z)|^p dm(Z)\right)^\frac2p\right] \left[ \sum_{i=k+1}^{j} \frac{2^{-i}}{\rho(x,2^{-i})} \right]   \, d\mu(x),
\end{split} 
\end{equation}
where we use Cauchy-Schwarz's inequality for the last line, and where $\rho$ is the function defined in \eqref{defrho}. 

From now on, let $\varepsilon = C_5^{-1}$ denote the same small constant as in (\hyperref[H5]{H5}). Then
\eqref{HH5a} says that $\displaystyle\frac{2^{-i(1-\varepsilon)}}{\rho(x,2^{-i})} 
\lesssim \frac{2^{-k^*(1-\varepsilon)}}{\rho(x,2^{-k^*})}$, and hence
\[ \sum_{i=k+1}^{j} \frac{2^{-i}}{\rho(x,2^{-i})}  
\lesssim \frac{2^{-k^*}}{\rho(x,2^{-k^*})}\sum_{i=k+1}^{j} \left(\frac{2^{-k^*}}{2^{-i}}\right)^{-\epsilon} 
\lesssim \frac{2^{-k^*}}{\rho(x,2^{-k^*})} \, 2^{-(k-k^*)\epsilon}.
\]
Moreover, thanks to \eqref{defrho} and 
Lemma \ref{lemmQ}, 
\begin{equation} \label{rhoUQ}
\rho(x,2^{-i}) = \frac{m(B(x,2^{-i})\cap \Omega)}{2^{-i}\mu(B(x,2^{-i}))} \approx \frac{m(U^*_{Q^i(x)})}{2^{-i} \mu(Q^i(x))}.
\end{equation}
We deduce from \eqref{Trthc} and the two last estimates
(including \eqref{rhoUQ} for $k^\ast$ for the second line) that
\[ 
\begin{split}
\int_{Q^*}  |\Tr_j u & - \Tr_k u|^2 d\mu  \\
& \lesssim 
\frac{2^{-k^*}}{\rho(x,2^{-k^*})} \, 2^{-(k-k^*)\epsilon} 
\int_{Q^*} \sum_{i=k+1}^{j}  \frac{m(U^*_{Q^i(x)})^{1-2/p}}{\mu(Q^i(x))} 
\left( \int_{U^*_{Q^i(x)}} |\nabla u(Z)|^p dm(Z)\right)^\frac2p \, d\mu(x)\\ 
& \lesssim \frac{2^{-2k^*} \mu(Q^*)}{m(U^*_{Q^*})} \, 2^{-(k-k^*)\epsilon} 
\int_{Q^*} \sum_{i=k+1}^{j}  \frac{m(U^*_{Q^i(x)})^{1-2/p}}{\mu(Q^i(x))} 
\left( \int_{U^*_{Q^i(x)}} |\nabla u(Z)|^p dm(Z)\right)^\frac2p \, d\mu(x)\\ 
& \quad =  \frac{2^{-2k^*} \mu(Q^*)}{m(U^*_{Q^*})} 2^{-(k-k^*)\epsilon} 
\sum_{i>k} \sum_{Q\in \D_i \atop Q \subset Q^*} m(U^*_Q)^{1-2/p} 
\left( \int_{U^*_Q} |\nabla u|^p \, dm \right)^\frac2p,
\end{split}
\]
where for the last line we decomposed $Q^\ast$ into cubes $Q = Q^i(x)$ for each $i$.
We write $\Omega_k$ for
\[\Omega_{Q^*,k} := 
\bigcup_{x\in Q^*} \gamma^*_{Q^k(x)} 
= \bigcup_{Q \subset Q^* \atop k(Q) \leq k} U^*_{Q},\]
(see \eqref{defgammaQ*}), so that the difference of traces is bounded by
\[ \int_{Q^*} \left|\Tr_j u - \Tr_k u\right|^2 d\mu \lesssim\frac{2^{-2k^*} \mu(Q^*)}{m(U^*_{Q^*})} 2^{-(k-k^*)\epsilon} \sum_{Q\in \D} m(U^*_Q)^{1-2/p} \left( \int_{U^*_Q} \1_{\Omega_{Q^*,k}}|\nabla u|^p \, dm \right)^\frac2p.\]
Now we use lemma \ref{HardyInequality} with $q=2/p$ and $g = \1_{\Omega_{Q^*,k}}|\nabla u|^p$ to obtain that for $j > k > k^\ast$
\begin{equation} \label{Trthd}
\int_{Q^*} \left|\Tr_j u - \Tr_k u\right|^2 d\mu 
\lesssim \frac{2^{-2k^*} \mu(Q^*)}{m(U^*_{Q^*})} 2^{-(k-k^*)\epsilon} \int_{\Omega_{Q^*,k}} |\nabla u|^2 dm.
\end{equation}

The last result is pretty nice, and 
just keeping the information that for each $Q^\ast$, 
$||\Tr_j u - \Tr_k u||^2_{L^2(Q^\ast,d\mu)} \leq C(u,Q^\ast) 2^{-k\varepsilon}$
for $j > k > k^\ast$, we get that the series $\sum_{k} \Tr_{k+1} u - \Tr_k u$ converges
normally in every $L^2(Q^\ast,d\mu)$, hence in $L^2_{loc}(\Gamma,\mu)$ and
$\mu$-almost everywhere. That is,
$$\Tr u(x) = \lim_{k\to +\infty} \Tr_k u(x) =  \lim_{X \in \gamma(x) \atop {\delta(X) \to 0}} \fint_{B(X,\delta(X)/2)} u(Z) \, dm(Z)$$
exists for $\mu$-almost every $x\in \Gamma$ and by \eqref{Trthd}
\begin{equation} \label{Trthe}
\int_{Q^*} \left|\Tr_k u - \Tr u\right|^2 d\mu \lesssim \frac{2^{-2k^*} \mu(Q^*)}{m(U^*_{Q^*})} 2^{-(k-k^*)\epsilon} \int_{\Omega_{Q^*,k}} |\nabla u|^2 dm.
\end{equation}

The estimate \eqref{Trthe} is not strong enough to imply the Lebesgue density property \eqref{defLebesgue}. However, observe that for $k<j$ and $X\in U_{Q_k(x)}$,
\[\fint_{B_x^k} |u - \Tr_j u(x)| \, dm \leq \fint_{B_x^k} |u - \Tr_k u(x)| \, dm + |\Tr_j u(x) - \Tr_k u(x)| \]
and, thanks to Lemma \ref{Poincare00} (improved Poincar\'e inequality) and the fact that $m(U^*_{Q^k_x}) \approx m(B_x^k)$ (by (\hyperref[H4]{H4}) and Lemma \ref{lemmQ}), 
\[ \fint_{B_x^k} |u - \Tr_k u(x)| \, dm  \lesssim 2^{-k} \left(\fint_{U^*_{Q^k(x)}} |\nabla u|^p \, dm\right)^\frac1p.\]
Together with \eqref{claimTrth1}, this implies that 
\[ \fint_{B_x^k} |u - \Tr_j u(x)| \, dm  \lesssim
\sum_{i=k}^{j}  
2^{-i} \left(\frac{1}{m(U^*_{Q^i(x)})} \int_{U^*_{Q^i(x)}} |\nabla u|^p dm\right)^\frac1p. \]
We integrate on $x\in \Gamma$ and 
invoke Cauchy-Schwarz's inequality and Lemma \ref{HardyInequality} to get, 
analogously to \eqref{Trthd},
\[
\int_{x\in Q^*} 
\left( 
\fint_{B_x^k}  
|u - \Tr_j u(x)| \, dm \right)^2 d\mu(x) 
\lesssim  \frac{2^{-2k^*} \mu(Q^*)}{m(U^*_{Q^*})} 2^{-(k-k^*)\epsilon} \int_{\Omega_{Q^*,k}} |\nabla u|^2 dm.\]
The right-hand side does not depend on $j$; we take 
the limit as $j$ approaches $+\infty$ and 
obtain
\begin{equation}
\int_{x\in Q^*} 
\left(\fint_{B_x^k} |u - \Tr u(x)| \, dm \right)^2 
d\mu(x)
\lesssim  \frac{2^{-2k^*} \mu(Q^*)}{m(U^*_{Q^*})} 2^{-(k-k^*)\epsilon} \int_{\Omega_{Q^*,k}} |\nabla u|^2 dm.
\end{equation}
It follows that $x\to 
\fint_{B_x^k}
|u - \Tr u(x)| \, dm$ converges to 0 in $L^2_{loc}(\Gamma,\mu)$ as $k\to +\infty$, and this
implies the Lebesgue density property \eqref{defLebesgue}.

\medskip

It remains to prove that $\Tr$ is a bounded operator from $W$ to $H$. 
If $x,y\in \Gamma$, we write $k(x,y)$ for the only integer $k$ that satisfies 
$2^{-k-1} \leq |x-y| < 2^{-k}$. We use
\eqref{defH} and decompose the integral as
\[\begin{split}
\|\Tr u\|_H & := \int_\Gamma \int_\Gamma  \dfrac{\rho(x,|x-y|)^2|\Tr u(x)-\Tr u(y)|^2}{m(B(x,|x-y|) \cap \Omega)} \, d\mu(y)\, d\mu(x) \\
 & \hspace{0cm} \lesssim \int_\Gamma\int_\Gamma \dfrac{\rho(x,|x-y|)^2 |\Tr u(x)-\Tr_{k(x,y)} u(x)|^2}{m(B(x,|x-y|) \cap \Omega)} \, d\mu(y)\, d\mu(x) \\
  & \hspace{1cm} + \int_\Gamma \int_\Gamma \dfrac{\rho(x,|x-y|)^2|\Tr u(y)-\Tr_{k(x,y)} u(y)|^2}{m(B(x,|x-y|) \cap \Omega)} \, d\mu(y)\, d\mu(x) \\
  & \hspace{2cm}  + \int_\Gamma \int_\Gamma \dfrac{\rho(x,|x-y|)^2 |\Tr_{k(x,y)} u(x) -\Tr_{k(x,y)} u(y)|^2}{m(B(x,|x-y|) \cap \Omega)} \, d\mu(y)\, d\mu(x) \\
 & \hspace{1cm}:= I_1 + I_2 + I_3.
\end{split}\]

Let us first treat the term $I_3$.
More precisely, 
we 
start with the difference $|\Tr_{k(x,y)} u(x) -\Tr_{k(x,y)} u(y)|$. To lighten the notation, 
write $k$ for $k(x,y)$. As before, denote by 
$B_x^k$ and $B_y^k$  
the balls used to define 
$\Tr_k u(x)$ and $\Tr_k u(y)$.
That is, 
$B_x^k$ and $B_y^k$ are such that
$$\Tr_k u(x) = \fint_{B_x^k} u \quad \text{ and } \quad  \Tr_k u(y) = \fint_{B_y^k} u.$$
Since the balls $B_x^k,B_y^k$ 
lie at distances at least $c2^{-k}$ from the boundary $\Gamma$ and at most 
$C2^{-k}$ from each other,
the Harnack chain condition (Proposition \ref{propHarnack}) says that we can find a chain of balls joining $B_x^k$ to $B_y^k$, with 
uniformly bounded 
length, staying at a distance at least $c2^{-k}$ from the boundary and  
at distance at most $C2^{-k}$ from both $x$ and $y$. We define $P_{x,y}^k$ as the union of the cubes in $\mathcal W$ that meet one of the balls of the chain. From what we just said, $P_{x,y}^k$ is a Whitney region associated to both $(x,2^{-k})$ and $(y,2^{-k})$, and so it
is the union of a bounded 
number of adjacent cubes in $\cW$. Therefore, similarly to the sets $U^*_Q$, $P_{x,y}^k$ satisfies the $C(\kappa,M)$ chain condition for some uniform $\kappa,M$ and is thus fitted for the Poincar\'e inequality.
These observations allow us to use Theorem \ref{Poincare1} and write
\[\begin{split}
|\Tr_{k} u(x) -\Tr_{k} u(y)| & := \left| \fint_{B_x^k} u\, dm - \fint_{B_y^k} u\, dm  \right| 
\lesssim \fint_{B_x^k} \Big|u - \fint_{B_y^k} u\, dm \Big|\, dm \\
& \lesssim \frac{1}{m(B_x^k)}  \int_{P_{x,y}^k} \Big|u - \fint_{B_y^r} u \, dm \Big| dm \\
& \lesssim \frac{2^{-k} m(P_{x,y}^k)^{\frac12}}{m(B_x^k)}  \left(\int_{P_{x,y}^k} \left|\nabla u \right|^2 dm\right)^\frac12.
\end{split}\]
Since both $P^k_{x,y}$ and $B_x^k$ are Whitney region associated to $y$ and $2^{-k}$ (i.e. there exists a large constant $C$ such that both sets are contained in $B(y,C2^{-k})$, contain a ball $B$ of radius $C^{-1}2^{-k}$, and are at distance at least $C^{-1}2^{-k}$ of $\Gamma$), the doubling measure condition (\hyperref[H4]{H4}) implies that $m(B_x^k) \approx m(P_{x,y}^k) \approx m(B(y,2^{-k-1}) \cap \Omega)$. Therefore,
\begin{equation}\label{Trthg}
|\Tr_{k} u(x) -\Tr_{k} u(y)|^2 \lesssim \frac{2^{-2k}}{m(B(y,2^{-k-1}) \cap \Omega)}   \int_{P_{x,y}^k} \left|\nabla u\right|^2 dm.
\end{equation}
We inject \eqref{Trthg} in $I_3$ and
observe that $m(B(x,|x-y|)\cap \Omega) \approx m(B(x,2^{-k(x,y)})\cap \Omega)$ and $\rho(x,|x-y|) \approx \rho(x,2^{-k(x,y)}) \approx \rho(y,2^{-k(x,y)})$ by (\hyperref[H3]{H3})--(\hyperref[H4]{H4}). 
Therefore, 
\begin{equation}\label{Trthh} \begin{split}
I_3 & \lesssim \int_\Gamma \int_\Gamma 
\frac{2^{-k(x,y)}\rho(x,2^{-k(x,y)})}{m(B(x,2^{-k(x,y)})\cap \Omega)} 
\frac{2^{-k(x,y)}\rho(y,2^{-k(x,y)})}{m(B(y,2^{-k(x,y)}) \cap \Omega)} \int_{P_{x,y}^{k(x,y)}} \left|\nabla u\right|^2 dm \, d\mu(x) \, d\mu(y) \\
& \quad = \int_\Gamma \int_\Gamma \frac{1}{\mu(B(x,2^{-k(x,y)}))} \frac{1}{\mu(B(y,2^{-k(x,y)}))} \int_{P_{x,y}^{k(x,y)}} \left|\nabla u\right|^2 dm \, d\mu(x) \, d\mu(y),
\end{split} \end{equation}
by the definition \eqref{defrho} of $\rho$.
Since $P_{x,y}^{k(x,y)}$ is a `Whitney region' for both $x$ and $y$,  we have that 
$x,y \in B(Z,C\delta(Z))$ for $Z\in P^{k(x,y)}_{x,y}$ 
where the constant $C \geq 2$ depends only on $n,C_1,C_2$,
and moreover 
$2^{-k(x,y)} \approx \delta(Z)$.  
Then by Fubini's lemma
\begin{equation}\label{Trthi0} \begin{split}
I_3&  \lesssim 
\int_\Omega |\nabla u(Z)|^2 dm(Z) \,
 \int_{x\in B(Z,C\delta(Z))} \frac{d\mu(x)}{\mu(B(x,c\delta(Z)))}  \int_{y\in B(Z,C\delta(Z))} \frac{d\mu(y)}{\mu(B(y,c\delta(Z)))}.
\end{split} \end{equation}
Yet
the doubling property (\hyperref[H3]{H3}) implies that for $z\in \Gamma \cap B(Z,C\delta(Z))$,
\[\mu(B(z,c\delta(Z))) \gtrsim \mu(B(z,C^2\delta(Z))) \geq \mu(B(Z,C\delta(Z))),\]
hence we can simply bound $I_3$ by
\begin{equation}\label{Trthi}
I_3 \lesssim \int_\Omega |\nabla u(Z)|^2 \, dm(Z)
\end{equation}
as desired.

\medskip

We turn now to the bound on $I_1$. 
Notice that the estimate for $I_2$ is the same as for $I_1$, either by symmetry or
since by (\hyperref[H3]{H3})--(\hyperref[H4]{H4}), $m(B(x,|x-y|) \cap \Omega) \approx \mu(B(y,|x-y|) \cap \Omega)$ and $\rho(x,|x-y|) \approx \rho(y,|x-y|)$.

Notice that $I_1$ depends on $y$ only via $|x-y|$, 
so, by the doubling property (\hyperref[H3]{H3}) again and then \eqref{defrho}
\begin{equation} \label{I1firstbd}
\begin{split}
I_1 & \lesssim \int_{x\in \Gamma} \sum_{k\in \mathbb Z} \frac{\rho(x,2^{-k})^2 |\Tr u(x) - \Tr_k u(x)|^2}{m(B(x,2^{-k}) \cap \Omega)}  \int_{y \in B(x,2^{-k})\setminus B(x,2^{-k-1})} d\mu(y) \, d\mu(x) \\
& \lesssim \int_{x\in \Gamma} 
\sum_{k\in \mathbb Z}
2^k \rho(x,2^{-k}) |\Tr u(x) - \Tr_k u(x)|^2 \, d\mu(x).
\end{split}\end{equation}
The trace operator is defined for $\mu$-almost every $x\in \Gamma$ by \eqref{defTr}. For such $x$, one get by letting $j$ tend to $+\infty$
in \eqref{claimTrth1} that 
\begin{equation} \label{Trthj}
\left|\Tr_k u(x) - \Tr u(x)\right| \lesssim \sum_{i>k} 2^{-i} \left( \frac1{m(U^*_{Q^i(x)})}\int_{U^*_{Q^i(x)}} |\nabla u|^p dm \right)^\frac1p.
\end{equation}
We use the above estimate in \eqref{I1firstbd} to obtain that
\[\begin{split}
I_1 & \lesssim \int_{x\in \Gamma} 
\sum_{k \in \mathbb Z}
2^k \rho(x,2^{-k}) \left( \sum_{i>k}  2^{-i} \left( \frac1{m(U^*_{Q^i(x)})}\int_{U^*_{Q^i(x)}} |\nabla u|^p dm \right)^\frac1p \right)^2 \, d\mu(x) \\
& \qquad = \int_{x\in \Gamma} 
\sum_{k \in \mathbb Z}
2^k \rho(x,2^{-k}) \left( \sum_{i>k}  \frac1{2^{i} \rho(x,2^{-i})} \, g_i(x)
\right)^2 \, d\mu(x), 
\end{split}\]
where 
\[g_i(x) = \rho(x,2^{-i}) \left( \frac1{m(U^*_{Q^i(x)})} \int_{U^*_{Q^i(x)}} |\nabla u|^p dm\right)^\frac1p.\] 
Thanks to (\hyperref[H5]{H5}), the sequence $\{s_i\}_{i \geq k_0}$ defined as 
$s_i := [2^i\rho(x,2^{-i})]^{-1}$ satisfies \eqref{condsi}.
As a consequence, Lemma \ref{rmkHardyIneq} 
(with $q=2$) gives that for each $x\in \Gamma$,
\[
\sum_{k \in \mathbb Z}
2^k \rho(x,2^{-k}) \left( \sum_{i>k}  
\frac1{2^{i} \rho(x,2^{-i})} \, g_i(x) \right)^2 
\lesssim 
\sum_{i \in \mathbb Z}
\frac1{2^i \rho(x,2^{-i})} \, |g_i(x)|^2.\] 
Thus the bound on $I_1$ becomes
\[\begin{split}
I_1 & \lesssim \int_\Gamma 
\sum_{i \in \mathbb Z}
\frac1{2^i \rho(x,2^{-i})} |g_i(x)|^2 \, d\mu(x)\\
& \qquad = \int_\Gamma 
\sum_{i \in \mathbb Z}
2^{-i} \rho(x,2^{-i}) \left( \frac1{m(U^*_{Q^i(x)})} \int_{U^*_{Q^i(x)}} |\nabla u|^p dm\right)^\frac2p d\mu(x).
\end{split}\]
We use \eqref{rhoUQ} to get rid of the function $\rho$, and then write the bound we obtained 
as a sum over
$\D$, which gives
\[\begin{split}
I_1 & \lesssim \int_\Gamma 
\sum_{i \in \mathbb Z}
\frac{m(U^*_{Q^i(x)})^{1-2/p}}{\mu(Q^i(x))} \left(\int_{U^*_{Q^i(x)}} |\nabla u|^p dm\right)^\frac2p d\mu(x) 
= \sum_{Q\in \D} m(U^*_Q)^{1-2/p}  \left(\int_{U^*_{Q}} |\nabla u|^p dm\right)^\frac2p.
\end{split}\]
We can now apply Lemma \ref{HardyInequality} 
with
$q=2/p$ and $g = |\nabla u|^p$. Recall recall that $q > 1$ because $p$ comes from 
Theorem \ref{Poincare1} and was chosen above \eqref{claimTrth1}, with $1 < p < 2$.
We get that 
\[I_1 \lesssim \int_\Omega |\nabla u|^2 \, dm ; \] 
Theorem \ref{TraceTh} follows.
\ep

We end this section with a useful result concerning the trace of a product.

\begin{lemma} \label{lmult}
Let $(\Omega,m,\mu)$ satisfy (\hyperref[H1]{H1})--(\hyperref[H6]{H6}). 
Suppose $u\in W$ and $\varphi \in C_0^\infty(\R^n)$. 
Then $u\varphi \in W$, with the product rule 
\begin{equation} \label{lem6.40z}
\nabla (u \varphi) = \varphi \nabla u + u \nabla \varphi
\end{equation}
for the gradient, and 
\begin{equation} \label{lem6.40y}
\Tr (u\varphi)(x) = \varphi(x) \Tr u(x)
\quad \hbox{for every point $x\in \Gamma$ satisfying \eqref{defLebesgue}.}
\end{equation}
\end{lemma}

\bp This result is the analogue of \cite[Lemma 5.4]{DFMprelim}. 
The proof is similar, so we only sketch it. 

We start with the simplest case is when we can see $u$ as a distribution on $\Omega$; 
this is the case when the stronger form (\hyperref[H7]{H6'}) of our assumption (\hyperref[H6]{H6}) holds. 
Then $\nabla(u\varphi) = \varphi \nabla u + u \nabla \varphi$ in the sense of distributions,
and we are about to check that $\nabla(u\varphi) \in L^2(\Omega, dm)$.

Choose $Q\in \mathbb D$ so large that $\supp\, \varphi \cap \Omega \subset T_{2Q}$.
Then, setting $\bar u = \fint_{U^*_Q} u\, dm$,
\begin{equation} \label{e6.21}
\|\nabla(u\varphi)\|_{L^2(\Omega,dm)}  
\leq \|\varphi\|_\infty \|u\|_{W} + \|\nabla \varphi\|_\infty \left(\int_{T_{2Q}} |u-\bar u|^2 \, dm \right)^\frac12 + \|\nabla \varphi\|_\infty m(T_{2Q})^{\frac12} |\bar u|. 
\end{equation}
All three terms in the right hand side are finite, since $\varphi$ is smooth, $u\in W \subset L^1(U^*_Q,m)$, and by Theorem \ref{Poincare1} (Poincar\'e's inequality). Consequently, $u\varphi \in W$ as desired.

As for the trace, observe that if $x$ satisfies \eqref{defLebesgue}, and if $B_x^k$ denotes $B(X,\delta(X)/2)$ for some $X \in U_{Q_k(x)}$, where 
as usual $Q_k(x)$ is the only cube of $\mathbb D_k$ that contains $x$, then
\[\fint_{B_x^k} |\varphi u - \varphi(x)\Tr u(x)| \, dm \leq \|\varphi\|_\infty \fint_{B_x^k} |u -\Tr u(x)| \, dm + |Tu(x)| \fint_{B_x^k} |\varphi -\varphi(x)| \, dm.\]
The first term converges to 0 as $k\to +\infty$ since $x$ is a Lebesgue point for $u$, and the second term also tends to 0 since $|Tu(x)|<+\infty$ ($x$ is a Lebesgue point) and $\varphi$ is continuous. Therefore
\[\lim_{\delta(X) \to 0, \atop X \in \gamma(x)} \fint_{B_x^k} |\varphi u - \varphi(x)\Tr u(x)| \, dm = 0,\]
which easily implies $\Tr (u\varphi)(x) = \varphi(x) \Tr u(x)$ by the definition of $\Tr$. 

This takes care of the lemma when $u$ and $\nabla u$ are taken as distributions.
In general, we used Definition \ref{defW} to define the space $W$ and the gradient $\nabla u$. 
Recall that we wrote $u$ as a limit in $L^1_{loc}(\Omega,m)$ of smooth functions 
$\varphi_i \in C^\infty(\ol\Omega)$, as in \eqref{defW1}, and required that \eqref{defW9}
and \eqref{defW2} hold for a suitable $v \in L^2(\Omega,dm)$ which is unique by (\hyperref[H6]{H6})
and which we also called $\nabla u$.

Now we want to show that $u \varphi \in W$, so we approximate it by the 
smooth functions $\varphi_i \varphi$. It is easy to see that the $\varphi_i \varphi$ 
converge to $u\varphi$ in $L^1_{loc}(\Omega,m)$ as in \eqref{defW1}, 
and that $\int_\Omega |\nabla(\varphi \varphi_i)|^2 < +\infty$
for every $i$, as in \eqref{defW9} (recall that $\varphi \in C_0^\infty(\R^n)$).
We try the gradient $w = \varphi v + u\nabla \varphi = \varphi \nabla u + u\nabla\varphi$
in the definition \eqref{defW2}. First observe that $w \in L^2(\Omega,dm)$ by the proof
of \eqref{e6.21} (and where Theorem \ref{Poincare1} is applied in the general context of $W$).
We claim that
\[
\lim_{i \to \infty} \int_\Omega |\nabla(\varphi \varphi_i) - w|^2 dm = 0,
\]
as needed for \eqref{defW2}. The first part of $\nabla (\varphi \varphi_i)$
is $\varphi \nabla\varphi_i$, which converges to $\varphi v$ in $L^2(\Omega,dm)$,
by \eqref{defW2}. Thus we are left with showing that $\varphi_i \nabla \varphi$
converges to $u \nabla \varphi$ in $L^2(\Omega,dm)$. Or, since $\varphi$ is bounded
and $\supp\, \varphi \cap \Omega \subset T_{2Q}$, that 
\begin{equation} \label{e6.22}
\lim_{i \to \infty} \int_{T_{2Q}} |\varphi_i - u|^2 dm = 0.
\end{equation}
Denote by $c_i$ the average of $\varphi_i - u$ on $2Q$. Then by Poincar\'e' inequality
(Theorem \ref{Poincare1}) and \eqref{defW2}, $\int_{T_{2Q}} |\varphi_i - u -c_i|^2 dm$ tends to $0$.
But also \eqref{defW1} says that $\int_{B} |\varphi_i - u| dm$ tends to $0$ for some small
ball $B \subset T_{2Q}$, so in fact $c_i$ tends to $0$, \eqref{e6.22} holds, and $u\varphi \in W$
with a derivative equal to $w$. The remaining estimates are as in the easier case, and 
Lemma \ref{lmult} follows.
\ep

\section{Poincar\'e inequalities on the boundary}

We are interested in 
a version of the Poincar\'e inequality for functions that 
have a vanishing trace at the boundary. 
The proofs shall use the tent sets $T_{2Q}$ that were constructed in Section \ref{Scones}, 
where we assumed $\Gamma$ and $\Omega$ unbounded. But as explained in Section \ref{SBounded}, 
the same construction works for bounded $\Gamma$ and/or $\Omega$ (with maybe a restriction on the size of possible $Q$), and the proofs in the section are directly adaptable to this case. The Poincar\'e inequalities that we prove here are a local results, so it makes sense anyway that they don't depend on the boundedness of $\Omega$.

\begin{theorem} \label{PoincareTh2} Let $(\Omega,m,\mu)$ satisfy (\hyperref[H1]{H1})--(\hyperref[H6]{H6}). 
There exists $p_1\in [1,2)$ and $k:= k(C_4)>1$ such that for $p\in [p_1,2]$, 
$Q\in \D$,  
and $u\in W$ such that
$\Tr u = 0$ on a set
$E \subset 2Q$ such that
$\mu(E) \geq c \mu(2Q)$, we have
\begin{equation} \label{PoincareTh2a}
\left( \frac1{m(T_{2Q})} \int_{T_{2Q}} |u|^{kp} \, dm \right)^{1/kp} \leq C \ell(Q) \left( \frac1{m(T_{2Q})} \int_{T_{2Q}} |\nabla u|^p \, dm \right)^{1/p},
\end{equation}
where $T_{2Q}$ is the same tent set over $2Q$ as in \eqref{defTQ}, and 
$C>0$ depends only on $n$, the constants $C_1$-$C_6$, and $c$.
\end{theorem}

\bp Take $z\in 2Q$. 
We start as in 
the proof of Lemma \ref{lemchain}, and 
for each $j\in \bN$ 
we define $Q_j^z$ as the dyadic cube in $\mathbb D_{j+k(Q)}$ that contains $z$. 
Let $X_j^z$ be
a Corkscrew point associated to $z$ and the scale
$\ell(Q^z_j)=2^{-j}\ell(Q)$; by construction of $\mathcal W_{Q_j^z}$, we can find a cube $I^{j,z} \in \mathcal W_{Q_j^z}$ containing $X_j^z$, and we 
denote by $Y_j^z$
the center of $I^{j,z}$. By construction of $\mathcal W_{Q_j^z}^*$,
\[B_j^z:= B(Y_j^z,\delta(Y_j^z)/2) \subset U^*_Q.\]
By Proposition \ref{propHarnack}, we can find a uniform integer $N=N(n,C_1,C_2)$ such that we can link $B_j^{z}$ to $B^{z}_{j+1}$ by a Harnack chain of length $N$. We construct 
a chain
of balls $\{\mathcal B_i^z\}$ as follows: for $j\in \mathbb N$, $\{\mathcal B_i^z\}_{jN \leq i \leq (j+1)N}$ is the chain linking $B_j^{z}$ to $B^{z}_{j+1}$ given by Proposition \ref{propHarnack} and used to built $\mathcal W^*_{Q^z_{j+1}}$.
The collection of balls $(\mathcal B_i^z)_{i\geq 0}$ 
that we just constructed has bounded overlap, is 
included in $\gamma^*_Q(z)$, and is such that 
$\diam \, \mathcal B_i^z \approx 2^{-i/N}\ell(Q)$; observe also that by construction of $\cW_Q$, 
we can choose $I^{0,z}$ (and thus 
$B_0^z = \mathcal B_0^z$) 
independent of $z$, hence we write $B_0$ for $B_0^z$.

\medskip

For any subset $S\subset T_{2Q}$, we write as before $u_{S}$ for $\fint_{S} u\, dm$. 
Theorem \ref{PoincareTh2} will follow from Theorem \ref{Poincare1} as soon as we prove 
that for some $p_1\in [1,2)$,
\begin{equation} \label{PoincareB}
|u_{B_0}| \lesssim \ell(Q) \left( \frac1{m(T_{2Q})} 
\int_{T_{2Q}} |\nabla u|^p \, dm \right)^{1/p} 
\end{equation}
holds for all $p \in [p_1,2]$.

\medskip

Let $q\in [p_0,2]$, where $p_0$ is the value provided by Lemma \ref{Poincare00}. 
Thanks to Theorem \ref{TraceTh}, for $\mu$-almost every $z\in E$, we have 
\[\lim_{j\to \infty} |u_{\mathcal B^z_{jN}}| = 0.\] 
In particular, 
\[\begin{split} 
|u_{B_0}| & \leq \lim_{j\to +\infty}
\big\{|u_{\mathcal B^z_{jN}}| + |u_{B_0}-u_{\mathcal B^z_{jN}}| \big\}
\leq \lim_{j\to \infty} |u_{\mathcal B^z_{jN}}| 
+ \sum_{0\leq i <jN} |u_{\mathcal B^z_{i}}-u_{\mathcal B^z_{i+1}}|  
\leq \sum_{i\in \bN} |u_{\mathcal B^z_{i}}-u_{\mathcal B^z_{i+1}}|.
\end{split}\]
Since $\mathcal B^z_{i} \cap \mathcal B^z_{i+1}$, $\mathcal B^z_{i}$, and $\mathcal B^z_{i+1}$ 
have comparable sizes, the Poincar\'e inequality (Lemma \ref{Poincare00}) gives that
\[\begin{split}
 |u_{\mathcal B^z_{i}}-u_{\mathcal B^z_{i+1}}| & \leq  |u_{\mathcal B^z_{i}}-u_{\mathcal B^z_{i+1}
 \cap \mathcal B^z_i}| + |u_{\mathcal B^z_{i+1}}-u_{\mathcal B^z_{i+1}\cap \mathcal B^z_i}| \\
 & \leq \fint_{\mathcal B^z_{i+1}\cap \mathcal B^z_i} \left( |u - u_{B^z_{i}}| 
 + |u - u_{\mathcal B^z_{i+1}}| \right) \, dm \\
 & \lesssim \fint_{\mathcal B^z_{i}} |u - u_{\mathcal B^z_{i}}| \, dm 
 + \fint_{\mathcal B^z_{i+1}} |u - u_{\mathcal B^z_{i+1}}| \, dm\\
 & \lesssim   2^{-i/N} \ell(Q) \left( \fint_{\mathcal B^z_{i}} |\nabla u|^q \, dm 
 + \fint_{\mathcal B^z_{i+1}} |\nabla u|^q \, dm \right)^\frac1q.\\
\end{split}\]
So the last two computations yield that, for $\mu$-almost every $z\in E$
\[\begin{split}
|u_{B_0}| & \lesssim \sum_{i\in \bN} 2^{-i/N} \ell(Q) 
\left(\fint_{\mathcal  B^z_{i}} |\nabla u|^q dm\right)^\frac1q  
\lesssim \ell(Q)  
\left(\sum_{i\in \bN}  2^{-iq \alpha/N}
\fint_{\mathcal  B^z_i} |\nabla u|^q dm\right)^\frac1q
\end{split}\]
where we applied H\"older's inequality for the last part, and the price to pay is that
we need to introduce $\alpha\in (0,1)$, close to 1, that will be fixed later on. Now
observe that, for 
$Z\in \mathcal B^z_i$, 
we have $\delta(Z) \approx  2^{-i/N}\ell(Q)$ and 
hence $2^{-iq\alpha /N} \approx \delta(Z)^{q\alpha} \ell(Q)^{-q\alpha}$.
So, by (\hyperref[H4]{H4}),
\[m(\mathcal B_i^z) \approx 
m(B(Z,2\delta(Z))) \quad 
\text{ for } 
Z\in \mathcal B_i^z.\] 
Hence
\[\begin{split}
|u_{B_0}| &\lesssim \ell(Q)^{1-\alpha}\left(\sum_{i\in \bN} \int_{\mathcal B^z_i} |\nabla u(Z)|^q \, \frac{\delta(Z)^{q\alpha}}{m(B(Z,2\delta(Z)))} \, dm(Z) \right)^\frac1q \\
& \lesssim \ell(Q)^{1-\alpha} \left(\int_{\gamma_{Q}^*(z)} |\nabla u(Z)|^q \, \frac{\delta(Z)^{q\alpha}}{m(B(Z,2\delta(Z)))} \, dm(Z) \right)^\frac1q \\
\end{split}\]
since the balls $\mathcal B_i^z$ have bounded overlap and are contained in $\gamma_{Q}^*(z)$
(see the beginning of the proof, slightly above \eqref{PoincareB})
We 
average over $z\in E$ this estimate for $|u_{B_0}|^q$ to obtain
\[\begin{split}
|u_{B_0}| & = \left(\fint_{z\in E} |u_{B_0}|^q d\mu(z)\right)^\frac1q \\
& \lesssim 
\ell(Q)^{1-\alpha}
\left(\fint_E \int_{\gamma_{Q}^*(z)} |\nabla u(Z)|^q 
\, \frac{\delta(Z)^{q\alpha}}{m(B(Z,2\delta(Z)))} \, dm(Z)\, d\mu(z)\right)^\frac1q \\
& \lesssim \ell(Q)^{1-\alpha} \left(\fint_{z\in 2Q} \int_{\gamma_{Q}^*(z)} |\nabla u(Z)|^q \, \frac{\delta(Z)^{q\alpha}}{m(B(Z,2\delta(Z)))} \, dm(Z)\, d\mu(z)\right)^\frac1q
\end{split}\]
because we assume $\mu(E) \geq c \mu(2Q)$. 
Notice that $Z \in \gamma_{Q}^*(z)$ implies that
$z\in B(Z,C\delta(Z))$.
Therefore, by Fubini's lemma and \eqref{defTQ},
\[\begin{split}
|u_{B_0}| & \lesssim \ell(Q)^{1-\alpha} \left(\frac{1}{\mu(2Q)} \int_{Z\in T_{2Q}} |\nabla u(Z)|^q \, \frac{\delta(Z)^{q\alpha}}{m(B(Z,C\delta(Z)))} \, dm(Z) \int_{z\in B(Z,C\delta(Z))} d\mu(z)\right)^\frac1q.
\end{split}\]
If $Z\in \Omega$, 
we pick $z_0 \in \Gamma$ such that $|Z-z_0| = \dist(Z,\Gamma)$, and 
define $\rho(Z)$ by $\rho(z_0,\delta(Z))$. The point $z_0$, and then $\rho(Z)$, are not
uniquely defined, 
but it is easy to check that two choice of $z_0$ will be equivalent (up to constants independent of $Z$),
or if the reader prefers, the estimates below do not depend on our choice of $z_0$.
The doubling conditions (\hyperref[H3]{H3}) and (\hyperref[H4]{H4}) and the definition \eqref{defrho} imply that
\[\int_{z\in B(Z,C\delta(Z))} d\mu(z) \approx \mu(B(z_0,\delta(Z))) \approx \frac{m(B(z_0,\delta(Z))\cap \Omega)}{\delta(Z) \rho(z_0,\delta(Z))} \approx \frac{m(B(Z,2\delta(Z)))}{\delta(Z) \rho(Z)},\]
and the estimate on $u_{B_0}$ becomes
\[ |u_{B_0}| \lesssim \ell(Q)^{1-\alpha} \left(\frac{1}{\mu(2Q)} \int_{Z\in T_{2Q}} |\nabla u(Z)|^q \, \frac{\delta(Z)^{q\alpha-1}}{ \rho(Z)} \, dm(Z) \right)^\frac1q. \]
By \eqref{defzQrQ}, \eqref{propTQ}, 
\eqref{defrho}, 
and (\hyperref[H4]{H4}), one can show that $m(T_{2Q}) \approx \ell(Q) \rho(z_Q,\ell(Q)) \mu(2Q)$, where $z_Q$ is a fixed point in $Q$. Therefore, for $p\in (q,2]$,
\[\begin{split}
|u_{B_0}| &  \lesssim \ell(Q) \left(\ell(Q)^{1-q\alpha } \rho(z_Q,\ell(Q)) \fint_{Z\in T_{2Q}} |\nabla u(Z)|^q \frac{\delta(Z)^{q\alpha-1}}{\rho(Z)} dm(Z) \right)^\frac1q \\
& \lesssim \ell(Q) \Big(\fint_{T_{2Q}} |\nabla u|^p dm  \Big)^\frac1p 
\left( \big[\ell(Q)^{1-\alpha q} \rho(z_Q,\ell(Q)) \big]^{p/(p-q)}\fint_{Z\in T_{2Q}} \left(\frac{\delta(Z)^{[q\alpha-1]}}{\rho(Z)}\right)^\frac{p}{p-q} dm(Z)  \right)^\frac{p-q}{qp}
\end{split}\]
by the H\"older inequality. The claim \eqref{PoincareB} will thus be proven if we can establish that, for some $\alpha\in (0,1)$ close to $1$,
\begin{equation} \label{PoincareC}
I_0:= \fint_{Z\in T_{2Q}} \left(\frac{\delta(Z)^{[q\alpha-1]}}{\rho(Z)}\right)^\frac{p}{p-q} dm(Z)  \lesssim [\ell(Q)^{1-\alpha q} \rho(z_Q,\ell(Q))]^{-p/(p-q)}
\end{equation}
when $p<2$ is close enough to $2$.

By construction of $T_{2Q}$ (see \eqref{defTQ} and \eqref{defgammaQ*}),
\[T_{2Q} \subset \bigcup_{R\subset 100Q} U_R^*\]
and the covering has a uniformly finite overlap.
Notice also that
\[\rho(Z) \approx \rho(R) \qquad \text{ for $Z \in U^*_R$ and $R\in \D$},\]
where $\rho(R)$ is defined in \eqref{defrhoQ} and where the constants are independent of $Z$ and $R$. We call $\rho(R)$ the value of $\rho(Z)$ for a $Z\in U^*_R$. The two last observation allow us to write
\begin{equation} \label{boundonI0a}
I_0 \lesssim \frac{1}{m(T_{2Q})} \sum_{R\in \D \atop R\subset 100Q} m(U^*_R) \left(\frac{\ell(R)^{[q\alpha-1]}}{\rho(R)}\right)^\frac{p}{p-q}.
\end{equation}
We let the reader check that by definition of $\rho$, and by arguments similar to the ones used to prove Lemma \ref{lemmQ},
\[m(T_{2Q}) \approx \ell(Q) \rho(Q) \mu(Q)\]
and
\[m(U^*_R) \approx \ell(R) \rho(R) \mu(R).\]
The bound \eqref{boundonI0a} becomes
\begin{equation} \label{boundonI0b}
I_0 \lesssim \frac{1}{\ell(Q) \rho(Q) \mu(Q)} \sum_{R\in \D \atop R\subset 100Q} 
\ell(R)^{[q\alpha-1] \frac{p}{p-q} + 1} \rho(R)^{1-\frac{p}{p-q}} \mu(R).
\end{equation}
By (\hyperref[H5]{H5}), 
\[\rho(R) \gtrsim \rho(Q) \left( \frac{\delta(R)}{\delta(Q)} \right)^{1-\epsilon}
\approx \rho(Q) \left( \frac{\ell(R)}{\ell(Q)} \right)^{1-\epsilon},\]
where $\epsilon = C_5^{-1}$ is the one given in (\hyperref[H5]{H5}).
We use this to replace $\rho(R)$ in \eqref{boundonI0b} by $\rho(Q)$; notice that the
inequality goes in the right direction because the exponent $1-\frac{p}{p-q} = \frac{-q}{p-q}$
is negative (recall that $p \in (q,2]$). We get that 
\begin{equation} \label{boundonI0c}
I_0 \lesssim \rho(Q)^{-p/(p-q)} \ell(Q)^{a} \mu(Q)^{-1} \sum_{R\in \D \atop R\subset 100Q} \ell(R)^{b} \mu(R),
\end{equation}
with the exponents $a= -1 + (\varepsilon-1)\big( 1-\frac{p}{p-q}\big)$
and $b= [q\alpha-1] \frac{p}{p-q} + 1 + (1-\varepsilon) \big( 1-\frac{p}{p-q}\big)
= [q\alpha-2+\epsilon]\frac{p}{p-q} + 2-\epsilon$.

If $2-p \leq \min\{\epsilon/2, (2-p_0)/2\}$, we can pick $\alpha\in (0,1)$ (small) and $q\in [p_0,2)$ such that $q\alpha -2 +\epsilon \geq 0$. With these values, 
we can still pick $p \in (q,2]$ as above, and since the power for $\ell(R)$ is $b > 0$,
we can bound $\ell(R)$ brutally by $\ell(Q)$, which gives
\begin{equation} \label{boundonI0d}
I_0 \lesssim \rho(Q)^{-p/(p-q)} \ell(Q)^{a+b} \mu(Q)^{-1} \sum_{R\in \D \atop R\subset 100Q} \mu(R) \lesssim \rho(Q)^{-p/(p-q)} \ell(Q)^{a+b}.
\end{equation}
Notice that $a+b= [q\alpha-1]\frac{p}{p-q}$; the
claim \eqref{PoincareC} follows from the observation that $\rho(Q) \approx \rho(z_Q,\ell(Q))$, and 
we have seen before that \eqref{PoincareB} and Theorem \ref{PoincareTh2} follow.
\ep

In the following corollary of Theorem \ref{PoincareTh2} we replace the tents $T_{2Q}$ by balls.  

\begin{corollary} \label{PoincareCor2}
Let $(\Omega,m,\mu)$ satisfy (\hyperref[H1]{H1})--(\hyperref[H6]{H6}).
There exists $p_1\in [1,2]$ and $k:= k(C_4)>1$ such that the following happens for $p\in [p_1,2]$.
Let $\lambda >1$ be given, and let $x\in \Gamma$, $r>0$, and $u\in W$ be such that $\Tr u = 0$ on $B(x, \lambda r) \cap \Gamma$; then
\begin{equation} \label{PoincareCor2a}
\left( \fint_{B(x,r) \cap \Omega} |u|^{kp} \, dm \right)^{1/kp} 
\leq C_\lambda r \left( \fint_{B(x,\lambda r) \cap \Omega} |\nabla u|^p \, dm \right)^{1/p},
\end{equation} 
where $C>0$ depends only on $n$, $C_1$ to $C_6$, and $\lambda$.
\end{corollary}

\smallskip

\bp
Let $x'\in \Gamma$ and $r'>0$ be given. Let $Q'\in \mathbb D$ be
the only dyadic cube such that $x'\in Q'$ and $r' \leq k(Q')< 2r'$. Then 
$B(x',r') \subset 2Q'$, and by Theorem \ref{PoincareTh2} and \eqref{propTQ}, 
there exists $K>1$ that depends only on $n$, $C_1$, and $C_2$ such that
\begin{equation} \label{PoincareCor2b}
\left( \fint_{B(x',r') \cap \Omega} |u|^{kp} \, dm \right)^{1/kp} \leq C r' \left( \fint_{B(x',Kr') \cap \Omega} |\nabla u|^p \, dm \right)^{1/p},
\end{equation}
provided that $\Tr u \equiv 0$ on $Q' \subset B(x',Kr') \cap \Gamma$.

This looks like the desired estimate, but the constant $K$ is too large; we will fix this with
a covering argument. Set $\tau = (\lambda - 1)r/100K$, with $\lambda$ as in the statement
and $K$ as above. Then let $x \in \Gamma$ and $r>0$ be given. Denote by
$(x_i)_{i\in I}$ a maximal collection of points of $\Gamma \cap B(x,(1+2\tau)r)$ such that 
$|x_i-x_j| \geq \tau r$ for $i\neq j$. Thus the balls $B_i = B(x_i, 2\tau r)$,  
cover $\Gamma \cap B(x,(1+2\tau) r)$, and the sets $D_i = \Omega \cap B(x_i, 4\tau r)$
cover
\[
H := \big\{X\in \Omega: \, \dist(X, \Gamma \cap B(x,(1+2\tau)r) \leq 2\tau r \big\}
\]
Notice that $I$ has at most $C$ elements, with a constant that depends also on $\lambda$
and $K$ through $\tau$, but this is all right. We can apply \eqref{PoincareCor2b} to each
$B(x_i,4\tau r)$, and we get that 
\begin{equation} \label{PC2b2}
\left( \fint_{D_i} |u|^{kp} \, dm \right)^{1/kp} \leq C \tau r 
\left( \fint_{B(x,\lambda r) \cap \Omega} |\nabla u|^p \, dm \right)^{1/p},  \end{equation}
because $B(x_i, 4K\tau r) \subset B(x,\lambda r)$ by choice of $\tau$.
We may sum over $i$ and get that 
\begin{equation} \label{PC2b3}
\left( \fint_{H} |u|^{kp} \, dm \right)^{1/kp} \leq C \tau r 
\left( \fint_{B(x,\lambda r) \cap \Omega} |\nabla u|^p \, dm \right)^{1/p}, \end{equation}
and now we just need to take care of $H_1 = \Omega \cap B(x,r) \sm H$.
Let $(y_j)_{j\in J}$ be a maximal collection of points of $H_1$, with
$|y_i-y_j| \geq \tau r$ for $i\neq j$. Thus $J$ has at most $C = C(\tau)$
points, and the set $B_j = B(y_j,2\tau r)$, $j\in J$, cover $H$. We want to control
each $\fint_{B_j} |u|^{kp} \, dm$, and then we'll sum.

Fix $j\in J$, and call $z_j$ the first point  of $[y_j,x]$ (starting from $y_j$)
that lies within $\tau r$ from $\Gamma$. Obviously $z_j \in B(x,r)$, and
$B(z_j,\tau r) \subset H$ because $\Gamma \sm B(x,(1+2\tau)r)$ is too far.
Now denote by $W_j$ the convex hull of $B(y_j,\tau r)$ and $B(z_j,\tau r)$
(a nice tube contained in $\Omega$) and set $\wt W_j = W_j \cup B_j$
(with a larger head around $y_j$, and still contained in $\Omega$). It is easy to see
that $\wt W_j$ satisfies the chain condition $C(\kappa,M)$ of Definition \ref{defPchain}, 
with any small $\kappa$ chosen in advance, and with an $M$ that depends only on $\kappa$ 
and $\tau$; we can take $B(z_j,\tau r/2)$ as the distinguished ball. This allows us to apply 
Theorem \ref{Poincare1}, and prove that
\begin{equation} \label{PC2b4}
\left( \fint_{B_j} |u- \bar u|^{kp} \, dm \right)^{1/kp} 
\leq C \left( \fint_{\wt W_j} |u - \bar u|^{kp} \, dm \right)^{1/kp} 
\leq C r \left( \fint_{\wt W_j} |\nabla u|^p \, dm \right)^{1/p}, 
\end{equation}
where $\bar u$ denotes the average of $u$ on $B(z_j,\tau r/2)$.
Now $|\bar u| \leq C \left( \fint_{B(x,\lambda r) \cap \Omega} |\nabla u|^p \, dm \right)^{1/p}$,
by \eqref{PC2b3} and because $B(z_j,\tau r) \subset H$, and $\wt W_j \subset B(x,(1+2\tau)r)
\subset B(x,\lambda r)$ by definition of $\tau$, we we may sum \eqref{PC2b4} over $j$
and get that 
\begin{equation} \label{PC2b5}
\left( \fint_{H_1} |u|^{kp} \, dm \right)^{1/kp} \leq C \tau r 
\left( \fint_{B(x,\lambda r) \cap \Omega} |\nabla u|^p \, dm \right)^{1/p}.
\end{equation}
We combine this with \eqref{PC2b3} and get \eqref{PoincareCor2b}, as needed for 
Corollary \ref{PoincareCor2}.
\ep

\section{The Extension Theorem}
\label{Sext}

The aim of this section is the construction of an extension operator $\Ext :\, H \to W$ such that the composition $\Tr \circ \Ext$ is the identity on $H$.
The section can be seen as the dual of Section \ref{STrace}. As in Section  Section \ref{STrace}, the results will be only proven when assuming $\Gamma$ and $\Omega$ unbounded, and the proof in the bounded case is very similar and is discussed in Section \ref{SBounded}.

We assume that $\Gamma$ and $\Omega$ are unbounded.
The beginning of this section is similar to \cite[Section 7]{DFMprelim}, 
but the proof of the density result Lemma \ref{densityH} is different.

We shall construct $\Ext$ with the help of a Whitney extension. But first, it is crucial to observe that for any $g\in L^1_{loc}(\Gamma,\mu)$ and $\mu$-almost every $x\in \Gamma$, one has
\begin{equation} \label{Lebesguemu}
\lim_{r\to 0} \fint_{B(x,r)} |g(y)-g(x)|\, d\mu(y) = 0.
\end{equation}
This is a consequence 
of the Lebesgue differentiation theorem in doubling spaces (see for instance \cite[Sections 2.8-2.9]{Federer}). It is easy to verify that $H \subset L^1_{loc}(\Gamma,\mu)$ 
(see \eqref{defH})
and so
\eqref{Lebesguemu} holds for any function $g\in H$.

\medskip

Our construction will rely on
the family $\cW$ of dyadic Whitney cubes already used 
in Section \ref{Scones}. We associate to 
$\mathcal W$ a partition of unity $\{\varphi_I\}_{I\in \cW}$ where 
the $\varphi_I$ are smooth functions supported in $2I$ that satisfy $0 \leq \varphi_I \leq 1$, $|\nabla \varphi_I| \leq C/\ell(I)$ and $\sum_{I\in \cW} \varphi_I = \1_\Omega$.

We record a few properties of $\cW$, that can be found in \cite[Chapter VI]{Stein}. 
If two dyadic cubes $I$ and $I'$ are such that $2I \cap 2I' \neq \emptyset$, then 
$\ell(I)/\ell(I') \in \{1/2,1,2\}$, and also $I' \subset 6I$. 

Hence, for a given $I$, 
\begin{equation} \label{Ext1}
\text{the number of cubes $I'\in \cW$ such that $2I' \cap 2I \neq 0$ is at most $2\cdot 12^n$,} 
\end{equation}
because each such $I'$ needs to be a dyadic cube in $6I$ such that $\ell(I) \geq \ell(I)/2$.

For each $I \in \mathcal W$, we write $\delta(I)=\dist(I,\Gamma)$, pick a point $\xi_I \in \Gamma$
such that $\dist(\xi_I,I) \leq 2 \delta(I)$, and set $B_I=B(\xi_I,\ell(I))$.

We define
the extension operator $\Ext$ on functions $g\in L^1_{loc}(\Gamma,\mu)$ by
\begin{equation} \label{defExt}
\Ext g(X) := \sum_{I\in \cW} \varphi_I(X) y_I,
\end{equation}
where 
\begin{equation} \label{defyQ}
y_I := \fint_{B_I} g(z) \, d\mu(z).
\end{equation}
If we wanted to have an extension operator on - for instance - Lipschitz function, we could take $y_I = g(\xi_I)$. However, since the function $g$ is 
not smooth (and maybe not even defined everywhere),
we need this extra average; a good way to see this 
is to notice that otherwise we would only use the values of $g$ on the countable set 
$\{\xi_I\}_{I\in \cW}$, 
which does not make sense for functions in $L^1_{loc}(\Gamma,\mu)$.

Notice that $\Ext g$ lies in $C^\infty(\Omega)$ because \eqref{Ext1} yields that the sum in \eqref{defExt} is locally finite. Moreover, if $g$ is continuous on $\Gamma$, then $\Ext g$ is continuous on $\overline{\Omega}$ (see \cite[Proposition VI.2.2]{Stein}).

\begin{theorem} \label{ThExt} Let $(\Omega,m,\mu)$ satisfies (\hyperref[H1]{H1})--(\hyperref[H6]{H6}).
For any $g\in L^1_{loc}(\Gamma,\mu)$
\begin{equation} \label{Ext2}
\Tr \circ \Ext g = g \qquad \text{ $\mu$-a.e. in $\Gamma$}.
\end{equation}
 Moreover, $\Ext$ is a bounded linear operator from $H$ to $W$, i.e. there exists $C:=C(C_3,C_4,C_5)>0$ such that for any $g\in H$,
\begin{equation} \label{Ext3}
\int_\Omega |\nabla \Ext g|^2 \, dm 
\leq C \| g \|_H^2
:= C \int_\Gamma \int_\Gamma \frac{\rho(x,|x-y|)^2 |g(x) - g(y)|^2}{m(B(x,|x-y|))}\,  d\mu(x) \, d\mu(y).
\end{equation}
\end{theorem}

\bp
Let $g \in L^1_{loc}(\Gamma,\mu)$ be given. 
We write $u$ for $\Ext g$ and we want to show that $\Tr u = g$, 
in the sense that \eqref{defTr} holds with $Tr u(x) = g(x)$ for $\mu$-almost
every $x\in \Gamma$, regardless of whether $g\in H$ or $u \in W$. 
We will actually prove the following stronger result, analogous to \eqref{defLebesgue}: 
for $\mu$-a.e. $x\in \Gamma$, one has
\begin{equation} \label{Ext4}
\lim_{X \in \gamma(x) \atop {\delta(X) \to 0}} \fint_{B(X,\delta(X)/2)} |u(Z)-g(x)| \, dm(Z) = 0.
\end{equation}
Since we only want to prove \eqref{Ext4} for $\mu$-a.e. point, we can restrict to 
the case when $x$ is a Lebesgue point of $g$, that is, when \eqref{Lebesguemu} is satisfied.

Fix such an $x\in \Gamma$ and $X \in \gamma(x)$. 
We write $B$ for $B(X,\delta(X)/2)$. Then 
\[ \fint_{B} |u(Z)-g(x)| \, dm(Z)  \leq 
\frac1{m(B)} \sum_{R \in \cW(B)} \int_{R} |u(Z) - g(x)|\, dm(Z), \]
where $\cW(B)$ is the set of dyadic cubes $I' \in \cW$ that meet $B$. It is easy to check that $\cW(B)$ contains a finite number of cubes $I'$ (the number is bounded uniformly in $X \in \Omega$),

for which $\ell(I') \approx \delta(X)$, and then,  by (\hyperref[H4]{H4}),  $m(I') \approx m(B)$.
So \eqref{Ext4} will be proven if we can establish that
\begin{equation} \label{Ext5}
\fint_{I'} |u(Z)-g(x)| \, dm(Z) \longrightarrow 0 \quad \text{ as $\delta(I') \to 0$},
\end{equation}

where we restrict to dyadic cubes $I'\in \cW$ such that $x\in KI'$ for some large enough constant $K:=K(n)$.
Recall from the definition \eqref{defExt} of $u=\Ext g$ that $u(Z) = \sum_{I\in \cW} \varphi_I(Z) y_I$.
We observed earlier that
the sum is locally finite on $I'$, and the cubes $I$ for which $\varphi_I$
does not vanish identically on $I$
are such that $I' \subset 6I$ and $\ell(I)/\ell(I') \in \{1/2,1,2\}$. We deduce that any such $I$ satisfies $B_I \subset K'B_{I'} \subset B(x,K''\delta(I'))$ and $B_{I'} \subset K'B_I \subset B(x,K''\delta(I'))$, and then by (\hyperref[H3]{H3}) that $\mu(B_I) \approx \mu(B(x,K''\delta(I')))$. The conclusion is that 
\[\begin{split}
\fint_{I'} |u(Z)-g(x)| \, dm(Z) & =  \fint_{I'} \bigg| \sum_{I \in \cW:\, 2I \cap 2I'\neq \emptyset} \varphi_I(Z) \fint_{B_I} [g(z)-g(x)] \, d\mu(z)   \bigg| \, dm(Z) \\
& \lesssim \sum_{I \in \cW: \, 2I \cap 2I'\neq \emptyset}  \fint_{B_I} |g(z)-g(x)| \, d\mu(z) \\
& \lesssim \fint_{B(x,K''\delta(I'))} |g(z)-g(x)| \, d\mu(z)
\end{split}\]
because the number of $I\in \cW$ that verify $2I \cap 2I'\neq \emptyset$ is uniformly 
bounded.
Thanks to \eqref{Lebesguemu}, the right-hand side above converges to 0 as $\delta(I') \to 0$. 
The claims \eqref{Ext5}, \eqref{Ext4}, and then \eqref{Ext2}, follow.

\medskip

Now, we want to show that for $g\in H$, $u\in W$ and even $\|u\|_W \lesssim \|g\|_H$. 
Recall that $u$ is smooth on $\Omega$ because the sum in \eqref{defExt} is locally
finite, so $u$ is locally integrable in $\Omega$, and its distribution derivative is locally
integrable too, and given by
\begin{equation}\label{Ext6}
\nabla u(X) = \sum_{I\in \cW} y_I \nabla \varphi_I(X) 
= \sum_{I\in \cW} [y_I - y_{I'}] \nabla \varphi_I(X),
\end{equation}
where $I'$ is any cube (that may depend on $X$ but not on $I$), and the identity holds because 
$\sum_I \nabla \varphi_I = \nabla (\sum_I \varphi_I) = \nabla 1 = 0$. 
So we only need to show that 
$\|u\|_W \leq C \|g\|_H < +\infty$.
First decompose $\|u\|^2_W$ as 
\begin{equation} \label{Ext7}
\|u\|^2_W
= \sum_{I' \in \cW} \int_{I'} |\nabla u|^2 \, dm.
\end{equation}
For the moment, we fix $I' \in \cW$ and $X\in I'$, and 
get a bound on $|\nabla u(X)|$. If $\cW(I')$ denotes the sets of dyadic cubes $I \in \cW$ such that $2I$ meets $I'$, then  
\[|\nabla u(X)| \leq \sum_{I \in \cW(I')} 
|y_I-y_{I'}| |\nabla \varphi_I(X)| \lesssim \ell(I')^{-1} \sum_{I \in \cW(I')} |y_I-y_{I'}| \]
because  $\nabla \varphi_I \lesssim \delta(I)^{-1} \approx \delta(I')^{-1}$. We use the definition of $y_I,y_{I'}$, the facts that $I\subset 6I'$ and $\delta(I) \approx \delta(I')$ to obtain that
\[\begin{split}
|y_I - y_{I'}| &\leq \fint_{B_I} \fint_{B_{I'}} |g(x)-g(y)| d\mu(x) \, d\mu(y) \\
& \leq \left( \fint_{B_I} \fint_{B_{I'}} |g(x)-g(y)|^2 d\mu(x) \, d\mu(y) \right)^\frac12 \\
& \lesssim \mu(B_{I'})^{-1} \left( \int_{100B_{I'}} \int_{B_{I'}} |g(x)-g(y)|^2 d\mu(x) \, d\mu(y) \right)^\frac12.
\end{split}\]
The combination of the last two computations gives
\[|\nabla u(X)| \lesssim \mu(B_{I'})^{-1}\ell(I')^{-1} \left( \int_{100B_{I'}} \int_{B_{I'}} |g(x)-g(y)|^2 d\mu(x) \, d\mu(y) \right)^\frac12\]
since $\cW(I')$ contains at most $2\cdot12^n$ elements, 
and then
\[\int_{I'} |\nabla u|^2 \, dm \lesssim \ell(I')^{-2} \mu(B_{I'})^{-2} m(I') \int_{100B_{I'}} \int_{B_{I'}} |g(x)-g(y)|^2 d\mu(x) \, d\mu(y).\]
We inject the above estimate in \eqref{Ext7} and  
obtain that
\[\begin{split}
\|u\|_W^2 & \lesssim \sum_{I'\in \cW} \ell(I')^{-2} \mu(B_{I'})^{-2} m(I') \int_{100B_{I'}} \int_{B_{I'}} |g(x)-g(y)|^2 d\mu(x) \, d\mu(y) \\
& \lesssim \int_{\Gamma} \int_{\Gamma} |g(x)-g(y)|^2 h(x,y) d\mu(x) \, d\mu(y), 
\end{split}\]
where 
\[h(x,y) := \sum_{I'\in \cW} \ell(I')^{-2} \mu(B_{I'})^{-2} m(I') \1_{100B_{I'}}(x) \1_{B_{I'}}(y).\]
Fix $x,y\in \Gamma$. Observe that if $I'$ satisfies $(x,y) \in 100B_{I'} \times B_{I'}$, 
then by (\hyperref[H3]{H3}), $\mu(B_{I'}) \approx \mu(B(x,\ell(I')))$ and by (\hyperref[H4]{H4}), $m(I') \approx m(B(x,\ell(I')) \cap \Omega)$. Hence by \eqref{defrho} 
\[ 
\frac{m(I')}{\ell(I')^2 \mu(B_{I'})^2} \approx \frac{\rho(x,\ell(I'))}{\ell(I') \mu(B(x,\ell(I')))} \, .
\]
Under the same assumption on $I'$, we also have $|x-y| \leq 101\ell(I')$, so by  
(\hyperref[H4]{H4}) again, 
$\mu(B(x,\ell(I')))^{-1} \lesssim \mu(B(x,|x-y|))^{-1}$. In addition, (\hyperref[H5]{H5}) gives that
\[\rho(x,\ell(I')) \lesssim \rho(x,|x-y|) \left( \frac{\ell(I')}{|x-y|} \right)^{1-\epsilon}\]
where $\epsilon := C_5^{-1} >0$
(notice that if $\ell(I') \leq |x-y| \leq 101\ell(I')$, we don't need
to use (\hyperref[H5]{H5}), just the doubling properties). All this yields
\[h(x,y) \lesssim \frac{\rho(x,|x-y|)}{\mu(B(x,|x-y|)) |x-y|^{1-\epsilon}} \sum_{I'\in \cW \atop \delta(I') \geq |x-y|/101} \ell(I')^{-\epsilon} \, \1_{B_{I'}}(y).\]
Since $B_{I'} \subset \kappa I'$
for some constant $\kappa:=\kappa(n) >1$ that does not depend on $I'$,
we see that for each $j\in \mathbb Z$, the number of dyadic cubes $I'$ such that $\ell(I') = 2^j$
and $y\in B_{I'}$ is uniformly bounded.
Together with the fact that $\delta(I') \approx \ell(I')$, this yields 
\[\sum_{I'\in \cW \atop \ell(I') 
 \geq |x-y|/101} \ell(I')^{-\epsilon} \, \1_{B_{I'}}(y)\lesssim \sum_{k\in \bN} (2^{k}|x-y|)^{-\epsilon} \lesssim |x-y|^{-\epsilon}.\] Altogether, 
\[h(x,y) \lesssim \frac{\rho(x,|x-y|)}{\mu(B(x,|x-y|)) |x-y|} \approx \frac{\rho(x,|x-y|)^2}{m(B(x,|x-y|)\cap \Omega)}\]
by \eqref{defrho} and thus
\[ \|u\|_W^2 \lesssim \int_{\Gamma} \int_{\Gamma} \frac{\rho(x,|x-y|)^2|g(x)-g(y)|^2}{m(B(x,|x-y|)\cap \Omega)}d\mu(x) \, d\mu(y):= \|g\|_H^2\]
as desired (see the definition \eqref{defH}). Theorem \ref{ThExt} follows.
\ep

\ms
\begin{lemma} \label{densityH}
For every $g\in H$, we can find a sequence $(g_k)_{k\in \bN}$ of functions in $C^\infty(\R^n)$ whose restrictions to $\Gamma$ (we still call them $g_k$) belong to $H$ and such that $(g_k)_k$ converges to $g$ in $H$, $L^2_{loc}(\Gamma,\mu)$ and $\mu$-a.e. pointwise.
\end{lemma}

\begin{remark}
The above density result (whose proof doesn't use Theorem \ref{ThExt}) actually entails the Lebesgue density result given as \eqref{Lebesguemu}. 
The proof of this implication uses 
maximal functions, is classical, and is left to the reader. 
\end{remark}

\bp
For the density of smooth functions, we are given $g\in H$ and we want to approximate it with smooth functions. The simplest way 
for us to construct functions $g_k$ will be to use our dyadic decompositions
$\D_k$ of $\Gamma$, but coverings of $\Gamma$ with balls of radius $2^{-k}$ would work as well.
We associate to $\D_k$ a collection of smooth functions $\{\varphi_Q\}_{Q\in \D_k}$ such that $\varphi_Q$ is supported in $2Q$, $\sum_{Q \in \D_k} \varphi_Q = 1$ near $\Gamma$, and $\|\nabla \varphi_Q\|_\infty \leq C 2^k$.
Finally we set 
\begin{equation}\label{3.18}
g_k(x) = \sum _{Q \in \D_k} \varphi_Q(x) y_Q
\end{equation}
for $x\in \Gamma$, where we take $y_Q = \fint_{2Q} g(y) d\mu(y)$. 
It is obvious that $g_k$ is a smooth function on $\R^n$ (the sum in \eqref{3.18} is locally finite). 
We shall prove now that
\begin{equation}\label{3.19}
\| g- g_k \|^2_H \leq C J(k),
\ \text{ where }
J(k) = \iint_{x,y \in \Gamma \, ; \, |x-y| \leq 
2^{3-k}} 
{\rho(x,|x-y|)^2|g(x)-g(y)|^2 \over m(B(x,|x-y|) \cap \Omega)} d\mu(x)d\mu(y).
\end{equation}
Notice that $J_k$ is a subintegral of $\| g\|_H^2$, where the domain of integration 
decreases to the empty set when $k$ tends to $+\infty$; thus
$\lim_{k \to +\infty} J(k) = 0$, 
and as soon as we prove \eqref{3.19}, we will get that $g_k \in H$ 
and $g_k$ tends to $g$ in $H$; the density of smooth functions in $H$ will follow.

\medskip

We need some notation. Set $h_k = g-g_k$ and for $Q\in \D$
\begin{equation}\label{3.20}
\cR_Q = \big\{ (x,y)\in Q \times 2Q \, ; |x-y| \geq \ell(Q)/2 \big\}.
\end{equation}
Every
pair 
of points 
$(x,y) \in \Gamma^2$ lies in at least one $\cR_Q$: choose $j$ 
such that $2^{-j-1} \leq |x-y| < 2^{-j}$, let $Q = Q^j(x)$ be the element of $\cD_j$,
and observe that $y\in B(x,2^{-j}) \subset 2Q$, hence $(x,y)\in \cR_{Q}$.
Because of this,
\begin{equation} \label{3.20b}
\|g- g_k \|^2_\H \leq \sum_{Q\in \D} T_Q^k,
\end{equation}
where 
\begin{eqnarray}\label{3.21}
T_Q^k &:=& \iint_{\cR_Q} {\rho(x,|x-y|)^2 |h_k(x)-h_k(y)|^2 \over m(B(x,|x-y|)\cap \Omega)} d\mu(x)d\mu(y)
\nn\\
& = & \iint_{\cR_Q} {|h_k(x)-h_k(y)|^2 \over \mu(B(x,|x-y|))}\frac{\rho(x,|x-y|)}{|x-y|} d\mu(x)d\mu(y) \\
&\lesssim& \frac{\rho(Q)}{\ell(Q)} \iint_{\cR_Q} \frac{|h_k(x)-h_k(y)|^2}{\mu(2Q)} d\mu(x)d\mu(y) \nn
\end{eqnarray}
because $|x-y| \approx \ell(Q)$, and by (\hyperref[H3]{H3})--(\hyperref[H4]{H4}) and the definitions\eqref{defrho} and \eqref{defrhoQ}.

We start the estimate of $T_Q^k$ with 
the large scales, 
where we shall merely estimate the size of $h_k$ on $\Gamma$.
Let us check that for any cube $Q^*$ such that 
$\ell(Q^*) \geq 2^{-k}$, 
\begin{equation}\label{3.22}
\int_{2Q^*} |h_k(x)|^2 d\mu(x) 
\leq C 
\int_{x\in 8Q^*} 
\frac1{\mu(B(x,2^{2-k}))} \int_{z\in B(x,2^{2-k})} |g(x)- g(z)|^2 d\mu(z) d\mu(x).
\end{equation}
We shall first estimate the contribution of a 
a given cube $Q_0 \in \D_k(Q^*)$, where
\begin{equation} \label{3.23}
\D_k(Q^*) := \{Q_0 \in \D_k, \, 2Q_0 \cap 2Q^* \neq \emptyset\},
\end{equation}
and then sum. So let $Q_0 \in \D_k(Q^*)$
be given. 
We estimate
\begin{eqnarray}\label{3.24}
a(Q_0) &:=& \int_{x\in  2Q_0} |h_k(x)|^2 d\mu(x)
= \int_{ 2Q_0} |g(x)-g_k(x)|^2 d\mu(x)
\nn\\
&=& \int_{x\in  2Q_0} |g(x)- \sum_{Q\in \D_k} \varphi_Q(x) y_Q|^2 d\mu(x)
= \int_{x\in  2Q_0} \big|\sum_{Q\in \D_k} \varphi_Q(x) [g(x)-  y_Q] \big|^2 d\mu(x)
\nn\\
&\leq& \int_{x\in  2Q_0} \sum_{Q\in \D_k} \varphi_Q(x) 
|g(x)-  y_Q|^2 d\mu(x)
\leq \int_{x\in  2Q_0} \sum_{Q\in \D_k(Q_0)} |g(x)-  y_Q|^2 d\mu(x)
\end{eqnarray}
by \eqref{3.18}, the fact that $\sum_Q \varphi_Q(x) = 1$, Cauchy-Schwarz for a finite average, and the fact that $\varphi_Q(x) = 0$ outside of $2Q$.
Notice that when $Q\in \D_k(Q_0)$ and $x\in 2Q_0$,
\begin{eqnarray}\label{3.25}
|g(x)-y_Q|^2  &=& \left|\fint_{z\in  2Q} [g(x)-g(z)]d\mu(z)\right|^2
\leq \fint_{z\in  2Q} |g(x)-g(z)|^2 d\mu(z)
\nn\\
&\leq& C \fint_{z\in \Gamma \cap B(x,2^{2-k})} |g(x)-g(z)|^2 d\mu(z)
\end{eqnarray}
because for $x\in 2Q_0$, the fact that $2Q_0 \cap 2Q \neq \emptyset$ implies that 
$2Q \subset B(x,2^{2-k})$, and $\mu$ is doubling by (\hyperref[H3]{H3}). 
Hence, since the number of element in $\D_k(Q_0)$ is bounded, 
\begin{eqnarray}\label{3.26}
a(Q_0) &\leq& 
C\int_{x\in  2Q_0} \sum_{Q\in \D_k(Q_0)} \fint_{z\in \Gamma \cap B(x,2^{2-k})} 
|g(x)-g(z)|^2 d\mu(z)d\mu(x)
\nn\\
&\leq&
C \int_{x\in  2Q_0} \fint_{z\in \Gamma \cap B(x,2^{2-k})} |g(x)-g(z)|^2 d\mu(z)d\mu(x)
\nn\\
&\leq& C  \int_{x\in  2Q_0} \frac 1{\mu(B(x,2^{2-k}))} 
\int_{z\in \Gamma \cap B(x,2^{2-k})} |g(x)-g(z)|^2 d\mu(z)d\mu(x).
\end{eqnarray}
Now, we sum on $Q_0 \in \D_k(Q^*)$ and get that
\begin{eqnarray}\label{3.27}
\int_{2Q^*} |h_k(x)|^2 d\mu(x) &\leq& \sum_{Q_0 \in \D_k(Q^*)} 
\int_{ 2Q_0} |h_k(x)|^2 d\mu(x) \leq \sum_{Q_0 \in \D_k} a(Q_0)
\nn\\&\lesssim &  \sum_{Q_0 \in \D_k(Q^*)} 
\int_{x\in  2Q_0}  \frac 1{\mu(B(x,2^{2-k}))} \int_{z\in \Gamma \cap B(x,2^{2-k})} |g(x)-g(z)|^2 d\mu(z)d\mu(x)
\nn\\&\lesssim&  
\int_{x\in 8Q^*}  \frac 1{\mu(B(x,2^{2-k}))} \int_{z\in \Gamma \cap B(x,2^{2-k})} |g(x)-g(z)|^2 d\mu(z)d\mu(x)
\end{eqnarray}
because the $2Q_0$ cover $Q^*$ 
(actually, they cover $2Q^*$), are contained in $8Q^\ast$,
and have bounded covering; the estimate \eqref{3.22} follows.

Recall that 
since $\ell(Q^\ast) \geq 2^{-k}$, 
(\hyperref[H5]{H5}) and Lemma \ref{lemmQ} imply that 
\[ \frac{2^{2-k}}{\rho(x,2^{2-k})} \lesssim \frac{\ell(Q^*)}{\rho(Q^*)} \left( \frac{\ell(Q^*)}{2^{-k}} \right)^{-\epsilon} \]
(with $\epsilon = C_5^{-1}$ as usual) and, for 
$z\in B(x,2^{2-k}) \cap \Gamma$,
\[  \frac{\rho(x,2^{2-k})}{2^{2-k}} \lesssim \frac{\rho(x,|x-z|)}{|x-z|} \left( \frac{2^{-k}}{|x-z|} \right)^{-\epsilon} \lesssim \frac{\rho(x,|x-z|)}{|x-z|}.\]
We return to \eqref{3.22}, use the two estimates above and the fact that
$\mu(B(x,|x-z|)) \leq C \mu(B(x,2^{2-k}))$ when $B(x,2^{2-k})$, and get that
\begin{equation} \label{3.27b} \begin{split}
\int_{2Q^*} |h_k|^2 d\mu 
& \lesssim \int_{x\in 8Q^*} \frac{\rho(x,2^{2-k})}{2^{2-k}\mu(B(x,2^{2-k}))} 
\frac{2^{2-k}}{\rho(x,2^{2-k})} \int_{z\in \Gamma \cap B(x,2^{2-k})} |g(x)-g(z)|^2 d\mu(z)d\mu(x) 
\\
& \lesssim \frac{\ell(Q^*)}{\rho(Q^*)} \left( \frac{\ell(Q^*)}{2^{-k}} \right)^{-\epsilon}  
\int_{x\in 8Q^*} \int_{z\in B(x,2^{2-k})} \frac{\rho(x,|x-z|) |g(x)-g(z)|^2}{\mu(B(x,|x-z|)) |x-z|} d\mu(z)d\mu(x) 
\\
& = \frac{\ell(Q^*)}{\rho(Q^*)} \left( \frac{\ell(Q^*)}{2^{-k}} \right)^{-\epsilon}  
\int_{x\in 8Q^*} \int_{z\in B(x,2^{2-k})} \frac{\rho(x,|x-z|)^2 |g(x)-g(z)|^2}{m(B(x,|x-z|) \cap \Omega)} d\mu(z)d\mu(x) 
\end{split}\end{equation} 
where the last
estimate comes from the definition of $\rho$. 
The right-hand side tends to $0$ (for any fixed $Q^\ast$) when $k$ tends to $+\infty$,
so \eqref{3.27b} 
means that $(g_k)$ converges to $g$ in $L^2_{loc}(\Gamma,\mu)$.

\medskip

Let us return to $T^k_{Q^*}$, 
starting with the case when $\ell(Q^*) \geq 2^{-k}$.
Observe that by \eqref{3.21} and \eqref{3.27b}
\begin{eqnarray}\label{3.28}
T^k_{Q^*} &\lesssim &  \frac{\rho(Q^*)}{\ell(Q^*)} \iint_{\cR_{Q^*}} \frac{|h_k(x)-h_k(y)|^2}{\mu(2Q^*)} d\mu(x)d\mu(y)
\nn\\
&\lesssim&\frac{\rho(Q^*)}{\ell(Q^*)} \iint_{\cR_{Q^*}} \frac{|h_k(x)|^2+|h_k(y)|^2}{\mu(2Q^*)} d\mu(x)d\mu(y)
\nn\\
&\lesssim& \frac{\rho(Q^*)}{\ell(Q^*)}  \int_{2Q^*} |h_k(x)|^2 d\mu(x) \\
& \lesssim & \left( \frac{\ell(Q^*)}{2^{-k}} \right)^{-\epsilon}
\int_{x\in 8Q^*} \int_{z \in B(x,2^{2-k})} 
\frac{\rho(x,|x-z|)^2 |g(x)-g(z)|^2}{m(B(x,|x-z|) \cap \Omega)} d\mu(z)d\mu(x). \nn
\end{eqnarray}
We now sum this over $Q^\ast$ such that $\ell(Q^*) \geq 2^{-k}$. Fix $x, z \in \Gamma$;
for each generation $j$, $j \leq k$, there are at most $C$ cubes $Q^\ast$, $Q^* \in \D_j$,
such that $x\in 8Q^*$. Therefore,
\begin{equation}\label{3.29}
\sum_{j\leq k} \sum_{Q^* \in \D_j} T_{Q^*}^k 
\lesssim \sum_{j\geq k} 2^{(j-k)\epsilon} J(k)  \lesssim  J(k),
\end{equation}
where $J(k)$ is as in \eqref{3.19}. This part fits with \eqref{3.19} (see \eqref{3.20b}).

\medskip

For the small scales we shall use the regularity of $g_k$. That is, for $\ell(Q) \leq 2^{-k}$, we recall that
$h_k= g-g_k$, hence, by the first part of \eqref{3.21}, $T_{Q}^k \leq 2U_Q^k + 2V_Q^k$, where 
\begin{equation}\label{3.30}
V_Q^k :=  \iint_{\cR_Q} {|g_k(x)-g_k(y)|^2 \over \mu(B(x,|x-y|))}\frac{\rho(x,|x-y|)}{|x-y|} d\mu(x)d\mu(y)
\end{equation}
and
\begin{equation}\label{3.31}
W_Q^k :=  \iint_{\cR_Q} {|g(x)-g(y)|^2 \over \mu(B(x,|x-y|))}\frac{\rho(x,|x-y|)}{|x-y|} d\mu(x)d\mu(y)
\end{equation}
are the analogues of $T_Q^k$ for $g_k$ and $g$. 

We deduce from the definition \eqref{3.20} that $\ell(Q)/2 \leq |x-y| \leq 2 \ell(Q)$ when $(x,y) \in \cR_Q$,
so a given pair $(x,y)$ cannot lie in more than $C$ sets $\cR_Q$, $Q \in \D$, so \eqref{3.31} yields
\begin{equation}\label{3.29a}
\sum_{j> k} \sum_{Q \in \D_j} W_{Q}^k \lesssim J(k).
\end{equation}
As for the $V_Q^k$, we decide to estimate $|g_k(x)-g_k(y)|$ rather brutally. Again we localize at the scale $2^{-k}$. Let $Q_0 \in \D_k$ be given,
and then pick $x\in Q_0$ and $y\in 2Q_0$.
We want to estimate $|g_k(x)-g_k(y)|$ in terms of
\begin{equation}\label{3.32}
b(Q_0) := \fint_{y\in  2Q_0}\fint_{z\in  4Q_0} 
|g(y)-g(z)|^2 d\mu(y)d\mu(z).
\end{equation}
By \eqref{3.18},
\begin{eqnarray}\label{3.33}
|g_k(x)-g_k(y)| &=& \Big| \sum_{Q\in \D_k} [\varphi_Q(x)-\varphi_Q(y)] y_Q \Big|
= \Big| \sum_{Q\in \D_k} [\varphi_Q(x)-\varphi_Q(y)] [y_Q-y_{Q_0}] \Big|
\nn\\
&\leq& \sum_{Q\in \D_k} |\varphi_Q(x)-\varphi_Q(y)| \, |y_Q-y_{Q_0}|
\end{eqnarray}
because $\sum_Q \varphi_Q(x) = \sum_Q \varphi_Q(y) = 1$. Notice that if  
$|\varphi_Q(x)-\varphi_Q(y)| > 0$,
then $x$ or $y$ lies in $2Q$,
hence $Q$ lies in the collection $\D_k(Q_0)$ of \eqref{3.23}. For such balls, 
\begin{eqnarray}\label{3.34}
|y_Q-y_{Q_0}| &=& \Big| \fint_{y\in 2Q_0}\fint_{z\in 2Q} [g(y)-g(z)] \,d\mu(z) \, d\mu(y)\Big|
\leq \fint_{y\in 2Q_0}\fint_{z\in 2Q} |g(y)-g(z)| \,d\mu(z) \, d\mu(y)
\nn\\
&\leq& \Big\{\fint_{y\in 2Q_0}\fint_{z\in 2Q} 
|g(y)-g(z)|^2 d\mu(y)d\mu(z)\Big\}^{1/2}
\nn\\
&\leq& C \Big\{\fint_{y\in 2Q_0}\fint_{z\in 4Q_0} 
|g(y)-g(z)|^2 d\mu(y)d\mu(z)\Big\}^{1/2} = C b(Q_0)^{1/2}, 
\end{eqnarray}
where $b(Q_0)$ is as in \eqref{3.32}, and 
because $2Q \subset 4Q_0$
for $Q \in \D_k(Q_0)$, and by (\hyperref[H3]{H3}).
Recall that $|\nabla \varphi_Q| \leq C {2^k}$. Since there are no more than $C$ cubes $Q\in \D_k(Q_0)$, 
\eqref{3.33} yields
\begin{equation}\label{3.35}
|g_k(x)-g_k(y)|^2 \lesssim |x-y|^2 2^{2k} \sup_{Q\in \D_k(Q_0)} |y_Q-y_{Q_0}|^2
\lesssim 
 |x-y|^2 2^{2k} b(Q_0).
\end{equation}
We shall use this soon, but for the moment let us estimate a given 
$V_{Q}^k$, $\ell(Q) \leq 2^{-k}$. 
Observe that by (\hyperref[H3]{H3})--(\hyperref[H4]{H4}) and the definitions\eqref{defrho} and \eqref{defrhoQ}, and as in the last part of \eqref{3.21},
\[V_Q^k \lesssim \frac{\rho(Q)}{\ell(Q)} \iint_{\cR_Q} \frac{|g_k(x)-g_k(y)|^2}{\mu(2Q)} d\mu(x)d\mu(y).\]
Then let $Q_0$ denote the cube of $\D_k$ that contains $x$, then
$x\in Q_0$ and $y \in 2Q_0$ when $(x,y) \in \cR_Q$, so we can use
\eqref{3.35} and get that
\begin{eqnarray}\label{3.36}
V_Q^k &\lesssim& \frac{\rho(Q)}{\ell(Q)} \iint_{\cR_Q} \frac{|x-y|^2 2^{2k} b(Q_0)}{\mu(2Q)} \, d\mu(y) \, d\mu(x).
\end{eqnarray}
Since 
$|x-y|\approx \ell(Q)$ when $(x,y)\in \cR_Q$, 
\begin{eqnarray}\label{3.36a}
V_Q^k &\lesssim & 
2^{2k}\ell(Q)\rho(Q)
\int_{x\in Q} \int_{y\in 2Q} \frac{b(Q_0)}{\mu(2Q)} d\mu(y)d\mu(x)
\nn
\\
&\lesssim & 2^{2k} \ell(Q)\rho(Q)\mu(Q) b(Q_0) 
\\
& \lesssim & 2^{2k} m(U^*_Q) b(Q_0) \approx m(U^*_Q) \ell(Q_0)^{-2} b(Q_0)
\end{eqnarray}
by \eqref{defrhoQ}. We sum over the cubes $Q\subset Q_0$ to obtain
\[\begin{split}
\sum_{Q\subset Q_0} V_Q^k \lesssim 2^{2k} \sum_{Q\subset Q_0} m(U^*_Q) b(Q_0) \approx 2^{2k} m(T_{Q_0})b(Q_0) \approx \frac{\rho(Q_0) \mu(Q_0)}{\ell(Q_0)} b(Q_0)
\end{split}\]
because the $U^*_Q$ are contained in $T_{Q_0}$ (see \eqref{defTQ} and \eqref{defgammaQ*}) 
and have bounded overlap.
We now use the definition \eqref{3.32} of $b(Q_0)$ and the doubling property (\hyperref[H3]{H3}) to get
\[
\sum_{Q\subset Q_0} V_Q^k  \lesssim \int_{2Q_0} \int_{4Q_0} \frac{\rho(Q_0) |g(y) - g(z)|^2}{\mu(4Q_0) \ell(Q_0)}  \, d\mu(y) \, d\mu(z). 
\]
Observe that $\mu(B(y,|y-z|)) \leq C \mu(4Q_0)$.
Besides, due to (\hyperref[H5]{H5}),
\[  \frac{\rho(Q_0)}{\ell(Q_0)} 
\lesssim \frac{\rho(y,|y-z|)}{|y-z|} \left( \frac{\ell(Q_0)}{|x-z|} \right)^{-\epsilon} 
\lesssim \frac{\rho(y,|y-z|)}{|y-z|} \qquad \text{ for $(y,z)\in 2Q_0 \times 4Q_0$}.\]
It follows that
\[ \begin{split}
\sum_{Q\subset Q_0} V_Q^k & \lesssim \int_{2Q_0} \int_{4Q_0} \frac{\rho(y,|y-z|) |g(y) - g(z)|^2}{\mu(B(y,|y-z|)) |y-z|}  \, d\mu(y) \, d\mu(z) 
\\
& \quad = \int_{2Q_0} \int_{4Q_0} \frac{\rho(y,|y-z|)^2 |g(y) - g(z)|^2}{m(B(y,|y-z|)\cap \Omega)}  \, d\mu(y) \, d\mu(z)
\end{split}\]
by \eqref{defrho}. Notice that $|y-z| \leq 8 \ell(Q_0) = 2^{3-k}$ when $y\in 2Q_0$
and $z\in 4Q_0$. Also, for a given pair $(y,z) \in \Gamma \times \Gamma$, with
$y \neq z$, the set of cubes $Q_0$ of generation $k$ for which $y\in Q_0$ and $z\in 2Q_0$
has less than $C$ elements; because of this, 
\begin{equation} \label{3.37} \begin{split}
\sum_{j > k} \sum_{Q\in \D_j} V_Q^k  
& = \sum_{Q_0 \in \D_k} \sum_{Q\subset Q_0} V_Q^k 
\\
&
\lesssim   \iint_{x,y \in \Gamma \, ; \, |x-y| \leq 2^{3-k}}  
\frac{\rho(y,|y-z|)^2 |g(y) - g(z)|^2}{m(B(y,|y-z|)\cap \Omega)}
\, d\mu(y) \, d\mu(z) = J(k)
\end{split}\end{equation}
The combination of \eqref{3.29}, \eqref{3.29a}, and \eqref{3.37} gives that 
$\sum_{j\geq 0} T_j^k \lesssim J(k)$, and hence by \eqref{3.20b}
$\| g- g_k \|^2_H = \sum_j T_j^k \leq CJ(k)$, as needed for \eqref{3.19}.

This completes our proof of the density of smooth functions 
in $H$.
The fact that $g_k$ converges to $g$ in $L^2_{loc}(\Gamma,\mu)$ has been shown in \eqref{3.27b}, and up to a subsequence, we can also assume that $g_k$ also converges to $g$ $\mu$-a.e. on $\Gamma$. Lemma \ref{densityH} follows.
\ep

\section{Completeness of $W$ and density of smooth functions}
\label{SComplete}

In all of this section, 
$(\Omega,m,\mu)$ satisfies (\hyperref[H1]{H1})--(\hyperref[H6]{H6}).

\medskip

First  
we talk about  
completeness. 
In fact $W$ cannot really be a Banach space, because $\|.\|_{W}$ is not a norm on $W$,
only a semi-norm. Thus
we need to quotient $W$ by the functions $u$ such that $\|u\|_W = 0$, that is,
thanks to Lemma \ref{Nablau}, by the constant functions. So we work
 with the homogeneous space $\dot W$ defined as the quotient space $W/\R$ - i.e. element of $\dot W$ are classes $\dot u := \{u+c\}_{c\in \R}$ - and outfitted with the quotient norm that we still call $\|.\|_{W}$ by notation abuse.

\begin{lemma} \label{Wcomplete}
Then the quotient space $\dot W = W/\R$, equipped with the quotient norm 
$\|\cdot\|_W$, is complete.

Also, , 
if a sequence  
$\{u_k\}_{k=1}^\infty$ in  
$W$ and a function $u_\infty\in W$ are such that 
$\lim_{k \to +\infty} \|u_k - u_\infty\|_W = 0$, 
then there exists constants $c_k \in \R$ such that $u_k - c_k \to u$ in $L^1_{loc}(\Omega,m)$.
\end{lemma}

\begin{remark}
The measure $\mu$ does not 
play any role in this lemma, and we might be able to remove the assumption (\hyperref[H3]{H3}). However, it will be convenient to use the dyadic decomposition (and Theorem \ref{Poincare1}) given in Section \ref{Scones}.
\end{remark}

\bp We follow the arguments of \cite[Lemma 5.1]{DFMprelim}. Let $\{\dot u_k\}_{k\in \bN}$ be a Cauchy sequence in $\dot W$. We need to show that
\begin{enumerate}[(i)]
\item for every sequence $\{u_k\}_{k\in \bN}$, with $u_k \in \dot u_k$, there exists $u\in W$ and $\{c_k\}_{k\in \bN}$ such that $u_k - c_k \to u$ in $L^1_{loc}(\Omega,m)$
and $\nabla u = \lim_{k \to +\infty} \nabla u_k$ in $L^2(\Omega;m)$;
\item
 if $u_\infty, u'_\infty \in W$ 
 are such that there exist $\{u_k\}_{k\in \bN}$ and $\{u'_k\}_{k\in \bN}$ such that 
 $u_k,u'_k\in \dot u_k$ for all $k\in \bN$ and
\[\lim_{k\to \infty} \|u_k - u_\infty\|_W = \lim_{k\to \infty} \|u'_k - u'_\infty\|_W = 0,\]
then $\dot u_\infty = \dot u'_\infty$.
\end{enumerate}

First we prove (ii). 
Let $u_\infty$, $u'_\infty$, $\{u_k\}_{k\in \bN}$, and 
$\{u'_k\}_{k\in \bN}$ be as in (ii).
Notice that $\nabla (u_k - u'_k) = 0$ (in the sense of $W$) because $u_k$ and $u'_k$
represent the same class; by Lemma~\ref{Nablau} this also means that $u_k - u'_k$ is a constant.
By assumption, $u_\infty \in W$ and $\nabla u_\infty$ (in the sense of $W$) is the limit of 
$\nabla u_k$ in $L^2(\Omega;m)$. Similarly, $\nabla u'_\infty$ is the limit of $\nabla u'_k$.
Hence $\nabla u_\infty=\nabla u'_\infty$ in $L^2$, and again this means that 
$u_\infty- u'_\infty$ is constant; (ii) follows.

\medskip
Let us turn to the proof of (i). 
Pick a central point $x_0\in \Gamma$ and, for $j\in \mathbb Z$, denote by $Q^{j}$ 
the only cube in $\D_j$ that contains $x_0$.
Observe that when $j$ tends to $-\infty$, 
the sets $T_{2Q^j}$ defined by \eqref{defTQ} grow to $\Omega$, i.e., eventually
contain any compact subset of $\Omega$. Thus the convergence in each $L^1(T_{2Q^j},m)$
implies the convergence in $L^1_{loc}(\Omega,m)$.

We shall restrict our attention to $j \leq 0$.
Observe that because of Theorem \ref{Poincare1} (Poincar\'e's inequality), applied 
with $p=2$ and $D = T_{2Q^j}$,
there exists constants $C_j := C(C_1,C_2,C_4,C_6,j)$ such that for any $f\in W$, 
\begin{equation} \label{Wcompl1}
\fint_{T_{2Q^j}} |f - f^0| \, dm \lesssim C_j \int_{\Omega} |\nabla f|^2 \, dm,
\end{equation}
where we can take $f^0 = \fint_{U^*_{Q^0}} f \, dm$, i.e., take the fixed set $E = U^*_{Q^0}$,
which is contained in $T_{2Q^j}$ because $j \geq 0$.

Now let $u_k \in \dot u_k$ be as in the statement, and set
\begin{equation} \label{9.4b}
c_k = \fint_{U^*_{Q^0}} u_k \, dm.
\end{equation}
By \eqref{Wcompl1}, $(u_k -c_k)_k$ is a Cauchy sequence in $L^1(T_{Q^j})$ for each integer $j\leq0$. Hence, there exists $u^j \in L^1(T_{Q^j})$ such that $u_k - c_k$ converges to $u^j$ in $L^1(T_{Q^j})$. 
By uniqueness of the limit, 
$u^j=u^{i}$ almost everywhere on $T_{Q^j} \cap T_{Q^{i}}$, 
so we can define $u\in L^1_{loc}(\Omega,m)$ such that $u=u^j$ a.e. on $T_{Q^j}$.

It remains to check that $u\in W$
and $u_k \to u$ in $W$. 
Since $u_k-c_k \in W$, Definition \ref{defW} gives us a smooth function 
$\varphi_k\in C^\infty(\ol \Omega) \cap W$
such that
\begin{equation} \label{Wcompl2}
\fint_{U^*_{Q^0}} |u_k - c_k - \varphi_k| \, dm \leq \frac1k
\end{equation} 
and
\begin{equation} \label{Wcompl3}
\int_{\Omega} |\nabla u_k - \nabla \varphi_k|^2 \, dm \leq \frac1k.
\end{equation} 
Set $d_k=\fint_{U^*_{Q_0}} \varphi_k \, dm$.
By \eqref{Wcompl2} and \eqref{9.4b},
\begin{equation} \label{Wcompl4}
|d_k| \leq \frac1k + \Big|\fint_{U^*_{Q^0}} (u_k - c_k)\, dm \Big| = \frac1k.
\end{equation} 
We are now ready to prove that $\varphi_k \to u$ in $L^1_{loc}(\Omega,m)$. Write
\[\begin{split}
\fint_{T_{2Q^j}} |u-\varphi_k| \, dm &\leq \fint_{T_{2Q^j}} |u-(u_k-c_k)| \, dm + \fint_{T_{2Q^j}} |(u_k-\varphi_k) - (c_k-d_k) | \, dm + |d_k| \\
& := T_1 + T_2 + T_3.
\end{split}\] 
The term $T_1$ tends to 0 as $k\to \infty$ because by construction 
$u_k - c_k \to u = u^j$ in $L^1(T_{Q^j})$. 
The term $T_2$ tends also to $0$, thanks to \eqref{Wcompl1} and \eqref{Wcompl3}. The term $T_3$ converges to 0 because of \eqref{Wcompl4}. We conclude, since $T_{2Q^j} \uparrow \Omega$, that $\varphi_k \to u$ in $L^1_{loc}(\Omega,m)$; in particular
\begin{equation} \label{Wcompl5}
\lim_{k\to \infty} \int_B |\varphi_k-u| \, dm = 0 
\quad \text{ for any ball $B$ satisfying $2B \subset \Omega$.} 
\end{equation}
In addition, $\{ \dot u_k\}_{k\in \bN}$
is a Cauchy sequence in $W$, and if we combine this fact with \eqref{Wcompl3}, we get that $\{\nabla \varphi_k\}_{k\in \bN}$ is a Cauchy sequence in $L^2(\Omega,m)$. 
We conclude that there exists $v$ such that
\begin{equation} \label{Wcompl6}
\lim_{k\to \infty} \int_\Omega |\nabla \varphi_k-v|^2 \, dm = 0.
\end{equation}
We may now use our smooth functions $\varphi_k \in  C^\infty(\ol \Omega) \cap W$
to check that $u \in W$, as in Definition~\ref{defW}, and with the gradient $v$; indeed
\eqref{defW1} comes from \eqref{Wcompl5}, and \eqref{defW2} comes from \eqref{Wcompl6}.
Also, $u_k \to u$ in $W$ because $\nabla u_k \to v$ in $L^2(m)$, by
 by \eqref{Wcompl3} and \eqref{Wcompl6}.
Lemma \ref{Wcomplete} follows.
\ep

\begin{lemma}  \label{lem6.20}
Let $\{u_i\}_{i\in \bN}$ be a Cauchy sequence 
in $W$, i.e., $\|u_i -u_j\|_W \to 0$ as $i,j\to \infty$. If $u_i$ converges to $u$ in $L^1_{loc}(\Omega,m)$, then $u\in W$ and $\|u_i-u\|_W \to 0$ as $i\to \infty$.
\end{lemma}

In connection with this result, we'll
say that $\{u_i\}_{i\in \bN}$ converges to $u$ in $W$ and $L^1_{loc}(\Omega,m)$ if $\{u_i\}_{i\in \bN}$ is a Cauchy sequence in $W$ as above and $u_i \to u$ in $L^1_{loc}(\Omega,m)$.

\smallskip
\bp
Keep the same sets $Q^j$ 
as in the proof of Lemma \ref{Wcomplete} and let $\{ u_i \}$ be as in the statement.
As before, by Definition of $u_i\in W$, we can find $\varphi_i \in C^\infty(\ol\Omega) \cap W$ 
such that 
\begin{equation} \label{Wcompl7}
\fint_{U^*_{Q^0}} |u_i - \varphi_i| \, dm \leq \frac1i
\end{equation} 
and
\begin{equation} \label{Wcompl8}
\int_{\Omega} |\nabla u_i - \nabla \varphi_i|^2 \, dm \leq \frac1i. 
\end{equation} 
Set $c_i=\fint_{U^*_{Q_0}} (u_i - \varphi_i) \, dm$; by \eqref{Wcompl7}, 
$|c_i| \leq \frac1i$. 
For each fixed $j \geq 0$, and as in the proof of Lemma \ref{Wcomplete},
\[\begin{split}
\fint_{T_{2Q^j}} |u-\varphi_i| \, dm &\leq \fint_{T_{2Q^j}} |u- u_i| \, dm + \fint_{T_{2Q^j}} |(u_i-\varphi_i) - c_i | \, dm + |c_i|,
\end{split}\] 
which tends to $0$ as $ i \to \infty$. This
proves that $\varphi_i \to u$ in $L^1_{loc}(\Omega,m)$. 
Moreover, by \eqref{Wcompl8} and the fact that $\{\nabla u_i\}_{i\in \bN}$ is a Cauchy sequence in $L^2(\Omega,m)$, $\{\nabla \varphi_i\}_{i\in \bN}$ is a Cauchy sequence in $L^2(\Omega,m)$, hence there exists $v\in L^2(\Omega,m)$ such that $\nabla \varphi_i \to v$ in $L^2(\Omega,m)$. These two convergences - the convergence in $L^1_{loc}(\Omega,m)$ and the convergence of the gradients in $L^2(\Omega,m)$ - entail by definition of $W$ that $u \in W$, and by uniqueness of the gradient that $v = \nabla u$. The lemma follows.
\ep

\begin{lemma}  \label{lem6.2}
Let $\{u_i\}_{i\in \bN}$ be a sequence of functions in $W$, and let $u\in W$. If $u_i$ converges to $u$ in both $W$ and $L^1_{loc}(\Omega,m)$, then $\Tr u_i$ converges to $\Tr u$ in $H$ and $L^2_{loc}(\Gamma,\mu)$. 
\end{lemma}

\bp The convergence of the traces in $H$ is a direct consequence of Theorem \ref{TraceTh}
and the convergence of the initial sequence in $W$,
and actually does not need the convergence in $L^1_{loc}(\Omega,m)$.

The convergence of the traces in 
$L^2_{loc}(\Omega,m)$ 
is the analogue of (5.16) in \cite{DFMprelim}. 
Let us write $g$ for $\Tr u$ and $g_i$ for $\Tr u_i$. 
Since the operator $\Tr$ is linear, without loss of generality, we can assume that $u \equiv 0$ 
and thus $g\equiv 0$. So we want to prove that $\{ g_i \}$ converges to 0 in 
$L^2_{loc}(\Gamma,\mu)$.
That is,  
 if $x_0$ is a fixed point in $\Gamma$ and $Q_{j}$ is the only set in 
$\mathbb D_{j}$ containing $x_0$, we want to show that for  $j\in \bN$ and $\epsilon >0$, there exists $i_0\in \bN$ such that 
\begin{equation} \label{claim6.3}
\int_{Q_{j}} |g_i|^2 d\mu \leq \epsilon \qquad \text{ for } i \geq i_0.
\end{equation}
We introduce $g_i^k:=\Tr_k u_i$, where $\Tr_k$ is defined in the proof of Theorem \ref{TraceTh}. 
Then
\[ \begin{split}
\int_{Q_{j}} |g_i|^2 d\mu & \leq 2 \int_{Q_{j}} |g_i-g_i^k|^2 d\mu + 2 \int_{Q_{j}} |g_i^k|^2 d\mu \\
& \lesssim  \frac{2^{-2j} \mu(Q_j)}{m(U^*_{Q_j})} 2^{-(k-j)\epsilon} \int_{\Omega} |\nabla u_i|^2 dm + \int_{x\in Q_{j}} \left| \fint_{B_x^k} u_i \right|^2 d\mu(x) 
\
:= T_1 + T_2.
\end{split}\]
where we invoke \eqref{Trthe} for the second line 
and $B_x^k$ is the ball used to define $\Tr_k u(x)$. 
The values of $\|u_i\|_{W}^2 = \int_{\Omega} |\nabla u_i|^2 \, dm$ are uniformly bounded, 
since $\{ u_i \}$ converges in $W$. 
So we can fix $k$, so large that 
that $T_1 \leq \epsilon/2$ uniformly in $i\in \bN$. 
As for $T_2$, observe that
\[T_2 \leq C_{j,k} \int_{E_{j,k}} |u_i|\]
where $E_{j,k} = \bigcup_{Q\in \mathbb D_{k}:\, Q \subset Q_{j}} U^*_Q$ is relatively compact 
in $\Omega$. Since the values of $j,k$ are fixed and $\{ u_i \}$ converges to $u \equiv 0$ in $L^1(E_{j,k})$, we can choose $i_0$ such that $T_2 \leq \epsilon/2$ for 
$i\geq i_0$. 
Lemma \ref{lem6.2} follows.
\ep

\begin{lemma} \label{lem6.3}
The space 
\begin{equation} \label{defW0}
W_0 := \{u\in W, \, \Tr u = 0\}
\end{equation}
equipped with the scalar product $\left<u,v \right>_{W} := \int_\Omega \nabla u\cdot \nabla v \, dm$ (and the norm $\|.\|_W$) is a Hilbert space.
\end{lemma}

\bp Notice that $\|.\|_W$ is indeed a norm for $W_0$, because the only constant that is allowed in $W_0$ is $0$. The proof will be 
similar to \cite[Lemma 5.2]{DFMprelim}, and use Lemmas \ref{Wcomplete} and \ref{lem6.2}. 

Let $\{u_i\}_{i\in \bN}$ be a Cauchy sequence  
in $W_0$. 
By 
the proof of 
Lemma \ref{Wcomplete}, there exists $\bar u\in W$ and a sequence of constants 
$c_i = \fint_{U^*_{Q^0}} u_i \, dm$ (see \eqref{9.4b})
such that 
\begin{equation} \label{lem6.3b}
u_i - c_i \to \bar u \quad \text{ in } L^1_{loc}(\Omega,m).
\end{equation} 
Let us prove that $\{ c_i \}$ is a Cauchy sequence in $\R$. For $i,j\geq 0$,
\[
|c_i-c_j| \leq 
\fint_{U^*_{Q^0}} 
|u_i-u_j| \, dm \leq 
C \Big\{ \fint_{U^*_{Q^0}}
|\nabla u_i - \nabla u_j|^2 \, dm \Big \}^{1/2}
\leq C
\|u_i-u_j\|_{W}, 
\]
where the second inequality is due to the Poincar\'e inequality (Theorem \ref{Poincare1}). 
So $\{ c_i \}_{i\in \bN}$ is indeed a Cauchy sequence, and thus converges to a constant 
$c\in \R$. Define $u \in W$ as $\bar u - c$;
then \eqref{lem6.3b} says that 
$u_i \to u$ in $L^1_{loc}(\Omega,m)$, but the convergence also holds in $W$ 
by definition of $\bar u$. Lemma \ref{lem6.2} implies now that $\Tr u$ is the limit in 
$L^2_{loc}(\Gamma,\mu)$ of $\Tr u_i \equiv 0$, that is $\Tr u \equiv 0$ and hence $u \in W_0$. 
The lemma follows.
\ep

Recall from Lemma \ref{lmult} that $u \varphi \in W$ when $u\in W$ and $\varphi \in C^\infty_0(\R^n)$,
and that we have the product rule $\nabla (u \varphi) = \varphi \nabla u + u \nabla \varphi$ for
its derivative. Also, the trace of $u\varphi$ is $\varphi Tr u$. We can use this to prove that
$C^\infty_0(\Omega)$ is dense in $W_0$, as in the following lemma.

\begin{lemma} \label{lem6.4}
The completion of $C^\infty_0(\Omega)$ for the norm $\|.\|_W$ is $W_0$.

Moreover, if $E \subset \R^n$ is an open set and $u\in W$ is compactly supported in $E \cap \overline \Omega$, then $u$ can be approximated in the norm $\|.\|_{W}$ by functions in $C^\infty_0(E \cap \Omega)$.
\end{lemma} 

\bp This result is entirely similar to \cite[Lemma 5.5]{DFMprelim} and we refer to it for a complete proof. The main steps are: 
\begin{enumerate}[(i)]
\item we use cut-off functions $\varphi_r$ to approach $u \in W_0$ by functions that are equal to 0 on $\Gamma_r := \{X\in \Omega, \, \delta(X) \leq r\}$,
\item we use cut-off functions $\phi_R$ to approach the functions $u\varphi_r$ obtained in (i) by functions compactly supported in $\Omega$,
\item we use a mollifier to smooth the functions $u\varphi_r\phi_R$.
\end{enumerate}
And obviously, in order to deal with the functions $u\varphi_r$ or $u\varphi_r \phi_R$, 
we use in a crucial manner 
the aforementioned Lemma \ref{lmult}.
\ep

\begin{lemma} \label{lem6.5}
The set $C^\infty(\Omega) \cap C^0(\overline{\Omega}) \cap W$ is dense in $W$.
That is, for any $u \in W$, there exists a sequence $\{u_i\}_{i\in \bN}$  
in $C^\infty(\Omega) \cap C^0(\overline{\Omega}) \cap W$ 
such that $\{ u_i \}$ converges to $u$ pointwise a.e. and in $L^1_{loc}(\Omega,m)$, and
\[\|u_i - u\|_{W} \longrightarrow 0 \quad \text{ as } i \to +\infty.\]
\end{lemma}

\bp In \cite{DFMprelim}, the analogue of this result is given by \cite[Lemma 5.3]{DFMprelim}, but we cannot follow the same approach here (in \cite{DFMprelim}, the functions we considered were in $L^1_{loc}(\R^n)$, and thus allowed us to simply use a mollifier).

\smallskip
However, most of the job is already done by Lemmas \ref{densityH} and \ref{lem6.4}. 
We take $u \in W$, and we want to find a smooth approximating sequence $\{ u_i \}$. 
First write $u=v+w$ where
\[w = \Ext \circ \Tr u \quad \text{ and } \quad v = u - w.\]
Theorems \ref{TraceTh} and \ref{ThExt} imply that
$w$ - and thus $v$ - lies in $W$. Moreover, thanks to \eqref{Ext2},  
$\Tr w = \Tr u$ and hence $\Tr v = 0$;
that is, $v \in W_0$. 

Thanks to Lemma \ref{lem6.4}, we can find 
a sequence $\{ v_i \}_{i\geq 0}$
in $C^\infty_0(\Omega) \subset C^\infty(\Omega) \cap C^0(\overline{\Omega})\cap W$ 
such that $\|v_i- v\|_{W}$ tends to $0$. 
And since $W_0$ continuously injects  
in $L^1_{loc}(\Omega,m)$ (
by Theorem~\ref{PoincareTh2}), the sequence $\{ v_i \}$ converges also in 
$L^1_{loc}(\Omega,m)$ and, up to a subsequence, 
pointwise a.e. 
This takes care of $v$. 

We use Lemma \ref{densityH} to approximate $\Tr u$ by some functions $(g_i)_{i\geq 0}$ in $C^\infty(\R^n) \cap H$. Then we construct $w_i \in C^\infty(\Omega)$ as $\Ext g_i$. 
Thanks to Theorem \ref{ThExt}, $\|w_i - w\|_{W}$ tends to $0$ 
as $i$ goes to $+\infty$. Besides it is easy to check from the definition of the extension operator $\Ext$ that the convergence of $g_i$ to $\Tr u$ in $L^1_{loc}(\Gamma,\mu)$ (also given by Lemma \ref{densityH}) implies that $w_i$ converges to $w$ uniformly on an compact subsets of $\Omega$, 
and thus also pointwise a.e. 
and in $L^1_{loc}(\Omega,m)$. 

If we set $u_i = v_i + w_i$, we showed above the right convergences 
(in $W$, $L^1_{loc}(\Omega)$, and a.e. pointwise) of $u_i$ to $u$. 
The only unproved fact is  
that $w_i$ is continuous  
up to the boundary, that is $w_i \in  C^0(\overline{\Omega})$. 
We skip this part because it  
is very classical (see for instance in Section VI.2.2 of 
\cite{Stein}).
\ep

The next result states some basic properties of 
the derivative of $f\circ u$ when $u\in W$ (chain rule), and the fact that 
$uv$ lies in $W\cap L^\infty(\Omega)$ as soon as 
$u$ and $v$ both 
lie in $W\cap L^\infty(\Omega)$.

\begin{lemma} \label{chainrule}
The following properties hold:
\begin{enumerate}[(i)]
\item Let $f\in C^1(\R)$ be such that $f'$ is bounded, and let $u\in W$. Then $f\circ u \in W$, 
\[\nabla (f\circ u) = f'(u) \nabla u, \quad \text{ and } \quad \Tr (f\circ u) = f\circ (\Tr u),\]
where the last two equalities hold in the $m$-a.e. and $\mu$-a.e. sense, respectively.
\item Let $u,v\in W$. Then $\max\{u,v\}$ and $\min\{u,v\}$ lie in 
$W$,
\[ \nabla \max\{u,v\}(x) = \left\{\begin{array}{ll} 
\nabla u(x) & \text{ if } u(x) \geq v(x) \\
\nabla v(x) & \text{ if } v(x) \geq u(x), \\
\end{array}\right.\]
\[ \nabla \min\{u,v\}(x) = \left\{\begin{array}{ll} 
\nabla u(x) & \text{ if } u(x) \leq v(x) \\
\nabla v(x) & \text{ if } v(x) \leq u(x), \\
\end{array}\right.\]
\[ \Tr (\max\{u,v\}) = \max\{ \Tr u, \Tr v\}, \quad \text{ and } \quad \Tr (\min\{u,v\}) = \min\{ \Tr u, \Tr v\},\]
where the first two equalities hold $m$-a.e., and the last two $\mu$-a.e. 
\item If $\{u_k\}_{k\in \bN}$, $\{v_k\}_{k\in \bN}$ are two sequences of functions in $W$ that 
converge
to $u$, $v\in W$ both in $L^1_{loc}(\Omega,m)$, pointwise a.e., 
and in $W$ (that is, $\|u_k-u\|_W + \|v_k-v\|_W \to 0$ as $k\to \infty$),
then $\max\{u_k,v_k\}$ and $\min\{u_k,v_k\}$ lie in $W$ and converge to 
$\max\{u,v\}$ and $\min\{u_k,v_k\}$ in $L^1_{loc}(\Omega,m)$, pointwise a.e., 
and in $W$. In addition,
$\Tr \max\{u_k,v_k\}$ tends to $\max\{\Tr u,\Tr v\}$ and 
$\Tr \min\{u_k,v_k\}$ tends to $\min\{\Tr u,\Tr v\}$, in both case in $L^2_{loc}(\Gamma,\mu)$.
\end{enumerate}
\end{lemma}

\bp Point (i) and (ii) are the analogues of Lemmas 6.1 in \cite{DFMprelim}. The proof is the same as in \cite{DFMprelim} (which is itself based on the proof of results 1.18 to 1.23 in \cite{NDbook}), and strongly relies on Lemma \ref{lem6.5} (the approximation of elements in $W$ by smooth functions) and 
Lemma \ref{lem6.2} (the convergence in $W$ implies the convergence of traces). The conclusion (iii) is an intermediate result for (ii), proved as (6.16) and (6.17) in \cite{DFMprelim}.
\ep

\begin{lemma} \label{productrule}
Let $u,v\in W \cap L^\infty(\Omega)$. Then $uv\in W \cap L^\infty(\Omega)$, with
$\nabla[uv] = v \nabla u + u \nabla v$, and $\Tr (uv) = \Tr u\cdot \Tr v$.
\end{lemma}

\bp
If $u$ or $v$ is the zero constant, there is nothing to prove. Otherwise, we can divide 
$u$ and $v$ by their respective $L^\infty$ norm, and thus, without loss of generality, 
we can assume that $\|u\|_\infty = \|v\|_\infty = 1$.

By Lemma \ref{lem6.5}, we can find two sequences $\{\tilde u_k\}_{k\in \bN}$ and 
$\{\tilde v_k\in \bN\}$
in $C^\infty(\Omega) \cap C^0(\overline\Omega) \cap W$ such that $\tilde u_k \to u$ 
and $\tilde v_k \to v$ in $W$, $L^1_{loc}(\Omega,m)$, and pointwise a.e. 
By  
(iii) of Lemma \ref{chainrule}, the truncated functions 
\[u_k:= \max\{-1,\min\{1,\tilde u_k\}\} \quad \text{ and } v_k:= \max\{-1,\min\{1,\tilde v_k\}\}\]
lie in $C^0(\overline\Omega) \cap W$, are locally Lipschitz, and converge to respectively $u$ and $v$ in $W$, $L^1_{loc}(\Omega,m)$, and pointwise a.e.

Since the derivative is a local object,  
we can use \eqref{NablaLip} and the classical product rule to say that 
\[\nabla[u_kv_k] = u_k \nabla v_k + v_k \nabla u_k.\]
We conclude by showing, as in the proof of \cite[Lemma 6.3]{DFMprelim} that $u_kv_k \to uv$ in $L^1_{loc}(\Omega,m)$, $u_k \nabla v_k + v_k \nabla u_k \to u\nabla v+ v\nabla u$ in $L^2(\Omega,m)$, and $u_kv_k = \Tr(u_kv_k) \to \Tr u \cdot \Tr v$ in $L^1_{loc}(\Gamma,\mu)$. The lemma follows then from Lemma \ref{lem6.20} and Theorem \ref{TraceTh}.
\ep

\section{The localized versions $\WW(E)$ of our energy space $W$}
\label{SLocal}

The aim of this short section is to define local versions of $W$, which will be useful to study
local solutions to our degenerate elliptic equations.
As in the previous section, we assume throughout that $(\Omega,m,\mu)$ 
satisfies (\hyperref[H1]{H1})--(\hyperref[H6]{H6}).

In general, we want to localize $W$ with an open set $E'$ of $\R^n$, we set 
\begin{equation} \label{defEitself}
E = E' \cap \ol\Omega,
\end{equation}
and define the space of functions $\WW(E)$ by
\begin{equation} \label{defWWE2}
\WW(E) := \{u \in L^1_{loc}(E \cap \Omega,m):\, \varphi u \in W \text{ for all } \varphi \in C^\infty_0(E')\}.
\end{equation}
It is natural to call this space $\WW(E)$, as opposed to $\WW(E')$, because it does not depend 
on the part of $E'$ that leaves away from $\ol\Omega$. But there is an important special case,
when $E' \subset \Omega$ and so $E = E'$ is an open subset of $\Omega$.
In this case, the information that $f\in \WW(E)$ does
not give any control on $f$ at the boundary $\dr E$ (which may intersect $\Gamma$),
and $\WW(E)$ 
will be mainly used to give interior estimates for weak solutions (that will be defined soon).
In the general case, $E$ may contain pieces of the boundary $\Gamma$, and then the fact that
$f \in \WW(E)$ gives some information on the behavior of $f$ near $E \cap \Gamma$, in the same
way as the fact that $f \in W$ gives a global information on $f$ near $\Gamma$.
For instance, we can can take for $E$ (the interior in $\ol\Omega$ of) the set
$T_Q \cup Q$, for some dyadic cube $Q\in \D$.
Obviously 
$\WW(E) \subset \WW(F)$ when 
$F\subset E$, and in particular $\WW(T_Q \cup Q) \subset \WW(T_Q)$. 
In addition, if $F \subsetneq E$, it is not very hard to find a function 
$u\in \WW(F) \setminus \WW(E)$ - just make $|\nabla u(X)|$ blows up when $X$ 
gets close to $E \setminus F$ - and thus the local spaces $\WW(E)$ are all different. 
Thus for instance 
\[ W \subsetneq \WW(\overline{\Omega}) \subsetneq \WW(\Omega),\]
smooth functions on $\Omega$ that possibly explode along $\Gamma$ lie in the last space,
while they only lie in $\WW(\overline{\Omega})$ when they are locally controlled 
near $\Gamma$, and they only lie in $W$ when in addition their gradient lies in $L^2(\Omega,m)$.

Functions in $\WW(E)$ are not necessarily in $L^1_{loc}(E)$ (see Section \ref{SecW} where we defined $W$). They still have the a notion of gradient - that may be different from the distributional gradient - inherited from $W$. Indeed, if $E'$ is an open subset of $\R^n$ such that $E = E' \cap \overline{\Omega}$, consider $K$ any compact subset of $E'$ and take $\varphi_K \in C^\infty(E')$ such that $\varphi_K = 1$ on $K$, then we construct the $W$-gradient of $u\in \WW(E)$ on $K$ as the $W$-gradient of $\varphi_K u$. As an easy consequence, for $u\in \WW(E)$, we have $\nabla u \in L^2_{loc}(E,m)$ 
(where in fact we just integrate on $E \cap \Omega$, 
but local means in terms of the open set $E'$, or $E= E' \cap \overline{\Omega}$) 
and then $u \in L^2_{loc}(E,m)$ by Theorem~\ref{Poincare1}.
These observations are summarized in the next lemma.

\begin{lemma} \label{WWElem}
Let $E = E' \cap \overline{\Omega}$, for some open set $E' \subset \R^n$. 
Then every function $u \in \WW(E)$ lies in 
 $L^1_{loc}(E, m)$, and its gradient lies in $L^2_{loc}(E,m)$.
\end{lemma}

\begin{remark} 
We don't have many doubts that the reverse inclusion
\begin{equation} \label{WWEa}
\WW(E) \supset \{ u \in L^1_{loc}(E,m):\, \nabla u \in L^2_{loc}(E,m)\}.
\end{equation}
also holds.
The idea of the proof of \eqref{WWEa} would be to take $u\in L^1_{loc}(E)$ 
that satisfies $\nabla u \in L^2_{loc}(E,m)$, 
and $\varphi \in C^\infty_0(E')$. We would set then $K':= \supp\, \varphi$ which has a smooth boundary, and we would say that $u|_{K'}$ can be extended to a function $\bar u \in W$ such that $\bar u = u$ a.e. on $K'$. Then we would use Lemma \ref{lmult} in order to show that $\varphi u = \varphi \bar u \in W$. The problem with this proof is that we don't know any reference for the extension theorem needed to built $\bar u$ in weighted Sobolev spaces (an analogue of \cite[Section 1.1.17]{Mazya11} in the unweighted case), and we do not want to spend time on something that we will not need.
\end{remark}

The next lemma allows us to speak about
traces for functions in the local Sobolev spaces 
$\WW(E)$.

\begin{lemma} \label{defTrWE}

Let $E'\subset \R^n$ be open, and set $E = E' \cap \overline{\Omega}$ as in \eqref{defEitself}.
For every function $u\in \WW(E)$, we can define the trace of $u$ on $\Gamma \cap E$ by
\begin{equation} \label{defTrWE2}
\Tr u (x) = \lim_{X \in \gamma(x) \atop \delta(X) \to 0}  \fint_{B(X,\delta(X)/2)} u
\qquad \text{ for $\mu$-almost every } x\in \Gamma \cap E,
\end{equation}
and $\Tr u \in L^2_{loc}(\Gamma \cap E',\mu)$.  
Moreover, for every choice of $f\in \WW(E)$ and 

$\varphi \in C^\infty(E')$, $\varphi u \in \WW(E)$ and
 \begin{equation} \label{defTrWE1}
\Tr (\varphi u)(x) = \varphi(x) \Tr u(x)  \ \text{ for $\mu$-almost every } x\in \Gamma \cap E.
 \end{equation}
\end{lemma}

\smallskip
\bp 
None of this is too surprising; the trace is a local notion, and $\WW(E)$ is designed to ba a local space.
Let $E'$, $E$, and $f$ be as in the statement, and let $B$ a compact ball in $E'$, and choose
$\psi \in C^\infty_0(E)$ such that $\psi \equiv 1$ near $B$. 
Then $\psi u \in W$ by \eqref{defWWE2}, and the analogue of \eqref{defTrWE2} 
for $\psi u$ comes
with the construction of the trace. This implies the existence
of the same limit for $u$, almost everywhere in $\Gamma \cap B$.

In addition, since $\psi u \in W$, Theorem \ref{TraceTh} says that
$Tr(\psi u) \in H$, and then $Tr(u) = Tr(\psi u) \in L^2(B,d\mu)$
(see the definition \eqref{defH}). Therefore $Tr(u) \in L^2_{loc}(E',d\mu)$, as announced.

The fact that  $\varphi u \in \WW(E)$ when $u\in \WW(E)$ and $\varphi \in C^\infty(E')$ comes
right from the definition \eqref{defWWE2} and Lemma \ref{lmult}, and \eqref{defTrWE1}
is immediate because when $B$ and $\psi$ are as above and $\mu$-almost everywhere on $B$,
\[
Tr(\varphi u) = Tr(\psi^2 \varphi u) = \psi \varphi Tr(\psi u)  = \varphi Tr(u)
\]
by \eqref{defTrWE2}, the formula for the trace of a product of $\psi u \in W$ and 
a $\psi \varphi \in C^\infty_0(\R^n)$, and the fact that in $B$, the formula \eqref{defTrWE2} does not
see the cut-off $\psi$.
\ep

\section{Definitions of solutions and their properties}
\label{SSolutions}

We now have all the functional analysis needed 
to deal with the main goal of this section, which is to
define weak solutions to appropriate degenerate elliptic operators, 
and give their first properties. We will
follow Section 8 in \cite{DFMprelim} (which itself
copies the frame of the first sections of \cite{KenigB}), and we will 
refer to \cite{DFMprelim} for most
of the proofs. As in the previous section, we systematically assume 
that $(\Omega,m,\mu)$ 
satisfies (\hyperref[H1]{H1})--(\hyperref[H6]{H6}).

Recall that we intend to work with the degenerate elliptic operators $L = - \div A \nabla$,
where the matrix $A: \Omega \to \mathbb M_n(\R)$ satisfies the following elliptic and boundedness conditions:
\begin{equation} \label{defEllip}
A(X) \xi \cdot \xi \geq C_A^{-1} w(x) |\xi|^2 \quad \text{ for } X\in \Omega \text{ and } \xi \in \R^n
\end{equation}
and
\begin{equation} \label{defBdd}
A(X) \xi \cdot \zeta \leq C_A w(x) |\xi||\zeta| \quad \text{ for } X\in \Omega \text{ and }  \xi,\zeta \in \R^n,
\end{equation}
where 
$w$ is the weight associated to the measure $m$ given in as part of (\hyperref[H4]{H4}). 
We shall also
use the ``normalized'' matrix $\A:= w^{-1}A$ which satisfies the unweighted ellipticity 
and boundedness conditions
\begin{equation} \label{defEllip2}
\A(X) \xi \cdot \xi \geq C_A^{-1} |\xi|^2 \quad \text{ for } X\in \Omega \text{ and } \xi \in \R^n
\end{equation}
and
\begin{equation} \label{defBdd2}
\A(X) \xi \cdot \zeta \leq C_A |\xi||\zeta| \quad \text{ for } X\in \Omega \text{ and }  \xi,\zeta \in \R^n.
\end{equation}
We introduce the bilinear form $a$ defined by
\begin{equation} \label{defauv} \begin{split}
a(u,v):= \int_\Omega A \nabla u \cdot \nabla v & = \int_\Omega \A \nabla u \cdot \nabla v \, dm \\
& \qquad \text{ for any $u,v$ that satisfies } \int_\Omega |\nabla u| |\nabla v| \, dm < +\infty.
\end{split} \end{equation}
The conditions \eqref{defEllip}--\eqref{defBdd} entail that $a$ is 
bounded on $\dot W\times \dot W$ (the homogeneous quotient space) and coercive on $\dot W$, i.e.,
\begin{equation}
a(u,u) \geq C_A^{-1} \|u\|_W^2 \quad \text{ and } \quad  a(u,v) \leq C_A \|u\|_W \|v\|_{W} \qquad \text{ for } u,v\in W.
\end{equation}
It is also coercive on $W_0$ (no need to take a quotient, because $W_0$ does not contain nontrivial
constant functions).

\begin{definition}
Let $E \subset \Omega$ be a open set. We say that $u \in \WW(E)$ is a (weak) 
solution to $Lu = 0$ in $E$ when 
\begin{equation} \label{defsol}
a(u,\varphi) = \int_\Omega A \nabla u \cdot \nabla \varphi = \int_\Omega \A \nabla u \cdot \nabla \varphi \, dm = 0 \qquad \text{ for any } \varphi \in C^\infty_0(\Omega).
\end{equation}
Similarly 
 $u \in \WW(E)$ is a subsolution (respectively supersolution) to $Lu = 0$ in $E$ when
\begin{equation} \label{defsubsol}
a(u,\varphi) \leq 0 \ \text{(resp. $\geq 0$)} \qquad \text{ for any } \varphi \in C^\infty_0(\Omega) \text{ that satisfies } \varphi \geq 0.
\end{equation}
\end{definition}

\ms
In the rest of the section, we present the analogues of the results in \cite[Section 8]{DFMprelim}, and we discuss the differences in the proofs when needed.

The first result enlarges the class of possible test functions.

\begin{lemma}  \label{rdefsol}
Let $E \subset \Omega$ be an open set and 
let $u\in \WW(E)$ be a solution to 
$Lu=0$ in $E$.  We write $E^\Gamma$ for $E \cup (\Gamma \cap \dr E)$, that is, $E^\Gamma$ is 
the union of $E$ with the part of its boundary that intersects $\Gamma$. The identity \eqref{defsol} holds: 
\begin{itemize}
\item when $\varphi \in W_0$ is compactly supported in $E$;
\item when $\varphi \in W_0$ is compactly supported in $E^\Gamma$ \ub{and} $u \in \WW(E^\Gamma)$;
\item when $E = \Omega$, $\varphi \in W_0$, \ub{and} $u\in W$.
\end{itemize}
In addition, 
\eqref{defsubsol} holds when $u$ is a subsolution (resp. supersolution) in $E$, $\varphi$ is a non-negative test function, and the couple $(u,\varphi)$ satisfies one of the above conditions.
\end{lemma}

\bp
See the proof of \cite[Lemma 8.3]{DFMprelim}. This lemma is a consequence of Lemma \ref{lem6.4}, that gives that the functions in $W_0$ can be approximated by smooth functions. 
\ep

The next result proves the stability of subsolutions/supersolutions under $\max$/$\min$.

\begin{lemma} \label{lstabsol}
Let $E\subset \Omega$ be an open set. 
\begin{itemize}
\item If $u,v\in \WW(E)$ are subsolutions (to $Lu=0$) in $E$, then $t = \max\{u,v\}$ is also a subsolution in $E$.
\item If $u,v\in \WW(E)$ are supersolutions in $E$, then $t = \min\{u,v\}$ is also a supersolution in $E$.
\end{itemize}

In particular if
$k\in \R$, then $(u-k)_+ := \max\{ u-k,0\}$ is a subsolution in $E$ whenever $u \in \WW(E)$ is a subsolution in $E$ and $\min\{u,k\}$ is a supersolution in $E$ whenever $u\in \WW(E)$ is a supersolution in $E$.
\end{lemma}

\bp
The proof is the same as the one of \cite[Lemma 8.23]{DFMprelim} and \cite[Theorem 3.5]{Stampacchia65}. Lemma 8.5 in \cite{DFMprelim} shows that the result can be localized into a relatively compact open subset $F$ of $E$. Theorem 3.5 in \cite{Stampacchia65} relies on the fact the bilinear form $a$ is coercive and continuous (on appropriate local spaces) and on convex analysis.   
\ep

In the sequel, the notation $\sup$ and $\inf$ are used for the essential supremum and essential infimum, 
since they are the definitions that makes sense for the functions in $W$ or in $\WW(E)$, for $E= E' \cap \overline{\Omega}$ and $E' \subset \R^n$ open. 

In addition, the expression ``$\Tr u = 0$ a.e. on $B$'', for a function $u\in \WW(B \cap \overline \Omega)$, 
means 
that $\Tr u$, which is defined on $\Gamma \cap B$ and lies in $L^1_{loc}(B\cap \Gamma,\mu)$ thanks to Lemma~\ref{defTrWE}, is equal to $0$ $\mu$-almost everywhere on $\Gamma \cap B$.
The expression ``$\Tr u \geq 0$ a.e. on $B$'' is defined similarly.

\medskip

We now state some classical regularity results inside the domain and at the boundary. 

\begin{lemma}[interior Caccioppoli inequality] \label{CaccioI} 
Let $E \subset \Omega$ be an open set, and let $u\in \WW(E)$ be a non-negative subsolution in $E$.  Then for any $\alpha \in C^\infty_0(E)$,
\begin{equation} \label{Caccio1}
\int_\Omega \alpha^2 |\nabla u|^2 dm \leq C  \int_{\Omega} |\nabla \alpha|^2 u^2 dm,
\end{equation}
where $C$ depends only upon the constant $C_A$. 

In particular, if $B$ is a ball of radius $r$ such that $2B \subset \Omega$ and $u\in \WW(2B)$ is a non-negative subsolution in $2B$, then
\begin{equation} \label{Caccio2}
\int_B |\nabla u|^2 dm \leq C r^{-2} \int_{2B} u^2 dm.
\end{equation}
\end{lemma}

\begin{lemma}[Caccioppoli inequality on the boundary] \label{CaccioB} 
Let $B \subset \R^n$ be a ball of radius $r$ centered on $\Gamma$, and let $u\in \WW(2B \cap \overline\Omega)$ be a non-negative subsolution in $2B \cap \Omega$ such that 
$\Tr u = 0$ a.e. on $2B$. Then for any $\alpha \in C^\infty_0(2B)$,
\begin{equation} \label{Caccio3}
\int_{2B \cap \Omega} \alpha^2 |\nabla u|^2 dm \leq C  \int_{2B \cap \Omega} |\nabla \alpha|^2 u^2 dm,
\end{equation}
where $C$ depends only on the constant $C_A$. In particular, we can take $\alpha \equiv 1$ on $B$ and $|\nabla \alpha| \leq \frac2r$, which gives
\begin{equation} \label{Caccio4}
\int_{B\cap \Omega} |\nabla u|^2 dm \leq C r^{-2} \int_{2B \cap \Omega} u^2 dm.
\end{equation}
\end{lemma}

\bp The proofs of the two lemmas are similar to Lemma 8.6 and Lemma 8.11 
in \cite{DFMprelim}. There is not any difficulty here, maybe it is worth saying that we use $\varphi = \alpha^2 u$, where $\alpha$ is an appropriate cut-off function; and $\varphi$ is a valid test function due to Lemma \ref{rdefsol} and, for the boundary version, Lemma \ref{defTrWE}.
\ep

Let us turn to the statement of the Moser estimates.

\begin{lemma}[interior Moser estimate] \label{MoserI}
Let $p>0$ and $B$ be a ball such that $2B \subset \Omega$. 
If $u \in \WW(2B)$ is a non-negative subsolution in $2B$, then
\begin{equation} \label{Moser1}
\sup_{B} u \leq C \left( \frac1{m(2B)} \int_{2B} u^p \, dm \right)^\frac1p,
\end{equation}
where $C$ depends on $n$, $C_4$, $C_6$, $C_A$, and $p$.
\end{lemma}

\begin{lemma}[Moser estimates on the boundary] \label{MoserB}
Let $p>0$, $B$ be a ball centered on $\Gamma$, and $u \in \WW(2B\cap \overline\Omega)$ be a non-negative 
subsolution in $2B \cap \Omega$ such that $\Tr u=0$ a.e. on $2B$. Then
\begin{equation} \label{Moser2}
\sup_{B \cap \Omega} u \leq C_p \left(m(2B)^{-1}\int_{2B \cap \Omega} |u|^p dm \right)^\frac1p,
\end{equation}
where $C_p$ depends only on $n$, $C_1$ to $C_6$, $C_A$, and $p$.
\end{lemma}

\bp The proofs for these two results are analogous to the ones of \cite[Lemmas 8.7 and 8.12]
{DFMprelim}, and relies on the so-called Moser iterations. 

What we need are Lemma \ref{lstabsol}, a Cacciopoli inequality (Lemma \ref{CaccioI} or Lemma \ref{CaccioB}, according to the version we want to prove), a Sobolev-Poincar\'e inequality (Theorem \ref{Poincare1} or Corollary \ref{PoincareCor2}, the balls in the right-hand side of \eqref{PoincareCor2a} are slightly bigger than the ones in the left-hand side, but the argument can easily be adapted), and the doubling property (\hyperref[H4]{H4}).  
\ep

The next step is the H\"older continuity of solutions. We shall give a 
few intermediate results, starting by the density lemmas.

\begin{lemma}[Density lemma inside the domain] \label{DensityI}
Let $B$ be a ball such that $4B \subset \Omega$ and $u \in \WW(4B)$ be a non-negative 
supersolution in $4B$ such that
\[m(\{X\in 2B, \, u(X) \geq 1\}) \geq \epsilon m(2B).\] 
Then 
\begin{equation} \label{Harnack33}
\inf_{B} u \geq C^{-1},
\end{equation}
where $C>0$ depends only on  $n$, $C_4$, $C_6$, $C_A$, and $\epsilon$.
\end{lemma}

\begin{lemma}[Density lemma on the boundary] \label{DensityB}
Let $B$ be a ball centered on $\Gamma$ and $u \in \WW(4B \cap \overline{\Omega})$ be a 
non-negative supersolution in 
$4B  \cap \Omega$ such that $\Tr u = 1$ a.e. on $4B$. 
Then 
\begin{equation} \label{Harnack3}
\inf_{B\cap \Omega} u \geq C^{-1},
\end{equation}
where $C>0$ depends only on  $n$, $C_1$ to $C_6$ and $C_A$.
\end{lemma}

\bp The proof of Lemma \ref{DensityI} can be copied from the one of Density Theorem (Section~4.3, Theorem 4.9) in \cite{HLbook}. The proof of Lemma \ref{DensityB} is similar to the one Lemma 8.14 
in \cite{DFMprelim} (which is itself inspired from the Density Theorem in \cite{HLbook}).

Formally, the ideas of the proof are to say that $v= -\ln u$ is a subsolution that satisfies $\Tr u = 0$ a.e. on $4B$ (if needed), and then to use Moser estimates (Lemma \ref{MoserI} or Lemma~\ref{MoserB}) and a Poincar\'e inequality (Theorem \ref{Poincare1} or Corollary \ref{PoincareCor2}) in an appropriate way. 

Of course, we need to be very careful: for instance when constructing $v$, we want to use Lemma \ref{chainrule} in order to verify that $v$ is indeed in $\WW(2B)$, yet the function $-\ln$ is not Lipschitz... But the pitfalls are the same as in the proof of \cite[Lemma 8.14]{DFMprelim}.
\ep

Next comes oscillation estimates.

\begin{lemma}[interior Oscillation estimates] \label{OscI}
Let $B$ be a ball such that $4B \subset \Omega$ and $u \in \WW(4B)$ be a solution in $4B$. Then, there exists $\eta \in (0,1)$ such that
\begin{equation} \label{Osc111}
\osc_{B} u \leq \eta \osc_{4B} u,
\end{equation}
where the constant $\eta$ depends only on on  $n$, $C_4$, $C_6$, and $C_A$.
\end{lemma}

\begin{lemma}[Oscillation estimates on the boundary] \label{OscB}
Let $B$ be a ball centered on $\Gamma$ and $u \in \WW(4B \cap \overline{\Omega})$ be a solution in $4B \cap \Omega$ such that $\Tr u $ is uniformly  bounded on $4B \cap \Gamma$. Then, there exists $\eta \in (0,1)$ such that
\begin{equation} \label{Osc1}
\osc_{B \cap \Omega} u \leq \eta \osc_{4B \cap \Omega} u + (1-\eta) \osc_{\Gamma \cap 4B} \Tr u.
\end{equation}
The constant $\eta$ depends only on on  $n$, $C_1$ to $C_6$, and $C_A$.
\end{lemma}

\bp
Lemma \ref{OscI} and Lemma \ref{OscB} can be proved respectively as Theorem 2.4 in \cite[Section 4.3]{HLbook} and as Lemma 8.15 in \cite{DFMprelim}. The proofs work as long as Lemma \ref{DensityI} or Lemma \ref{DensityB} is true. 
\ep

We shall now present the H\"older regularity of solutions.

\begin{lemma}[interior H\"older continuity] \label{HolderI}
Let $x\in \Omega$ and $R>0$ be such that $B(x,2R) \subset \Omega$, and let $u\in \WW(B(x,2R))$ be a solution to $Lu=0$ in $B(x,2R)$. Set 
\[\osc_B u:= \sup_B u - \inf_B u.\]
Then there exists $\alpha\in (0,1]$ and $C>0$ such that for any $0<r<R$,
\begin{equation} \label{Holder1}
\osc_{B(x,r)} u \leq C \left( \frac rR \right)^\alpha \left( \frac{1}{m(B(x,R))} \int_{B(x,R)} u^2 \, dm \right)^\frac12,
\end{equation}
where $\alpha$ and $C$ depend only on $n$, $C_4$, $C_6$, and $C_A$.
Hence $u$
is (possibly after modifying it on a set of measure $0$)
locally H\"older continuous with exponent $\alpha$.
\end{lemma}

\begin{lemma} \label{HolderB}
Let $B=B(x,r)$ be a ball centered on $\Gamma$ and $u \in \WW(B \cap \overline{\Omega})$ be a solution in $B \cap \Omega$ 
such that $\Tr u$ is continuous and bounded on $B \cap \Gamma$. 
There exists $\alpha >0$ such that for
$0<s <r$,
\begin{equation} \label{Holder2}
\osc_{B(x,s) \cap \Omega} u \leq C \left(\frac{s}r\right)^\alpha \osc_{B(x,r) \cap \Omega} u + C \osc_{B(x,\sqrt{sr}) \cap \Gamma} \Tr u
\end{equation}
 where the constants $\alpha,C$ depend only on $n$, $C_1$ to $C_6$, and $C_A$. In particular, $u$ (possibly after modification on a set of measure 0) is continuous on $B \cap \Omega$, can be extended by continuity on $B \cap \Gamma$, and the values of this extension on $B \cap \Gamma$ are $\Tr u$.
 
If, in addition, $\Tr u\equiv 0$ on $B$, then for any $0<s<r/2$
\begin{equation} \label{Holder3}
\osc_{B(x,s) \cap \Omega} u \leq C \left(\frac{s}r\right)^\alpha \left( \frac1{m(B)} \int_{B \cap \Omega} |u|^2 dm \right)^\frac12.
\end{equation}
\end{lemma}

\bp The proof of the two last lemmas are the same as the ones of Theorem 2.5 in  \cite[Section 4.3]{HLbook} and Lemma 8.16 in \cite{DFMprelim}. 
\ep

It remains to treat the Harnack inequality.

\begin{lemma}[Harnack Inequality] \label{HarnackI}
Let $B$ be a ball such that $2B\subset \Omega$, and let $u\in \WW(2B)$ be a non-negative solution to $Lu=0$ in $2B$. Then 
\begin{equation} \label{Harnack1}
\sup_B u \leq C \inf_B u,
\end{equation}
where $C$ depends only on  $n$, $C_4$, $C_6$, and $C_A$.
\end{lemma} 

\bp The proof in \cite{DFMprelim} uses, roughtly speaking, the condition (\hyperref[H7]{H6'}) to say that the Harnack inequality can be proved using the classical theory of uniformly elliptic operators in divergence form.

If we were to have (\hyperref[H7]{H6'}) instead of (\hyperref[H6]{H6}), we could proceed in a similar manner. 
Fortunately for us, the proof in the classical theory can easily adapted to our setting. This observation was already made in \cite{FKS}, but since our conditions are slightly weaker than \cite{FKS}, we sketch the proof to check that we don't have any extra difficulties.

\ms

\noindent {\bf Step 1:} The John-Nirenberg lemma.

Let $O$ be an open subset of $\Omega$. Suppose that $w\in L^1(O,m)$ lies 
in $BMO(O,m)$, in the sense that 
that for every ball $B \subset O$
\begin{equation} \label{JNL1}
\fint_B |w-w_B| \, dm \leq C_{JN}
\end{equation}
for a constant $C_{JN}$ independent of $B$, and where $w_B$ denotes $\fint w \, dm$. Then we claim that for any $B \subset O$, 
\begin{equation} \label{JNL2}
\fint_B \exp\left( \frac{p_0}{C_{JN}} |w-w_B| \right) \, dm \leq C,
\end{equation}
where $p_0$ and $C$ depend
only on $C_4$ (and $n$).

The claim is the John-Nirenberg lemma, whose proof uses only the Calder\'on-Zygmund decomposition 
(see for instance \cite[Chapter 3, Theorem 1.5]{HLbook}).

\ms

\noindent {\bf Step 2:} The weak Harnack inequality.

Suppose $2B \subset \Omega$ and let 
$u \in W(2B)$ be a non-negative supersolution to $Lu = 0$. Then we claim that 
there exists $p_1>0$ such that 
\begin{equation} \label{wHI1}
\inf_{B} u \geq C^{-1} \left( \fint_{2B} u^{p_1} \, dm \right)^\frac1{p_1},
\end{equation}
where $C^{-1}$ depends only on $C_4$, $C_6$ $C_A$, and $n$.

For any $\epsilon>0$, we consider
the supersolution $\bar u = u+\epsilon>0$ and then $v = \bar u^{-1}$.  
For any $\varphi \in C^\infty_0(2B)$, the function $v^2\varphi$ belongs to $W_0$ thanks to Lemmas \ref{lmult}
and \ref{chainrule}, and is compactly supported in $\Omega$. So $v^2\varphi$ can be used as a test function, by
Lemma \ref{rdefsol}, hence
\[\int_{2B} \A\nabla \bar u \cdot  \nabla[v^2\varphi]\, dm \geq 0, \] 
that is, 
\[ \int_{2B} \bar u^{-2} (\A\nabla \bar u \cdot \nabla \varphi) \, dm 
\geq 2 \int_{2B} \bar u^{-3}\varphi (\A\nabla \bar u \cdot \nabla u) \, dm\]
hence, by the positivity of $\A$ 
\begin{equation} \label{wHI2}
\int_{2B} (\A\nabla v \cdot \nabla \varphi) \, dm 
\leq -2 \int_{2B} \bar u^{-3}\varphi (\A\nabla \bar u \cdot \nabla u) \, dm \leq 0.
\end{equation}
We deduce that $v$ is a non-negative subsolution in $2B$, so using Moser's estimate (Lemma \ref{MoserI}), we get that for any $p>0$
\[\sup_{B} v \leq C_p \left( \fint_{\frac32B} v^{p} \, dm \right)^\frac1p\]
where $C_p$ depends on $C_4$, $C_6$, $C_A$, $n$, and $p$. 
Using the fact that $v = \bar u^{-1}$, 
we deduce that 
\begin{equation} \label{wHI3}
\inf_{B} \bar u \geq C_p \left( \fint_{\frac32B} \bar u^{-p} \, dm \right)^{-\frac1p}.
\end{equation}
The claim \eqref{wHI1} will be established as soon as we prove 
that for some $p_1>0$, one has
\begin{equation} \label{wHI4}
\left( \fint_{\frac32B} \bar u^{-p_1} \, dm \right) \left( \fint_{\frac32B} \bar u^{p_1} \, dm \right) 
\leq C,
\end{equation}
with a bound $C$ independent of 
 $u$ and 

the $\epsilon$ used to define $\bar u$, and we shall now 
prove \eqref{wHI4}
using the John-Nirenberg inequality.

Take $w = \log \bar u$;
we want to check that $w \in BMO(\frac32B)$. 
We test $\bar u$ against the test function $ \bar u ^{-1} \varphi^2$, 
where $\varphi \in C^\infty_0(2B)$ to obtain
\[ 2 \int_{2B} \varphi \bar u^{-1}(\A\nabla \bar u \cdot \nabla \varphi )\, dm -  \int_{2B} \varphi^2 \bar u^{-2}(\A\nabla \bar u \cdot \nabla \bar u)\, dm \geq 0. \]
We use the fact that $\nabla w = \bar u^{-1} \nabla u$ and the ellipticity conditions \eqref{defEllip2}--\eqref{defBdd2} to obtain
\[\int_{2B} \varphi^2 |\nabla w|^2 \, dm \leq C \int_{2B} \varphi |\nabla w| |\nabla \varphi| \, dm,\]
which implies, by the Cauchy-Schwarz inequality, that 
\begin{equation} \label{wHI5}
\int_{2B} \varphi^2 |\nabla w|^2 \, dm \leq C \int_{2B} |\nabla \varphi|^2 \, dm.
\end{equation}
For any ball $B'\subset \frac32B$ of radius $r'$, we can built a smooth function $\varphi$ such that $\varphi \equiv 1$ on $B'$, $\varphi \equiv 0$ on $\frac{9}{8}B'$, and $|\nabla \varphi| \leq 10/r'$. Using those test functions in \eqref{wHI5} gives that for any $B' \subset \frac32B$,
\begin{equation} \label{wHI6}
\fint_{B'} |\nabla w|^2 \, dm \leq C (r')^{-2},
\end{equation}
where $C$ depends only on $C_A$. The assumption (\hyperref[H6]{H6}), i.e. the Poincar\'e inequality, 
infers now that
\[\fint_{B'} |w-w_{B'}| \, dm \leq C,\]
as in \eqref{JNL1}. From step 1, the inequality \eqref{JNL2} thus holds, that is we can find a $p_1>0$ such that 
\begin{equation} \label{wHI7}
\fint_{\frac32 B} \exp\left( p_1 |w-w_{\frac32 B}| \right) \, dm \leq C.
\end{equation}
We are now ready for the proof of \eqref{wHI4}. Indeed, just observe that
\[\begin{split}
\left( \fint_{\frac32B} u^{-p_1} \, dm \right) & \left( \fint_{\frac32B} u^{p_1} \, dm \right) 
 = \left( \fint_{\frac32B} \exp(-p_1 w) \, dm \right) \left( \fint_{\frac32B} \exp(p_1 w) \, dm \right) \\
& = \left( \fint_{\frac32B} \exp(-p_1 [w-w_{\frac32 B}]) \, dm \right) \left( \fint_{\frac32B} \exp(p_1 [w-w_{\frac32 B}]) \, dm \right) \\
& \leq \left( \fint_{\frac32B} \exp(p_1 |w-w_{\frac32 B}|) \, dm \right)^2 \lesssim 1,
\end{split}\]
by \eqref{wHI7}.

\ms

\noindent {\bf Step 3:} Conclusion.

We combine \eqref{wHI1} with Lemma \ref{MoserI} - the Moser inequality inside the domain, 
to get the desired Harnack inequality
\[\sup_B u \leq C \inf_B u\]
We should require $B$ to satisfy $4B \subset \Omega$, but we can easily solve 
this issue by covering $B$ by balls $B'$ of smaller radius that satisfy
$4B' \subset 2B \subset \Omega$.
\ep

We shall also need the following version of the Harnack inequality, which will be useful to define the harmonic measure.

\begin{lemma} \label{HarnackI2}
Let $K$ be a compact subset of $\Omega$ and let $u\in \WW(\Omega)$ be a non-negative 
solution in $\Omega$. Then 
\begin{equation} \label{Harnack2}
\sup_K u \leq C_K \inf_K u,
\end{equation}
where $C_K$ depends only on 
$n$, $C_1$, $C_2$, $C_4$, $C_6$, $C_A$, and $\diam\, K /\dist(K,\Gamma)$.
\end{lemma} 

\bp The proof is the same as the one of \cite[Lemma 8.10]{DFMprelim}. The topological conditions (\hyperref[H1]{H1})--(\hyperref[H2]{H2}) allow us to connect any couple of points in $K$ by a chain of balls that says away from the boundary (see Proposition \ref{propHarnack}). 
The length of the chain can be bounded by a constant depending only on $\diam K/\dist(K,\dr \Omega)$.  We then use the Harnack inequality above on those balls.
\ep

\section{Construction of the harmonic measure}
\label{SHarm}

We follow Section 9 in \cite{DFMprelim} and, as in the previous section, we will refer to \cite{DFMprelim} when the proofs do not require any new argument.

\medskip

The objective of the section is, as the title suggests, to construct a harmonic measure associated to our degenerate
operator $L = - \div A \nabla$ that still satisfies \eqref{defEllip}--\eqref{defBdd}. By harmonic measure, 
we mean
a family of measures $\omega^X_L$, where $X\in \Omega$ is called pole of the harmonic measure, such that for any Borel set $E \subset \Gamma$, the function $u_E$ defined as $u_E(X) = \omega^X_L$ solves
the Dirichlet problem
\begin{equation} \label{Dphm}
\left\{
\begin{array}{ll}
Lu_E=0 & \text{ in } \Omega \\
u_E= \1_E & \text{ on } \Gamma.
\end{array}
\right.
\end{equation}
But
\eqref{Dphm} does not make a lot of sense for the moment. The part ``$Lu_E = 0$ in $\Omega$'' is easy to interpret: we want $u_E$ to lie in $\WW(\Omega)$ and to be a solution to $Lu=0$ in $\Omega$. 
The part ``$u_E = \1_E$ on $\Gamma$'' is harder to understand: we could hope that the meaning is 
$\Tr u_E = \1_E$ $\mu$-a.e. on $\Gamma$, but it is unclear that it is possible at this point.

Another issue
is the uniqueness: just take $\Omega = \R^n_+$ - i.e. $\Gamma = \R^{n-1}$ - and $E = \emptyset$ - i.e. $\1_E \equiv 0$ - and we can find at least two solutions ($u_1\equiv 0$ and $u_2 = \delta$) to $Lu=0$ that satisfy both $\Tr u = 0$. Imposing that $u$ lies in $W$ is not immediately
possible, since characteristic functions of non-trivial sets do not always lie 
in $H$.

Our salvation will come  
from the maximum principle. And instead of \eqref{Dphm}, we shall say that the harmonic measure $\omega_L^X$ is built such that for any $g\in C^\infty_0(\R^n)$, the function defined by
\[u(X) = \int_\Gamma g(y) d\omega^X_L(y)\]
lies in $W$, is a solution to $Lu = 0$, and satisfies $\Tr u = g$. 
Let us now give a full presentation. 

\medskip

We say that $f\in W^{-1}$ if $f$ is a linear form on $W_0$ that satisfies 
\[|\left< f, v \right>_{W^{-1},W_0}| \leq C_f \|v\|_{W},\] 
where we anticipate slightly and denote by $\left< f, v \right>_{W^{-1},W_0}$
the effect of $f$ on $v$.
The best constant $C_f$ in the inequality above is denoted $\|f\|_{W^{-1}}$.

Let us give first the existence and uniqueness of solutions $u \in W$ to 
$Lu= f$ and $\Tr u = g$, where $f\in W^{-1}$ and $g\in H$ are given. 

\begin{lemma} \label{lLM}
For any $f\in W^{-1}$ and any $g\in H$, there exists a unique $u \in W$ such that 
\begin{equation} \label{DpL1}
\int_\Omega \A \nabla u \cdot \nabla v \, dm = \left< f,v \right>_{W^{-1},W_0} \qquad \text{ for all } v\in W_0,
\end{equation}
and
\begin{equation} \label{DpL2}
\Tr u=g  \text{ a.e. on } \Gamma.
\end{equation}
Moreover, there exists $C>0$,
independent of $f$ and $g$,
 such that
\begin{equation} \label{ubygf}
\|u\|_W \leq C (\|g\|_H + \|f\|_{W^{-1}}).
\end{equation}
\end{lemma}

\bp The lemma follows from the extension theorem (Theorem \ref{ThExt}) and the Lax-Milgram theorem. Details are given in the proof of \cite[Lemma 9.1]{DFMprelim}.
\ep

The next result needed is a maximum principle. In its weak form, the maximum principle is as follows.

\begin{lemma} \label{lMP}
Let $u \in W$ be a supersolution in $\Omega$ satisfying $\Tr u \geq 0$ $\mu$-a.e. on $\Gamma$. Then $u\geq 0$ a.e. in $\Omega$.
\end{lemma} 

\bp
Take $v:= \min\{u,0\} \leq 0$;
we want 
to prove that $v \equiv 0$. 
Lemma \ref{chainrule} allows us to say that $v\in W$, $\nabla v = \nabla u \1_{\{u<0\}}$, and $\Tr v = 0$ a.e. in $\Gamma$. In particular $v \in W_0$, which makes $v$ a valid test function to be tested against the supersolution $u$ (see Lemma \ref{rdefsol}). This gives
\begin{equation} \label{lMP3}
0 \geq \int_\Omega \A \nabla u \cdot \nabla v \, dm = \int_{\{u<0\}} \A \nabla u \cdot \nabla u \, dm = \int_\Omega \A \nabla v \cdot \nabla v \, dm \geq C_A^{-1} \|\nabla v\|_{W}^2 \geq 0.
\end{equation}
that is, $\|\nabla v\|_{W} = 0$.
Yet, $\|.\|_W$ is a norm on $W_0 \ni v$, hence $v = 0$ a.e. in $\Omega$.  
\ep

Here is a stronger form of the maximum principle.

\begin{lemma}[Maximum principle] \label{lMPs}
Let $u \in W$ be a solution to $Lu=0$ in $\Omega$. Then
\begin{equation} \label{lMP6}
\sup_\Omega u \leq \sup_\Gamma \Tr u
\end{equation}
and
\begin{equation} \label{lMP7}
\inf_\Omega u \geq \inf_\Gamma \Tr u,
\end{equation}
where we recall that $\sup$ and $\inf$ actually essential supremum and infimum. 
In particular, if $\Tr u$ is essentially bounded (for the measure $\mu$), then 
\begin{equation} \label{lMP5}
\sup_\Omega |u| \leq \sup_\Gamma |\Tr u|.
\end{equation}
\end{lemma} 

\bp
Let us prove \eqref{lMP6}. Write $M$ for the essential supremum of $\Tr u$ on $\Gamma$;
we may assume that $M < +\infty$, because otherwise \eqref{lMP6} is trivial. 
Then $M-u \in W$ and $\Tr (M-u) \geq 0$ a.e. on $\Gamma$. 
Lemma~\ref{lMP} yields $M-u \geq 0$ a.e. in $\Omega$, that is
\begin{equation} \label{lMP6a}
\sup_\Omega u \leq \sup_\Gamma Tu.
\end{equation}
The lower bound \eqref{lMP7} is similar 
and \eqref{lMP5} follows.
\ep

The harmonic measure will be defined with the help of the Riesz representation theorem (for measures), so we need a linear form on $C^0_0(\Gamma)$, the space of compactly supported continuous functions on $\Gamma$. 
We also write $C^0_b(\overline{\Omega})$ for the space of continuous bounded functions on $\overline{\Omega}$.

\begin{lemma} \label{ldefhm}
There exists a unique linear operator
\begin{equation} \label{defhm4}
U : \, C^0_0(\Gamma) \to  C^0_b(\overline{\Omega}) 
\end{equation}
such that, for every every $g\in C^0_0(\Gamma)$,
\begin{enumerate}[(i)]
\item if $g \in C^0_0(\Gamma) \cap H$, then $Ug \in W$, and it is the solution of \eqref{DpL1}--\eqref{DpL2}, with $f=0$, provided by Lemma \ref{lLM};
\item $\ds \sup_{\Omega} Ug = \sup_\Gamma g$ and $\ds \inf_{\Omega} Ug = \inf_\Gamma g$;
\end{enumerate}
In addition, $U$ enjoys the following properties:
\begin{enumerate}[(i)] \addtocounter{enumi}{2}
\item the restriction of $Ug$ to $\Gamma$ is $g$;
\item $Ug \in \WW(\Omega)$ and is a solution to $Lu = 0$ in $\Omega$;
\item if $B$ is a ball centered on $\Gamma$ and $g \equiv 0$ on $B$, then
$Ug$ lies in $\WW(B \cap \overline{\Omega})$;
\end{enumerate}
\end{lemma}

\bp
The proof of the existence of $U$ and its properties is similar to the one of Lemma 9.4 in \cite{DFMprelim}. 
Still, let us give a sketch of  
the proof of 
existence for 
$U$. 

\medskip

First, we use $(i)$ to define $U$ on $C^0_0(\Gamma) \cap H$. Lemma \ref{lMPs} proves that $(ii)$ is satisfied for any $g \in C^0_0(\Gamma) \cap H$, and in particular 
\[U:\, C^0_0(\Gamma) \cap H \to  C^0_b(\overline{\Omega}) \]
is a continuous operator if we equip both $C^0_0(\Gamma) \cap H$ and $C^0_b(\overline{\Omega})$ with the norm $\|.\|_\infty$. Then, observe that the space $\C_0^0(\Gamma) \cap H$ contains all the restrictions to $\Gamma$ of functions in $C^\infty_0(\R^n)$, and hence $\C_0^0(\Gamma) \cap H$ is dense in $C_0^0(\Gamma)$ (equipped with the norm $\|.\|_\infty$). We define $U:\, C^0_0(\Gamma) \to  C^0_b(\overline{\Omega})$ as the only bounded extension of $U:\, C^0_0(\Gamma) \cap H \to  C^0_b(\overline{\Omega})$, and in particular $(ii)$ is preserved.

The property $(iii)$ is true when $g\in C^0_0(\Gamma) \cap H$ thanks to Lemma \ref{HolderB}, and the property is kept by the extension. The property $(iv)$ is true when $g\in C^0_0(\Gamma) \cap H$ by Lemma \ref{lLM}, and can be extended for all $g\in C^0_0(\Gamma)$ with the help of Cacciopoli's inequality (Lemma \ref{CaccioI}). As for the property $(v)$ - which is immediate by construction for any $g\in C^0_0(\Gamma) \cap H$ - we prove it by approaching $g \in C^0_0(\Gamma)$ by functions in $ C^0_0(\Gamma) \cap H$ and we use Lemma \ref{CaccioB} to control the estimate on the gradient when we take the limit.

\medskip

The uniqueness of $U$ is also immediate, since $(ii)$ forces $U: \, C^0_0(\Gamma) \to C^0_b(\overline{\Omega})$ to be the continuous extension of $U:\,  C^0_0(\Gamma) \cap H \to C^0_b(\overline{\Omega})$ given by $(i)$. \ep

\begin{lemma} \label{ldefhm2}
For any $X\in \Omega$, there exists a unique positive regular Borel measure $\omega^X:= \omega_L^X$ on $\Gamma$ such that 
\begin{equation} \label{defhm}
Ug(X) = 
\int_\Gamma g(y) d\omega^X(y)
\end{equation}
for any $g \in C^0_0(\Gamma)$. Besides, for any Borel set $E \subset \Gamma$,
\begin{equation} \label{defhm6}
\omega^X(E) = \sup\{ \omega^X(K): \, E \supset K, \, K \text{ compact}\} 
= \inf\{ \omega^X(V): \, E \subset V, \, V \text{ open}\}. 
\end{equation}

In addition, the harmonic measure is a probability measure, that is
\begin{equation} \label{defhm12}
\omega^X(\Gamma) = 1. 
\end{equation}
\end{lemma}

\bp
The first part, that is the existence of a positive regular Borel measure satisfying \eqref{defhm}, and the property \eqref{defhm6}, is immediate by applying the Riesz representation theorem (see for instance \cite[Theorem 2.14]{Rudin}) to $U$. The positivity of the harmonic measure comes from $\ds \inf_{\Omega} Ug = \inf_\Gamma g$ given by Lemma \ref{ldefhm} $(ii)$.

\medskip

The fact that $\omega^X(\Gamma) \leq 1$ comes from the fact that $\ds \sup_{\Omega} Ug = \sup_\Gamma g$. We can prove that 
$\omega^X(\Gamma) \geq 1$ by using the H\"older regularity at the boundary (Lemma \ref{HolderB}). See the proof of Lemma~9.6 in \cite{DFMprelim} for details. 
\ep

\begin{lemma} \label{lprhm}
Let $E \subset \Gamma$ be a Borel set and define the function $u_E$ on $\Omega$ by $u_E(X) = \omega^X(E)$. 
Then
\begin{enumerate}[(i)]
\item if there exists $X\in \Omega$ such that $u_E(X) = 0$, then $u_E \equiv 0$;
\item the function $u_E$ lies 
in $\WW(\Omega)$ and is a solution in $\Omega$;
\item if $B\subset \R^n$ is a ball such that
$E\cap B = \emptyset$, then $u_E \in \WW(B \cap \overline{\Omega})$ and $\Tr u_E = 0$ a.e. on $B$.
\end{enumerate}
\end{lemma}

\bp The proof of this result is analogous to the one of Lemma 9.7 in \cite{DFMprelim}. 
Here are some of the main ideas. 
The proof of $(i)$ is quite easy. 
We approach $\1_E$ by $g\in C^0_0(\Gamma)$, and compare $u_E$ with $u_g = Ug$. 
We get that 
$\left|u_g(X) - u_E(X) \right| = u_g(X)$ 
is as small as we want. Then we use  
Lemma \ref{HarnackI2} to say that $0 \leq u_E(Y) \lesssim u_g(Y) \lesssim u_g(X) \leq \epsilon$, and we  let $\epsilon$ tend to $0$.  

The proofs of $(ii)$ and $(iii)$ are longer. They consist in approaching $u_E$ by functions $u_g = Ug$ that have all the desired properties by Lemma \ref{ldefhm}, then controlling $\nabla u_g$ uniformly 
with the help of 
Lemma \ref{CaccioI} and Lemma \ref{CaccioB}. We eventually use Lemma~\ref{HolderI}, 
Lemma~\ref{HolderB}, and Lemma~\ref{lMPs} to ensure that the $u_g$'s 
are nice functions that satisfy $g \leq h \Rightarrow u_g \leq u_h$.
\ep

\section{Bounded boundaries}
\label{SBounded}

So far, in Section \ref{Scones} and hence in all the sections following it, 
we have been working with a boundary set $\Gamma$ which is unbounded, 
and hence an unbounded domain $\Omega$ too. But it is some times interesting, 
and not too difficult, to deal with bounded sets $\Gamma$. In this section we describe how to modify our assumptions, and some times the proofs, to extend the results of this paper to the case of bounded $\Gamma$. So let us assume now that 
$\Gamma = \d \Omega$ is bounded, and (to normalize things) that
\begin{equation} \label{Bd1}
0 \in \Gamma \ \text{ and } \ \diam(\Gamma) = R_0 > 0.
\end{equation}
There will be two slightly different cases to consider, Case 1 when $\Omega$ also is bounded (and connected - due to (\hyperref[H1]{H1})),
and Case~2 when $\Omega$ is the unbounded component of $\R^n \sm \Gamma$. 
When the dimension of $\Gamma$ is smaller than $n-1$, we are in Case 2, 
but Case 1 is interesting too, especially in the context of mixed co-dimensions, 
where we may do it on purpose to add pieces of boundary 
that isolate some parts of a domain. That is, even if we start with the unbounded
component of $\R^n \sm \Gamma$, we could for instance add to $\Gamma$ a large sphere 
like $S = \d B(0,2R_0)$ to $\Gamma$, and restrict our attention to the bounded component
of $\R^n \sm (\Gamma \cup S)$ that touches $S$ because this is simpler.

Most of the results above
are local, in the sense that they rely on computations that do not go too far.
The only difference that it makes on our assumptions is that - if $\Omega$ is bounded - we need 
to take the $r$ in (\hyperref[H1]{H1}) not too large, for instance not bigger than $\diam \Omega$, 
while the unbounded case will require that $r$ to be taken in the full range $(0,+\infty)$.

Observe also that in the case where $\Gamma$ is bounded, we just need 
the analogue of (\hyperref[H1]{H1})-(\hyperref[H6]{H6}), where 
we keep the same statement as before but only ask 
(\hyperref[H1]{H1}) and (\hyperref[H3]{H3}) to hold when $B(x,r) \subset B_0  = B(0,2R_0)$, and (\hyperref[H2]{H2}) to hold for points $X, Y \in B_0$ 
(the other case would follow anyway). When $\Omega$ is also bounded, 
we can restrict to $B(0,2R_0)$ in the definition of (\hyperref[H4]{H4})
(the absolute continuity and doubling property  for $m$) and (\hyperref[H6]{H6})
(the density and weak Poicar\'e inequality for $m$).
One could see such apparent weakening as an improvement, but it is easy to check that the cases that we dropped are automatically true for bounded $\Gamma$ and/or $\Omega$.
However the truth is that we are not really interested in studying wild weights $w$ far from 
$\Gamma$, and a simple monomial equivalent at infinity would probably be enough.

With these assumptions, most of our local estimates still hold, with very little changes in the proofs.
Let us be a little more specific.

We keep the definition of $W$ as it was. Notice that constant functions still lie in $W$ (with a vanishing
norm); depending on the behavior of $m$ far from $\Gamma$, the functions $u\in W$ may have a more or less rich behavior near infinity, but let us not bother yet. Section \ref{SecW} goes through without modification
(we kept the same assumptions on $m$ alone).

The definition of dyadic pseudocubes has to be changed a little bit: we only use the 
partition $\Gamma = \bigcup_{j \in {\mathcal J}_k} Q_j^k$ for $k \geq k_0$, where 
$k_0$ is such that $2^{-k_0} \sim 4 R_0$, and also it is customary to take for $k = k_0$
the trivial decomposition into a unique cube $Q_0 = \Gamma$. 
 Of course, all the subsequent sums in $k$ will be restricted to $k \geq k_0$. 

Then, even in the definition of the access regions $\gamma(x)$ (as in \eqref{defgamma})
for unbounded domains $\Omega$, we will only consider cubes of size at most $CR_0$, and so our 
access regions will be bounded. We are not shocked because for many of our results we already 
considered the truncated regions of \eqref{defgammaQ*}. 
The results of Section \ref{Scones}, and in particular the improved Poincar\'e inequality in 
Theorem \ref{Poincare1a}, are local and stay the same, but we only consider sets that are contained in $C B_0$. 
Recall however that the case of balls $B$ such that $2B \subset \Omega$, even when $B$ is large,
is taken care of in Lemma \ref{Poincare00}, so we will never be in real trouble anywhere.

Our definition \eqref{defH} of the Hilbert space $H$ on $\Gamma$ stays the same; as before
constants lie in $H$, with a vanishing norm. Theorem \ref{TraceTh} on the existence of a trace 
operator is still valid with the same proof. The proof does not use the values of $u\in W$ at distance
more than $CR_0$ from $\Gamma$, so we may even forget the corresponding part of $||f||_W$ in the
estimate for $\|\Tr(f)\|_H$. That is,
\begin{equation}\label{Bd2}
\|\Tr u\|_H^2 \lesssim \int_{\{Z\in \Omega \cap B(0,CR_0) \}} |\nabla u(Z)|^2 \, dm(Z).
\end{equation}
Another way (softer but just a bit more complicated technically)
to check this is to notice that, when $\Omega$ is unbounded, we may always 
truncate any $u\in W$ in the following way. We select a smooth cut-off function $\varphi$ such that 
$\varphi = 1$ in $B_0$ and $\varphi = 0$ outside of $2B_0$, pick a ball  $B_1$ of radius $R_0$
such that $B_1 \subset 2B_0 \sm B_0$ and $B_1$ touches $\d B_0$, let $m_1$ denote the average
of $u$ on $B_1$, and consider the ``truncated'' function $\wt u = \varphi u + (1-\varphi) m_1$.
Obviously $\wt u$ has the same trace, and is would be easy to see, using the extension of Lemma \ref{Poincare00} 
to a $(2,2)$-Poincar\'e inequality, as in Theorem \ref{Poincare1} with $p=2$, that 
$\wt f \in W$, with $||\wt f||_W^2 \leq C \int_{2B_0} |\nabla f|^2 dm$. Of course all this is much easier if
$w$  is reasonably smooth on $2B_0 \sm B_0$.

The  product rule for the trace and gradient (Lemma \ref{lmult}), as well as all the local algebraic formulas,
go through. Similarly, the Poincar\'e inequalities on the boundary
(Theorem \ref{PoincareTh2} and Corollary \ref{PoincareCor2}) stay the same, 
except that we restrict to balls of radius at most $10R_0$, say. 

Our extension theorem (Theorem \ref{ThExt}) is still true; the construction also easily gives that
$\Ext(f)$ is constant on $\R^n \sm CB_0$ (when $\Omega$ is unbounded), and we can take the constant
equal to the average $m_0=\fint_\Gamma f d\mu$ . The simplest way to see this is to consider $f-m_0$
and use the formula \eqref{defExt}, but restrict the sum to Whitney cubes of size at most $CR_0$. Or said
differently, for the function $f-m_0$ we can use $y_I = 0$ for all the large Whitney cubes.

There is no difficulty with the density or algebraic results of Section \ref{SComplete}, and the local
spaces of Section \ref{SLocal} are (just a bit) simpler. 
The definition of solutions is local, and all the regularity theorems for solutions found in Section \ref{SSolutions} stay the same. 
This statement may look obvious, we are saying since the beginning of the section that all the results are exactly the same for bounded $\Gamma$ and unbounded $\Gamma$, but let us observe the following interesting fact. 
The boundary regularity results, such as Lemmas \ref{MoserB} and \ref{HolderB}, hold for all balls $B$ centered at the boundary even when the radius of $B$ is way bigger than the diameter of $\Gamma$, and so can be applied for instance to the Green functions - that we shall introduce in a 
next paper.

Let us now review the basic features of the harmonic measure. Its construction given in Section~\ref{SHarm} still goes through; that is, the existence of solutions as in Lemma \ref{lLM}
is still valid, by Lax Milgram, and so is the maximum principle that allows one to solve the 
Dirichlet problem for $f \in C^0_0(\Gamma) = C^0(\Gamma)$. Thus $\omega^X$ is defined
(by the Riesz representation theorem), and is again a probability measure because the best extension
of the function $1$ is $1$. 

\medskip

If $\Gamma$ is bounded and $\Omega$ is the unbounded component of $\R^n \setminus \Gamma$, then Brownian paths leaving from $X \in \Omega$ 
have a nonzero probability of never touching $\Gamma$ before going to infinity. It means that the classical harmonic measure - defined from the Laplacian - is not a probability measure. 
This simple case is however not included in our theory; indeed the assumption (\hyperref[H5]{H5}) fails for large $r$ when we take $\mu$ as the surface measure on the bounded set $\Gamma$ and $m$ as the Lebesgue measure on $\Omega$.
On the contrary, our theory roughly says that a modified Brownian motion, that imposes a drift in the direction of $\Gamma$ when we are far from it, is sufficient to guarantee to touch $\Gamma$ with probability one.

For the two last sections, where we study the Green functions and the harmonic measure, leading to a comparison principle, we do not assume that $\Gamma$ is unbounded anymore. So both $\Gamma$ and $\Omega$ 
can be bounded or unbounded, and we believe that the case where $\Gamma$ is bounded and $\Omega$ is the unbounded component of $\R^n \sm \Gamma$ is the most tricky one.

\section{Green functions}

\label{SGreen}

We associate Green functions to the degenerate elliptic operator $L$. A Green function is, formally, a function $g$ defined on $\Omega \times \Omega$ and such that for any $y \in \Omega$, the function $g(.,y)$ satisfies \eqref{DpL1} and \eqref{DpL2} with $f = \delta_y$ - the Dirac distribution at the point $y$ - and $g\equiv 0$.

The harmonic measure can be seen as a fundamental tool to solve the problem $Lu = 0$ in $\Omega$ with $\Tr u = g$ on $\Gamma$, while the Green function is a key ingredient to be able to solve $Lu=f$ in $\Omega$ with $\Tr u = 0$ on $\Gamma$. Their properties are actually related, as we shall see in Section \ref{SCP}. 

Let us recall that, as in the previous sections, we assume (\hyperref[H1]{H1})--(\hyperref[H6]{H6}), and \eqref{defEllip}--\eqref{defBdd}.

\medskip

In order to define the Green function in our context, we will follow closely the proof of Gr\"uter and Widman \cite{GW} (as in \cite{DFMprelim}). In the article \cite{GW}, the authors proved the existence of the Green functions $g(.,y)$ by taking a weak limit of some $g^\rho(.,y)$ that solves $Lu = f^\rho$ and $\Tr u \equiv 0$ for some $f^\rho$ that `approximates' the delta distribution $\delta_y$.

Some difficulties appears when we try to get `local' estimates, i.e. when the distance between $x$ and $y$ is small compared to the distance of both points to the boundary. Those estimates are needed to show that our $g^\rho(.,y)$ are uniformly bounded in some good space. We solve those difficulties by using methods inspired from \cite{FJK}, where the authors deal with degenerate elliptic operators but they define Green function via another method.

For short, we claim here that Gr\"uter and Widman's method can be applied - up to few changes - to obtain Green functions in a large varieties of situations, and for instance doesn't require to have a global Sobolev inequality.

Instead of giving a big theorem for the start, as in \cite{DFMprelim}, we choose here to divide the work, and prove plenty of small lemmas, whose proofs are sometimes omitted because they are the same as in \cite{DFMprelim}. The important results are gathered at the end of the section, in Theorem \ref{GreenEx}.

\begin{definition} \label{Defgrho}
Let $y\in \Omega$ and $\rho>0$ satisfy 
$100\rho < \delta(y)$. The function $g^\rho(.,y)$ is the function in $W_0$ that satisfies 
\begin{equation} \label{defgrho2}
\int_\Omega \A \nabla g^\rho(.,y) \cdot \nabla v \, dm = \fint_{B(y,\rho)} v \, dm \qquad \text{ for all } v\in W_0,
\end{equation}
as given by Lemma \ref{lLM}.

Notice that the definition makes sense, because $v \to \fint_{B(y,\rho)} v \, dm$ is a 
bounded linear form on $W_0$ (and hence an element $f^\rho \in W^{-1}$
to which we apply the lemma), by the doubling condition (\hyperref[H4]{H4}) and the Poincar\'e inequality 
Corollary \ref{PoincareCor2}. The norm of $f^\rho$ in $W^{-1}$ depends on $y$ and $\rho$, 
but it doesn't matter.

Since $y$ will be fixed for a long part of our section, we write in the sequel 
$\g^\rho$ for $g^\rho(.,y)$ and $B_\rho$ for $B(y,\rho)$.
Then the condition \eqref{defgrho2} in the definition becomes
 \begin{equation} \label{defgrho}
\int_\Omega \A \nabla \g^\rho \cdot \nabla v \, dm = \fint_{B_\rho} v \, dm \qquad \text{ for all } v\in W_0,
\end{equation}
\end{definition}

We deduce at once from the definition that
\begin{equation} \label{Green4}
\text{$\g^\rho \in W_0$ is a solution to $Lu = 0$ in $\Omega \setminus \overline{B_\rho}$.}
\end{equation}
In particular, by Lemmas \ref{HolderI} and \ref{HolderB}, the function $\g^\rho$ is continuous on $\overline \Omega \setminus \overline{B_\rho}$.

\begin{lemma} \label{grho>0}
For all $y\in \Omega$, the function $\g^\rho = g^\rho(.,y)$
is nonnegative.
\end{lemma}

\bp
The proof of this fact 
is identical to the one given for 
\cite[Lemma~10.1]{DFMprelim},  
and relies only on the stability of $W_0$ provided by Lemma \ref{chainrule}.
\ep

We 
now prove pointwise estimates on $\g^\rho$ and start with the case 
when $x$ is far from $y$.

\begin{lemma} \label{grhoxyfar}
If $x,y\in \Omega$ are such that $10|x-y| \geq \delta(y)$, then
\[\g^\rho(x) \leq C\frac{|x-y|^2}{m(B(y,|x-y|) \cap \Omega)},\]
where $C$ depends on $C_1$ to $C_6$, $C_A$, and $n$.
\end{lemma}

\bp Let $R\geq \delta(y)> 100\rho>0$, and write $B_R$ for $B(y,R)$. Let $p$ be 
in the range given by Corollary~\ref{PoincareCor2}; that is, any $p\in [1,2k)$, where $k$ is a constant that depends on the geometry, will do.
We want to prove that for all $t>0$,
\begin{equation} \label{Green8}
\frac{m(\{x\in B_R, \, \g^\rho(x) > t\})}{m(B_R \cap \Omega)} \leq C \left(\frac{t \, m(B_R \cap \Omega)}{R^2}\right)^{-\frac p2}
\end{equation}
with a constant $C$ independent of $\rho$, $t$ and $R$. The proof of the claim is analogous to the one 
in \cite{DFMprelim}, but we repeat it because we will use similar computations later on. We use \eqref{defgrho} with the test function 
\begin{equation} \label{Green9}
\varphi(z) := \left(\frac 2{t}-\frac 1{\g^\rho(z)}\right)^+ = \max\left\{0, \frac 2{t}-\frac 1{\g^\rho(z)}\right\}
\end{equation}
(and $\varphi(z) = 0$ if $\g^{\rho}(z) = 0$), which lies in $W_0$ by Lemma~\ref{chainrule}. 
Set  
$\Omega_{s} := \big\{z\in \Omega, \, \g^\rho(z) > s \big\}$ 
and observe that $\varphi$ is supported in $\Omega_{t/2}$. Hence 
\begin{equation} \label{GreenA}
a(\g^\rho,\varphi) = \int_{\Omega_{t/2}} \frac{\A \nabla \g^\rho \cdot  \nabla \g^\rho}{(\g^\rho)^2} \, dm 
= \fint_{B_\rho} \varphi \, dm \leq \frac2{t}
\end{equation}
and then, 
thanks to the ellipticity condition \eqref{defEllip2},
\begin{equation} \label{GreenB}
\int_{\Omega_{t/2}} \frac{|\nabla \g^\rho|^2}{(\g^\rho)^2} \, dm \leq \frac Ct.
\end{equation}
Pick a point $y_0 \in \Gamma$ such that
$|y-y_0| = \delta(y)$. Set $\wt B_R$ for $B(y_0,2R) \supset B_R$. Define $v$ by
$v(z) : = (\ln(\g^\rho(z))-\ln t + \ln 2)^+$; then $v \in W_0$ too, 
thanks to Lemma~\ref{chainrule}, and 
$|\nabla v|^2 = |\nabla \g^\rho|^2/(g^\rho)^2$.
Corollary \ref{PoincareCor2} (the 
Sobolev-Poincar\'e inequality at the boundary) implies that 
\begin{equation} \label{GreenC} \begin{split}
\left(\int_{\Omega_{t/2} \cap \wt B_R} |v|^p \, dm\right)^\frac1p & \leq C R \,  m(\wt B_R \cap \Omega)^{\frac1p-\frac12} \left(\int_{\Omega_{t/2} \cap 2\wt B_R} |\nabla v|^2 \, dm\right)^\frac12 \\
& \leq C R \, m(B_R \cap \Omega)^{\frac1p-\frac12} t^{-\frac12}
\end{split}\end{equation}
by  
\eqref{GreenB} and (\hyperref[H4]{H4}). Yet, $v \geq \ln(2)$ on $\Omega_t$, and thus the above inequality gives that
\begin{equation} \label{GreenE}
m(\Omega_t \cap B_R) \leq m(\Omega_t \cap \wt B_R) \lesssim R^p [m(B_R \cap \Omega)]^{1-\frac p2} t^{-\frac p2} \leq m(B_R \cap \Omega) \left( \frac{m(B_R \cap \Omega) t}{R^2}\right)^{-\frac p2}.
\end{equation}
The claim \eqref{Green8} follows. 
We are now ready to establish pointwise estimates on $\g^\rho$ when $x$ is far from $y$. 
We now 
aim to prove \eqref{Green8} with a constant independent of $\rho$.
Set $R=10|x-y| \geq \delta(y)$. By \eqref{Green4}, $\g^\rho \in W_0$ is a solution to $Lu=0$ in $\Omega \setminus \overline{B_\rho}$, so we can use Moser's estimates; 
we claim that we get that 
\begin{equation} \label{GreenG0}
\g^\rho(x) \leq  \frac{C}{m(B(x,R/2) \cap \Omega)} \int_{B(x,R/2) \cap \Omega} \g^\rho \, dm.
\end{equation}
When $\delta(x) \geq R/50$ we apply Lemma~\ref{MoserI} in the ball $B(x,R/100)$), and 
when $\delta(x) \leq R/50$ we apply Lemma~\ref{MoserB} in the ball $B(x_0,R/25)$, 
where $x_0$ is such that $|x-x_0| = \delta(x)$. We can use (\hyperref[H4]{H4}) to replace
the measure of the ball by $m(B(x,R/2) \cap \Omega)$.

\noindent We can use now the fact that $B(x,R/2) \subset B_R$ and 
Cavalieri's formula (see for instance \cite[p. 28, Proposition 2.3]{Duoandi}) to get that
 \begin{equation} \label{GreenG}
\g^\rho(x) \lesssim \int_0^{+\infty}  \frac{m(\Omega_{t} \cap B_R)}{m(\Omega \cap B_R)}dt. 
\end{equation}
Take $s>0$, to be chosen later. 
We bound the interior of the integral above by $1$ when $t<s$, and for $t\geq s$ we use
\eqref{Green8}, which we apply with some $p>2$ (this is possible); we get that
 \begin{equation} \label{GreenH0} \begin{split}
\g^\rho(x) & \lesssim \int_0^{s}  \frac{m(\Omega_{t} \cap B_R)}{m(\Omega \cap B_R)}dt +  \int_s^{+\infty}  \frac{m(\Omega_{t} \cap B_R)}{m(\Omega \cap B_R)}dt \\
& \lesssim \int_0^s 1\, dt + \left(\frac{m(B_R \cap \Omega))}{R^2}\right)^{-\frac p2} \int_s^{+\infty} t^{-\frac p2} dt \\
& \lesssim s\left[ 1 + \left(\frac{s\, m(B_R \cap \Omega)}{R^2}\right)^{-\frac p2} \right].
\end{split} \end{equation}
Now we 
minimize the right-hand side in $s$. 
We find $s \approx R^{2}/m(B_R \cap \Omega)$ and then $\g^\rho(x) \lesssim R^{2}/m(B_R \cap \Omega)$. Since $R = 10|x-y|$, the lemma follows from (\hyperref[H4]{H4}).
\ep

The next result deals with the case when $x$ and $y$ are close to each other.

\begin{lemma} \label{grhoxyclose}
If $x,y\in \Omega$ are such that $40\rho \leq 2|x-y| \leq \delta(y)$, then
\[\g^\rho(x) \leq C\int_{|x-y|}^{\delta(y)} \frac{r^2}{m(B(y,r))} \, \frac{dr}{r},\]
where $C$ depends on $C_1$ to $C_6$, $C_A$, and $n$.
\end{lemma}

\bp The proof of this result in the classical case, at least the one in \cite{GW}, uses a global Sobolev equality. In our setting, we don't have Sobolev embeddings, only a Sobolev-Poinca\'e inequality 
(Theorem \ref{Poincare00}), and the $L^q$ norm given in the right-hand side of our  
Sobolev-Poincar\'e inequality may just be $L^{2+\epsilon}$. 
In particular, we have no reason to get close to the desired $L^{2n/(n-2)}$.

Fortunately, the slight improvement in the exponent of the $L^p$ space given by
Theorem~\ref{Poincare00} 
will be - as for Lemma \ref{grhoxyfar} - sufficient. 
Even better, the proof will follow the same ideas as Lemma \ref{grhoxyfar}. 

Let $j_0 \geq 0$
be the biggest integer such that $2^{j_0+1}|x-y| \leq \delta(y)$. 
To lighten the notation, we write $B^j$ for $B(y,2^j|x-y|)$.
 We shall prove that for any $j$ between $0$ and $j_0 - 1$, 
\begin{equation} \label{diffsupgrho}
\sup_{B^{j+1} \setminus B^{j}} \g^\rho - \sup_{B^{j+2} \setminus B^{j+1}} \g^\rho 
\leq C \, \frac{(2^{j}|x-y|)^2}{m(B_j)}. 
\end{equation}
We write $\g_j^\rho$ for $\g^\rho - \sup_{B^{j+2} \setminus B^{j+1}} \g^\rho$. 
We also write $\Omega_{s,j}$ for $\{x\in \Omega, \, \g^\rho_j >s \}$. 
Notice that $B^{j+2} \subset \Omega$, and $\sup_{B^{j+2} \setminus B^{j+1}}
= \sup_{\Omega \setminus B^{j+1}}$ by the maximum principle. 
Hence, $\Omega_{s,j} \subset B^{j+1}$ when $s > 0$.
Let $t > 0$ be given; we 
use again \eqref{defgrho},
but with the test function 
\[\varphi_j(z) := \left( \frac{2}{t} - \frac1{(\g^\rho_j(z))^+} \right)^+,\]
to get, as we did for \eqref{GreenA},
\[a(\g^\rho,\varphi) = \int_{\Omega_{t/2,j}} \frac{\A \nabla \g^\rho_j \cdot \nabla \g^\rho_j}{(\g^\rho_j)^2} \, dm \leq \fint_{B_\rho} \varphi_j \, dm \leq \frac2t,\]
and by 
\eqref{defEllip2}
\begin{equation} \label{GreenBj}
\int_{\Omega_{t/2,j}} \frac{|\nabla \g^\rho_j|^2}{(\g^\rho_j)^2} \, dm \leq \frac Ct.
\end{equation}
Set $v_j(z):= (\ln(\g^\rho_j(z)) - \ln t + \ln 2)^+$;
as before $v_j \in W$, it is  supported in $\Omega_{t/2,j} \subset B^{j+1}$, and
 $|\nabla v_j| = |\nabla \g^\rho_j|/\g^\rho_j$ on $\Omega_{t/2}$.
 Since $v_j = 0$ on $B^{j+2} \sm B^{j+1} \subset \Omega$, we can apply 
 Theorem \ref{Poincare1} and Remark \ref{rmkL2g}, and we get that
 \begin{equation} \label{GreenCj} \begin{split}
\left(\int_{\Omega_{t/2,j}} |v_j|^p \, dm\right)^\frac1p & \leq C 2^{j+1}|x-y| \, m(B^{j+1})^{\frac1p-\frac12} \left(\int_{\Omega_{t/2,j}} |\nabla v_j|^2 \, dm\right)^\frac12 \\
& \leq C 2^{j}|x-y| \, m(B^j)^{\frac1p-\frac12} t^{-\frac12}
\end{split}\end{equation}
by \eqref{GreenBj} and (\hyperref[H4]{H4}), 
where $p\in [k,2k]$
plays the role of $kp$ in Theorem \ref{Poincare1} and Remark~\ref{rmkL2g}, 
and we are mostly interested in $p=2$ there which yields $p=2k$ here. Of course $C$ is independent of $j$.
Since $|v_j| > \ln 2$ on $\Omega_{t,j}$
and $\Omega_{t,j} \subset \Omega_{t/2,j}$, 
\eqref{GreenCj} implies that
\begin{equation} \label{Green8j}
\frac{m(\Omega_{t,j})}{m(B^{j+2})} 
\leq C \left(\frac{t \, m(B^j)}{(2^j|x-y|)^2}\right)^{-\frac p2}.
\end{equation}
The rest of the proof of \eqref{diffsupgrho} is 
similar to what we did for 
Lemma \ref{grhoxyfar}.  
Since $\g^\rho_j$ is a solution in $B^{j+2} \setminus B^{j-1}$, we can use the Moser inequality inside 
$B^{j+2} \setminus B^{j-1}$ to get that for $z\in B^{j+1} \setminus B^j$,
$\g^\rho_j(z)$ is smaller - up to a constant - than its average on $B(z,2^{j-10}|x-y|)$.
The measure of this last ball
is equivalent, by (\hyperref[H4]{H4}), to the measure of $B_{j+2}$ and thus for any $z\in B^{j+1} \setminus B^j$
 \begin{equation} \label{GreenG2} 
 \g^\rho_j(z) \lesssim \fint_{B^{j+2}} \g^{\rho}_j \, dm =  \int_0^{+\infty}  \frac{m(\Omega_{t,j})}{m(B^{j+2})}dt 
\end{equation}
(compare with \eqref{GreenG}).
Again we split the last integral  
into two pieces, and for the second one
we use \eqref{Green8j}; we
obtain that for all $z\in B^{j+1} \setminus B^j$ 
 \begin{equation} \label{GreenH0j} \begin{split}
\g^\rho_j(z) & \lesssim \int_0^s 1\, dt + \left(\frac{m(B^{j}))}{(2^j|x-y|)^2}\right)^{-\frac p2} \int_s^{+\infty} t^{-\frac p2} dt \\
& \lesssim s\left[ 1 + \left(\frac{s\, m(B^j)}{(2^j|x-y|)^2}\right)^{-\frac p2} \right].
\end{split} \end{equation}
where we can choose $p$, which comes from Theorem \ref{Poincare1} (for instance applied with
the exponent $2$), 
strictly bigger than 2. We optimize in $s$ 
and take the supremum in $z$ to get the desired estimate \eqref{diffsupgrho}.

We can use (\hyperref[H4]{H4}) to rewrite 
\eqref{diffsupgrho} as
\begin{equation} \label{diffsupgrho2}
\sup_{B^{j+1} \setminus B^{j}} \g^\rho - \sup_{B^{j+2} \setminus B^{j+1}} \g^\rho \leq C \int_{2^j|x-y|}^{2^{j+1}|x-y|} \frac{r^2}{m(B(y,r))} \frac{dr}{r}. 
\end{equation}
This was for $j < j_0$, but for $j=j_0$ we will be able to apply Lemma \ref{grhoxyfar}.
Recall that $2^{j_0+2}|x-y| >\delta(y)$ by definition of $j_0$; hence for 
$z\in B^{j_0+1} \setminus B^{j_0}$, 
$|z-y| \geq 2^{j_0}|x-y| > \delta(y)/10$ and by Lemma \ref{grhoxyfar} and the same trick
as for \eqref{diffsupgrho2},
\[
\sup_{B^{j_0+1} \setminus B^{j_0}} \g^\rho 
\leq C \int_{2^{j_0}|x-y|}^{2^{j_0+1}|x-y|} \frac{r^2}{m(B(y,r))} \frac{dr}{r} 
\leq C \int_{2^{j_0}|x-y|}^{\delta(y)} \frac{r^2}{m(B(y,r))} \frac{dr}{r}.
\]
Now, 
since $\g_\rho$ is continuous around $x$ (by the interior H\"older estimates and \eqref{Green4}), 
\[\begin{split}
\g^{\rho}(x) & \leq \sup_{B^1\setminus B^0} \g^\rho 
\leq \sup_{B^{j_0+1}\setminus B^{j_0}} \g + \sum_{j=0}^{j_0-1} \left(\sup_{B^{j+1} \setminus B^{j}} \g^\rho - \sup_{B^{j+2} \setminus B^{j+1}} \g^\rho \right) \\
& \lesssim \int_{2^{j_0}|x-y|}^{\delta(y)} \frac{r^2}{m(B(y,r))} \frac{dr}{r} + \sum_{j=0}^{j_0-1} \int_{2^j|x-y|}^{2^{j+1}|x-y|} \frac{r^2}{m(B(y,r))} \frac{dr}{r} \\
& \lesssim \int_{|x-y|}^{\delta(y)} \frac{r^2}{m(B(y,r))} \frac{dr}{r} \, ; 
\end{split}\]
Lemma \ref{grhoxyclose} follows.\ep

Before we continue to prove estimate about Green functions, we take a little time to talk about cut-off functions. Pick 
$\phi \in C^\infty(\R_+)$ such that $0 \leq \phi \leq 1$, $\phi \equiv 0$ on $(2,+\infty)$, $\phi \equiv 1$ on $[0,1)$, and $|\phi'| \leq 2$.  If we want a cut-off function adapted to the ball $B(x,r)$, the first choice will be
\begin{equation} \label{alpha1}
\alpha^1(y) 
:= \phi\left( \frac{|x-y|}{r} \right).
\end{equation}
If the above cut-off function fails to work, we might try to use a cut-off function 
that involves logarithms, 
in the spirit of the one used by Sobolev. For instance, if we work on the 
Green function for the unit disk (with the classical Lebesgue measure) on $\R^2$, the good cut-off 
may be
\begin{equation} \label{alpha2}
\alpha^2(y) 
:= \phi\left( \frac{\ln(\delta(y)/r)}{\ln(\delta(y)/|x-y|)} \right).
\end{equation}
In the classical theory, where the domains are equipped 
with the usual Lebesgue measure, we would use $\alpha^1$ when $n\geq 3$ and $\alpha^2$ when $n=2$. 
In our article, $\alpha^1$ or $\alpha^2$ may be needed, or something different. The cut-off functions that we shall need depend on the measure $m$ and the purpose of the next lines is to define them.

We define the function $\gamma = \gamma_y$
on $(0,\delta(y))$ by 
\begin{equation} \label{defgg}
\gamma(s) := \int_{s}^{\delta(y)} \frac{t^2}{m(B(y,t))} \, \frac{dt}{t}.
\end{equation}
The function $\gamma$ is well defined (because $m(B(y,t)) > 0$ since $m$ is doubling on $\Omega$),
and decreasing. In addition, $t \mapsto  m(B(y,t))$ is continuous, because $m$ is absolutely
continuous with respect to the Lebesgue measure 
(and not because $m$ is doubling),
and so $\gamma$ is of class $C^1$, with a derivative equal to 
\begin{equation} \label{defgg'}
\gamma'(s) := -\frac{s}{m(B(y,s))}.
\end{equation}
Next we use our function $\gamma$ and set 
\begin{equation} \label{alpha}
\alpha_r(s) := \phi\left( \frac{\gamma(r)}{\gamma(s)} \right)
\hbox{ for } 0 < r, s < \gamma(y).
\end{equation}
By construction, 
\begin{equation} \label{alpha22}
\alpha_r \equiv 1 \text{ on }[0,r),
\end{equation}
\begin{equation} \label{alpha3}
\text{$\alpha_r(s) = 0$ when $\gamma(s) < \frac12 \gamma(r)$,}
\end{equation}
and $\alpha_r$ is of class $C^1$ on $(0,\delta(y))$, with a derivative equal to 
$- \frac{\gamma'(s) \gamma(r)}{\gamma(s)^2}\,
\phi'\left( \frac{\gamma(r)}{\gamma(s)} \right)$.
Thus
\begin{equation} \label{equivgg}
\alpha'_r  \text{ is supported on the interval where } \ 
\frac12 \gamma(r) \leq \gamma(s) \leq \gamma(r)
\end{equation}
(recall that $\gamma$ is decreasing), and
\begin{equation} \label{bdalpha'}
|\alpha'_r(s)| \leq 8 \,\frac{\gamma'(s)}{\gamma(r)}.
\end{equation}

We shall also need a variation of the maximum principle (Lemma \ref{lMP}).

\begin{lemma} \label{lMPg}
Let $F \subset \R^n$ be a closed set and $E \subset \R^n$ an open set such that
$F \subset E \subset \R^n$ and $\dist(F,\R^n \setminus E)>0$. 
Let $u \in \WW(E \cap \overline \Omega)$ be a supersolution 
for $L$ in $\Omega \cap E$ such that 
\begin{enumerate}[(i)]
\item $\ds \int_{E \cap \Omega} |\nabla u|^2 \, dm < +\infty$,

\smallskip
\item $\Tr u \geq 0$ a.e. on $\Gamma \cap E$,
\item $u \geq 0$ a.e. in $(E \setminus F) \cap \Omega$.
\end{enumerate}
Then $u\geq 0$ a.e. in $E \cap \Omega$.
\end{lemma}

\bp The proof of this result is the same as the one of \cite[Lemma 11.3]{DFMprelim}. \ep

We are now ready to establish a lower bound on $\g^\rho(x)$ when $x$ and $y$ are close. Those lower bounds are not necessary in our article to prove the existence of the Green function, and a reader 
who is only interested in existence can skip the next lemma.

\begin{lemma} \label{grhoxylb}
If $x,y\in \Omega$ are such that $40\rho \leq 2|x-y| \leq \delta(y)$, then
\[\g^\rho(x) \geq C^{-1} \int_{|x-y|}^{\delta(y)} \frac{t^2}{m(B(y,t))} \, \frac{dt}{t},\]
where $C>0$ depends on $C_4$, $C_6$, $C_A$, and $n$.
\end{lemma}

\bp The first point that we need to verify is that $\g^\rho(x)$ is increasing when $x\to y$, at least in a weak sense, and when $|x-y|$ is way bigger than $\rho$. Pick $r>10\rho$. The function $v_r:= (\sup_{B_r \setminus B_{r/2}} \g^\rho) - \g^\rho$ - when $B_s$ denotes as usual $B(y,s)$ - is a solution to $Lu = 0$ in $\Omega \setminus \overline{B_{r/2}}$. Moreover we can easily observe that $\Tr v_r \geq 0$ and $v_r \geq 0$ on $B_r \setminus B_{r/2}$. 
We can apply 
Lemma~\ref{lMPg} with $E = \R^n \setminus \overline{B_{r/2}}$ and $F = \R^n \setminus B_r$, 
and in particular (i) holds because $\g^\rho \in W_0$; 
we deduce that 
$v_r\geq 0$ in $\Omega \sm \overline{B_{r/2}}$, that is,
\begin{equation} \label{gincrease}
\sup_{\Omega \setminus B_{r/2}} \g^\rho \leq \sup_{B_r \setminus B_{r/2}} \g^\rho.
\end{equation}

For the rest of the proof, we write $r$ for $|x-y|$, and for $i\geq 0$, we set $r_i$ as the only value such that 
\begin{equation} \label{defri}
\gamma(r_i) = 2^{-i}\gamma(r).
\end{equation}
Such a point exists, because $\gamma$ is a (strictly) decreasing continuous function with 
$\gamma(\delta(y)) = 0$. Notice that $r_0 = r$, and $\{ r_i \}$ is an increasing sequence whoes limit 
is $\delta(y)$ (but we won't go that far).
First, we use the test function on $\eta_1(x):= \alpha_{r_1}(|x-y|)$ in \eqref{defgrho}. Thanks to \eqref{equivgg}, $\nabla \eta_1$ is supported in $B_{r_2} \setminus B_{r_1}$ and one obtains
\begin{align} \label{pre99}
1 & = \int_{B_{r_2} \setminus B_{r_1}} \A \nabla \g^\rho \cdot \nabla \eta_1 \, dm \leq \frac{C}{\gamma(r_1)} \int_{B_{r_2} \setminus B_{r_1}} |\nabla \g^\rho| |\gamma'(|x-y|)| \, dm 
\nonumber\\
& \leq  \frac{C}{\gamma(r_1)} \left(\int_{B_{r_2} \setminus B_{r_1}} |\nabla \g^\rho|^2 \, dm\right)^\frac12 \left(\int_{B_{r_2} \setminus B_{r_1}} \frac{|x-y|^2}{m(B(y,|x-y|))^2} \, dm\right)^\frac12
\end{align}
We now want to prove that
\begin{equation} \label{Green99}
\int_{B_{r_2} \setminus B_{r_1}} \frac{|x-y|^2}{m(B(y,|x-y|))^2} \, dm \leq C \int_{r_1}^{r_2} \frac{s^2}{m(B_s)} \frac{ds}{s} \leq C \gamma(r_1).
\end{equation}
The second part is just the definition of $\gamma(r_1)$ (we integrate further); 
the first inequality will be a little longer to prove, because we want to avoid
the unpleasant situation where $B_{r_2} \setminus B_{r_1}$ is a very thin annulus.

Let $C_4$ denote the doubling constant for $m$, as in \eqref{HH4a}, then set
$C'_4 = 2C_4+4$, and let $\tau > 0$ be small, to be chosen soon, depending on $C_4$.

First assume that $(1+\tau) r_1\leq r_2 \leq 2 r_1$.
Then the integrand on the left is comparable to $r_1^2 m(B(y,r_2)^{-2}$ and
$m(B_{r_2} \setminus B_{r_1})$ is comparable to $m(B(y,r_2)$; the desired inequality follows,
with a constant that depends on $\tau$, because $\int_{r_1}^{r_2} \frac{ds}{s} \geq C \tau^{-1}$.
When $r_2 > 2 r_1$, this is also easy: cut $B_{r_2} \setminus B_{r_1}$ into annuli
of modulus comparable to $1$, and prove the inequality separately on each one as we just did.

We may now assume that $r_2 \leq (1+\tau) r_1$; our defense will be that this does not happen
in the present circumstances. We claim that if $\tau$ is chosen so small that $C'_4\tau < 1$,
$r_2 \leq (1+\tau) r_1$ implies that $(1+C'_4 \tau) r_1 \geq \delta(y)$.

Indeed, suppose instead that $(1+C'_4 \tau) r_1 < \delta(y)$.
By definition, $\gamma(r_1) = 2 \gamma(r_1)$, hence by \eqref{defgg}
\begin{equation} \label{addedg36}
\int_{r_1}^{r_2} \frac{t^2}{m(B(y,t))} \, \frac{dt}{t} 
= \int_{r_2}^{\delta(y)} \frac{t^2}{m(B(y,t))} \, \frac{dt}{t}
\geq \int_{(1+\tau) r_1}^{(1+C'_4\tau) r_1} \frac{t^2}{m(B(y,t))} \, \frac{dt}{t}.
\end{equation}
The left-hand side is at most $(r_2-r_1) r_2 m(B_{r_1})^{-1} \leq 2\tau r_1^2 m(B_{r_1})^{-1}$
because $\tau < 1$, and the left-hand side is at least
$[(C'_4-1) \tau r_1] \, [(1+C'_4\tau) r_1] \, m(B_{2r_1})^{-1}$
because $C'_4\tau < 1$; since $m(B_{2r_1}) \leq C_4 m(B_{r_1})$ because $m$ is doubling,
this is at least $C_4^{-1} (C'_4-1) \tau r_1^2 m(B_{r_1})^{-1}$. We chose $C'_4 = 2C+4$, 
and the ensuing contradiction proves the claim. 

For the purposes of this lemma, we can take $\tau = (100 C'_4)^{-1}$, and then we just proved that
\eqref{Green99} holds as soon as $(1+10^{-2}) r_1 \leq \delta(y)$. Similarly,  
\begin{equation} \label{Green99bis}
\int_{B_{r_{i+1}} \setminus B_{r_i}} \frac{|x-y|^2}{m(B(y,|x-y|))^2} \, dm 
\leq C \int_{r_i}^{r_{i+1}} \frac{s^2}{m(B_s)} \frac{ds}{s} \leq C \gamma(r_i)
\end{equation}
as long as $(1+10^{-2}) r_i \leq \delta(y)$. We want to show that this does not
happen for $i \leq 2$, so we need a control on the variations of $\delta(y)-r_i$
along our sequence. Let us check that 
\begin{equation} \label{add37}
\delta(y) - r_{i} \leq 3 (\delta(y) - r_{i+1}) 
\, \text{ for $i \geq 0$ such that } r_{i+1}\geq \frac56 \delta(y).
\end{equation}
Suppose not, set $R = \delta(y) - 3 (\delta(y) - r_{i+1}) > r_i$, and observe that (as in \eqref{addedg36})
\begin{equation} \label{add38}
\int_{R}^{r_{i+1}} \frac{t^2}{m(B(y,t))} \, \frac{dt}{t} 
< \int_{r_i}^{r_{i+1}} \frac{t^2}{m(B(y,t))} \, \frac{dt}{t} 
= \int_{r_{i+1}}^{\delta(y)} \frac{t^2}{m(B(y,t))} \, \frac{dt}{t}.
\end{equation}
When we replace $m(B(y,t))$ by the larger number $m(B(y,r_{i+1}))$ on the left-hand side,
we get a smaller integral; similarly, when we replace $m(B(y,t))$ by the smaller number $m(B(y,r_{i+1}))$ 
on the right-hand side, we get a larger integral. Hence
$\int_{R}^{r_{i+1}} t dt < \int_{r_{i+1}}^{\delta(y)} t dt$.
Notice that $R \geq \delta(y)/2$ because $r_{i+1}\geq \frac56 \delta(y)$
and the interval on the left is twice as long as on the right; this gives a contradiction, and 
\eqref{add37} follows.

We may now prove that in the present circumstances, \eqref{Green99} holds, 
and even \eqref{Green99bis} for $0 \leq i \leq 2$.
Indeed, we start from $r_0 = r = |x-y| \leq \delta(y)/2$, so $\delta(y) - r_0 \geq \delta(y)/2$,
and it follows from a short iteration of \eqref{add37} that $\delta(y) - r_2 \geq \delta(y)/24$,
and so $(1+10^{-2}) r_2 \leq \delta(y)$.

We may now return to \eqref{pre99}. Since $\gamma(r_1) \approx \gamma(r_0) = \gamma(r)$ 
by \eqref{defri},
\eqref{pre99} and \eqref{Green99} imply that
\begin{equation} \label{1lbgxy}
1 \leq \frac{C}{\gamma(r)} \int_{B_{r_2} \setminus B_{r_1}} |\nabla \g^\rho|^2 \, dm.
\end{equation}

Next, 
since $\g^\rho$ is a solution to $Lu = 0$ in $B_{r_3} \setminus B_{r}$, we can use the Cacciopoli inequality (Lemma \ref{CaccioI}) with the test function $\eta_2(x):= \alpha_{r_2}(|x-y|)[1-\alpha_{r_0}(|x-y|)]$, 
which lies 
in $C^\infty_0(B_{r_3} \setminus B_{r})$ and satisfies $\eta_2 \equiv 1$ on $B_{r_2} \setminus B_{r_1}$. 
This yields
\[\int_{B_{r_2} \setminus B_{r_1}} |\nabla \g|^2 \, dm \leq C \int_{B_{r_3} \setminus B_{r_0}} |\g^\rho|^2 |\nabla \eta_2|^2 \, dm \leq C \sup_{B_{r_3} \setminus B_{r_0}} (\g^\rho)^2 \int_{B_{r_3} \setminus B_{r_0}} |\nabla \eta_2|^2 \, dm.\]
We bound the gradient of $\eta_2$ with the help of \eqref{bdalpha'}, \eqref{defgg'}, and \eqref{defri}. 
We find 
that $|\nabla \eta_2(x)| \leq \frac{C}{\gamma(r)} \frac{|x-y|}{m(B_{|x-y|})}$.
So by \eqref{Green99bis}
\[\int_{B_{r_3} \setminus B_{r_0}} |\nabla \eta_2|^2 \, dm \leq \frac{C}{\gamma(r)^2} 
\int_{B_{r_3} \setminus B_{r_0}} \frac{|x-y|^2}{m(B(y,|x-y|))^2} \, dm \leq \frac{C}{\gamma(r)}.\] 
As a consequence, 
\[\int_{B_{r_2} \setminus B_{r_1}} |\nabla \g^\rho|^2 \, dm \leq \frac{C}{\gamma(r)} \sup_{B_{r_3} \setminus B_{r_0}} (\g^\rho)^2 \]
and, 
together with \eqref{1lbgxy}, 
\begin{equation} \label{2lbgxy}
1 \leq \frac{C}{\gamma(r)^2} \sup_{B_{r_3} \setminus B_{r_0}} (\g^\rho)^2.
\end{equation}

To conclude, we invoke \eqref{gincrease} and then the interior Harnack inequality (Lemma \ref{HarnackI}), 
to get (since $r_0 = r$)
\[ \sup_{B_{r_3} \setminus B_{r_0}} \g^\rho \leq \sup_{B_{r} \setminus B_{r/2}} \g^\rho \leq C \g^\rho(x).\]
We use the above estimate in \eqref{2lbgxy} and notice $\gamma(r)$ is exactly the bound required 
for Lemma \ref{grhoxylb} (see \eqref{defgg}); the lemma follows.
\ep

\begin{lemma} \label{NgrhoL2}
Suppose that $\rho \leq 10^{-2}\delta(y)$ and $r < \delta(y)/2$, and set $B_r=B(y,r)$ as usual. Then
\[\int_{\Omega\setminus B_r} |\nabla \g^\rho|^2 \, dm \leq C\int_{r}^{\delta(y)} \frac{s^2}{m(B(y,s))} \, \frac{ds}{s},\]
where $C>0$ depends on $C_1$ to $C_6$, $C_A$, and $n$.
\end{lemma}

\bp By \eqref{defgrho}
for the test function $v := \g^\rho$ (and the elliptic condition \eqref{defEllip2}),
\begin{equation} \label{GreenN0}
\int_\Omega |\nabla \g^\rho|^2 dm 
\lesssim \int_\Omega \A\nabla \g^\rho \cdot \nabla \g^\rho \, dm 
= \fint_{B_\rho} \g^\rho \, dm.
\end{equation}
Then cut the last integral in two; this yields
\[\begin{split}
\int_\Omega |\nabla \g^\rho|^2 dm & \lesssim  \fint_{B_\rho} \left( \g^\rho - \fint_{B_{50\rho} \setminus B_{25\rho}} \g^\rho\, dm \right)  \, dm + \fint_{B_{50\rho} \setminus B_{25\rho}} \g^\rho\, dm \\
& \lesssim  \fint_{B_{50\rho}} \left|\g^\rho - \fint_{B_{50\rho} \setminus B_{25\rho}} \g^\rho\, dm \right|  \, dm + \fint_{B_{50\rho} \setminus B_{25\rho}} \g^\rho\, dm. 
\end{split}\]
We use the Poincar\'e inequality (Theorem \ref{Poincare1}) to bound the 
first integral, and Lemma \ref{grhoxyclose} to bound the last one. 
Notice in particular that $40\rho \leq 2|x-y| < \delta(y)$ for $x\in B_{50\rho} \setminus B_{25\rho}$,
so the lemma applies. This gives
\[\begin{split}
\int_\Omega |\nabla \g^\rho|^2 dm & \lesssim \rho \left(\fint_{B_{50\rho}} \left|\nabla \g^\rho \right|^2  \, dm\right)^\frac12 + \int_{25\rho}^{\delta(y)} \frac{s^2}{m(B_s)} \frac{ds}{s}\\
& \lesssim \frac{25\rho}{m(B_{25\rho})^{1/2}} \left(\int_{\Omega} \left|\nabla \g^\rho \right|^2  \, dm\right)^\frac12 + \int_{25\rho}^{\delta(y)} \frac{s^2}{m(B_s)} \frac{ds}{s}.
\end{split}\]
Since $m$ is doubling, $m(B_t) \leq C_4 m(B_{25\rho})$ for $25\rho \leq t \leq 50\rho$, so
\[
\frac{(25\rho)^2}{m(B_{25\rho})} \leq C_4 \int_{25\rho}^{50\rho} \frac{t^2}{m(B_t)} \frac{dt}{t}
\leq C_4 \gamma(25\rho)
\]
(because $50\rho < \delta(y)$ and by the definition \eqref{defgg}).
So the above estimate can be written
\[\int_\Omega |\nabla \g^\rho|^2 dm \lesssim \gamma(25\rho)^\frac12 \left(\int_{\Omega} \left|\nabla \g^\rho \right|^2  \, dm\right)^\frac12 + \gamma(25\rho),\]
or
\begin{equation} \label{L2bdgrho}
\int_\Omega |\nabla \g^\rho|^2 dm \lesssim \gamma(25\rho):= \int_{25\rho}^{\delta(y)} \frac{s^2}{m(B_s)} \frac{ds}{s}.
\end{equation}

The estimate \eqref{L2bdgrho} is 
already good, and it proves
Lemma \ref{NgrhoL2} for any $r>0$ such that
$\gamma(r)\geq \gamma(25\rho)/2$. Assume now that $\gamma(r) \leq \gamma(25\rho)/2$. 
Since $\gamma$ is decreasing, there is a unique $R \in [25\rho, r)$ such that 
$\gamma(R) = 2\gamma(r)$.
Then $\g^\rho$ is a solution to $Lu=0$ in $\Omega \setminus B_{R}$.
By Lemma~\ref{grhoxyclose} and the proof of \eqref{gincrease}, 
$\g^\rho(x) \leq C \gamma(R) = 2C\gamma(|x-y|)$ 
for any $x\in \Omega \setminus B_{R}$. 
We claim that if $\eta_3$ is a nonnegative smooth function that satisfies 
$\eta_3 \equiv 0$ on $B_{R}$, 
we have the following Cacciopoli-type inequality:
\begin{equation} \label{GreenR} 
\int_{\Omega} |\nabla \g^\rho|^2 |\eta_3|^2 \, dm 
\leq C \int_{\Omega} (\g^\rho)^2 |\nabla \eta_3|^2 dm.
 \end{equation}
The above bound is not an application of the Cacciopoli inequalities stated in Section \ref{SSolutions}, 
because the test function $\eta_3$ is not 
contained in a ball $2B$ such that $\g^\rho$ is a solution in $2B$. 
However, the proof of \eqref{GreenR} is 
very similar to the proof of the usual Cacciopoli inequality, and we leave it to the reader.
It relies on
the fact that $\g^\rho\eta_3 \in W_0$ even though $\eta_3$ is not compactly supported.

We use \eqref{GreenR} with $\eta_3(x) := 1-\alpha_{R}(|x|)$, 
where $\alpha_{R}$ 
is the function in \eqref{alpha}. 
Notice that $\eta_3(x) = 0$ on $B_R$ and $\eta_3(x) = 1$ on $\R^n \sm B_r$ 
(because $\gamma(|x|) < \gamma(r) = \gamma(R)/2$).
So we obtain
\begin{equation} \label{GreenRZ} 
\int_{\Omega\setminus B_r} |\nabla \g^\rho|^2 \, dm 
\leq \frac{C}{\gamma(r)^2}  \int_{B_{r}\setminus B_{R}} 
|\g^\rho(x)|^2 \frac{|x-y|^2}{m(B(y,|x-y|))^2} dm(x).
 \end{equation}
by
\eqref{bdalpha'} and
\eqref{defgg'}. But for $x\in B_{r}\setminus B_{R}$, Lemma \ref{grhoxyclose} says that
$\g^\rho(x) \leq C \gamma(R) = 2C \gamma(r)$; hence by the proof of \eqref{Green99},
\begin{equation} \label{GreenRY} 
\int_{\Omega\setminus B_r} |\nabla \g^\rho|^2 \, dm 
\lesssim \int_{B_{r}\setminus B_{R}}  
\frac{|x-y|^2}{m(B(y,|x-y|))^2} dm(x) 
\lesssim \gamma(R) 
\lesssim \gamma(r).
 \end{equation}
 Lemma \ref{NgrhoL2} 
 follows.
 \ep

\begin{lemma} \label{NgLq} There exists  
$q>1$, 
that depends only on $C_4$, 
 such that for any $y\in \Omega$  
 and any $\rho\leq \delta(y)/100$,
 \[\int_{B(y_0,2\delta(y))} |\nabla \g^\rho|^q \, dm \leq C  m(B_{\delta(y)}) \left( \frac{\delta(y)}{m(B_{\delta(y)})}\right)^q,\]
 where $y_0\in \Gamma$ is such that $|y-y_0| = \delta(y)$, 
 and where $C>0$ depends only on $C_1$ to $C_6$, $C_A$, and $n$.
\end{lemma}

\begin{remark}
If
$\gamma(r) := \int_{r}^\delta(y) \frac{s^2}{m(B_s)} \frac{ds}{s}$ is uniformly bounded, 
or in other words 
if $\int_{0}^\delta(y) \frac{s^2}{m(B_s)}$ is finite, then 
by \eqref{L2bdgrho} we can take $q=2$ in Lemma \ref{NgLq}.
\end{remark}

\bp We first start by proving general results, which are only consequences of 
the doubling property
(\hyperref[H4]{H4}). 
There exists $\alpha\in (0,1)$ such that 
\begin{equation} \label{lbdoubling}
 m(B) \leq \alpha m(2B)
\end{equation} 
for every ball $B\subset \R^n$ such that $2B\subset \Omega$.
Indeed, if $r$ denotes
the radius of $B$, then we can find a ball $B_0$ of radius $r/2$ in $2B\setminus B$. 
Then $B \subset 3B_0$, hence $m(B) \leq C_4^2 m(B_0)$ by (\hyperref[H4]{H4}),
and now $m(2B) - m(B) \geq m(B_0) \geq C_4^{-2} m(B)$, and \eqref{lbdoubling}
holds with $\alpha = (1+C_4^{-2})^{-1}$.

Similarly to \eqref{mdoublingG}, the estimate \eqref{lbdoubling} can be improved into
\begin{equation} \label{lbdoublingG}
m(B(y,r)) \leq C\left(\frac{r}{s}\right)^{2\epsilon} m(B(y,s)) \qquad \text{ for } r\leq s\leq\delta(y),
\end{equation} 
where $C,\epsilon>0$ depends only on $C_4$,
and we use $2\epsilon$ instead of $\epsilon$ to simplify the later computations. Indeed, let $k$ be the integer such that $2^{-k-1} < r/s \leq 2^{-k}$. 
Then by \eqref{lbdoubling}
\[\begin{split}
m(B(y,r)) & \leq m(B(y,2^{-k}s)) \leq \alpha^{k} m(B(y,s)) \leq (2^{-k})^{\ln_2(1/\alpha)} m(B(y,s)) \\
& \leq \frac1\alpha \left( \frac rs \right)^{\ln_2(1/\alpha)} m(B(y,s)).
\end{split}.\]
The claim \eqref{lbdoublingG} follows by taking $\epsilon = \frac12 \ln_2(1/\alpha)>0$.

Let us use again $B_s$ for $B(y,s)$. The inequality \eqref{lbdoublingG} implies in particular that 
\[\frac{r^\epsilon}{m(B_r)} \geq C\left(\frac{r}{s}\right)^{-\epsilon} \frac{s^{\epsilon}}{ m(B_s)} \qquad \text{ for } r \leq s \leq \delta(y),\]
which proves that the function $r \to \delta(y)^{1-\epsilon}r^\epsilon/m(B_r)$ reaches all the values between $\delta(y)/m(B_{\delta(y)})$ and $+\infty$. Moreover if $t$ is in the given range, the values of $t$ 
that satisfy$t= \delta(y)^{1-\epsilon}r^\epsilon/m(B_r)$ are all the same up to a harmless constant.

\medskip

For the next step we 
prove weak $L^q$ estimates on the gradient of $\g^\rho$. 
Set $\wh \Omega_t := \{x\in \Omega, \, |\nabla g^\rho(x)|>t\}$. 
Thanks to Lemma \ref{NgrhoL2}, for all $r\in (0,\delta(y)/2)$, we have
\[m(\wh \Omega_t) \leq m(B_r) + \frac C{t^2} \int_{r}^{\delta(y)} \frac{s^2}{m(B_s)} \frac{ds}{s}.\] 
Then by \eqref{lbdoublingG}
\[m(\wh \Omega_t) 
\leq m(B_r) +\frac C{t^2}  \frac{1}{m(B_r)} \int_{r}^{\delta(y)} \left( \frac{r}{s} \right)^{2\epsilon} s\, ds .\] 
Or,
since we can always chose $\epsilon <1$,
\begin{equation} \label{mhO1}
m(\wh \Omega_t) \leq  m(B_r) + \frac C{t^2}  \frac{\delta(y)^{2(1-\epsilon)} r^{2\epsilon}}{m(B_r)}.
\end{equation} 
We aim to optimize the above expression in $r$. But we shall only care about big values of $t$,
so let us only consider $t \geq \delta(y)/m(B_{\delta(y)})$ for the moment.
First assume that $2^{-\epsilon}\delta(y)/m(B_{\delta(y)/2}) \geq t \geq \delta(y)/m(B_{\delta(y)})$. 
Then we choose $r = \delta(y)/2$ in \eqref{mhO1}, and it is easy to see, using (\hyperref[H4]{H4})
and \eqref{mdoublingG} in particular, that 
\[m(\wh \Omega_t) \leq C m(B_{\delta(y)})^{-\frac \epsilon{d-\epsilon}} \left( \frac{\delta(y)}{t}\right)^{\frac{d}{d-\epsilon}},\]
where 
$d = d_m > 0$ is the (possibly large) exponent of \eqref{mdoublingG}.
Notice that we may always replace $d_m$ with a larger exponent in \eqref{mdoublingG}, so
we may assume that $d\geq 2\epsilon$, and this way the exponent 
$-\frac \epsilon{d-\epsilon}$ is rather small and negative.
We strive for the same bound when $t \geq 2^{-\epsilon}\delta(y)m(B_{\delta(y)/2})$. We then take
$r$ such that 
\begin{equation} \label{defr(t)}
t = \delta(y)^{1-\epsilon}r^\epsilon/m(B_r), 
\end{equation}
and we have seen in the previous paragraph that, even if we may have different choices for $r$, they are all the same up to a constant. Using $r$ as in \eqref{defr(t)} in \eqref{mhO1}, we obtain that 
\begin{equation} \label{mhO2}
m(\wh \Omega_t) \lesssim m(B_r) = \frac{\delta(y)^{1-\epsilon}r^\epsilon}{t}.
\end{equation}
Yet, by \eqref{defr(t)} and \eqref{mdoublingG},
\[\begin{split}
m(B_r) & = \frac{\delta(y)}{t} \left(\frac{r}{\delta(y)}\right)^\epsilon \leq C \frac{\delta(y)}{t} \left(\frac{m(B_r)}{m(B_{\delta(y)})}\right)^\frac\epsilon d,
\end{split}\]
or equivalently
\[m(B_r)^{\frac{d-\epsilon}{d}} \leq C \, \frac{\delta(y)}{t} m(B_{\delta(y)})^{-\frac\epsilon d}.\]
Using this bound on $m(B_r)$ in \eqref{mhO2}, we obtain that 
\begin{equation} \label{mhO3}
m(\wh \Omega_t) \leq C m(B_r) \leq C' m(B_{\delta(y)})^{-\frac \epsilon{d-\epsilon}} \left( \frac{\delta(y)}{t}\right)^{\frac{d}{d-\epsilon}}
\end{equation}
as desired.

We are ready to conclude. 
We write $q_0$ for $\frac{d}{d-\epsilon} >1$. The bound \eqref{mhO3} becomes 
\begin{equation} \label{mhO4}
m(\wh \Omega_t) \leq C m(B_{\delta(y)})^{1-q_0} \left( \frac{\delta(y)}{t}\right)^{q_0}.
\end{equation}
We take $q = (1+q_0)/2 >1$. Then 
\[\begin{split}
\frac1q 
\int_{B(y_0,2\delta(y))} |\nabla \g^\rho|^q \, dm 
& = \int_0^\infty t^{q-1} m(\wh \Omega_t \cap B(y_0,2\delta(y)))  dt 
\\
& \leq m(B(y_0,2\delta(y))) \int_0^{\delta(y)/m(B_{\delta(y)})} t^{q-1} dt + \int_{\delta(y)/m(B_{\delta(y)})}^\infty t^{q-1} m(\wh \Omega_t) dt.
\end{split}\]
Then by (\hyperref[H4]{H4}) and \eqref{mhO4},
\[\begin{split} 
\int_{B(y_0,2\delta(y))} |\nabla \g^\rho|^q \, dm & \lesssim m(B_{\delta(y)}) \left( \frac{\delta(y)}{m(B_{\delta(y)})}\right)^q + m(B_{\delta(y)}) \left( \frac{\delta(y)}{m(B_{\delta(y)})}\right)^{q_0} \int_{\delta(y)/m(B_{\delta(y)})}^\infty t^{q-q_0-1} dt \\
& \lesssim m(B_{\delta(y)}) \left( \frac{\delta(y)}{m(B_{\delta(y)})}\right)^q
\end{split}\]
since $q<q_0$. Lemma \ref{NgLq}
follows. \ep

We are now ready for the big theorem.

\begin{theorem} \label{GreenEx}
There exists a non-negative function $g:\ \Omega \times \Omega \to \R \cup \{+\infty\}$ 
with the following properties.
\begin{enumerate}[(i)]
\item For any $y\in \Omega$ and 
any function $\alpha \in C_0^\infty(\R^n)$ 
such that $\alpha \equiv 1$ in a neighborhood of $y$
\begin{equation} \label{GreenE1}
(1-\alpha) g(.,y) \in W_0.
\end{equation}
In particular, $g(.,y) \in \WW
(\overline \Omega \setminus \{y\}) \subset L^1_{loc}(\overline \Omega \setminus \{y\},dm)$ and $\Tr [g(.,y)] = 0$ on $\Gamma$.
\smallskip

\item There exists $q>1$ that depends only on $C_4$ such that for every choice of $y\in \Omega$, 
\begin{equation} \label{GreenE4}
\nabla g(.,y) \in L^q(B(y,\delta(y)),dm).
\end{equation}
\smallskip

\item For 
$y\in \Omega$ and  
$\varphi \in C^\infty_0(\Omega)$,
\begin{equation} \label{GreenE6}
\int_\Omega A \nabla_x g(x,y)\cdot \nabla \varphi(x) dx = \varphi(y).
\end{equation}
In particular, $g(.,y)$ is a solution 
of $Lu=0$ 
in $\Omega \setminus \{y\}$. 
\smallskip

\end{enumerate}
In addition, the following bounds hold.
\begin{enumerate}[(i)] \setcounter{enumi}{3}
\item For $x,y\in \Omega$ such that $|x-y| \geq \delta(y)/10$, 
\begin{equation} \label{GreenE7}
0 \leq g(x,y) \leq C\frac{|x-y|^2}{m(B(y,|x-y|) \cap \Omega)},
\end{equation}
where $C>0$ depends on $C_1$ to $C_6$, $C_A$, and $n$.

\smallskip

\item For $x,y\in \Omega$ such that $|x-y|\leq \delta(y)/2$,
\begin{equation} \label{GreenE2}
c \int_{|x-y|}^{\delta(y)} \frac{s^2}{m(B(y,s))} \frac{ds}{s} \leq g(x,y) \leq C \int_{|x-y|}^{\delta(y)} \frac{s^2}{m(B(y,s))} \frac{ds}{s},
\end{equation}
where $C>0$ depends on $C_1$ to $C_6$, $C_A$, and $n$; and where $c>0$ depends on $C_4$, $C_6$, $C_A$, and $n$.

\smallskip

\item For $r \in (0,\delta(y)/2)$ and $y\in \Omega$,
\begin{equation} \label{GreenE3}
\int_{\Omega \setminus B(y,r)} |\nabla_x g(x,y)|^2 dm(x) \leq C \int_{r}^{\delta(y)} \frac{s^2}{m(B(y,s))} \frac{ds}{s},
\end{equation}
where $C>0$ depends again on $C_1$ to $C_6$, $C_A$, and $n$.

\smallskip

\item If $q>1$ is the exponent in \eqref{GreenE4}
\begin{equation} \label{GreenE5}
\left(\fint_{B(y,\delta(y))} |\nabla_x g(x,y)|^{q} dm(x)\right)^\frac1q \leq C \frac{\delta(y)}{m(B(y,\delta(y))},
\end{equation}
where $C>0$ depends as usual on $C_1$ to $C_6$, $C_A$, and $n$.
\end{enumerate}
\end{theorem}

\bp As we shall see, we already have all the desired estimates on the $\g^\rho := g^\rho(x,y)$; the proof will mainly
consist in choosing a right limit to those $\g^\rho$.

We start with a standard exercise on compactness. 
For every compact set $K$ in $\overline{\Omega} \setminus \{y\}$, 
Lemmas \ref{grho>0}, \ref{grhoxyfar}, and \ref{grhoxyclose} prove that
the set $F_K:= \big\{\g^\rho(x), \, ; \,  x\in K \text{ and }  0<\rho<\dist(y,K)/100 \big\}$ is bounded;
then by Lemma \ref{HolderB} the functions $\g^\rho$, $\rho <  \dist(y,K)/100$ are H\"older continuous
on $K$ (on a slightly smaller compact set), with uniform bounds. 
In particular, for every compact set $K \subset \overline{\Omega} \setminus \{y\}$
the set $A_K:= \big\{\g^\rho, \, ; \,  0<\rho \leq \dist(y,K)/200 \big\}$ - seen as a subset of the continuous functions on $K$ - is equicontinuous. 
Ascoli's theorem entails that $A_K$ is relatively compact in $C^0(K)$, that is we can find a sequence
of radii $\rho$, that tends to $0$,  such that the corresponding $\g^\rho$ converge, uniformly on $K$,
to a (continuous) function written $\g_K$ .
We take a sequence of compacts sets $K_i$ such that $K_i \subset K_{i+1}$ and 
$\bigcup_i K_i = \overline \Omega \setminus \{y\}$, and by a diagonal process, we can find a sequence $(\rho_\eta)_{\eta \in \bN}$ and a continuous function $\g$ on $\overline \Omega \setminus \{y\}$ 
such that $\rho_\eta \to 0$ and
\begin{equation} \label{grhocv1}
\text{$\g^{\rho_\eta}$ converges to $\g$ uniformly on every
compact set of $\overline \Omega \setminus \{y\}$.}
\end{equation}

We shall use again the cut-off functions $\alpha_r$ defined in \eqref{alpha} and their properties. Set $\wt \alpha_r(x) = \alpha_r(|x-y|)$; 
we want to prove that the $\{\g^{\rho_\eta}(1-\wt\alpha_r)\}_{\eta \in \bN}$ 
form a Cauchy sequence in $W_0$.
For any $r < \delta(y)/2$, define $r_1 \in (r,\delta(y))$ 
as the only value such that $\gamma(r_1) = \gamma(r)/2$; then 
for $\eta,\nu \in \bN$,
\begin{align} \label{grhoggy}
\int_\Omega |\nabla [(\g^{\rho_\eta}-\g^{\rho_\nu})& (1-\wt\alpha_r)]|^2
\, dm \nonumber\\
& \leq 2 \int_\Omega |\nabla [\g^{\rho_\eta} -\g^{\rho_\nu}]|^2 |1-\wt\alpha_r|^2
\, dm + 2 \int_\Omega |\g^{\rho_\eta} -\g^{\rho_\nu}|^2 |\nabla \wt\alpha_r|^2
\, dm \nonumber\\
& \leq C \int_\Omega |\g^{\rho_\eta} -\g^{\rho_\nu}|^2 |\nabla \wt\alpha_r|^2
\, dm
\end{align} 
where, for the last line, we take $\eta,\nu$ big enough so that $\rho_\eta, \rho_\nu$ 
are way smaller than $r$ and we use the Cacciopoli-type inequality \eqref{GreenR}. 
Since $\nabla \wt\alpha_r$ 
is supported in $B_{r_1} \setminus B_r$, and the later is a compact set in 
$\Omega \setminus \{y\}$, the convergence \eqref{grhocv1} forces the right-hand side of \eqref{grhoggy} to tend to $0$.
In addition, all the $\g^\rho_\eta$ have a vanishing trace,
and so do the $\g^{\rho_\eta} (1-\wt\alpha_r)$
(see Lemma~\ref{lmult}). 
We deduce that $\{\g^{\rho_\eta}(1-\wt\alpha_r)\}_{\eta \in \bN}$ 
is indeed a Cauchy sequence in $W_0$, so it converges {\bf strongly} in $W_0$ to a function $g^{(r)}$. By uniqueness of the limit, 
 $\g^{(r)} = \g (1-\wt\alpha_r)$. 
 In short, we proved that for $0 < r < \delta(y)/2$
\begin{equation} \label{grhocv2}
\text{$\g^{\rho_\eta}(1-\wt\alpha_r)$ 
converges strongly to $\g(1-\wt\alpha_r)$ 
in $W_0$. }
\end{equation}
Notice that $\g$ has a gradient in $L^2_{loc}(\Omega \setminus \{y\}, dm)$ defined as 
\begin{equation} \label{grhocv3}
\text{$\nabla \g(x) = \nabla [\g(1-\wt\alpha_r)]$ if 
$\wt\alpha_r(x) = 0$.} 
\end{equation}

We still need a last convergence, one that goes across the pole $\{y\}$. Lemma \ref{NgLq} provides us with the uniform bound
\[\left( \fint_{B_{\delta(y)}} |\nabla g^{\rho_\eta}|^q \, dm \right)^\frac1q \leq C \frac{\delta(y)}{m(B_{\delta(y)})}.\]
So, 
 up to a subsequence, the quantities $\nabla \g^{\rho_\eta}$ converges {\bf weakly} to a function $h \in L^q(B_{\delta(y)}, dm)$. But since $\nabla \g^{\rho_\eta}$ already converges to $\nabla \g$ in $L^2_{loc}(B_{\delta(y)} \setminus \{y\}, dm)$, it forces $\nabla \g = h$ except maybe at the point $y$, but it has no importance because $m(\{y\}) = 0$. To sum up,
\begin{equation} \label{grhocv4}
\text{$\nabla \g^{\rho_\eta}$ converges weakly to $\nabla \g$ in $L^q(B_{\delta(y)}, dm)$.}
\end{equation}

\medskip

Now
let us show $(i)$--$(vii)$ of the theorem.
For the first statement $(i)$, let us start with the more likely situation where 
$\lim_{r \to 0} \gamma(r) = +\infty$. Since $\alpha = 1$ near $y$, 
we can we can find $s > 0$ such that $|y-x| > s$ when $\alpha(x) \neq 1$.
Choose $r$ so small that $\gamma(r) > 2\gamma(s)$; then for $x$ such that $\alpha(x) \neq 1$,
$\gamma(|x-y|) \leq \gamma(s) < \frac12 \gamma(r)$, so 
$\wt\alpha_r(x) = \alpha_r(|x-y|) = 0$ by \eqref{alpha3}.
Because of this, $(1-\alpha)\g = (1-\alpha)(1-\wt\alpha_r)\g$. 
This function lies in $W_0$, as needed,  by \eqref{grhocv2} and Lemma~\ref{lmult}.

In the other case when $\lim_{r \to 0} \gamma(r) < +\infty$, we are in the happy situation where \eqref{L2bdgrho}
says that $\int_\Omega |\nabla \g^\rho|^2 dm \leq C$, with a constant that depends on $y$, 
but not on $\rho$; then the almost everywhere pointwise limit $\g^\rho$ satisfies 
$\int_\Omega |\nabla \g|^2 dm \leq C$ too, and its trace is still $0$ on $\Gamma$.
Finally $(1-\alpha)\g$ does the same; see for instance  the proof of Lemma \ref{lem6.3} for the limit, and Lemma \ref{defTrWE} for the product. This takes care of $(i)$.

\smallskip
The statement $(ii)$ is part of \eqref{grhocv4}.
\smallskip 

For the identity $(iii)$, we take $r_0$ so that 
$\gamma(r_0) = \frac12 \gamma(\delta(y)/2)$. 
We then write $\varphi = \varphi_1 + \varphi_2$ where 
$\varphi_1 = \varphi \wt\alpha_{r_0}$ 
and $\varphi_2 = \varphi (1-\wt\alpha_{r_0})$. 
The function $\varphi_1$ is continuous and smooth enough
for $\nabla \varphi_1$ to lie in $L^{q'}(B_{\delta(y)},dm)$,
and so by \eqref{grhocv4}
and then definition \eqref{defgrho},
\begin{equation} \label{grhocv5}
\int_\Omega A \nabla \g\cdot \nabla \varphi_1 \, dx = \lim_{\eta \to \infty} \int_\Omega \A \nabla \g^{\rho_\eta} \cdot \nabla \varphi_1 \, dm = \lim_{\eta \to \infty} \fint_{B_{\rho_\eta}} \varphi_1 \, dm = \varphi_1(x) = \varphi(x).
\end{equation}
When $x \in B_{r_0}$, $\wt\alpha_{r_0}(x) = \alpha_{r_0}(|x-y|)= 1$ by \eqref{alpha22},
and hence $\varphi_2(x)=0$. But otherwise $\wt\alpha_{\delta(y)/2}(x) = 
\alpha_{\delta(y)/2}(|x-y|) = 0$ because $\gamma(|x-y|) \leq \gamma(r_0) = \frac12 \gamma(\delta(y)/2)$, and by \eqref{alpha3}. Hence 
$\g(1-\alpha_{\delta(y)/2}) = \g$ on the support of $\nabla \varphi_2$, and so
 \begin{equation} \label{grhocv6} \begin{split}
\int_\Omega A \nabla \g\cdot \nabla \varphi_2 \, dx 
& = \int_\Omega A \nabla [\g(1-\wt\alpha_{\delta(y)/2})] 
\cdot \nabla \varphi_2 \, dx 
=  \lim_{\eta \to \infty} \int_\Omega 
\A \nabla [\g^{\rho_\eta}(1-\wt\alpha_{\delta(y)/2})] 
\cdot \nabla \varphi_2 \, dm \\
& = \lim_{\eta \to \infty} \int_\Omega \A \nabla \g^{\rho_\eta} \cdot \nabla \varphi_2 \, dm
= \lim_{\eta \to \infty} \fint_{B_{\rho_\eta}} \varphi_2 \, dm = \varphi_2(x) = 0
\end{split} \end{equation}
where we used \eqref{grhocv2} for the second equality
and then returned by the same path. 
The combination of \eqref{grhocv5} and \eqref{grhocv6} infers $(iii)$.

The estimates given in $(iv)$ and $(v)$ are direct consequences of Lemmas \ref{grho>0}, \ref{grhoxyfar}, \ref{grhoxyclose}, \ref{grhoxylb}, and the convergence \eqref{grhocv1}. The bound found in $(vi)$ is due to Lemma \ref{NgrhoL2} and \eqref{grhocv2}, while $(vii)$ comes from Lemma \ref{NgLq} and \eqref{grhocv4}.
Theorem \ref{GreenEx}
follows. \ep

\begin{remark} \label{rmkL1g}
Before stating the next result, let us comment a bit on 
Theorem \ref{GreenEx}.
One can easily see that $g(.,y)$ lies 
in $L^1_{loc}(\overline \Omega\setminus \{y\}, dm)$, since the latter is bigger than 
the space of 
continuous functions on $\overline \Omega\setminus \{y\}$ (and $g(.,y)$ is continuous on $\overline \Omega\setminus \{y\}$ thanks to \eqref{grhocv1}). However, we said nothing about the integration of $g(.,y)$ on a neighborhood of $\{y\}$.
The fact that $g(.,y)$ can be integrated over a bounded region that covers $\{y\}$ is a simple consequence of \eqref{GreenE2}. Indeed, if $y\in \Omega$ and $r\leq \delta(y)/2$, 
first observe that
\[
\int_{B(y,r)} g(x,y) \, dm(x) 
\leq C \int_{B(y,r)} \int_{|x-y|}^{\delta(y)} \frac{s^2}{m(B(y,s))} \frac{ds}{s} \, dm(x) 
\leq C \int_{t=0}^r \int_{s=t}^{\delta(y)} \frac{m(B(y,t))}{m(B(y,s))} \, s\, ds \frac{dt}{t}.
\] 
Let $\epsilon\in (0,1)$ be the constant in \eqref{lbdoublingG}; then
\[\int_{B(y,r)} g(x,y) \, dm(x) \leq C  \int_{0}^r \int_{t}^{\delta(y)} \, s^{1-2\epsilon}ds \, t^{2\epsilon - 1}dt \leq C \delta(y)^{2-2\epsilon} r^{2\epsilon}.\]
In addition, by \eqref{GreenE7} 
\[\int_{B(y,\delta(y)) \setminus B(y,\delta(y)/2)} g(x,y) \, dm \leq C m(B(y,\delta(y)) \sup_{x\in B(y,\delta(y)) \setminus B(y,\delta(y)/2)} g(x,y) \leq C \delta(y)^2.\]
The combination of the last two 
estimates implies 
that
\begin{equation} \label{L1g}
\int_{B(y,\delta(y))} g(x,y) \, dm \leq C \delta(y)^2.
\end{equation}
Due to \eqref{grhocv1}, the functions $g^{\rho_\eta}(.,y)$ converges pointwise a.e. to 
$g(.,y)$ on $B(y,\delta(y))$. So by 
the Lebesgue domination theorem
(and the fact that bounds above are also valid for the $g^{\rho_\eta}$),
we even have
\[
\lim_{\eta \to \infty}\int_{B(y,\delta(y))} |g^{\rho_\eta}(x,y)- g(x,y)|\, dm(x) = 0 .\]
Together with \eqref{grhocv1}, we proved that 
\begin{equation} \label{L1gcv}
\text{ 
the functions $g^{\rho_\eta}(.,y)$ converge 
to $g(.,y)$ in $L^1_{loc}(\overline \Omega)$.}
\end{equation}
 \end{remark}

For the next lemma, we need some additional notation. We write $A^T$ for the transpose 
matrix of $A$, i.e. $(A^T)_{ij}(x) = A_{ji}(x)$ for all $1\leq i,j\leq n$ and $x\in \Omega$. Obviously, $A^T$ satisfies the ellipticity and boundedness conditions \eqref{defEllip}--\eqref{defBdd} with the same constant as $A$. The elliptic operator $L_T:= -\div A^T \nabla$ enjoys the very same properties as $L$, in particular, Theorem~\ref{GreenEx} 
yields the existence of $g^T :\, \Omega \cap \Omega \to \R \cup \{+\infty\}$ 
with the same properties as $g$ (except for \eqref{GreenE6}, where $A$ is replaced by $A^T$).

\begin{lemma} \label{GreenSym}
With the notation above,
\begin{equation}
g(x,y) = g_T(x,y) \quad \text{ for } x,y \in \Omega.
\end{equation}
In particular, the functions 
$x\to g(y,x)$ satisfy
the estimates in Theorem \ref{GreenEx}.
\end{lemma}

\bp
The result is the same as the one of \cite[Theorem 1.3]{GW} (or Lemma 10.6 in \cite{DFMprelim}). Yet, the limits we took in the proof of Theorem \ref{GreenEx} is a bit different to the one in \cite{GW} and \cite{DFMprelim}. So our result deserves a proof.

Actually, the convergence property \eqref{grhocv4} will make the proof very easy for us. Let $x,y\in \Omega$ be such that $x\neq y$. 
By our
construction of the Green functions, 
there exist two sequences $(\rho_\eta)_{\eta\in \bN}$ and $(\sigma_\nu)_{\nu \in \bN}$ 
such that $\rho_\eta, \sigma_\nu \to 0$ and 
\begin{equation} \label{grhocv8}
\text{$\g^{\rho_\eta}$ converges to $\g$ uniformly on any compact set of $\overline \Omega \setminus \{y\}$}
\end{equation}
and
\begin{equation} \label{grhocv9}
\text{$\g_T^{\sigma_\nu}$ converges to $\g_T$ uniformly on any compact set of $\overline \Omega \setminus \{x\}$.}
\end{equation}

Using \eqref{defgrho2}  and  \eqref{Green4} for both $g(.,y)$ and $g_T(.,x)$, we see 
that for any $\eta,\nu\in \bN$
\begin{equation} \label{grhocv0}
\int_\Omega \A \nabla g^{\rho_\eta}(z,y) \cdot \nabla g_T^{\sigma_\nu} (z,x) \, dm(z) = \fint_{B(y,\rho_\eta)} g_T^{\sigma_\nu}(z,x) \, dm(z) = \fint_{B(x,\sigma_\nu)} g^{\rho_\eta}(z,y) \, dm(z)
\end{equation}
We use the uniform convergence of $\g_T^{\sigma_\nu}$ on $\overline {B(y,|x-y|/2)} \subset \Omega \setminus \{x\}$ and the uniform convergence of $g(.,y)^{\rho_\eta}$ on $\overline {B(x,|x-y|/2)} \subset \Omega \setminus \{y\}$ given by \eqref{grhocv8}--\eqref{grhocv9} in the last equality of \eqref{grhocv0}. We get that 
$g_T(y,x) = g(x,y)$ as desired. \ep

\begin{lemma} \label{GreenDi}
The Green function satisfies
\begin{equation} \label{GreenDi1}
 g(x,y) \leq C \delta(x)^\alpha \frac{|x-y|^{2-\alpha}}{m(B(x,|x-y|) \cap \Omega)} \ \  \text{ for
 $x,y \in \Omega$ such that }   |x-y | \geq 4\delta(x),
\end{equation}
where $C>0$ and $\alpha>0$ depend only on $C_1$ to $C_6$, $C_A$, and $n$.
\end{lemma}

\bp
See the proof of \cite[Lemma 10.9]{DFMprelim}. The arguments are based on the pointwise bounds \eqref{GreenE7} and the H\"older regularity at the boundary (Lemma \ref{HolderB}). Actually, the coefficient $\alpha$ is the one of Lemma \ref{HolderB}. 
\ep

The next result that we wanted is the representation of solutions by Green functions. More precisely, we wanted to take a smooth function $f \in C^\infty_0(\Omega)$ and construct $u(x)$ for $x\in \Omega$ as
\begin{equation} \label{GreenFS1}
u(x) = \int_\Omega g(x,y) f(y) dm(y).
\end{equation}
We have seen in Remark \ref{rmkL1g} that $g(.,y)$ 
lies in $L^1_{loc}(\overline \Omega, dm)$. Moreover, due to 
Lemma~\ref{GreenSym}, 
we also have that $g(x,.)$ is in $L^1_{loc}(\overline \Omega, dm)$.
Yet, 
in the case of a general weights $w$,
we do not know if $g(x,.)$ lies  
in the unweighted space $L^1_{loc}(\overline \Omega)$. That is why, in the definition \eqref{GreenFS1}, the function $u$ has to be defined as an integral over the measure $m$.

In doing so, the formal identity satisfied by $u$ is not $Lu = f$ but $Lu = fw$, where $w$ is the weight used to define the measure $m$. Another way to see it, maybe more relevant, is to say that $w^{-1}Lu = f$. That is, we are solving the Dirichlet problem $\mathcal L u = f$ for an elliptic operator $\mathcal L = -w^{-1} \div [\A w\nabla ]$ where $\A$ satisfies the classical elliptic conditions \eqref{defEllip2}--\eqref{defBdd2}. 

At last, by using $L$ instead of $\mathcal L:= w^{-1}L$, we are somehow linking the measure $m$ to the plain Lebesgue measure on $\R^n$. So some readers may want to use $\mathcal L$ all the time. The theory is identical to what we have done until now, since we only worked with solutions to $Lu=0$ before the Green functions, and $Lu = 0$ is obviously equivalent to $\mathcal L u = 0$.

We expect from the Green representation of solutions that the function $u=u_f$ constructed in \eqref{GreenFS1} lies in $W_0$ and is a weak solution to $Lu=fw$ in the sense that 
\begin{equation} \label{GreenFS2}
\int_\Omega A\nabla u \cdot \nabla \varphi  =  \int_\Omega f \varphi w = \int_\Omega f \varphi \, dm \ \ \text{ for every }
\varphi \in C^\infty_0(\R^n). 
\end{equation}
Our assumptions (\hyperref[H1]{H1})--(\hyperref[H6]{H6}) are enough to have \eqref{GreenFS2} and the fact that $\Tr u \equiv 0$ (for the former, we still need to be careful about our weird definition of the gradient, and for the later, just use Lemma \ref{GreenDi}). However, we did not succeed to prove that $u \in W$. That is why our next results will be restricted 
to the case where the weight $w$ is nice enough, that is when (\hyperref[H7]{H6'}) is satisfied instead of (\hyperref[H6]{H6})

\begin{lemma} \label{GreenFS} 
Assume that $(\Omega,m,\mu)$ satisfies (\hyperref[H1]{H1})--(\hyperref[H5]{H5}) and (\hyperref[H7]{H6'}). Let $g: \Omega \times \Omega \to \R \cup \{+\infty\}$ 
be the non-negative Green function constructed in Lemma~\ref{GreenEx}. Take $f \in C^\infty_0(\Omega)$ and construct $u(x)$ for $x\in \Omega$ as
\begin{equation} \label{GreenFS3}
u(x) = \int_\Omega g(x,y) f(y) dm(y).
\end{equation}
Then $u$ belongs to $W_0$ and is the solution to $Lu=fw$ (given by Lemma \ref{lLM}) in the sense that
\begin{equation} \label{GreenFS4}
\int_\Omega A\nabla u \cdot \nabla \varphi  =  \int_\Omega f \varphi w = \int_\Omega f \varphi \, dm \ \ \text{ for every }
\varphi \in W_0. 
\end{equation}
\end{lemma}

\begin{remark} In \eqref{GreenFS3} and \eqref{GreenFS4}, we can replace $dm$ by the classical $n$-dimensional Lebesgue measure $dx$.
\end{remark}

\bp
The proof is the same as the one of Lemma 10.7 in \cite{DFMprelim}. It relies on the fact that the solutions to $Lu=fw$ are continuous inside $\Omega$, because as long as we consider inside estimates, (\hyperref[H7]{H6'}) implies that the classical unweighted elliptic theory can be applied. See Theorem 8.22 in \cite{GT} for the theorem in the classical case.

We also need the fact that the approximations $g^\rho(.,y)$ converges in $L^1_{loc}(\Omega)$ to $g(.,y)$. Under (\hyperref[H7]{H6'}), this result is a consequence of the weak $L^q$ convergence of the gradients and the $L^1$-Poincar\'e inequality for inside balls, the latter is true because inside estimates works here exactly like the classical unweighted case. With only  (\hyperref[H6]{H6}), we can only use $L^2$ - or $L^{2-\epsilon}$ - Poincar\'e inequalities.
\ep

\begin{lemma} \label{GreenUn}
Assume that $(\Omega,m,\mu)$ satisfies (\hyperref[H1]{H1})--(\hyperref[H5]{H5}) and (\hyperref[H7]{H6'}).
There exists a unique function $g: \Omega \times \Omega \mapsto \R \cup \{ + \infty\}$ such that 
$g(x,.)$ is continuous on $\Omega \setminus \{x\}$
and locally integrable in $\Omega$ for every
$x\in \Omega$, and such that for every
$f \in C^\infty_0(\Omega)$
the function $u$ given by
 \begin{equation} \label{GreenUn1}
u(x) : = \int_\Omega g(x,y) f(y) dm(y)
 \end{equation}
 belongs to $W_0$ and is a solution of $Lu = f$ in the sense that
 \begin{equation} \label{GreenUn2}
\int_\Omega A\nabla u \cdot \nabla \varphi  = \int_\Omega \A\nabla u \cdot \nabla \varphi \, dm 
=  \int_\Omega f \varphi dm \qquad \text{for every }
\varphi \in W_0. 
\end{equation}
\end{lemma}

\bp
See the proof of Lemma 10.8 in \cite{DFMprelim}. In short, 
the existence is due to Theorem~\ref{GreenEx}, 
and Lemmas \ref{GreenSym} and \ref{GreenFS},  
while the uniqueness of $g$ comes from the uniqueness of $u\in W_0$ satisfying $Lu=fw$, and the latter is due to Lemma \ref{lLM}.
\ep

\section{Comparison principle}

\label{SCP}

First, let us state the non-degeneracy of the harmonic measure.

\begin{lemma} \label{ltcp4}
Let $\alpha >1$, $B:=B(x_0,r)$ be a ball. Take $X_0 \in \Omega$ be any corkscrew point associated to $x_0$ and $r$ given by the assumption (\hyperref[H1]{H1}).
Then 
\begin{equation} \label{tcp33a}
\omega^X(B \cap \Gamma) \geq C_\alpha^{-1} \ \ \ \text{ for } 
X \in \frac1\alpha B
\end{equation}
and
\begin{equation} \label{tcp33}
\omega^X(B \cap \Gamma) \geq C^{-1}_\alpha \ \ \ \text{ for }  
X \in B(X_0,\delta(X_0)/\alpha),
\end{equation}
\begin{equation} \label{tcp33b}
\omega^X(\Gamma \setminus B) + \frac{m(B\cap \Omega)}{r^2} g(X,X_0) \geq C_{\alpha}^{-1} \ \ \ \text{ for } 
X \in \Omega \sm \alpha B,
\end{equation}
and
\begin{equation} \label{tcp33c}
\omega^X(\Gamma \sm B) + \frac{m(B\cap \Omega)}{r^2} g(X,X_0) \geq C^{-1}_{\alpha} \ \ \ \text{ for }  
X \in B(X_0,\delta(X_0)/\alpha),
\end{equation}
where in the four estimates, $C_{\alpha}$ depends on $C_1$ to $C_6$, $C_A$, $n$, and $\alpha$.
\end{lemma}

\begin{remark} \label{needofG}
The estimates \eqref{tcp33a}--\eqref{tcp33} are classical results about the non-degeneracy of the harmonic measure. However, the reader can be at first surprised by the appearance of the Green functions in \eqref{tcp33b}--\eqref{tcp33c}. 
The terms that involves the Green functions are yet necessary. Indeed, none of our assumptions stops the boundary $\Gamma$ to be a bounded and $\Omega$ to be still infinite. Simply take for instance $\Omega = \R^n \setminus \{0\}$ and $\Gamma = \{0\}$ with appropriate measure $\mu$ and $m$.
Under those conditions, we can actually have $\Gamma \setminus B = \emptyset$, which leads to $\omega^X(\Gamma\setminus B) = 0$ for all $X \in \Omega$.

We claim -  without proof but we pretend that there are real difficulties to it - that the estimates 
\begin{equation} \label{tcp33bb}
\omega^X(\Gamma \setminus B) \geq C_{\alpha}^{-1} \ \ \ \text{ for } 
X \in \Omega \sm \alpha B
\end{equation}
holds whenever $\Omega$ is bounded (since the $\Omega \setminus \alpha B$ would be empty when $\Gamma \setminus B = \emptyset$) or when we can find a point in $\Gamma$ close to $\alpha B$ yet outside of $\alpha B$, i.e. whenever $[100B \setminus \alpha B] \cap \Gamma \neq \emptyset $.

At last, the estimate \eqref{tcp33c} is given to make it similar to \eqref{tcp33b}. The harmonic measure is actually unnecessary in \eqref{tcp33c}.
\end{remark}

\bp 
The proof of \eqref{tcp33a}--\eqref{tcp33} is the same as the one of \cite[Lemma 11.10]{DFMprelim}, and relies on the H\"older continuity at the boundary (Lemma \ref{HolderB}), the existence of Harnack chains (Proposition \ref{propHarnack}), and the Harnack inequality (Lemma \ref{HarnackI}).

Rapidly, the H\"older inequality at the boundary will imply that $\omega^X(B\cap \Gamma)$ is bigger than $1/2$ for any points ``close'' to $\Gamma \setminus \frac1\alpha B$. Then we use Harnack chains of balls to link any point in $\frac1\alpha B$ to one of the previous points, and the Harnack inequality repeatedly on the balls of the Harnack chain.

\medskip

Let us make the proof of \eqref{tcp33b}--\eqref{tcp33c}.
First, let us prove \eqref{tcp33c}. Thanks to \eqref{GreenE2}, we have
\begin{equation} \label{tcp33e1}
g(X,X_0) \geq C^{-1} \int_{\delta(X_0)/2}^{\delta(X_0)} \frac{s^2}{m(B(X_0,s))} \frac{ds}{s} \qquad \text{ for all } X \in B(X_0,\delta(X_0)/2).
\end{equation}
Using the doubling property, since $\delta(X_0) \approx r$, we have $m(B(X_0,s))/s^2 \approx m(B \cap \Omega)/r^2$ for all $s\in (\delta(X_0)/2,\delta(X_0))$. The estimate \eqref{tcp33e1} becomes 
\begin{equation} \label{tcp33e2}
\frac{m(B\cap \Omega)}{r^2} g(X,X_0) \geq C^{-1} \qquad \text{ for all } X \in B(X_0,\delta(X_0)/2).
\end{equation}
We let the reader check that $\Omega \sm \{X_0\}$, obtained from $\Omega$ by removing a single point, will still satisfy (\hyperref[H1]{H1})--(\hyperref[H2]{H2}), maybe with some constant $C'_1,C'_2$ smaller than $C_1,C_2$. 
Indeed, if $X_0$ is close to a Corkscrew point for $\Omega$ associated to $(x,r)$, then the Corkscrew point for $\Omega$ associated to $(x,C_1^{-1}r)$ will be far from $X_0$ and so Corkscrew point for $\Omega \setminus \{X_0\}$ with a constant $C'_1 = (C_1)^2$. 
The Harnack chains in $\Omega \setminus \{X_0\}$ are the same as in $\Omega$, except if they got close to $X_0$. In this case, we consider smaller balls, and we avoid $X_0$ by taking balls in $B(X_0,\delta(X_0)/2) \setminus B(X_0,\delta(X_0)/4) \subset \Omega$.
As a consequence, we can link any point from 
\[[\{X \in \Omega, \, \dist(X,\Gamma) > \eta r\} \cap 4B] \setminus B(X_0,\delta(X_0)/4)\] 
to a point in $B(X_0,\delta(X_0)/2) \setminus B(X_0,\delta(X_0)/4)$ by a Harnack chain of ball with uniformly finite length (the length of the chain is bounded by a constant that depends only on $\eta$ and $n$). For the sequel, we write $\Omega_\eta$ for $\{X \in \Omega, \, \dist(X,\Gamma) > \eta r\}$. We use the fact that $g(.,X_0)$ is a solution to $Lu=0$ on $\Omega \setminus \{X_0\}$ and the Harnack inequality (Lemma \ref{HarnackI}) on each balls of those Harnack chain to improve \eqref{tcp33e2} into
\begin{equation} \label{tcp33e3}
\frac{m(B\cap \Omega)}{r^2} g(X,X_0) \geq C_\eta^{-1} \qquad \text{ for all } X \in \Omega_\eta \cap 4B.
\end{equation}
In particular, if $\eta = (\alpha-1)/C_1$, we get \eqref{tcp33c} without the harmonic measure, so we get \eqref{tcp33c} since $\omega^X$ is non-negative.

The proof of \eqref{tcp33b} needs additional computations.
We write $h$ for the smooth function in $C^\infty_0((\alpha+1)B/2))$ satisfying $0\leq h \leq 1$ and $h \equiv 1$ on $B$. We set $u_h \in W$ for the solution to $Lu_h = 0$ with $\Tr u_h = 1- \Tr h$. By positivity of the harmonic measure, 
\begin{equation} \label{tcp33g}
 \omega^X(\Gamma\sm B) \geq u_h(X) \geq \omega^X(\Gamma\sm (\alpha+1)B/2) \geq 0 \qquad \text{ for } X \in \Omega.
\end{equation}
We prefer $u_h$ to $\omega^X(\Gamma \sm B)$ because $u_h$ is in $W$, which makes him suitable for the use of Lemma \ref{lMPg} (our maximum principle). The first estimate that we state comes from \eqref{tcp33e3} without difficulty:
 \begin{equation} \label{tcp33e4}
u_h(X) + \frac{m(B\cap \Omega)}{r^2} g(X,X_0) \geq C_\eta^{-1} \qquad \text{ for all } X \in \Omega_\eta \cap [4B \setminus \alpha B].
\end{equation}
We want the estimate on the larger set $\Omega \cap [4B \setminus \alpha B]$, so we need to prove that \eqref{tcp33e4} is also true when $X$ is close to $\Gamma$. Let $\eta >0$ be small and to be fixed. Let $X \in [\Omega \sm \Omega_\eta] \cap [4B \setminus \alpha B]$. Take $x\in \Gamma$ so that $|X-x| \leq \eta r$, which is possible since $X \in \Omega \sm \Omega_\eta$. Due to the fact that $X$ is also in $2B \setminus \alpha B$, it forces $x$ to be in $\Gamma \setminus (\alpha - \eta)B$. We chose then $\eta = (\alpha-1)/8$, which makes $x\in \Gamma \setminus \frac{7\alpha + 1}8B$. Consequently, for $X \in [\Omega \sm \Omega_\eta] \cap [4B \setminus \alpha B]$,
\[u_h(X) \geq \omega^X(\Gamma\sm (\alpha+1)B/2) \geq \omega^X(B(x,3\eta r) \geq  C^{-1} \]
by \eqref{tcp33g}, the construction of $\eta$ and $x$, and \eqref{tcp33a}. The combination of the last estimate with \eqref{tcp33e4} entails
 \begin{equation} \label{tcp33e5}
u_h(X) + \frac{m(B\cap \Omega)}{r^2} g(X,X_0) \geq C_\alpha^{-1} \qquad \text{ for all } X \in \Omega \cap [4B \setminus \alpha B]
\end{equation}
since Green functions are non-negative. We finish the proof with the maximum principle given by Lemma \ref{lMPg}, which will become our favorite tool for the section. Indeed, we define $v$ as
\[v:= u_h(X) + \frac{m(B\cap \Omega)}{r^2} g(X,X_0) - C_\alpha^{-1}\]
where $C_\alpha^{-1}$ is the constant in the right-hand side of \eqref{tcp33e5}, and we aim to apply Lemma \eqref{lMPg} for the solution $v$ with the sets $E = \R^n \setminus \alpha B$ and $F = 4B \setminus \alpha B$. Recall that the term $u_h(X)$ lies in $W$. Together with \eqref{GreenE5}, we deduce that assumption (i) of Lemma \ref{lMPg} is true. The other assumptions required by the lemma are given by \eqref{tcp33e5} and the fact that $\Tr v = 1-C_\alpha^{-1} > 0$ on $\Gamma \setminus \alpha B$. We deduce that $v \geq 0$ on $E$, which is exactly the desired estimate \eqref{tcp33b}. The lemma follows.
\ep

If $B$ is a ball centered on the boundary $\Gamma$, we bound the values in $B \cap \Omega$ of a solution $u$ (to $Lu=0$ in $KB \cap \Omega$) by the value of $u$ at a Corkscrew point associated to the ball $B$.

\begin{lemma} \label{ltcp2}
There exists $K:= K(C_1,C_2,n)$ such that the following holds.

Let $B = B(x_0,r)$ be a ball centered on $\Gamma$ and let $X_0$ be a Corkscrew point associated to $x_0$ and $r$ given by (\hyperref[H1]{H1}). Let $u \in W_r(KB\cap \overline\Omega)$ be a non-negative,
non identically zero,
solution of $Lu=0$ in $KB\cap \Omega$,
such that 
$\Tr u \equiv 0$ on $KB\cap \Gamma$. Then
\begin{equation} \label{tcp6}
u(X) \leq C u(X_0)   \ \ \ \ \text{ for } X \in B \cap \Omega,
\end{equation}
where $C>0$ depends on $C_1$ to $C_6$, $C_A$, and $n$.
\end{lemma}

\bp
We get inspiration from the proof of \cite[Lemma~4.4]{KJ} (see also \cite[Lemma 11.8]{DFMprelim}). Lemma 11.8 in \cite{DFMprelim} deals with balls centered at the boundary, and $K=2$. However, in \cite{DFMprelim}, the connectedness is not a issue, while we need to be careful here. Indeed, taking the universal constant $K = 2$ in Lemma \ref{ltcp2} is not possible, since nothing garantees that we can link 2 points in $B \cap \Omega$ by a path that stays in $2B$. 

We solve this problem by taking the tent sets constructed in Section \ref{Scones}, which can be seen as connected substitute of the sets $B \cap \Omega$. Then we use the property \eqref{propTQ} to conclude.

\medskip

First, let us recall the following fact. Let $x\in \Gamma$ and $s>0$ such that $\Tr u \equiv 0$ on $B(x,s) \cap \Gamma$. Then the H\"older continuity of solutions given by Lemma~\ref{Holder1} proves the existence of $\epsilon>0$ (that depends on $C_1$ to $C_6$, $C_A$, and $n$) such that
\begin{equation} \label{tcp7}
\sup_{B(x,\epsilon s)} u \leq \frac12 \sup_{B(x,s)} u.
\end{equation}
Without loss of generality, we can choose $\epsilon < 1/1000C_1 < \frac12$.

A rough idea of the proof of \eqref{tcp6} is that $u(X)$ should not be near the maximum of $u$ when $X$ lies close to $B\cap \Gamma$, because of \eqref{tcp7}. Then we are left with
points $x$ that lie far from the boundary, and we can use the Harnack inequality to control $u(x)$.
The difficulty is that when $X \in B\cap \Gamma$ lies close to $\Gamma$, $u(x)$ can be bounded by values 
of $u$ \underline{inside the domain}, and not by values of $u$ near $\Gamma$ but from the exterior 
of $B$.
We will prove this latter fact by contradiction: we show that if $\sup_{B} u$ exceeds a certain bound, 
then we can construct a sequence of points $X_k \in \frac K2B$, where $K$ is large enough, such that $\delta(X_k) \to 0$ and 
$u(X_k) \to +\infty$, and hence we 
contradict the H\"older continuity of solutions at the boundary.

\medskip

As said in the beginning of the proof, the quantities $\lambda B \cap \Omega$ lack connectedness, and it will be more convenient to work with a tent set $T_{2Q^*}$, which has all the desired connectedness by Lemma  \ref{lemchain}. The cube $Q^*$ and the constant $K$ are defined as follow. Let $k \in \mathbb Z$ be such that $2^{-k-2} \leq r \leq 2^{-k-1}$, and we write $Q$ for the unique cube in $\D_k$ containing $x$. Notice that $2Q \supset B \cap \Gamma$, but $T_{2Q}$ is not necessarily bigger than $B$, and so we take $Q^*$ the first ancestor of $Q$ such that $\dist(T_{2Q},T_{2Q^*}^c) >\ell(Q)$. Check that the difference of generations between $Q$ and $Q^*$ is uniformly finite, and so combined with \eqref{propTQ}, we obtain that we can find $K$ that depends only on $n$, $c_1$, and $C_2$ such that
\begin{equation} \label{tcp7a}
2B \cap \Omega \subset T_{2Q^*} \subset \frac12 B^*:= \frac K2B.
\end{equation}

\medskip

We can link any couple of points in $T_{2Q^*}$ by a chain of balls $B_i$ that satisfies $2B_i \subset 2B^*\cap \Omega$. The proof of this fact is similar to the proof of fact that $T_{2Q^*}$ satisfies the chain condition $C(\kappa,M)$ for some $\kappa,M$ (see Lemma \ref{lemchain}) and thus will be omitted.
Therefore, the fact that $u(X) > 0$ somewhere, that $T_{2Q^*}$ is connected, and the Harnack inequality (Lemma~\ref{HarnackI}), maybe applied a few times,  yield $u(X_0) > 0$. We can rescale $u$ and assume that $u(X_0) = 1$. 

We claim 
that there exists
$M>0$ such that for any integer $N \geq 1$ and
$Y \in T_{2Q^*}$,
\begin{equation} \label{tcp8}
\delta(Y) \geq \epsilon^N r \Longrightarrow u(Y) \leq M^N  ,
\end{equation}
where $\epsilon$ 
comes from 
\eqref{tcp7} and the constant $M$ depends only upon $n$, $C_1$ to $C_6$, and $C_A$.
We will prove the claim by induction.
The base case is given by the following. We want to show the existence of $M_1 \geq 1$ such that
\begin{equation} \label{tcp9}
u(Y) \leq M_1 \ \ \ \ \ \text{ for every }
Y \in T_{2Q^*} \text{ such that } \delta(Y) \geq \epsilon^2 r.
\end{equation}
Indeed, if $Y \in T_{2Q^*}$ satisfies $\delta(Y) \geq \epsilon^2 r$, Proposition \ref{propHarnack} implies that we can link $Y$ to $X_0$ by a chain of balls that stay away from the boundary and with length uniformly bounded by $C(\epsilon)$. We can construct the chain such that it stays also far from the boundary of $B^*$\footnote{This fact is true because by construction of $T_{2Q^*}$, one can see that the center of the balls constituting the chain can be taken in $T_{2Q^*}$, and \eqref{tcp7a} let us a bit of freedom, but if may be easier for the reader to think that this statement would be also true by taking a larger $K$}.
Together with the Harnack inequality (Lemma~\ref{HarnackI}), we obtain \eqref{tcp9}, and hence \eqref{tcp8} for $N=1,2$ as long as $M$ is chosen bigger than $M_1$.

Now, let $Y\in T_{2Q^*}$ such that $\delta(Y) \leq \epsilon^2 r$. By construction of $T_{2Q^*}$, $Y$ belongs to some $\gamma^*_{Q^*}(z)$ for some $z\in 2Q^*$. We take $Z$ a Corkscrew point associated to $z$ and $C_1\delta(Y)/\epsilon$. Since $\epsilon < C_1^{-1}$, $Z \in B(z,r)$, and so $Z$ stays in $T_{2Q^*}$. In addition, by construction, $\delta(Z) \geq \delta(Y)/\epsilon$ and $|Z-Y| \leq \delta(Y)/\epsilon^2$; these two estimates can be combined to Proposition \ref{propHarnack} (existence of Harnack chains, as before the chain can stay far from $\dr B^*$) and Lemma~\ref{HarnackI} (Harnack inequality) to get we  the existence of $M_2 \geq 1$  such that $u(Y) \leq M_2 u(Z)$. 
So we just proved that
\begin{equation} \label{tcp10} \begin{array}{l}
\text{for any $Y \in T_{Q^*}$ such that 
$\delta(Y) \leq \epsilon^{-2} r$, } \\
\qquad \text{there exists $Z\in T_{2Q^*}$ such that $\delta(Z) \geq  \delta(Y)/\epsilon$ and $u(Y) \leq M_2 u(Z)$.}
\end{array} \end{equation}

We turn
to the main
induction step. Set $M = \max\{M_1, M_2\} \geq 1$ and let $N \geq 2$ be given. 
Assume, by induction hypothesis, that for any $Z \in T_{2Q^*}$ satisfying $\delta(Z) \geq \epsilon^{N} \ell(Q)$, we have $u(Z) \leq M^N$. 
Let 
$Y \in T_{2Q^*}$ be such that 
$\delta(Y) \geq \epsilon^{N+1} \ell(Q)$ . 
The assertion \eqref{tcp10} yields the existence of $Z\in T_{2Q^*}$ such that $\delta(Z) \geq \delta(Y)/\epsilon \geq \epsilon^N \ell(Q)$ and $u(Y) \leq M_2 u(Z) \leq Mu(Z)$. 
By the induction hypothesis, $u(Y) \leq M^{N+1}$. 
This completes our induction step, and the proof of \eqref{tcp8} for every $N \geq 1$.

\medskip

Choose an integer $i$ such that $2^i \geq M$, where $M$ is the 
constant of \eqref{tcp8} that we just found, and then set $M' = M^{i+3}$.
We want to prove by contradiction that
\begin{equation} \label{tcp11}
u(X) \leq M' u(X_0) = M' \ \ \ \text{ for every }
X \in B(x_0,r).
\end{equation}
So we assume that  
\begin{equation} \label{tcp12}
\text{there exists $X_1 \in B(x_0,r)$ such that $u(X_1) > M'$}
\end{equation}
and we want to prove by induction that for every 
integer $k \geq 1$,
\begin{equation} \label{tcp13}
\text{there exists $X_k \in T_{2Q^*}$ such that $u(X_k) > M^{i+2+k}$ and $\dist(X_k,B) \leq (1-2^{1-k})r$.}
\end{equation}
The base step of the induction is given by \eqref{tcp12} and we want to do  
the induction step. Let $k\geq 1$ be given and assume that 
\eqref{tcp13} holds.
From the contraposition of \eqref{tcp8}, we deduce that $\delta(X_k) < \epsilon^{i+2+k}r$. 
Choose 
$x_k \in \Gamma$ such that $|X_k-x_k| = \delta(X_k) < \epsilon^{i+2+k}r$.
By the induction hypothesis,
\begin{equation} \label{tcp14}
\dist(x_k,B) \leq |x_k-X_k| + \dist(X_k,B) \leq (1-2^{1-k})r + \epsilon^{i+2+k} r 
\end{equation}
and, 
since $\epsilon \leq \frac12$,
\begin{equation} \label{tcp15}
|x_k - x_0| \leq (1-2^{1-k} - 2^{-2-k})r. 
\end{equation}
Now, due to \eqref{tcp7}, we can find $X_{k+1} \in B(x_k,\epsilon^{2+k} r)$ such that 
\begin{equation} \label{tcp16}
u(X_{k+1}) \geq 2^i \sup_{X\in B(x_k,\epsilon^{i+2+k} r)} u(X) \geq 2^i u(X_k) \geq M^{i+2+(k+1)}.
\end{equation}
The induction step will be complete if we can prove that $\dist(X_{k+1},T_{2Q}) \leq (1-2^{-k})r$. Indeed, 
\begin{equation} \label{tcp17} \begin{split}
\dist(X_{k+1},B) & \leq |X_{k+1} - x_k| +\dist(x_{k},B) \leq \epsilon^{2+k} r + (1-2^{1-k} - 2^{-2-k})r \\
& \leq (1-2^{-k})r
\end{split}\end{equation}
by \eqref{tcp15} and because $\epsilon \leq \frac12$.

\medskip
Let us sum up. We assumed the existence of  $X_1 \in B$ such that $u(X_1)>M'$ and we end up 
with \eqref{tcp13}, that is a sequence $X_k$ of values in $2B$ such that $u(X_k)$ 
increases to $+\infty$. 
Up to a subsequence, we can thus find a point in $\overline{2B} \subset B^*$ where $u$ is not continuous, 
which contradicts Lemma~\ref{HolderB}. Hence 
$u(X) \leq M' = M' u(X_0)$ for
$X \in B$. Lemma~\ref{ltcp2} follows.
\ep

We can now compare the harmonic measure with the Green function, that can be seen as a weak version of the comparison principle.

\begin{lemma} \label{ltcp3}
Let $B:=B(x_0,r)$ be a ball centered on $\Gamma$. Let $X_0$ is a corkscrew point associated to $x_0$ and $r$. Then one has
\begin{equation} \label{tcp18}
C^{-1} \frac{m(B \cap \Omega)}{r^2} g(X,X_0) \leq \omega^X(B \cap \Gamma) \leq C \frac{m(B \cap \Omega)}{r^2} g(X,X_0) 
\ \ \ \text{ for }
X\in \Omega \setminus 2B,
\end{equation}
and
\begin{equation} \label{tcp18a}
\omega^X(\Gamma \setminus \frac54B) 
\leq C \frac{m(B \cap \Omega)}{r^2} g(X,X_0) 
\ \ \ \text{ for }
X\in [B \cap \Omega] \setminus B(X_0,\delta(X_0)/4),
\end{equation}
where $C>0$ depends only upon $n$, $C_1$ to $C_6$, and $C_A$.
\end{lemma}

\begin{remark}
The bound \eqref{tcp18a} may look a bit artificial. 
There is nothing deep about the constant $\frac54$ in the left-hand side.
We could have used $2$ instead, and obtain a statement which looks a little weaker
but is actually equivalent (we leave the proof of this to the reader); 
we simply reproduced $\frac54$ in the form given by our
general comparison principle (Theorem \ref{lCP3}).

Observe also that we do not necessarily have  
the lower bound in \eqref{tcp18a}. 
See Remark \ref{needofG}. 
\end{remark}

\bp This lemma is the analogue of \cite[Lemmas 11.9 and 11.11]{DFMprelim}.

\medskip

First, we quickly prove the first inequality in \eqref{tcp18}. The 
upper bound \eqref{GreenE7} for the Green function,
together with (\hyperref[H4]{H4}), implies that
\begin{equation} \label{tcp18b}
0 \leq \frac{m(B \cap \Omega)}{r^2} g(X,X_0) \leq C \qquad \text{ for } X \in B(X_0,\delta(X_0)/2) \sm B(X_0,\delta(X_0)/4).
\end{equation}
As in
the proof of Lemma \ref{ltcp4}, we take 
$h \in C^\infty_0(B)$ such that
$0\leq h \leq 1$ and $h \equiv 1$ on $\frac12B$. We set $u_h \in W$ for the solution to $Lu_h = 0$ with $\Tr u_h = \Tr h$. By the
positivity of the harmonic measure, 
\begin{equation} \label{tcp18c}
 \omega^X(\Gamma \cap \frac12 B) \leq u_h(X) \leq \omega^X(\Gamma \cap B) \qquad \text{ for } X \in \Omega.
\end{equation}
We combine \eqref{tcp18c} with the non-degeneracy of the harmonic measure \eqref{tcp33} to get that 
\begin{equation} \label{tcp18d}
u_h \geq C^{-1} \qquad \text{ for } X \in B(X_0,\delta(X_0)/2).
\end{equation}
The estimates \eqref{tcp18b} and \eqref{tcp18d} easily infer the existence of a constant $\kappa$ such that 
\[v(X):= \kappa u_h(X) - \frac{m(B\cap \Omega)}{r^2} g(X,X_0) \geq 0 \qquad \text { for } X \in B(X_0,\delta(X_0)/2) \sm B(X_0,\delta(X_0)/4).\]
The assumptions of Lemma \ref{lMPg} for the function $v$ and the sets 
$E = \R^n \setminus B(X_0,\delta(X_0)/4)$ and 
$F = B(X_0,\delta(X_0)/2) \sm B(X_0,\delta(X_0)/4)$ are satisfied, which implies that 
$v\geq 0$ on $E$, i.e.,
\begin{equation} \label{tcp18z}
\frac{m(B \cap \Omega)}{r^2} g(X,X_0) \leq \kappa u_h \leq \omega^X(\Gamma \cap B) \qquad \text{ for } X \in \Omega \setminus B(X_0,\delta(X_0)/4).
\end{equation}
This is stronger than 
the first inequality in \eqref{tcp18}.

\medskip
 
We shall also use the following result on Green functions: 
for $\phi \in C^\infty(\R^n) \cap W$ and $X \notin \supp \, \phi$,
\begin{equation} \label{tcp20}
u_\phi(X) = - \int_\Omega A \nabla \phi(Y) \cdot \nabla_y g(X,Y) dY,
\end{equation}
where $u_\phi \in W$ is the solution ro $Lu = 0$, with the Dirichlet condition $\Tr u_\phi = \Tr \phi$ on $\Gamma$, given by Lemma \ref{lLM}. 
The identity \eqref{tcp20} is the same as (11.70) in \cite{DFMprelim}, 
and its proof - which only relies on the properties on the Green functions given in Section \ref{SGreen} - is the same as in \cite{DFMprelim}.

\medskip

We turn to the proof of the upper bound in \eqref{tcp18}, that is,  
\begin{equation} \label{tcp19}
\omega^X(B(x_0,r) \cap \Gamma) \leq C \, \frac{m(B \cap \Omega)}{r^2}\, 
g(X,X_0) \ \ \ \text{ for }
X\in \Omega \setminus 2B.
\end{equation}
For the rest of the proof, $K$ is the constant in Lemma \ref{ltcp2}. Let $X\in \Omega \setminus 2B$ be given, and choose $\phi \in C^\infty_0(\R^n)$ 
such that $0\leq \phi \leq 1$, $\phi \equiv 1$ on $B\cap \Gamma$, 
$\supp \, \phi  \subset E_B:= \{Y\in \Omega, \, \dist(Y,B\cap \Gamma) \geq (100K)^{-1}r \}$, and $|\nabla \phi| \leq 200K/r$. 
We get that 
\begin{equation} \label{tcp30}
u_\phi(X) \leq \frac{C}r \int_{E_B}  
|\nabla_y g(X,Y)| dm(Y)
\end{equation}
by \eqref{tcp20} and \eqref{defBdd}, and since $\omega^X(B \cap \Gamma) \leq u_\phi(X)$
by the positivity of the harmonic measure,
\begin{equation} \label{tcp30bis} \begin{split}
\omega^X(B\cap \Gamma) & \leq \frac{C}r \int_{E_B} |\nabla_y g(X,Y)| dm(Y).
\end{split} \end{equation}
We cover $E_B$ 
by a finitely overlapping collection of balls $(B_i)_{i\in \mathcal I}$ centered on $B\cap \Gamma$ and of radius $(10K)^{-1}r$. Then  
\begin{equation} \label{tcp30ter} \begin{split}
\omega^X(B\cap \Gamma) & \leq \frac{C}r \sum_{i \in \mathcal I} \int_{B_i \cap \Omega} |\nabla_y g(X,Y)| dm(Y) \\
& \lesssim \sum_{i\in \mathcal I} \frac{m(B_i \cap \Omega)^{1/2}}r  \left( \int_{B_i \cap \Omega} |\nabla_y g(X,Y)|^2 dm(Y)\right)^\frac12 \\
& \lesssim \sum_{i\in \mathcal I} \frac{m(B_i \cap \Omega)^{1/2}}{r^2}  \left( \int_{2B_i \cap \Omega} |g(X,Y)|^2 dm(Y)\right)^\frac12,
\end{split} \end{equation}
where we use successively the Cauchy-Schwarz 
inequality and Cacciopoli's inequality at the boundary (Lemma \ref{CaccioB}); 
the use of Cacciopoli's inequality is indeed allowed 
because $Y \to g(X,Y)$ is a solution of $L_T u := -\div A^T \nabla u$ in $2B_i \cap \Omega$ by Lemmas \ref{GreenSym} and \ref{GreenEx} (iii).
Observe that $Y \to g(X,Y)$ is more generally a solution of $L_T u := -\div A^T \nabla u$ 
in each set  
$2KB_i \cap \Omega$, because the radius of $2KB_i$ is less than $r/2$ and hence $2KB_i \subset 2B \not\ni X$. As a consequence, Lemma \ref{ltcp2} yields that 
\begin{equation} \label{tcp32} 
\omega^X(B\cap \Gamma) 
\lesssim \sum_{i\in \mathcal I} \frac{m(B_i)}{r^2} \, g(X,X_i), 
 \end{equation}
where $X_i$ is a corkscrew point associated to the ball $B_i$. Hence 
\begin{equation} \label{tcp32z} 
\omega^X(B\cap \Gamma) \lesssim \frac{m(B)}{r^2} g(X,X_0) 
 \end{equation}
because of the finite overlapping of the $(B_i)_i$, the Harnack inequality, and the fact that we can easily find Harnack chains of balls that link $X_i$ to $X_0$ and that avoids $X$. The bounds \eqref{tcp19} and then \eqref{tcp18} follow.

\medskip

The proof of \eqref{tcp18a} can be treated in a similar manner, and we refer to  
\cite[Lemma~11.11]{DFMprelim} 
for additional ideas on the proof.
\ep

\begin{lemma}[Doubling volume property for the harmonic measure] \label{ldphm}
Let $\alpha >1$, and take a ball $B:=B(x_0,r)$ in $\R^n$ centered on $\Gamma$. One has
\begin{equation} \label{dphm1} 
\omega^X(2B\cap \Gamma) \leq C_\alpha \omega^X (B \cap \Gamma) 
 \ \ \ \text{ for }
X \in \Omega \setminus 2\alpha B,
 \end{equation}
where $C_\alpha>0$ depends only upon $n$, $C_1$ to $C_6$, $C_A$, and $\alpha$.
\end{lemma}

\bp The proof is the same as the one of \cite[Lemma 11.102]{DFMprelim}. Here are some ideas.

\medskip

When $\alpha = 2$, we use Lemma \ref{ltcp3}, the doubling property (\hyperref[H4]{H4}), the Harnack inequality, and the existence of Harnack chains of balls to write
\begin{equation} \label{dphm2} 
\omega^X(2B\cap \Gamma) \lesssim \frac{m(2B)}{(2r)^2} g(X,X_2) \lesssim \frac{m(B)}{r^2} g(X,X_1) \lesssim \omega^X(B\cap \Gamma),
\end{equation}
where $X_1$ and $X_2$ are corkscrew points associated to respectively $(x_0,r)$ and $(x_0,2r)$.

When $1 < \alpha < 2$, we cover $2B \cap \Gamma$ by a collection of finitely overlapping balls $2B_i$ of radius $2r_\alpha:= (\alpha-1)r$ and centered on $B(x_0,2r-\frac32 r_\alpha) \cap \Gamma$. In this case, for any $X \in \Omega \setminus 2\alpha B$
\begin{equation} \label{dphm3} 
\omega^X(2B\cap \Gamma) \leq \sum_i \omega^X(2B_i\cap \Gamma) \lesssim  \sum_i \omega^X(B_i\cap \Gamma) \lesssim \omega^X(B(x_0,2r-\frac12 r_\alpha) \cap \Gamma),
\end{equation}
where the second estimate is due to \eqref{dphm2}. We repeat the argument a finite number of time (depending in $\alpha-1>0$) to get \eqref{dphm1}. The lemma follows.
\ep

\begin{lemma}[Comparison principle for global solutions] \label{lCP}
Let $B:=B(x_0,r)$ be a ball centered on $\Gamma$, and let 
$X_0\in \Omega$ be a corkscrew point associated to $(x_0,r)$.
Let $u,v\in W$ be two non-negative, non identically zero, 
solutions to $Lu = Lv = 0$ in $\Omega$ such that
$\Tr u = \Tr v = 0$ on $\Gamma \setminus B(x_0,r)$.
Then
 \begin{equation} \label{CP1} 
C^{-1} \frac{u(X_0)}{v(X_0)} \leq \frac{u(X)}{v(X)} \leq C \frac{u(X_0)}{v(X_0)}
\ \ \ \text{ for } X \in \Omega \setminus 2B,
 \end{equation}
where $C>0$ depends only on $n$, $C_1$ to $C_6$, and $C_A$. 
\end{lemma}

\begin{remark}
We also have \eqref{CP1} for any $X\in \Omega \setminus B(x_0,\alpha r)$, where $\alpha >1$. In this case, the constant $C$ depends also on $\alpha$. We let the reader check that the proof below can be easily adapted to prove this too.
\end{remark}

\bp The proof is very similar to the one of \cite[Lemma 11.14]{DFMprelim}. 
Let us recall the main steps and show the differences.

\medskip

By symmetry of the roles of $u$ and $v$, we only need to show the upper bound
  \begin{equation} \label{CP2} 
\frac{u(X)}{v(X)} \leq C \frac{u(X_0)}{v(X_0)}  \ \ \ \text{ for }
X \in \Omega \setminus 2B.
 \end{equation}
 Notice also that thanks to the Harnack inequality (Lemma \ref{HarnackI}), $v(X)>0$ on the whole $\Omega$, so we we don't need to be careful when we divide by $v(X)$.
 
 \medskip
 
We introduce some notation for two boundary balls: set $\Gamma_1 := \Gamma \cap B$ and $ \Gamma_2 : = \Gamma \cap  \frac{15}8 B$. 
The proof of the lemma is composed of three steps :
\begin{enumerate}[(a)]
\item we prove
the lower bound
  \begin{equation} \label{CP4} 
v(X) \geq C^{-1} \omega^X(\Gamma_1) v(X_0) \ \ \ \text{ for }
X \in \Omega \setminus 2B ;
 \end{equation}
 \item we prove
 the upper bound
 \begin{equation} \label{CP8} 
u(X) \leq Cu(X_0)\omega^X(\Gamma_2) \ \ \ \text{ for }
X \in \Omega \setminus 2B ; 
\end{equation}
\item we 
conclude by using the fact that the harmonic measure is doubling (Lemma \ref{ldphm}).
\end{enumerate}

The proof of \eqref{CP4} is can be done exactly as in 
\cite[Lemma 11.14]{DFMprelim}. Let us quickly sketch it. 
By the Harnack inequality (Lemma \ref{HarnackI}), for all $X\in B(X_0,\delta(X_0)/2)$, we have $v(X) \gtrsim v(X_0)$. Together with the upper bound \eqref{GreenE7}, we get the existence of a constant $K_1$ such that the function
\[ v_1(X) : =  
K_1 v(X) - \frac{m(B \cap \Omega)}{r^2} v(X_0) g(X,X_0)\]
satisfies all the assumption of the maximum principle (Lemma \ref{lMPg}) with $E = \R^n \sm B(X_0,\delta(X_0)/4)$ and $F = B(X_0,\delta(X_0)/2) \sm B(X_0,\delta(X_0)/4)$. Indeed, since $v$ is non-negative everywhere, it forces $\Tr v_1 = \Tr v \geq 0$. We deduce that $v_1 \geq 0$ on $E$, i.e.
\[ \frac{m(B \cap \Omega)}{r^2} v(X_0) g(X,X_0) \leq K_1 v(X) \qquad \text{ for } X \in \Omega \sm B(X_0,\delta(X_0)/4).\]
The claim \eqref{CP4} is now an immediate consequence of Lemma \ref{ltcp3}.
 
We turn to the proof of \eqref{CP8}. We first check that  
\begin{equation} \label{CP9} 
u(X) \leq Cu(X_0) \ \ \ \text{ for }
X \in \frac{13}8B \setminus \frac{11}8B.
\end{equation}
Let $K$ as in Lemma \ref{ltcp2}. We want to establish \eqref{CP9} in the two sets:
\begin{equation} \label{CP10} 
\Omega_1 : = \Omega \cap \{X\in B(x_0,\frac{13}{8}r) \setminus B(x_0,\frac{11}{8}r), \, \delta(X) < \frac1{8K} r\}
\end{equation}
and
\begin{equation} \label{CP11} 
\Omega_2 : =  \{X\in B(x_0,\frac{13}{8}r) \setminus B(x_0,\frac{11}{8}r), \, \delta(X) \geq \frac1{8K} r\}.
\end{equation}
The proof of \eqref{CP9} on $\Omega_2$ is easy, it is only a consequence of the Harnack inequality (Lemma \ref{HarnackI}) and the existence of Harnack chains.

Then, we prove \eqref{CP9} for
$X\in \Omega_1$. Let thus $X \in \Omega_1$ be given. 
Take $x\in \Gamma$ such that $|X-x| = \delta(X)$;  in particular,  particular, 
$|X-x| \leq \frac r{8}$, and hence  
$x \in \frac74B$. Now let $X_1$ be 
a Corkscrew point for 
$(x,\frac{r}{4K})$. 
Since $u$ is a non-negative solution of $Lu=0$ in $B(x,\frac r{4}) \cap \Omega$ satisfying $\Tr u = 0$  on $B(x,\frac r4 ) \cap \Gamma$, Lemma \ref{ltcp2} gives that $u(Y) \leq C u(X_1)$ for 
$Y \in B(x,\frac r{4K})$ and thus  
$u(X) \leq C u(X_1)$. 
By the existence of Harnack chains (Proposition \ref{propHarnack}) and the Harnack inequality (Lemma \ref{HarnackI}), $u(X_1) \leq Cu(X_0)$. Hence $u(X) \leq u(X_1)$, which completes  
the proof of \eqref{CP9} on $\Omega_1$

\medskip

The end of the proof is as in 
in \cite[Lemma 11.117]{DFMprelim}, but let us recall it. We proved \eqref{CP9} and now we want to get \eqref{CP8}. 
Recall from Lemma \ref{ltcp4} that 
$\omega^X(\frac74B \cap \Gamma) \geq C^{-1}$ for 
$X \in \frac{13}{8}B \cap \Omega$.
Hence, by \eqref{CP9},
\begin{equation} \label{CP12} 
u(X) \leq Cu(X_0) \omega^X(\frac74B\cap \Gamma) \ \ \ \text{ for }
X \in \left[\frac{13}{8}B \setminus \frac{11}{8}B\right] \cap \Omega.
\end{equation}
Let $h\in C^\infty_0(B(x_0,\frac{15}{8}r))$ be such that $0\leq h \leq 1$ and $h \equiv 1$ on $B(x_0,\frac74 r)$. Then let $u_h \in W$ be the solution of $Lu_h = 0$ with the Dirichlet condition $\Tr u_h = \Tr h$. 
By the positivity of the harmonic measure,
\begin{equation} \label{CP13} 
u(X) \leq Cu(X_0) u_h(X) \ \ \ \text{ for } 
X \in \left[\frac{13}{8}B \setminus \frac{11}{8}B\right] \cap \Omega.
\end{equation}
The maximum principle given by Lemma \ref{lMPg} - where we take $E = \R^n \setminus \overline{\frac{11}{8}B}$ and $F = \R^n \setminus \frac{13}{8}B$ - yields
\begin{equation} \label{CP14} 
u(X) \leq Cu(X_0) u_h(X) \ \ \ \text{ for }
X \in \Omega \setminus \frac{13}{8}B
\end{equation}
and hence 
\begin{equation} \label{CP15} 
u(X) \leq Cu(X_0) \omega^X(\Gamma_2) \ \ \ \text{ for }
X \in \Omega \setminus \frac{13}{8}B,
\end{equation}
where we use again the positivity of the harmonic measure. The assertion \eqref{CP8} is now proven.

\medskip

We conclude as follows. Because of 
\eqref{CP4} and \eqref{CP8},
\begin{equation} \label{CP16} 
\frac{u(X)}{v(X)} \leq C\frac{u(X_0)}{v(X_0)} \frac{\omega^X(\Gamma_2)}{\omega^X(\Gamma_1)} 
\ \ \ \text{ for }
X \in \Omega \setminus 2B.
\end{equation}
The bound \eqref{CP2} - and thus the lemma - follows then from the fact that $\omega^X(\Gamma_2) \lesssim \omega^X(\Gamma_1)$, which is given by Lemma \ref{ldphm}.
 \ep
 
 \begin{lemma}[Comparison principle for harmonic measures / Change of poles] \label{lCP2}
Let $B:=B(x_0,r)$ be a ball centered on $\Gamma$ and 
let $X_0$ be 
a corkscrew point associated to $(x_0,r)$. Let $E,F \subset \Gamma \cap B$ be two Borel subsets of $\Gamma$ such that $\omega^{X_0}(E)$ and $\omega^{X_0}(F)$ are positive. Then 
  \begin{equation} \label{CP17} 
C^{-1} \frac{\omega^{X_0}(E)}{\omega^{X_0}(F)} 
\leq  \frac{\omega^{X}(E)}{\omega^{X}(F)} \leq C  \frac{\omega^{X_0}(E)}{\omega^{X_0}(F)}
\ \ \ \text{ for } X \in \Omega \setminus 2B,
 \end{equation}
where $C>0$ depends only on $n$, $C_1$ to $C_6$, and $C_A$. 
In particular, with the choice $F = B \cap \Gamma$, 
  \begin{equation} \label{CP18} 
C^{-1} \omega^{X_0}(E) \leq  \frac{\omega^{X}(E)}{\omega^{X}(B \cap \Gamma)} 
\leq C \omega^{X_0}(E) \ \ \ \text{ for } X \in \Omega \setminus 2B,
 \end{equation}
where $C>0$ depends on the same quantity as for \eqref{CP17}.
\end{lemma}

\bp
This result can be deduced from Lemma \ref{lCP} with the same proof we obtained \cite[Lemma 11.135]{DFMprelim} from \cite[Lemma 11.117]{DFMprelim}. It relies on approximations of 
$X\mapsto \omega^{X}(E)$ and $X\mapsto \omega^{X}(F)$ 
by solutions \underline{in $W$} to $Lu=0$ in $\Omega$.
\ep

\begin{theorem}[Comparison principle for locally defined functions] \label{lCP3}
There exists $K\geq 2$ depending only on $n$, $C_1$, and $C_2$ such that the following holds.

Let $B:=B(x_0,r)$ be a ball centered on $\Gamma$, and let $X_0 \in \Omega$ be a corkscrew point associated to $(x_0,r)$. 
Take $u,v\in W_r(KB \cap \overline\Omega)$ 
to be two non-negative,
not identically zero, solutions to $Lu = Lv = 0$ in $KB \cap \Omega$,
such that $\Tr u = \Tr v = 0$ on $KB \cap \Gamma$. Then
  \begin{equation} \label{CP26} 
C^{-1} \frac{u(X_0)}{v(X_0)} \leq \frac{u(X)}{v(X)} \leq C \frac{u(X_0)}{v(X_0)}
\ \ \ \text{ for } X \in \Omega \cap B,
 \end{equation}
 where $C>0$ depends only on depends only on $n$, $C_1$ to $C_6$, and $C_A$. 
\end{theorem}

\bp Two strategies can be used to prove this theorem:
\begin{enumerate}[(i)]
\item If we mimic the classical proof from the codimension 1 case, we  
need to find a (good enough) domain $D$ such that $B\cap \Omega \subset D \subset KB\cap \Omega$, and we work with the harmonic measure on $\dr D$ of the operator $L$ restricted to $D$. One might think that for instance $D=2B \cap \Omega$ will work out, but this choice of $D$ will not be ``good enough'' if it is not connected (and it can easily happen).

The difficulty is first to construct such a $D$ that satisfies the corkscrew point condition and the Harnack chain condition, but this part could be possibly done by considering the tents sets constructed in Section \ref{Scones}. Yet, even with such good $D$, we still need to build a measure $\mu_D$ on $\dr D$ which is suitable, in particular satisfies (\hyperref[H5]{H5}) for this particular domain. Well, at the present moment, we don't even know if building such $\mu_D$ is possible with our assumptions.

\item The second strategy, that we shall apply, is to follow the ideas used in \cite[Theorem 11.146]{DFMprelim}. In this strategy, we are not allowed to consider a harmonic measure different from the one we defined on $\Gamma$. The main pitfall in the present theory which did not exist in \cite{DFMprelim} is the fact that $\omega^X(\Gamma \sm B)$ can be 
null because $\Gamma \sm B$ is empty, 
and so we do not necessary have the non-degeneracy of the harmonic measure, and we shall use the estimate \eqref{tcp33b} involving the Green function instead.

We may have chosen to restrict our attention to the balls $B$ that do not cover entirely $\Gamma$.
Here we decided to allow more balls, but then when we take $r$ large in our theorem - when $\Omega$ is unbounded and $\Gamma$ is bounded - we need to impose stronger conditions
(and we get a stronger conclusion).
\end{enumerate}

\medskip

\noindent {\bf Step 1:} Construction of a  
function $f_{y_0,s}$.

\smallskip

Let $y_0\in \Gamma$ and $s>0$. We write $Y_0$ for a corkscrew point associated to $y_0$ and $s$. Roughly speaking, we would like 
to use the function $f_{y_0,s}(X)$ defined by 
\begin{equation}\label{a148}
f_{y_0,s}(X) := \frac{m(B(y_0,s) \cap \Omega)}{s^2} g(X,Y_0) - K_2 \left[ \omega^X(\Gamma \setminus B(y_0,K_1 s)) + \frac{m(B(y_0,K_1s) \cap \Omega)}{(K_1s)^2} g(X,Y_{K_1}) \right]
\end{equation}
where $Y_{K_1}$ is a Corkscrew point associated to $y_0$ and $K_1s$, 
and where $K_1,K_2>0$ are some large constants 
that depend only on $n$, $C_1$ to $C_6$, and$C_A$. We could show that with large enough choices of $K_1$ and $K_2$,
$f_{y_0,s}$ is positive in $B(y_0,s)$ and negative outside of a big ball $B(y_0,2K_1s)$. 
However, we want to use the maximum principle to extend such inequalities to larger regions,
and with this definition involving the harmonic measure, our $f_{y_0,s}$ is not smooth enough 
to be used in Lemma \ref{lMPg}. So we shall first replace 
$\omega^X(\Gamma \setminus B(y_0,K_1s))$ in \eqref{a148}
by some solution of $L u=0$ with smooth Dirichlet condition.

\medskip

Let $h\in C^\infty(\R^n)$ be such that $0\leq h \leq 1$, $h\equiv 0$ on $B(0,1/2)$ and 
$h\equiv 1$ on the complement of  $B(0,1)$. 
For $\beta > 1$ 
(which will be chosen large), we define $h_{\beta}$ by $h_\beta(x) = h(\frac{x-y_0}{\beta s})$. 
Let $u_\beta$ be the solution, given by Lemma \ref{lLM}, 
of $L u_\beta=0$ with the Dirichlet condition $\Tr u_\beta = \Tr h_\beta$. 
Notice that $u_\beta \in W$ because $1-u_\beta$ is the solution of $L$ with the smooth and compactly supported trace of 
$1-h_\beta$.
By the positivity of the harmonic measure, 
it holds that 
for any $X\in \Omega$ and $\gamma>0$,
\begin{equation} \label{CPA7}
 \omega^X(\Gamma \setminus B(y_0,\beta s)) \leq u_\beta(X) \leq  \omega^X(\Gamma \setminus B(y_0,\beta s/2)),
\end{equation}
and we can see here that $u_\beta$ will be used as a smooth substitute of the harmonic measure.

\medskip

Similarly to \eqref{tcp18b}, using the Green function upper bounds and (\hyperref[H4]{H4}), we have
\begin{equation} \label{CPA1}
 \frac{m(B(y_0,s) \cap \Omega)}{s^2} g(X,Y_0) \leq C  \qquad \ \text{ for } 
X \in B(Y_0,\delta(Y_0)/2) \setminus B(Y_0,\delta(Y_0)/4).
\end{equation}
Then by 
Lemma \ref{lMPg} with $E=\R^n \sm B(Y_0,\delta(Y_0)/4)$ and $F = \R^n \sm \overline{B(Y_0,\delta(Y_0)/2)}$, 
\begin{equation} \label{CPA3b}
 \frac{m(B(y_0,s) \cap \Omega)}{s^2} g(X,Y_0) \leq C \ \ \ \text{ for }
X \in \Omega \sm B(Y_0,\delta(Y_0)/4).
\end{equation}
From this and
the non-degeneracy of the harmonic measure (Lemma \ref{ltcp4}), we deduce
that for $\beta \in (1,\infty)$ and $X \in \Omega \setminus B(y_0,2\beta s)$
\begin{equation} \label{CPA3} \begin{split}
 \frac{m(B(y_0,s) \cap \Omega)}{s^2} g(X,Y_0) & 
 \leq K_2 \left[\omega^X (\Gamma \setminus B(y_0,\beta s))
 +  \frac{m(B(y_0,\beta s) \cap \Omega)}{(\beta s)^2} g(X,Y_\beta)) \right] \\
 & \leq K_2 \left[ u_\beta(X) + \frac{m(B(y_0,\beta s) \cap \Omega)}{(\beta s)^2} g(X,Y_\beta))\right],
\end{split} \end{equation}
where $Y_\beta$ is a corkscrew point associated to $(y_0,\beta s)$, and where the constant $K_2>0$ depends only on $n$, $C_1$ to $C_6$, and $C_A$; in particular, $K_2$ does not depend on $\beta$.

\medskip

Our aim now is to find $K_1\geq 20C_1$ such that, for 
$X\in \Omega \cap [B(y_0,s) \setminus B(Y_0,\delta(Y_0)/4)]$, 
\begin{equation} \label{CPA6}
K_2 \left[ u_{K_1}(X)  + \frac{m(B(y_0,K_1 s) \cap \Omega)}{(K_1 s)^2} g(X,Y_{K_1}) \right] \leq \frac12 \frac{m(B(y_0,s) \cap \Omega)}{s^2} g(X,Y_0).
\end{equation}

According to the H\"older continuity at the boundary (Lemma \ref{HolderB}), we have 
\begin{equation} \label{CPA4}
\sup_{B(y_0,10 s)} u_\beta \leq C \beta^{-\alpha}. 
\end{equation}
Moreover, Lemma \ref{GreenDi} 
applied to $g(.,Y_\beta)$ implies, that
\begin{equation} \label{CPA4z}
\frac{m(B(y_0,\beta s) \cap \Omega)}{(\beta s)^2} g(X,Y_\beta) 
\leq C \beta^{-\alpha}
\qquad \text{ for } X \in B(y_0,10s).
\end{equation}
If both cases, we need $\beta$ to be big enough, for instance $\beta \geq 20C_1$, and the constants $C$ and $\alpha>0$ depend only 
on $n$, $C_1$ to $C_6$, and $C_A$. 
Due to the non-degeneracy of the harmonic measure given by \eqref{tcp33b} - and (\hyperref[H4]{H4}) - we have 
\[u_4 + \frac{m(B(y_0,s) \cap \Omega)}{s^2} g(X,Y_0) \geq C^{-1} \quad \text{ for } X \in \Omega \setminus B(y_0,5s).\]
By  
the last estimate in \eqref{CPA4}--\eqref{CPA4z},
there exists 
$K_3>0$ such that 
\begin{equation} \label{CPA4b}
u_\beta(X) + \frac{m(B(y_0,\beta s) \cap \Omega)}{(\beta s)^2} g(X,Y_\beta) \leq K_3 \beta^{-\alpha} \left[ u_4(X) + \frac{m(B(y_0,s) \cap \Omega)}{s^2} g(X,Y_0) \right]
\end{equation}
for  
$X \in \Omega \cap [B(y_0,10s) \sm B(y_0,5s)]$. In addition, by increasing $K_3$ if needed, 
the estimates \eqref{CPA4}--\eqref{CPA4z} and the 
lower bound in \eqref{GreenE2} for the Green function imply that
\begin{equation} \label{CPA4d}
u_\beta(X) + \frac{m(B(y_0,\beta s) \cap \Omega)}{(\beta s)^2} g(X,Y_\beta) \leq K_3 \beta^{-\alpha} \frac{m(B(y_0,s) \cap \Omega)}{s^2} g(X,Y_0) 
\end{equation}
when $X\in B(Y_0,\delta(Y_0)/2) \sm B(Y_0,\delta(Y_0)/4)$. We invoke then Lemma \ref{lMPg}, used on the function 
\[X \to K_3 \beta^{-\alpha} \left[ u_4(X) + \frac{m(B(y_0,s) \cap \Omega)}{s^2} g(X,Y_0)  \right]
- \left[ u_\beta(X) + \frac{m(B(y_0,\beta s) \cap \Omega)}{(\beta s)^2} g(X,Y_\beta) \right]\]
and the sets $E = B(y_0,10s) \setminus \overline B(Y_0,\delta(Y_0)/4)$ and 
$F = \overline {B(y_0,5s)} \sm B(Y_0,\delta(Y_0)/2)$, to deduce that 
for $X \in \Omega \cap B(y_0,10s) \sm B(Y_0,\delta(Y_0)/4)$,
\[ \label{CPA4c} \begin{split}
u_\beta(X) + \frac{m(B(y_0,\beta s) \cap \Omega)}{(\beta s)^2} g(X,Y_\beta) & \leq K_3 \beta^{-\alpha} \left[ u_4(X) + \frac{m(B(y_0,s) \cap \Omega)}{s^2} g(X,Y_0) \right] \\
& \leq K_3 \beta^{-\alpha} \left[ \omega^X(B(y_0,4s)) + \frac{m(B(y_0,s) \cap \Omega)}{s^2} g(X,Y_0) \right] \\
& \leq C \beta^{-\alpha} \frac{m(B(y_0,s) \cap \Omega)}{s^2} g(X,Y_0)
\end{split} \]
by Lemma \ref{ltcp3}. Therefore, \eqref{CPA6} can be indeed achieved for some large $K_1$, that depends only on $n$, $C_1$ to $C_6$, and $C_A$ (recall that $K_2$ is already fixed, and depends on the same parameters).

\medskip

Define the function $f_{y_0,s}$ on $\Omega \setminus \{Y_0\}$ by
\begin{equation} \label{CPA8}
f(X): =  \frac{m(B(y_0,s) \cap \Omega)}{s^2} g(X,Y_0) - K_2 \left[ u_{K_1}(X) + \frac{m(B(y_0,K_1 s) \cap \Omega)}{(K_1 s)^2} g(X,Y_{K_1}))\right].
\end{equation}
The inequality \eqref{CPA3} gives
\begin{equation} \label{CPA9}
f_{y_0,s}(X) \leq 0 \ \ \ \text{ for }
X \in \Omega  \setminus B(y_0,2K_1s),
\end{equation}
and the estimate \eqref{CPA6} proves
that
\begin{equation} \label{CPA0}
f_{y_0,s}(X) \geq \frac12  \frac{m(B(y_0,s) \cap \Omega)}{s^2} g(X,Y_0) 
\end{equation}
for $X \in \Omega \cap [B(y_0,s) \setminus B(Y_0,\delta(Y_0)/4)]$.
\medskip

\noindent {\bf Step 2:} End of the proof

\medskip

Let us turn to the proof of the comparison principle. 
By symmetry and as
in Lemma~\ref{lCP}, it suffices to prove the upper bound in \eqref{CP26}, that is
  \begin{equation} \label{CP32} 
\frac{u(X)}{v(X)} \leq C \frac{u(X_0)}{v(X_0)} \ \ \ \text{ for } 
X \in \Omega \cap B,  \text{ where }  B= B(x_0,r).
 \end{equation}

\medskip

We claim that
  \begin{equation} \label{CP33} 
v(X) \geq C^{-1} \frac{m(B)}{r^2} v(X_0) g(X,X_0) 
\qquad \text{ for } X \in [\Omega \cap B] \sm B(X_0,\delta(X_0)/4),
 \end{equation}
 where $C>0$ depends only on $n$, $C_1$ to $C_6$, and $C_A$.
So let $X \in \Omega \cap B$ be given. Two cases may happen. 
If $\delta(X) \geq \frac r{8K_1}$, where $K_1$ comes from \eqref{CPA6} 
and is the same as in the definition of $f_{y_0,s}$, the
existence of Harnack chains (Proposition \ref{propHarnack}) and the Harnack inequality (Lemma \ref{HarnackI}) give that
\[v(X) \approx v(X_0) 
\]
For the above inequality to hold, we need the Harnack chains to stay in the area where $v$ is a solution; we take $K$ big enough to make sure that it happens, and by Proposition \ref{propHarnack}, $K$ need to depend only on $C_1$ and $C_2$.
Similarly, Proposition \ref{propHarnack} and Lemma \ref{HarnackI}, 
together with the bound \eqref{GreenE2} on the Green function,
give that
\[\frac{m(B)}{r^2} g(X,X_0) \approx 1 \qquad \text{ on } 
[\Omega \cap B(x_0,r)] \sm B(X_0,\delta(X_0)/4).\]
We conclude that for all $X \in [\Omega \cap B] \sm B(X_0,\delta(X_0)/4)$ satisfying $\delta(X) \geq \frac r{8K_1}$.
   \begin{equation} \label{CP34} 
v(X) \approx  v(X_0) \frac{m(B)}{r^2} g(X,X_0).
 \end{equation} 
The more interesting remaining case is when 
 $\delta(X) < \frac{r}{8K_1}$. 
 Take $y_0 \in \Gamma$ such that $|X-y_0| = \delta(X)$. Set 
 $s:= \frac r{8K_1}$ and $Y_0$ a corkscrew point associated to $(y_0,s)$.  
 The ball $B(y_0, \frac12 r) = B(y_0,8K_1s)$ is contained
 in $B(x_0,\frac 74 r)$. The following points hold :  
 \begin{itemize}
 \item 
 The quantity  $\int_{B(y_0,4K_1s) \setminus B(Y_0,\delta(Y_0)/4)} |\nabla v|^2 dm$ is finite 
 because $v \in \WW(B(x_0,2r))$. The fact that 
 $\int_{B(y_0,4K_1s) \setminus B(Y_0,\delta(Y_0)/4)} |\nabla f_{y_0,s}|^2 dm$ is finite as well 
 follows from the property \eqref{GreenE1} of the Green function. 
 \item 
 There exists $K_4>0$ such that 
  \begin{equation} \label{CP35} 
K_4 v(Y) - v(Y_0) f_{y_0,s}(Y) \geq 0 \ \ \ \text{ for }
Y \in B(Y_0,\delta(Y_0)/2) \setminus B(Y_0,\delta(Y_0)/4). \end{equation}
This latter inequality is due to the following two bounds: the fact that
 \begin{equation} \label{CP36} 
f_{y_0,s}(Y) \leq \frac{m(B(y_0,s) \cap \Omega)}{s^2}g(Y,Y_0) \leq C \ \ \ \text{ for }
Y \in B(Y_0,\delta(Y_0)/2) \setminus B(Y_0,\delta(Y_0)/4),
 \end{equation}
which is a consequence of
the definition \eqref{CPA8} and 
\eqref{GreenE7}, and the bound
  \begin{equation} \label{CP37} 
v(Y) \geq C^{-1}v(Y_0) \ \ \ \text{ for }
Y \in B(Y_0,\delta(Y_0)/2),
 \end{equation}
which comes from the Harnack inequality (Lemma \ref{HarnackI}).
 \item 
 The function $K_4 v - v(Y_0) f_{y_0,s}$
is nonnegative on 
 $\Omega \cap [B(y_0,4K_1s) \setminus B(y_0,2K_1s)]$. Indeed, $v\geq 0$ on $\Omega \cap B(y_0,4K_1s) $ and, thanks to \eqref{CPA9}, $f_{y_0,s} \leq 0$ on $\Omega \setminus B(y_0,2K_1s)$.
 \item 
 The trace of $K_4 v- v(Y_0) f_{y_0,s}$ is non-negative on 
$B(y_0,4K_1s) \cap \Gamma$
 again because $\Tr v = 0$ on 
$B(y_0,4K_1s) \cap \Gamma$ 
 and $\Tr [f_{y_0,s}] \leq 0$ on 
$B(y_0,4K_1s) \cap \Gamma$ 
 by construction.
 \end{itemize} 
 The previous points prove that $K_4 v- v(Y_0) f_{y_0,s}$ satisfies the assumptions of 
 Lemma \ref{lMPg} with
 $E = B(y_0,4K_1s) \setminus \overline{B(Y_0,\delta(Y_0)/4)}$ and $F = \overline{B(y_0,2K_1s)} \setminus B(Y_0,\delta(Y_0)/2)$. As a consequence, for any $Y\in B(y_0,4K_1s) \setminus B(Y_0,\delta(Y_0)/4)$
  \begin{equation} \label{CP38} 
K_4 v(Y) - v(Y_0) f_{y_0,s}(Y) \geq 0,
 \end{equation} 
 and hence, 
 for any $Y \in B(y_0,s) \setminus B(Y_0,\delta(Y_0)/4)$
   \begin{equation} \label{CP39} 
v(Y) \geq  (K_4)^{-1} v(Y_0) f_{y_0,s}(Y) \geq C^{-1} \frac{m(B(y_0,s) \cap \Omega)}{s^2} v(Y_0) g(Y,Y_0)
 \end{equation} 
 by  
 \eqref{CPA0}. The points $X_0$ and $Y_0$ are both corkscrew points, and they can be linked by a Harnack chain of balls whose length depends only on $r/s = 8K_1$, that is with uniformly bounded length. So using the 
 Harnack inequality on every ball of the chain,  
 we deduce that $v(Y_0) \approx v(X_0)$ and $g(Y,Y_0) \approx g(Y,X_0)$ whenever $Y \in B(y_0,s)$ is far from $Y_0$ and $X_0$ (but it cannot be close to $X_0$ since $K_1 \geq 20C_1$). Moreover, $s^{-2} m(B(y_0,s) \cap \Omega) \approx r^{-2} m(B \cap \Omega)$ by (\hyperref[H4]{H4}). Therefore \eqref{CP38} becomes 
\[v(Y) \geq C^{-1} \frac{m(B \cap \Omega)}{r^2} v(X_0) g(Y,X_0) \qquad \text{ for } Y \in [\Omega \cap B(y_0,s)] \sm B(Y_0,\delta(Y_0)/4).\]
Since the two functions of $Y$ in the inequality above are solutions in $\Omega \cap B(y_0,2s)$, the Harnack inequality yields the following improvement:
 \[v(Y) \geq C^{-1} \frac{m(B \cap \Omega)}{r^2} v(X_0) g(Y,X_0) \qquad \text{ for } Y \in \Omega \cap B(y_0,s).\]
 Recall that $X \in B(y_0,s)$ by construction of $y_0$ and $s$. We conclude, at last, that even when $X \in \Omega \cap B$ is such that $\delta(X) < \frac r{8K_1}$, we still have
 \[ v(X) \geq C^{-1} \frac{m(B)}{r^2} v(X_0) g(X,X_0).\]
 The claim \eqref{CP33} follows.
 
 \medskip
 
Now we want to prove that,  for all $X\in [\Omega \cap B] \sm B(X_0,\delta(X_0)/4)$,
\begin{equation} \label{CP42} 
u(X) \leq C u(X_0) \left[ \omega^X(\Gamma \setminus \frac54 B) + \frac{m(B \cap \Omega)}{r^2} g(.,X_{0})\right],
 \end{equation} 
 where the constant $C>0$ depends only on $n$, $C_1$ to $C_6$, and $C_A$.
By Lemma \ref{ltcp2},
\begin{equation} \label{CP43} 
u(X) \leq C u(X_0)  \ \ \ \text{ for } 
X \in \Omega \cap \frac74 B,
 \end{equation} 
as long as $K$ is large enough so that Lemma \ref{ltcp2} can be applied. 
But again, $K$ does not need to depend on anything else than 
$n$, $C_1$ and $C_2$. 
Pick  $h'\in C^\infty(\R^n)$ such that $0 \leq h'\leq 1$, $h'\equiv 1$ outside of $\frac32 B$, 
and $h' \equiv 0$ on $\frac54 B$. 
Let $u_{h'} = U(h')$ be the solution of $Lu_{h'} = 0$ with 
the data $\Tr u_{h'} = \Tr h'$ 
(given by Lemma~\ref{lLM}). 
As before, $u_{h'} \in W$ because $1-u_{h'} = U(1-h)$ and $1-h$ is a test function.
Also, $u_{h'}(X) \geq \omega^X(\Gamma \setminus \frac32 B)$ by monotonicity.
So \eqref{tcp33b}, which 
states the non-degeneracy of the harmonic measure, gives
\begin{equation} \label{CP44} 
u_{h'}(X) + \frac{m(\frac32B \cap \Omega)}{(\frac 32 r)^2} g(X,X_{3/2})\geq C^{-1} \ \ \ \text{ for } 
X\in \Omega \setminus \frac{13}{8}B,
 \end{equation} 
where $X_{3/2}$ is a corkscrew point associated to $(x_0,\frac32 r)$. The 
doubling property (\hyperref[H4]{H4}) and Harnack inequality for the function $g(X,.) = g_T(.,X)$ entail now that
\begin{equation} \label{CP44b} 
u_{h'}(X) + \frac{m(B \cap \Omega)}{r^2} g(X,X_{0})\geq C^{-1} \ \ \ \text{ for } 
X\in \Omega \setminus \frac{13}{8}B.
 \end{equation} 
The combination of \eqref{CP43} and \eqref{CP44b} yields the existence of $K_5>0$ such that 
\[\tilde u:= K_5 u(X_0) \left[ u_{h'} + \frac{m(B \cap \Omega)}{r^2} g(.,X_{0})\right] - u \geq 0 \quad \text{ on } \Omega \cap [\frac74 B \setminus \frac{13}{8}B].\]
Moreover, using the Harnack inequality and the Green function lower bounds, by increasing slightly $K_5$ if needed, we also have $\tilde u \geq 0$ in $B(X_0,\delta(X_0)/2) \sm B(X_0,\delta(X_0)/4)$. Now, it is not very hard to see that $\tilde u$ satisfies all the assumptions of Lemma \ref{lMPg}, with the sets $E = \frac74 B \sm \overline{B(X_0,\delta(X_0)/4)}$ and $F =  \overline{\frac{13}{8}B} \sm B(X_0,\delta(X_0)/2)$.
Observe in particular that $u_{h'} \in W$, $T u_{h'} \geq 0$, and as long as we choose $K\geq 2$, $u \in \WW(2B)$ and $Tu = 0$ on $\Gamma \cap 2B$.
Then by Lemma \ref{lMPg},
\begin{equation} \label{CP45} 
u(X) \leq K_5 u(X_0) \left[ u_{h'} (X)+ \frac{m(B \cap \Omega)}{r^2} g(X,X_{0})\right] \ \ \ \text{ for }
X\in [\Omega \cap \frac74B] \sm B(X_0,\delta(X_0)/4),
 \end{equation} 
and since $u_{h'}(X) \leq \omega^X(\Gamma \setminus \frac54 B)$
for all $X\in \Omega$, 
\[ 
u(X) \leq C u(X_0) \left[ \omega^X(\Gamma \setminus \frac54 B) 
+ \frac{m(B \cap \Omega)}{r^2} g(X,X_{0})\right] 
\ \ \ \text{ for } 
X\in [\Omega \cap \frac74B] \sm B(X_0,\delta(X_0)/4).\]
The claim \eqref{CP42} follows.

\medskip
 
The bounds \eqref{CP33} and \eqref{CP42} imply that
\begin{equation} \label{CP47} 
\frac{u(X)}{v(X)} \leq C \frac{u(X_0)}{v(X_0)} 
\left[\frac{ r^2 \omega^X(\Gamma \setminus \frac54B)}{m(B \cap \Omega) g(X,X_0)} +1 \right] 
\ \ \ \text{ for }
X\in [\Omega \cap B] \setminus B(X_0, \delta(X_0)/4).
 \end{equation} 
 The bound \eqref{CP32} in the set $[\Omega \cap B] \setminus B(X_0, \delta(X_0)/4)$ is now a consequence of the above inequality and \eqref{tcp18a}. The bound \eqref{CP32} in the full domain $\Omega \cap B$ is then an easy consequence of the Harnack inequality (Lemma \ref{HarnackI}).
\ep

\end{document}